%% file: main.tex
\documentclass[11pt, oneside]{book}
\pdfoutput=1
\usepackage{subfiles}

\usepackage[bib]{preamble/packages}
\usepackage{preamble/commands}
\usepackage{preamble/formatting}

\usepackage{hyperxmp}
\usepackage[
    type={CC},
    modifier={by-nc-nd},
    version={4.0},
    imageposition={left}
]{doclicense}

\usepackage{fancyhdr}

\title{An introduction to filtered and synthetic spectra}
\author{Sven van Nigtevecht}
\date{23 December 2025}

\begin{document}
\frontmatter
\maketitle

\chapter*{Preface}

These are expository notes about spectral sequences, filtered spectra, and synthetic spectra.
Feedback is always appreciated, so please do not hesitate to contact me in case you find typos, corrections, or have other comments.
I can make no claim to originality of the contents of these notes, except for the arrangement and presentation of the results.
The baseline for these notes is a mini-course I gave at Utrecht University in the fall of 2024.
The current version of these notes is essentially a lengthy excerpt from my PhD thesis, but unfortunately does not yet contain all the material I covered in the mini-course.
The main subjects that are currently missing are: an in-depth example computation, Goerss--Hopkins obstruction theory, and a look at the construction of the $\infty$-category of synthetic spectra (and variants thereof).
I hope to include these in later versions.
For the moment, I would direct readers to \cite{CDvN_part2} for a sample computation, and to \cite{pstragowski_vankoughnett_obstruction_theory} for an introduction to obstruction theories.

\section*{Acknowledgements}

First and foremost, I would like to thank my PhD advisor, Lennart Meier, for all of the extensive feedback he provided on everything I wrote, for the problems he suggested me to work on, and for the way he shaped my mathematical thinking.
Much of the material in these notes was born out of my collaborations with Christian Carrick and Jack Davies.
I would like to thank Shaul Barkan, Robert Burklund, Irakli Patchkoria, and Piotr Pstr\k{a}gowski for teaching me so much synthetic mathematics, for answering my questions patiently, and for the interest they showed in my work, which was incredibly encouraging.
Large portions of these notes grew out of my attempts to better understand the Omnibus Theorem of Burklund--Hahn--Senger; their work has been very influential.
Many thanks to John Rognes for his many and detailed comments and corrections, and thank you as well to Miguel Barata, Maite Carli, Bastiaan Cnossen, Marius Nielsen, and Ryan Quinn for catching typos and feedback on earlier versions.

These notes were written as part of my PhD, which was supported by the NWO grant \texttt{VI.Vidi.193.111}.
Parts of these notes were written during my stay at the Isaac Newton Institute for Mathematical Sciences, Cambridge, during the programme \emph{Equivariant homotopy theory in context} (supported by EPSRC grant \texttt{EP/Z000580/1}).

\tableofcontents

\mainmatter

\chapter{Introduction}
\label{ch:intro_to_part1}

\subfile{intro_fil_syn}


\chapter{Filtered spectra and spectral sequences}
\label{ch:filtered_spectra_sseqs}

\subfile{fil/basics_abridged}

\chapter{The \texorpdfstring{$\tau$}{tau}-formalism}
\label{ch:tau_formalism_filtered}

\subfile{fil/tau_abridged}

\subfile{fil/taubss}

\subfile{fil/omnibus}

\subfile{fil/deform2}

\chapter{Synthetic spectra}
\label{ch:synthetic_spectra}

\subfile{syn/syn}

\chapter{Variants of synthetic spectra}
\label{ch:variants_Syn}

\subfile{syn/syn_models}

\appendix

\subfile{fil/sseq_informal}

\printbibliography[heading=bibintoc]

\end{document}

%% file: intro_fil_syn.tex
It is no exaggeration to say that the road to modern homotopy theory is paved with spectral sequences.
Unfortunately, one big barrier to using spectral sequences is the abundance of indices, maps, and diagrams hiding in their definition.
This might lead one to regard the construction of a spectral sequence as a black box, only to be opened in the most dire of circumstances.
This is particularly unfortunate because there is much power to be had in working directly with the filtrations, much like how working with spectra has proved to be a better approach than working directly with homology theories.

We believe that this need not be so.
With these notes, we set out to show that working directly with filtrations is not only possible, but practical.
An important reason for this is the \emph{$\tau$-formalism}, which provides a way to off-load much of the notational headache.
More than merely being convenient for working with filtrations, we show that taking $\tau$ seriously also makes it possible to understand more exotic categories of filtrations, such as \emph{synthetic spectra}.

Our goal, then, is to give a mostly self-contained introduction to spectral sequences, filtered spectra and synthetic spectra, all through the lens of the $\tau$-formalism.
The essence of this part may be summarised by the following slogan.
\begin{center}
    \emph{Filtered spectra are to spectral sequences,\\
    as synthetic spectra are to Adams spectral sequences,\\
    as spectra are to homology theories.}
\end{center}
We mean this in the following sense.
In the second clause of each analogy, we are referring to the purely algebraic objects.
These lack good categorical properties (e.g., the category of spectral sequences is not abelian, is not monoidal, etc.).
Working with the objects of the first clause offers a way to remedy these problems, as these give rise to the corresponding algebraic objects, but constitute a good homotopy theory (having homotopy limits and colimits, a symmetric monoidal structure, etc.).

In the first analogy, there are some slight caveats.
First, as spectra are inherently stable objects, one can only hope to capture spectral sequences arising in the stable setting, so we ignore for these purposes spectral sequences coming from, e.g., towers of spaces.
However, even within this more specialised setting, there are spectral sequences that do not arise from filtered spectra.
In practise, most spectral sequences do come from a filtered spectrum, so we view this as more of a technicality.
In this sense, their relationship is akin to the one between stable $\infty$-categories and triangulated categories.

For synthetic spectra, the situation is the opposite: for every spectrum, its Adams spectral sequence comes from a preferred synthetic spectrum, but not every synthetic spectrum captures an Adams spectral sequence.
However, far from being a technical point, this is a key feature of synthetic spectra.
A general synthetic spectrum can be thought of as a \emph{modified Adams spectral sequence}: the mere fact that it lives in the synthetic category means that it has a much closer relationship to Adams spectral sequences than a plain filtered spectrum would have.

Almost all results in this part are well known or folklore.
The value, we believe, lies in having all of these results in one place.
In this introduction, we give an overview of the main results, discuss some of the history of synthetic spectra, and end with a more detailed outline of these notes.

\section{Filtered spectra}

A \emph{filtered spectrum} is a functor $X \colon \Z^\op \to \Sp$, where $\Z$ denotes the poset of the integers under the usual ordering.
We usually depict this as a sequence
\[
    \dotsb \to X^1 \to X^0 \to X^{-1} \to \dotsb.
\]
We refer to the maps between the individual spectra as the \emph{transition maps}.
This comes with a notion of bigraded homotopy groups, being given by
\[
    \oppi_{n,s}X = \pi_n(X^s).
\]
We write $\FilSp$ for the $\infty$-category of filtered spectra.
Note, however, that this use of $\infty$-categories is merely preferential.
In fact, nearly all results on filtered spectra in these notes are not inherently modern in any way, and could have been obtained long ago, even right alongside the introduction of spectral sequences arising from filtered chain complexes.

\begin{remark}
    Although we have written everything in terms of filtered spectra, readers who are more familiar with chain complexes can instead work with \emph{filtered chain complexes}.
    The definition is the same: it is a diagram
    \[
        \dotsb \to C^1 \to C^0 \to C^{-1} \to \dotsb
    \]
    where each $C^s$ is a chain complex, and each transition map is a map of chain complexes.
    (Beware, then, that the upper index is \emph{not} the chain complex degree.
    Alternatively, one might write $C^s_\bullet$ for the chain complex in position $s$.)
    There is some further change in terminology:
    \begin{itemize}
        \item Instead of working with homotopy groups, one should work with the homology groups of chain complexes.
        \item Instead of working with suspensions $\Sigma$, one should work with the shift operator $[1]$.
        \item Instead of working with cofibres, one should work with the mapping cone of chain complexes.
        If the map of chain complexes is injective, then this is (quasi-isomorphic to) the quotient in the usual sense.
        Up to quasi-isomorphism, one can replace a filtered chain complex by one where all maps to be injective, so this is not much of a restriction.
    \end{itemize}
\end{remark}

If $X$ is a filtered spectrum, then we let $X^{-\infty}$ denote its colimit.
We think of $X$ as a tool to understand its colimit.
More precisely, the homotopy groups of $X^{-\infty}$ are given by the colimit over the transition maps:
\[
    \pi_n X^{-\infty} \cong \colim_s \oppi_{n,s} X.
\]
Knowing the homotopy groups $\oppi_{n,s} X$ for all $s$ is more information than knowing the group $\pi_n X^{-\infty}$.
Even if one is only interested in knowing $\pi_n X^{-\infty}$, it is nevertheless a good idea to remember $\oppi_{n,s}X$ for every $s$, along with all transition maps between them.

However, in practice, the difficulty of understanding the homotopy groups of the spectra $X^s$ is on par with those of $X^{-\infty}$, even for cleverly chosen filtrations.
Usually, the homotopy groups of the cofibres of the transition maps are much easier to compute; we refer to these cofibres as the \emph{associated graded} spectra, and write
\[
    \Gr^s X \defeq \cofib(X^{s+1} \to X^s).
\]
The (attempted) passage from the homotopy groups $\pi_n \Gr^s X$ to the homotopy groups $\oppi_{n,s}X$ is exactly the structure of a \emph{spectral sequence}.

We give an informal but in-depth introduction to how the spectral sequence arises from a filtration in \cref{ch:informal_sseq}.
For the moment, we give a rough indication.
Fixing an integer $n$ and applying $\pi_n$ to the diagram $X$, we obtain a filtered abelian group.
For a fixed $s$, certain elements in $\pi_n X^s$ may not be in the image of the map $\pi_n X^{s+1}\to\pi_n X^s$; we think of these as `being born' at stage $s$.
An element might be sent to zero after an application of a number of transition maps $\pi_n X^s \to \pi_n X^{s-r}$ (for some $r\geq 1$); we think of these as `dying' at some later point.
A spectral sequence encodes the event of an element being born and dying $r$ steps to the right by a \emph{differential} of length $r$.

At this point, it becomes useful to introduce some notation.
We reserve the formal symbol $\tau$, and let it act on the homotopy groups $\set{\oppi_{*,*}X}_{n,s}$ via the transition maps: if $\alpha \in \oppi_{n,s}X = \pi_n(X^s)$ is an element, then we define $\tau \cdot \alpha$ to be the image of $\alpha$ under $X^s \to X^{s-1}$.
This turns $\oppi_{*,*}X$ into a bigraded $\Z[\tau]$-module, where $\tau$ has bidegree $(0,-1)$.
By the previous discussion, this means that $\tau^r$-torsion elements (i.e., elements that are annihilated by multiplication by $\tau^r$) correspond to differentials of length $r$ or shorter in the spectral sequence.
One can make the translation between the spectral sequence and the $\Z[\tau]$-module $\oppi_{*,*}X$ more precise; the resulting Rosetta stone is known as the \emph{Omnibus Theorem}.
Its proof essentially revolves around making the diagram chase of \cref{ch:informal_sseq} very precise, and also taking into account potential convergence issues.
We summarise this usage of $\tau$ and the surrounding comparison results by calling it the \emph{$\tau$-formalism}.

There are many advantages to remembering the filtration that gives rise to this spectral sequence, rather than only remembering the latter.
The boundary maps $\Gr^s X \to \Sigma X^{s+1}$ encode information about differentials of \emph{all} lengths.
We refer to these are \emph{total differentials}.
Using these instead of the ordinary differentials leads, for instance, to a significantly strengthened version of the Leibniz rule, which we refer to as the \emph{filtered Leibniz rule}.
None of this would be possible when working with bare spectral sequences.
The $\tau$-formalism helps us in keeping track of total differentials, but is again not inherently necessary.

More than help us with internal computations, the $\tau$-formalism also governs the structure of the $\infty$-category of filtered spectra.
This is particularly important for exporting the $\tau$-formalism to other contexts, where it can become significantly more powerful.
The map $\tau$ can be lifted to be a map in the $\infty$-category $\FilSp$.
We can perform homotopical constructions with it, such as forming its cofibre $C\tau$.
It turns out that $C\tau$ admits an $\E_\infty$-algebra structure in $\FilSp$.
Moreover, the associated graded functor
\[
    \FilSp \to \grSp, \quad X \mapsto \Gr X
\]
can be identified with the functor
\[
    \FilSp \to \Mod_{C\tau}(\FilSp), \quad X \mapsto C\tau \otimes X.
\]

In the main text, we discuss the general theory of exporting this formalism, through what has become known as \emph{deformations}.
For the purposes of this introduction, we will focus on our main example: the case of \emph{synthetic spectra}.
This will also allow us to discuss the differences between our account of the $\tau$-formalism and the existing literature.

\section{Synthetic spectra}


The generality of filtered spectra makes them very useful: proving something about filtered spectra yields applications to all spectral sequences.
However, their generality can also make them a little unwieldy.
Say, for instance, that we are working with a particular type of spectral sequence where the first page has more structure than merely the homotopy groups of a spectrum.
It would then be desirable to work in a modification of filtered spectra where this additional structure exists in the category itself.
This additional structure should make the category easier to work with, making it a more powerful tool for studying that specific type of spectral sequence.
In the case of Adams spectral sequences, this is exactly what the $\infty$-category of \emph{synthetic spectra} is.

Let $E$ be a multiplicative homology theory.
The $E$-based Adams spectral sequence tries to approximate maps between spectra by maps between their $E$-homology.
Taking one spectrum to be a sphere, we obtain a spectral sequence that tries to compute homotopy groups.
Under hypotheses on $E$, this is of the form
\[
    \uE_2^{n,s} \cong \Ext_{E_*E}^{s,\, n+s}(E_*,\, E_*X) \implies \pi_n X,
\]
where the Ext groups refer to Ext groups of $E_*E$-comodules, which roughly speaking is remembering $E$-homology (co)operations present on the $E$-homology of a spectrum.
The original case introduced by Adams \cite{Adams_Adams_sseq} is the one where $E$ is $\F_p$-homology, which to this day remains the main tool for computing stable homotopy groups of spheres.
Another popular, more chromatic flavour is the case where $E=\MU$, which is referred to as the \emph{Adams--Novikov spectral sequence} (abbreviated ANSS).
For general $E$, both the computation of the Ext groups as well as its differentials are highly nontrivial tasks.
In the case $E=\F_2$, the state of the art in terms of computing the Ext groups is for for $n+s\leq 200$ by Bruner--Rognes \cite{bruner_rognes_cohomology_steenrod_algebra}, and more recently for $n+s\leq 261$ by Lin \cite{lin_noncommutative_grobner}, and the state of the art in computing differentials is with almost complete information up until dimension~$90$ by Isaksen--Wang--Xu \cite{isaksen_wang_xu_dimension_90}, and with very recent further information going up until dimension~$126$ by Lin--Wang--Xu \cite{lin_wang_xu_last_kervaire_dimension}.
In these notes, our focus is mostly on how to compute differentials.



Rather than thinking of Adams spectral sequences as filtered spectra arising in a particular way, synthetic spectra let us picture them as living in their own category.
There is a stable $\infty$-category $\Syn_E(\Sp)$ of \emph{$E$-synthetic spectra}, along with functors
\[
    \nu_E \colon \Sp \to \Syn_E(\Sp) \qquad \text{and} \qquad \sigma \colon \Syn_E(\Sp)\to \FilSp,
\]
called the \emph{$E$-synthetic analogue} and \emph{signature} functor, respectively.
A series of computations in $\Syn_E(\Sp)$ shows that the composite $\sigma \circ \nu_E$ is, in a precise sense, the $E$-Adams spectral sequence functor.
We think of $\sigma$ as a forgetful functor, sending a synthetic spectrum to its `underlying spectral sequence'.



Not only can $\sigma$ be thought of as an underlying spectral sequence functor, it is also the mechanism through which we import the $\tau$-formalism.
Namely, $\sigma$ is the right adjoint in an adjunction
\[
    \begin{tikzcd}[column sep=3.5em]
        \FilSp \rar[shift left,"\rho"] & \Syn_E(\Sp). \lar[shift left, "\sigma"]
    \end{tikzcd}
\]
Synthetic spectra come with their own notion of bigraded spheres, and these are the image under $\rho$ of filtered bigraded spheres.
By adjunction therefore, if $X$ is a synthetic spectrum, then $\pi_{*,*}X$ is captured by the filtered homotopy groups $\pi_{*,*}(\sigma X)$.
Moreover, the functor $\rho$ sends $\tau$ in $\FilSp$ to a map that is normally called $\tau$ in $\Syn_E(\Sp)$.
In this way, all results in the $\tau$-formalism directly apply to $\Syn_E(\Sp)$ as well.
For instance, the total differentials and the Omnibus Theorem are available in this context too.

This does not mean that the use of filtered spectra does away with synthetic spectra.
Rather, synthetic spectra are a natural home for Adams spectral sequences, and the $\tau$-formalism therein is significantly more powerful.
The main aspects in which the structure of $\Syn_E(\Sp)$ is simpler than $\FilSp$ are the following.
\begin{itemize}
    \item There is a t-structure on $\Syn_E(\Sp)$, called the \emph{homological t-structure}, whose heart is equivalent to the abelian category of (graded) $E_*E$-comodules.
    In this t-structure, for \emph{every} spectrum $X$, the synthetic spectrum $\nu_E X$ is connective.
    \item The $\infty$-category of $C\tau$-modules in $\Syn_E(\Sp)$ is equivalent to (a version of) the derived $\infty$-category of (graded) $E_*E$-comodules.
\end{itemize}
These aspects are closely related, but not the same.
Together, they can be regarded as the reason that $\Syn_E(\Sp)$ is closely related to $E$-Adams spectral sequences: for instance, the Ext groups on the $\uE_2$-page page are precisely the mapping objects in the derived $\infty$-category of comodules.
More precisely, it is with these properties that we can compute $\sigma \circ \nu_E$ to be the Adams spectral sequence functor.

Neither of these features is present in $\FilSp$.
For instance, while $\nu_E X$ is connective in the homological t-structure, the filtered spectrum $\sigma(\nu_E X)$ is usually not connective in the standard t-structure on $\FilSp$.
This makes it easier to manipulate and study $\nu_E X$ in the synthetic context, again showing that $\Syn_E(\Sp)$ is the natural home for Adams spectral sequences.

As another example, we have a commutative diagram
\[
    \begin{tikzcd}
        \Syn_E(\Sp) \rar["\sigma"] \dar["C\tau\otimes \blank"'] & \FilSp \dar["C\tau\otimes\blank"] \\
        \Mod_{C\tau}(\Syn_E(\Sp)) \rar["\sigma"] & \Mod_{C\tau}(\FilSp).
    \end{tikzcd}
\]
As explained before, the vertical functor on the right can be identified with the associated graded functor.
The diagram now tells us that its composite with $\sigma$ factors through $C\tau$-modules in $\Syn_E(\Sp)$.
Because this is an $\infty$-category of an algebraic nature, we can much more effectively compute there, making the associated graded of the signature more accessible.

One can use this additional structure in another way, which has become known as the \emph{$C\tau$-method}\footnote{We warn that these terms should not be confused with what we call the \emph{$\tau$-formalism}. The latter is an overarching term, while $C\tau$-method is a specific technique that is particularly powerful in the synthetic $\tau$-formalism.} or \emph{$C\tau$-philosophy} of Gheorghe, Isaksen, Wang and Xu \cite{gheorghe_wang_xu_special_fibre,isaksen_wang_xu_dimension_90}.
It has been one of the landmark advances in computational stable homotopy theory.
It works as follows.
Fixing a spectrum $F$, one can set up the $\nu_E F$-Adams spectral sequence internal to $\Syn_E(\Sp)$ (where $F$ can be different from $E$).
We can push this spectral sequence along either one of the two functors
\[
    \Sp \from \Syn_E(\Sp) \to \Mod_{C\tau}(\Syn_E(\Sp)).
\]
Upon mapping it to spectra, we recover the ordinary $F$-Adams spectral sequence, while upon mapping it to $C\tau$-modules, we obtain a purely algebraic spectral sequence, where we may essentially compute differentials by hand (or by computer).
Because these two are now related to a spectral sequence in $\Syn_E(\Sp)$, we can thus deduce differentials in spectra from differentials in the algebraic realm.
Isaksen--Wang--Xu \cite{isaksen_wang_xu_dimension_90} work in the case $E=\MU$ and $F = \F_2$; see \cref{sec:history_tau_formalism} below for a further discussion.

A related application of this structure on synthetic $C\tau$-modules, is to the computation of Toda brackets.
Toda brackets in $\Syn_E(\Sp)$ map to Toda brackets under the functor $C\tau \otimes \blank$, whereupon they become Massey products, which can again be computed algebraically by hand.
This leads to a synthetic version of Moss's Theorem; see \cite[Section~3]{CDvN_part2}.

\section{History}
\label{sec:history_tau_formalism}

Our introduction of the $\tau$-formalism is quite different from the way it arose historically.
Whereas we presented it first of all as a notational device, its origins are much more complicated.

It first appeared when \emph{$\bC$-motivic spectra} began to be used to study the stable homotopy groups of spheres \cite{isaksen_stable_stems,gheorghe_wang_xu_special_fibre,isaksen_wang_xu_stable_stems,isaksen_wang_xu_dimension_90}.
In motivic spectra, one can run the motivic $\F_2$-Adams spectral sequence for the motivic sphere.
There is a functor $\Sp_\bC \to \Sp$ from $\bC$-motivic spectra to ordinary spectra called \emph{Betti realisation}, and this turns the motivic $\F_2$-ASS into the ordinary $\F_2$-ASS.
It turns out that the motivic version is remarkably similar to the ordinary one.
First off, the motivic version is trigraded, due to motivic homotopy groups being bigraded.
There is an endomorphism $\tau$ of the motivic sphere, and the motivic dual Steenrod algebra is almost the ordinary one tensored with $\F_2[\tau]$.
Isaksen \cite{isaksen_stable_stems} realised that this additional grading and the presence of $\tau$ introduce constraints, and used this to deduce motivic Adams differentials, which upon Betti realisation yield new differentials in the ordinary Adams spectral sequence.
One can get further information out of this approach by combining it with the movitic Adams--Novikov spectral sequence.

Later, Gheorge--Wang--Xu \cite{gheorghe_wang_xu_special_fibre} showed that modules over $C\tau$ in (cellular, $p$-complete) $\bC$-motivic spectra are equivalent to the derived $\infty$-category of $\BP_*\BP$-comodules.
Combined with the observation that differentials in the motivic Adams--Novikov spectral sequence correspond to differentials in the ordinary Adams--Novikov spectral sequence, this led to the \emph{$C\tau$-method} (explained above in synthetic terms), and with it a great advancement in our understanding of the sphere spectrum; see \cite{isaksen_wang_xu_dimension_90}.

There is a certain unreasonable effectiveness of these motivic methods, since motivic spectra are inherently algebro-geometric objects, and it is not a priori clear why this algebraic geometry is connected to the Adams--Novikov spectral sequence.
This was explained when different models were given for the subcategory of cellular motivic spectra (the subcategory in which these applications take place).
\begin{itemize}
    \item Pstr\k{a}gowski \cite{pstragowski_synthetic} defined $E$-synthetic spectra, and proved that if $E=\MU$, this gives a model for cellular motivic spectra (at least after $p$-completion).
    \item Gheorge--Isaksen--Krause--Ricka \cite{mmf} gave a model for ($p$-complete) cellular $\bC$-motivic spectra in terms of modules in filtered spectra.
    \item Burklund--Hahn--Senger \cite{burklund_hahn_senger_Rmotivic} directly compare these synthetic and filtered models.
\end{itemize}
This explains why $\bC$-motivic homotopy theory was so successfully applied: synthetic spectra are by nature designed to be a good category of Adams filtrations.
See also \cite[Remark~4.62]{pstragowski_synthetic} for a further discussion comparing synthetic spectra with \cite{isaksen_wang_xu_dimension_90}.

With the construction of synthetic spectra, it became possible to use the same type of techniques for other $E$ as well, not just in the case $E=\MU$.
Burklund--Hahn--Senger \cite[Theorem~9.19]{burklund_hahn_senger_manifolds} prove the \emph{Omnibus Theorem} to formalise the idea that $E$-synthetic analogues capture the $E$-Adams spectral sequence.
Later, Patchkoria--Pstr\k{a}gowski \cite{patchkoria_pstragowski_derived_inftycats} define a more general setting for defining synthetic categories, allowing for much more general $E$ and for a much more general stable $\infty$-category in the place of spectra.

If we only care about the applications to computations and spectral sequences rather than motivic homotopy theory as a whole, then synthetic spectra offer a more light-weight technical setup to do these computations with.
In these notes, we further argue that most synthetic techniques originate in the even more light-weight context of filtered spectra.
Not only does this make proofs and constructions more concrete, but it also allows for straightforward generalisations.
For instance, while the Omnibus Theorem of Burklund--Hahn--Senger only applies to synthetic analogues, the version we will deduce (using the adjunction $\rho \dashv \sigma$ above) from our filtered version applies to \emph{all} (convergent) synthetic spectra, and moreover applies in the same fashion to any (good enough) deformation.
We give a more detailed comparison of the proofs of these two versions of the Omnibus Theorem in \cref{rmk:comparison_omnibus_proofs}.

\section{Outline}

We begin by reviewing the theory of filtered abelian groups, filtered spectra and spectral sequences without the use of $\tau$.
This is the content of \cref{ch:filtered_spectra_sseqs}.
At this point, we do not yet introduce the $\tau$-formalism: rather, this chapter is aimed at setting up the basic concepts and terminology to be used later on, and to make various conventions and subtleties explicit.
In particular, we include a short discussion on the Adams spectral sequence in a non-synthetic sense, to provide all the necessary background for the later chapters.

Next, in \cref{ch:tau_formalism_filtered} we introduce the $\tau$-formalism in the filtered setting, starting with filtered abelian groups, and afterwards in filtered spectra.
Aside from discussing \emph{total differentials}, the main goal of this chapter is to prove the \emph{Omnibus Theorem} in the filtered setting.
The device for proving this is the \emph{$\tau$-Bockstein spectral sequence}, which we introduce in this chapter as well.
Finally, we end with a general discussion of \emph{deformations}, showing how to export the $\tau$-formalism to other $\infty$-categories.

With this preparation in hand, in \cref{ch:synthetic_spectra} we come to our other main topic, which is \emph{synthetic spectra}.
After reviewing the basic categorical properties, we show how the theory of deformations applies to synthetic spectra, and work this out in detail.
The main goal of this chapter is to compute the signature of a synthetic analogue.
After this, the synthetic Omnibus Theorem follows as a corollary.
Finally, in \cref{ch:variants_Syn} we discuss certain variants and properties of synthetic spectra, in particular the comparison between synthetic and motivic spectra.

\section{Conventions}

These notes are written mainly for a homotopical audience.
We have in mind a reader who has seen spectral sequences before, but is not necessarily intimately familiar with their construction.
In certain sections, familiarity with the Adams spectral sequence from a practical perspective will be useful, but is not strictly speaking required.

Throughout, we assume a working knowledge of $\infty$-categories in the sense of Joyal and Lurie; the standard references are \cite{HTT} and \cite{HA}.
We distinguish between categories and $\infty$-categories: by the term \emph{category}, we mean a $(1,1)$-category, while by the term \emph{$\infty$-category}, we mean an $(\infty,1)$-category.

When a morphism in an $\infty$-category admits a two-sided inverse up to homotopy, we refer to it as being an \emph{isomorphism}, rather than an \emph{equivalence}.
This should not cause much confusion, as we do not compare the $\infty$-categories we work in with a model category giving rise to them.
The only exception is that we speak of an \emph{equivalence} of $\infty$-categories rather than an \emph{isomorphism} of $\infty$-categories (which, we admit, is not a fully consistent choice of terminology).

We use the term \emph{space} for what is also referred to as an \emph{$\infty$-groupoid} or an \emph{anima}.
The $\infty$-category of spaces is denoted by $\Spaces$.
The $\infty$-category of spectra is denoted by $\Sp$, and we write $\otimes$ for the smash product of spectra.
By the term \emph{$\E_\infty$-ring}, we mean an $\E_\infty$ ring spectrum.

We do not distinguish notationally between an abelian group and its corresponding Eilenberg--MacLane spectrum.
For example, we do not write $\uH\Z$, but simply write~$\Z$.
The context will clarify whether we are dealing with a spectrum or an abelian group.

If $A$ is a graded abelian group and $n$ is an integer, then we write $A[n]$ for the graded abelian group given by
\[
    (A[n])_k = A_{k-n}.
\]
We will use the same formula for graded modules, comodules, etc.

We use Adams indexing for all of our spectral sequences, both in our formulas and in depicting spectral sequences.
This means a $d_r$-differential has bidegree $(-1,r)$.
Meanwhile, we use a homological-algebra indexing for Ext groups: for integers $s$ and $t$, we write
\[
    \Ext^{s,t}(M,N) = \Ext^s(M[t],\, N).
\]
In the case of Adams spectral sequences, this means we will often have the expression
\[
    \uE_2^{n,s} \cong \Ext_{E_*E}^{s,\, n+s}(E_*,\, E_*X).
\]
Our reasoning for this is that Adams indexing is most useful for working with spectral sequences, while the homological-algebra indexing on Ext is what one uses when computing these Ext groups.

Often, we refer to the $\uE_r$-page of a spectral sequence as its \emph{$r$-th page} or as \emph{page $r$}.

For ease of reference, we include here a list of the places where we make or clarify indexing conventions, sorted roughly by theme.
\begin{itemize}
    \item Homological vs.\ cohomological indexing, decreasing vs.\ increasing filtrations, towers vs.\ filtrations: \cref{rmk:decreasing_filtr_is_cohomological,rmk:filtrations_vs_tower,rmk:homological_vs_cohom_filab,rmk:motivation_Adams_indexing_good,constr:underlying_sseq,def:cycles_boundaries_and_infty_page,rmk:n_s_placement_sseq_notation,rmk:convenient_indexing_left_right_concentrated}.
    \item First- vs.\ second-page indexing: \cref{rmk:E2_indexing_FilSp,rmk:second_page_indexing_tau_BSS}.
    \item Filtered $\tau$-formalism: \cref{def:fil_tau,not:tau_Bocksteins,not:two_different_meanings_multiplication_by_tau,constr:trigraded_sseq_bifiltered_spectrum,not:tau_tilde,rmk:second_page_indexing_tau_BSS}.
    \item Synthetic indexing: \cref{def:synthetic_bigraded_spheres,rmk:motivic_grading_Syn,rmk:cohom_vs_hom_indexing_t_structure_Syn,rmk:why_second_page_indexing_Syn,var:tau_BSS_synthetic,thm:synthetic_omnibus,rmk:synthetic_vs_motivic_sphere_indexing}.
\end{itemize}

%% file: fil/basics_abridged.tex
As is well-known, filtrations give rise to spectral sequences.
The goal of this chapter is to review the theory of filtrations in the stable setting and the resulting spectral sequences, as well as dealing with more subtle issues like convergence.

We begin by studying filtered abelian groups in \cref{sec:fil_ab}.
Partially we do this as a warm-up, but mainly because it is the natural structure on the homotopy groups of a filtered spectrum.
After introducing and reviewing their basic category theory in \cref{sec:fil_sp}, we discuss their relation to spectral sequences in \cref{sec:spectral_sequences}.
For a more relaxed introduction to how spectral sequences arise from filtered spectra, we refer to \cref{ch:informal_sseq}.
Next, our goal is to describe our main example: the Adams spectral sequence.
This is the topic of \cref{sec:Adams_and_Tot_sseqs}, where we also include some background on spectral sequences arising from cosimplicial objects.
To prepare for this, we include a short digression on a duality between filtrations and towers of spectra, which we refer to as \emph{reflection}, in \cref{sec:reflecting}.

There is a myriad of sources on spectral sequences, which would be impossible to list here.
We learned much of this chapter from \cite{antieau_decalage}, \cite{boardman_convergence}, \cite{hedenlund_phd}, \cite{rognes_ASS_lecture_notes,rognes_sseq_lecture_notes}.

\section{Filtered abelian groups}
\label{sec:fil_ab}

We begin with an elementary algebraic concept.
We will use the adjective \emph{strict} (which we borrow from \cite{antieau_decalage}) to distinguish it from the later, more general concept of \cref{def:filtered_abelian_group}.

\begin{definition}
    \label{def:classical_filtered_abelian_group}
    Let $A$ be an abelian group.
    \begin{numberenum}
        \item A \defi{strict filtration} on $A$ is a sequence of subgroups
        \[
            \dotsb \subseteq F^1 \subseteq F^0 \subseteq F^{-1} \subseteq \dotsb \subseteq A.
        \]
    \end{numberenum}
    Let $\set{F^s}$ be a strict filtration on $A$.
    \begin{numberenum}[resume]
        \item If $a \in A$ is an element, then the \defi{filtration} of $a$ is the integer $s$ such that
        \[
            a \notin F^{s+1} \quad \text{but} \quad a \in F^s.
        \]
        We say that $a$ has filtration $\infty$ if it lies in all the $F^s$, and that it has filtration $-\infty$ if it lies in none of the $F^s$.
        \item The \defi{associated graded} of $\set{F^s}$ is the graded abelian group $\Gr F$ given by
        \[
            \Gr^s F = F^s/F^{s+1}.
        \]
    \end{numberenum}
\end{definition}

By definition, the subgroup $F^s$ is the subgroup of elements of filtration at least $s$.
It might therefore be helpful to think of $F^s$ as $F^{\geq s}$.
Note that $F^\infty$ is the limit $\lim_s F^s$, while $F^{-\infty}$ is the colimit $\colim_s F^s$.

We regard a strict filtration on $A$ is a tool to help us understand the group $A$.
One can think of it as starting with the elements of filtration $+\infty$ and moving down in filtration, where at each step the associated graded is measuring how many elements we `add'.
In the end, this procedure allows us to see all the elements that do not have filtration $-\infty$.
In practice, the associated graded is what one has the most control over.
As a result, we think of elements of filtration $\pm \infty$ as bad, and hope to find ourselves in situations where they do not exist.

An example of a result that formalises this idea is the following.
For a further discussion and other results in this direction, we refer to \cite[Section~2]{boardman_convergence}.

\begin{proposition}
    \label{prop:map_on_fil_ab_is_iso_under_conditions}
    Let $A$ and $B$ be abelian groups equipped with strict filtrations $\set{F^sA}$ and $\set{F^sB}$, respectively.
    Let $f\colon A \to B$ be a map that respects these filtrations.
    Suppose that
    \begin{numberenum}
        \item the map $f$ induces an isomorphism $F^\infty A \congto F^\infty B$;
        \item the first derived limit $\derlim{1}_s F^sA$ vanishes;
        \item the map $f$ induces an isomorphism on associated graded $\Gr^s A \congto \Gr^s B$ for all $s$;
        \item both $A$ and $B$ have no elements of filtration $-\infty$.
    \end{numberenum}
    Then $f$ is an isomorphism of abelian groups, and moreover restricts to an isomorphism $F^s A \congto F^s B$ for every $s$.
\end{proposition}
\begin{proof}
    See \cite[Theorem~2.6]{boardman_convergence}.
\end{proof}

\begin{remark}
    \label{rmk:decreasing_filtr_is_cohomological}
    In the above definition, we used a \emph{decreasing indexing} on the filtration.
    One should think of this as \emph{cohomological indexing} for filtrations.
    We follow this convention because most filtrations we consider (for example, the Adams filtration) are of the form
    \[
        \dotsb \subseteq F^2 \subseteq F^1 \subseteq F^0 = A.
    \]
\end{remark}

\begin{remark}
    There is an obvious variant of \cref{def:classical_filtered_abelian_group} for graded abelian groups.
    In this case, the associated graded is naturally a bigraded abelian group.
\end{remark}

\begin{remark}
    \label{rmk:Boardman_notation}
    In \cite[Section~2]{boardman_convergence}, the following terminology is introduced.
    \begin{itemize}
        \item If a filtration has no elements of filtration $-\infty$ (i.e., every element of $A$ appears in one of the $F^s$, or equivalently, if $\colim_s F^s = A$), then the filtration is said to be \emph{exhaustive}.
        \item If there are no elements of filtration $+\infty$ (i.e., if the limit $\lim_s F^s$ vanishes), then the filtration is said to be \emph{Hausdorff}.
        \item If the first-derived limit $\derlim{1}_s F^s$ vanishes, then the filtration is said to be \emph{complete}.
        (Note that a filtration can be complete without being Hausdorff.
        In other words, the limit of a ``Cauchy sequence'' need not be unique.)
    \end{itemize}
\end{remark}

\begin{warning}
In these notes, we will deviate from Boardman's terminology recalled in the previous remark: see \cref{def:derived_complete_filtered_ab} below.
\end{warning}

By definition, a strict filtration only grows as we move down in filtration.
It turns out to be useful to allow for a more general concept, one where we allow the groups to shrink as well.

\begin{definition}
    \label{def:filtered_abelian_group}
    \mbox{}
    \begin{numberenum}
        \item A \defi{filtered abelian group} is a functor $\Z^\op \to \Ab$, where we view $\Z$ as a poset under the usual ordering.
        We write
        \[
            \FilAb \defeq \Fun(\Z^\op,\Ab)
        \]
        for the (presentable, abelian) category of filtered abelian groups.
        \item If $A \colon \Z^\op \to \Ab$ is a filtered abelian group, then we write $A^\infty$ and $A^{-\infty}$ for its limit and colimit, respectively.
        \item The tensor product of abelian groups induces a presentably symmetric monoidal structure on $\FilAb$ via Day convolution, viewing $\Z^\op$ as a symmetric monoidal category under addition.
        A \defi{filtered commutative ring} is a commutative algebra object in $\FilAb$.
        \item If $A$ is a filtered abelian group, then its \defi{associated graded} is the graded abelian group $\Gr A$ given by
        \[
            \Gr^s A \defeq \coker(A^{s+1} \to A^s).
        \]
    \end{numberenum}
\end{definition}

In diagrams, a filtered abelian group $A$ consists of abelian groups $A^s$ for $s\in \Z$, together with maps
\[
    \dotsb \to A^1 \to A^0 \to A^{-1} \to \dotsb.
\]
We refer to these maps as \defi{transition maps}.

Let us explain our (perhaps slightly nonstandard) terminology.

\begin{remark}[Filtrations vs.\ towers]
    \label{rmk:filtrations_vs_tower}
    We deliberately use the term \emph{filtration} instead of \emph{tower} in the above.
    Throughout, we use the word \emph{filtration} to indicate that we think of the colimit as the underlying object, and the limit as an error term.
    When we use the word \emph{tower}, we instead regard the limit as the underlying object and the colimit as the error term.
    For an example of the difference, see \cref{ex:reflection_p_adic_filtr}.
\end{remark}

\begin{remark}[Homological vs.\ cohomological grading]
    \label{rmk:homological_vs_cohom_filab}
    Fitting with \cref{rmk:decreasing_filtr_is_cohomological}, we regard the usage of $\Z^\op$ to index filtered objects as \emph{cohomological} indexing of filtered objects.
    This is also why we use superscripts to indicate the index.
    If we instead think of these objects as towers in the sense of the previous remark, then the usage of $\Z^\op$ is a \emph{homological} indexing convention.
\end{remark}

Next, let us compare the notion of a filtered abelian group with that of a strict filtration as in \cref{def:classical_filtered_abelian_group}.
\begin{itemize}
    \item A strict filtration is a special case of a filtered abelian group, namely one whose transition maps are injective.
    The only difference is that the ambient abelian group from \cref{def:classical_filtered_abelian_group} is no longer present in \cref{def:filtered_abelian_group}.
    We will instead regard the colimit $A^{-\infty}$ as the ambient abelian group.
    Said differently, giving a strict filtration on an abelian group $B$ in the sense of \cref{def:classical_filtered_abelian_group} consists of providing a filtered abelian group $A$ in the sense of \cref{def:filtered_abelian_group} with injective transition maps, together with a map $A^{-\infty} \to B$.
    
    Going forward, we will usually use the term \emph{strict filtration} to refer to a filtered abelian group with injective transition maps.
    When we use the version of \cref{def:classical_filtered_abelian_group}, we will always check that there are no elements of filtration $-\infty$, to be consistent with the previous story.
    
    \item Conversely, a filtered abelian group $A$ gives rise to an \defi{induced strict filtration} $\set{F^s}$ on its colimit $A^{-\infty}$, via
    \begin{equation}
        F^s \defeq \im (A^s \to A^{-\infty}) \subseteq A^{-\infty}. \label{eq:induced_classical_filtration_filt_ab}
    \end{equation}
    This filtration has, essentially by definition, no elements of filtration $-\infty$.
    Note however that the assignment $A \mapsto \set{F^s}$ loses information: the transition maps in the filtered abelian group need not be injective.
\end{itemize}

We still need to deal with the potential presence of elements of filtration $+\infty$; we will use the following terminology.

\begin{definition}
    \label{def:derived_complete_filtered_ab}
    We say a filtered abelian group $A$ is \defi{derived complete} if
    \[
        \lim A = 0 \qquad \text{and} \qquad \derlim{1} A = 0.
    \]
\end{definition}

We will revisit this later in \cref{sec:tau_in_fil_ab}, where it is called \emph{$\tau$-completeness}.
For a discussion without the language of $\tau$, see \cite[Definition~2.7, Proposition~2.8]{boardman_convergence}.

\begin{remark}
    The reason we call this \emph{derived complete} is that it matches the notion of completeness for filtered spectra to be introduced in \cref{def:complete_filtered_spectrum} below.
    More specifically, by post-composing with the inclusion, a filtered abelian group $A$ determines a functor $\Z^\op \to \D(\Ab)$ to the derived $\infty$-category of abelian groups.
    Then $A$ is derived complete if and only if the (derived) limit of $\Z^\op \to \D(\Ab)$ is zero.
    Because the forgetful functor $\D(\Ab) \to \Sp$ preserves limits, this is equivalent to viewing $A$ as a functor $\Z^\op \to \Sp$ landing in discrete (a.k.a.\ Eilenberg--MacLane) spectra, and asking for the (homotopy) limit of this functor to vanish.
    See \cref{rmk:htpy_of_sequential_limit} for an elaboration on this point.
\end{remark}

\begin{remark}[Filtered tensor product]
Concretely, the tensor product of $A,B\in \FilAb$ is given levelwise by
\[
    (A \otimes B)^s = \colim_{i + j \geq s} A^i \otimes B^j,
\]
with the natural transition maps between them.
A filtered commutative ring is a filtered abelian group $A$ together with pairings
\[
    A^s \otimes A^t \to A^{s+t}
\]
for every $s,t\in\Z^\op$, being compatible with the transition maps and satisfying the obvious commutative ring diagrams.
The unit for this monoidal structure is
\[
    \begin{tikzcd}
        \dotsb \ar[r] & 0 \ar[r] & \Z \ar[r,equals] & \Z \ar[r,equals] & \dotsb,
    \end{tikzcd}
\]
with the first $\Z$ appearing in filtration $0$.
\end{remark}

We leave it to the reader to verify that the associated graded assembles to a symmetric monoidal functor
\[
    \Gr \colon \FilAb \to \grAb.
\]

\begin{remark}[Filtration is subadditive]
\label{rmk:subadditive_filtration}
Suppose that we have a filtered ring structure on a strict filtration $\set{F^s}$.
Then this structure is the same as a commutative ring structure on the ambient abelian group such that for all $s$ and $t$, we have
\[
    F^s \cdot F^t \subseteq F^{s+t}.
\]
Note that this means that the filtration might ``jump'': a product in $F^{s+t}$ might land in the subgroup $F^N$ for $N > s+t$.
In other words, filtration is \emph{subadditive} under multiplication.
Products that jump in filtration become zero in the associated graded.
\end{remark}

We introduced a strict filtration as a tool to better understand its ambient abelian group.
It is suggestive then that the only purpose of a filtered abelian group as in \cref{def:filtered_abelian_group} is to give rise to its induced strict filtration via \eqref{eq:induced_classical_filtration_filt_ab}.
Said differently, one might consider the kernels of the maps $A^s \to F^s$ to be an anomaly, because they determine the zero element in $A^{-\infty}$.
This is decidedly \emph{not} the perspective we will take: the entire filtration is the object of interest.
Many of the benefits from the synthetic perspective come from remembering the filtration as a whole.

However, there is one downside to working with filtered abelian groups: the associated graded cannot measure the kernels of the transition maps.
In order to take these into account as well, we have to move to the derived setting.

\section{Filtered spectra}
\label{sec:fil_sp}

Instead of moving to the derived setting by using derived abelian groups, we immediately go to the more universal case of spectra.

\begin{definition}
    \leavevmode
    \begin{numberenum}
        \item A \defi{filtered spectrum} is a functor $\Z^\op \to \Sp$, where we view $\Z$ as a poset under the usual ordering.
        We write
        \[
            \FilSp \defeq \Fun(\Z^\op,\, \Sp)
        \]
        for the (presentable, stable) $\infty$-category of filtered spectra.
        \item If $X \colon \Z^\op \to \Sp$ is a filtered spectrum, then we write $X^\infty$ and $X^{-\infty}$ for its limit and colimit, respectively.
        \item The smash product of spectra induces a presentably symmetric monoidal structure on $\FilSp$ via Day convolution.
        A \defi{filtered $\E_\infty$-ring} is an $\E_\infty$-algebra object in $\FilSp$.
        \item A \defi{graded spectrum} is a functor $\Z^\discr \to \Sp$, where $\Z^\discr$ is the discrete category with objects $\Z$.
        We write
        \[
            \grSp \defeq \prod_\Z \Sp = \Fun(\Z^\discr,\, \Sp).
        \]
        for the (presentable, stable) $\infty$-category of graded spectra.
        \item If $X$ is a filtered spectrum, then its \defi{associated graded} is the graded spectrum $\Gr X$ given by
        \[
            \Gr^s X \defeq \cofib(X^{s+1}\to X^s).
        \]
    \end{numberenum}
\end{definition}

We often depict a filtered spectrum $X \colon \Z^\op \to \Sp$ as a diagram
\[
    \dotsb \to X^1 \to X^0 \to X^{-1} \to \dotsb.
\]
The abuse of this notation is not great: a filtered spectrum is, up to contractible choice, determined by the spectra $\set{X^n}$ together with their transition maps; see \cite[Proposition~3.3, Corollary~3.4]{ariotta_cochain_complexes}.

\begin{remark}
    The role of spectra in the above definition is not special.
    Most of the results in this section apply to a suitable stable $\infty$-category in the place of spectra.
    For concreteness, and to prevent this chapter from becoming needlessly long, we stick to the case of spectra.
    For a discussion in this greater generality, we refer to \cite{antieau_decalage}, or parts of \cite[Part~II]{hedenlund_phd}.
\end{remark}

Like with filtered abelian groups, if $X$ is a filtered spectrum, then we think of $X^{-\infty}$ as the `underlying spectrum' of $X$.

Using that cofibres are functorial, the associated graded assembles into a functor
\[
    \Gr \colon \FilSp \to \grSp.
\]
We will see later that this can be upgraded to a symmetric monoidal functor; see \cref{rmk:Ctau_gives_Gr_symm_mon_str}.
For an alternative proof, see \cite[Section~II.1.3]{hedenlund_phd}.

The functor $\pi_* \colon \Sp \to \grAb$ induces a functor
\[
    \FilSp \to \Fil(\grAb).
\]
We can find corepresenting objects for the individual abelian groups of this functor.
By Yoneda, each transition map will then be corepresented by a map between these objects, but we defer a discussion of the resulting maps to the next chapter.
Although we could define these corepresenting objects by hand, it will be convenient to introduce the following functor.

\begin{definition}
    \label{def:functor_Z_to_FilSp}
    We write $i \colon \Z \to \FilSp$ for the functor given by
    \[
        s \mapsto \opSigma^\infty_+ \Hom_\Z(\blank,\, s).
    \]
\end{definition}

By the properties of Day convolution, the functor $i$ is naturally symmetric monoidal.

\begin{remark}
    The functor $i$ gives the $\infty$-category $\FilSp$ a universal property, which says that colimit-preserving functors out of $\FilSp$ correspond to functors out of~$\Z$.
    We discuss this, as well as the structure that such a functor puts on the target category, in detail later in \cref{sec:deformations}.
\end{remark}

\begin{definition}[Filtered bigraded spheres]
    \label{def:filtered_bigraded_spheres}
    Let $n$ and $s$ be integers.
    \begin{numberenum}
        \item The \defi{filtered $(n,s)$-sphere} is
        \[
            \S^{n,s} \defeq \opSigma^n i(s).
        \]
        We refer to $n$ as the \defi{stem}, and to $s$ as the \defi{filtration}.
        \item We write $\Sigma^{n,s} \colon \FilSp \to \FilSp$ for the functor given by tensoring with $\S^{n,s}$ on the left.
        \item We write $\pi_{n,s}\colon \FilSp \to \Ab$ for the functor
        \[
            \pi_{n,s}(\blank) \defeq [\S^{n,s},\, \blank].
        \]
    \end{numberenum}
\end{definition}

Unwinding definitions, we see that $\S^{n,s}$ is the filtered spectrum given by
\[
    \begin{tikzcd}
        \dotsb \ar[r] & 0 \ar[r] & \S^n \ar[r,equals] & \S^n \ar[r,equals] & \dotsb,
    \end{tikzcd}
\]
where the first $\S^n$ appears in position $s$.
In diagrams therefore, $\Sigma^{n,s}$ is given by applying $\Sigma^n$ levelwise, and by shifting the filtered spectrum $s$ units to the left.
Likewise, we see that for any filtered spectrum $X$, we have a natural isomorphism
\[
    \oppi_{n,s}X \cong \pi_n(X^s).
\]

\begin{remark}
    This definition of the filtered bigraded spheres is designed to be compatible with \emph{first-page indexing} of filtered spectra, i.e., the indexing that makes the underlying spectral sequence start on the first page.
    It is possible to change this to a second-page indexing (or any page); see \cref{rmk:E2_indexing_FilSp} for a further discussion.
\end{remark}

The filtered spheres are, in fact, generators for $\FilSp$.

\begin{proposition}
    \label{prop:properties_of_FilSp}
    \leavevmode
    \begin{numberenum}
        \item For every $s\in\Z$, the filtered sphere $\S^{0,s}$ is a compact and invertible object in $\FilSp$, with inverse $\S^{0,-s}$.
        In particular, the monoidal unit of $\FilSp$ is compact.
        \item The bigraded homotopy groups $\pi_{*,*}$ detect isomorphisms of filtered spectra.\label{item:fil_htpy_detects_isos}
        \item \label{item:FilSp_generated_by_spheres} As a stable $\infty$-category, $\FilSp$ is generated under colimits by the spheres $\S^{0,s}$ for $s\in \Z$.
        That is, the objects $\S^{n,s}$ for $n,s\in\Z$ generate $\FilSp$ under colimits.
        In particular, $\FilSp$ is compactly generated by dualisables.
        \item The monoidal $\infty$-category $\FilSp$ is \emph{rigid} in the sense that an object is compact if and only if it is dualisable.\label{item:FilSp_ridigly_compactly_gen}
    \end{numberenum}
\end{proposition}
\begin{proof}
    The first property follows because $\Map(\S^{0,s},X)\cong \Omega^\infty X^s$, and the second is evident.
    \Cref{item:FilSp_generated_by_spheres} follows from \cref{item:fil_htpy_detects_isos} using \cite[Corollary~2.5]{yanovski_monadic_tower}.
    To prove \cref{item:FilSp_ridigly_compactly_gen}, first notice that the unit $\S^{0,0}$ is compact, so that all dualisable objects are compact.
    As it is also generated by compact dualisable objects, it follows that every compact object is dualisable too; see, e.g., (the footnote to) \cite[Terminology~4.8]{naumann_pol_ramzi_monoidal_fracture}.
\end{proof}

By default, we will equip $\FilSp$ with the following t-structure.
We borrow the name from \cite{barkan_monoidal_algebraicity}, though there appears to be no agreed-upon name for it.

\begin{definition}
    The \defi{diagonal t-structure} on $\FilSp$ is the t-structure where a filtered spectrum $X$ is connective if and only if
    \[
        \oppi_{n,s}X = 0 \qquad \text{whenever } n<s.
    \]
    We write $\tau_{\geq n}^\diag$ and $\tau_{\leq n}^\diag$ for the $n$-connective cover and $n$-truncation functors with respect to this t-structure, respectively.
\end{definition}

As this specified class of connective objects is presentable and closed under colimits, this determines a unique accessible t-structure on $\FilSp$ by \cite[Proposition~1.4.4.11]{HA}.
It enjoys the following properties.

\begin{proposition}
    \label{prop:diagonal_t_structure_filsp}
    \leavevmode
    \begin{letterenum}
        \item \label{prop:connective_fil_sp} A filtered spectrum $X$ is connective if and only if
        \[
            \oppi_{n,s}X = 0 \qquad \text{whenever } n<s.
        \]
        \item \label{prop:truncated_fil_sp} A filtered spectrum $X$ is $0$-truncated if and only if
        \[
            \oppi_{n,s}X = 0 \qquad \text{whenever } n>s.
        \]
        \item The connective cover $\tau_{\geq0}^\diag X \to X$ induces an isomorphism
        \[
            \pi_{n,s}(\tau_{\geq 0}^\diag X) \congto \pi_{n,s}(X) \qquad \text{whenever }n\geq s.
        \]
        Likewise, the $0$-truncation $X \to \tau_{\leq 0}^\diag X$ induces an isomorphism
        \[
            \pi_{n,s}(X) \congto \pi_{n,s}(\tau_{\leq 0}^\diag X) \qquad \text{whenever }n \leq s.
        \]
        \item The functor $\pi_{*,*}$ induces a symmetric monoidal equivalence
        \[
            \FilSp^\heart \simeqto \grAb, \quad X \mapsto (\oppi_{n,n}X)_n.
        \]
        \item The diagonal t-structure is complete.
        \item The diagonal t-structure is compatible with filtered colimits.
        \item \label{item:diag_t_str_is_monoidal} The diagonal t-structure is compatible with the monoidal structure.
    \end{letterenum}
\end{proposition}

\begin{proof}
    See, e.g., \cite[Propositions~II.1.22--II.1.24]{hedenlund_phd}.
\end{proof}

This t-structure is a convenient device for giving functorial definitions of the Whitehead filtration and Postnikov tower.

\begin{definition}
    \label{def:whitehead_postnikov}
    \mbox{}
    \begin{numberenum}
        \item The \defi{constant filtration} is the functor $\Const \colon \Sp \to \FilSp$ given by precomposition with $\Z^\op \to \Delta^0$.
        \item The \defi{Whitehead filtration} is the functor $\Wh \colon \Sp \to \FilSp$ given by $\tau_{\geq 0}^\diag \circ \Const$.
        \item The \defi{Postnikov tower} is the functor $\Post \colon \Sp \to \FilSp$ given by $\tau_{\leq 0}^\diag \circ \Const$.
    \end{numberenum}
\end{definition}
From this definition, the functors $\Wh$ and $\Post$ come with natural transformations
\[
    \Wh \to \Const \qquad \text{and} \qquad \Const \to \Post.
\]

\begin{remark}
\label{rmk:whitehead_lax_mon}
    Note that $\Const$ is naturally a (strong) symmetric monoidal functor.
    Since the diagonal t-structure is monoidal, it follows that $\Wh\colon \Sp \to \FilSp$ is naturally a lax symmetric monoidal functor.
\end{remark}

We now check that these are the expected filtrations.

\begin{lemma}
    Let $X$ be a spectrum.
    Then $\Wh X$ is of the form
    \[
        \dotsb \to \tau_{\geq 1} X \to \tau_{\geq 0} X \to \tau_{\geq -1} X \to \dotsb
    \]
    with the natural maps as transition maps, and $\Post X$ is of the form
    \[
        \dotsb \to \tau_{\leq 1} X \to \tau_{\leq 0} X \to \tau_{\leq -1} X \to \dotsb
    \]
    with the natural maps as transition maps.
\end{lemma}
\begin{proof}
    If $F$ is a filtered spectrum that is $t$-truncated in the diagonal t-structure for some $t$, then it follows that $\pi_n F^{-\infty}=0$ for all $n$, so that $F^{-\infty}=0$.
    Consequently, using the cofibre sequence
    \[
        \Wh X \to \Const X \to \tau_{\leq -1}^\diag (\Const X)
    \]
    and that taking colimits is exact, we find that the natural map $\Wh X \to \Const X$ induces an isomorphism on colimits.
    In other words, the colimit of $\Wh X$ is $X$.
    Because $\Wh X$ is connective in the diagonal t-structure, the spectrum $(\Wh X)^s$ is $s$-connective.
    The cofibre of $(\Wh X)^s \to X$ is an $(s-1)$-truncated spectrum, because the cofibre of $\Wh X \to \Const X$ is $(-1)$-truncated in the diagonal t-structure.
    This means that the map $(\Wh X)^s \to X$ exhibits the source as the connective cover of the target.
    This applies to all $s$, implying that for every $s$, the connecting map $(\Wh X)^{s+1} \to (\Wh X)^{s}$ is given by the $(s+1)$-connective cover.
    
    The same reasoning applies to $\Post X$, using the limit instead of the colimit.
\end{proof}

The functor $\Const$ restricts to an equivalence from $\Sp$ to the full subcategory of $\FilSp$ on the \emph{constant} filtered spectra, i.e., those for which all transition maps are isomorphisms.
We will revisit this later under the guise of \emph{$\tau$-invertible filtered spectra} in \cref{ssec:filtered_inverting_tau}.

This definition of the Whitehead tower functor is also helpful to prove that $\Wh$ is fully faithful, and to characterise when a filtered spectrum is isomorphic to the Whitehead filtration on a spectrum.
We follow the proof given by \cite[Theorem~A.1]{annala_pstragowski_weight_filtration}.

\begin{theorem}
    The functor $\Wh \colon \Sp \to \FilSp$ is fully faithful.
    Its essential image consists of those filtered spectra $F \colon \Z^\op \to \Sp$ such that for every $s\in\Z^\op$, the natural map $F^s \to F^{-\infty}$ exhibits the source as the $s$-connective cover of the target.
\end{theorem}
\begin{proof}
    Let $X$ and $Y$ be spectra.
    The natural transformation $\Wh \to \Const$ induces the second map in the composite
    \[
        \Map_\Sp(X,Y) \to \Map_\FilSp(\Wh X,\, \Wh Y) \to \Map_\FilSp(\Const X,\, \Const Y).
    \]
    This composite is an isomorphism because $\Const$ is fully faithful.
    It therefore suffices to show that the second map in this composite is an isomorphism.
    We use that the natural transformation $\Wh \to \Const$ has two characterisations: first, it exhibits the source as the connective cover of the target; and second, it exhibits the target as the colimit (i.e., underlying spectrum) of the source.
    The first is by definition of the functor $\Wh$, and the second follows from completeness of the t-structure on spectra.
    The map under consideration therefore factors as
    \[
        \begin{tikzcd}
            \Map_\FilSp(\Wh X,\, \Wh Y) \rar["\cong"] \drar & \Map_\FilSp(\Wh X,\, \Const Y) \dar["\cong"] \\
             & \Map_\FilSp(\Const X,\, \Const Y),
        \end{tikzcd}
    \]
    where the horizontal map is an isomorphism by the universal property of the connective cover, and the second since $\Const$ is right adjoint to taking the colimit.

    To prove the essential image is of the claimed form, suppose $F$ is a filtered spectrum of the described form.
    Write $X$ for the spectrum $F^{-\infty}$.
    We claim that $F \to \Const X$ exhibits $F$ as the connective cover (in the diagonal t-structure) of the target.
    Since $F$ is connective in the diagonal t-structure, the natural map $F \to \Const X$ lifts to $F \to \Wh X$.
    Now both $F^s \to X$ and $(\Wh X)^s \to X$ exhibit the source as the $s$-connective cover of the target, so the map between them is an isomorphism.
\end{proof}

\begin{definition}
    \label{def:complete_filtered_spectrum}
    Let $X$ be a filtered spectrum.
    We say that $X$ is \defi{complete} if the limit $X^\infty$ of $X$ vanishes.
    We write $\Filhat\Sp$ for the full subcategory of $\FilSp$ on the complete filtered spectra.
\end{definition}

The inclusion $\Filhat\Sp \subseteq \FilSp$ admits a left adjoint, called the \emph{completion functor}
\[
    \FilSp \to \Filhat\Sp, \quad X \mapsto \widehat{X} \defeq \cofib(X^\infty \shortto X)
\]
We will explore this concept in more detail later under the name of \emph{$\tau$-completion} of filtered spectra in \cref{ssec:filtered_tau_completion}.
In particular, we will see that a map of filtered spectra is an isomorphism after completion if and only if it is an isomorphism on associated graded; see \cref{prop:tau_completion_filsp_concretely}.
For an alternative discussion without the language of~$\tau$, see \cite[Section~II.1.2]{hedenlund_phd}.

Finally, we turn to the relationship between filtered spectra and filtered abelian groups.

\begin{definition}
    \label{def:filtration_on_pistar}
    Let $X$ be a filtered spectrum.
    The \defi{induced strict filtration} on the abelian group $\pi_n X^{-\infty}$ is the (strict) filtration given by
    \[
        F^s\, \pi_n X^{-\infty} \defeq \im(\pi_n X^s \to \pi_n X^{-\infty}).
    \]
\end{definition}

Note that, because $\pi_* \colon \Sp \to \grAb$ preserves filtered colimits, this is the same as the strict filtration on $\pi_n X^{-\infty}$ induced by the filtered abelian group $\pi_n \circ X$.
In particular, the induced strict filtration on $\pi_n X^{-\infty}$ has no elements of filtration $-\infty$.

In practice, it is not easy to compute the homotopy groups $\pi_n X^s$ directly, so we should not compute this filtration from the definition.
What is usually much more accessible is the associated graded of the filtered spectrum, but this carries considerably less information.
One might try and invest the associated graded with as much structure as possible, so that it starts to remember the homotopy of the filtered spectrum itself.
This is precisely what a \emph{spectral sequence} does.

\section{Spectral sequences}
\label{sec:spectral_sequences}

So far, we have not yet delived on our promise that the derived setting is able to measure both the kernel and cokernel of transition maps.
The notion of a spectral sequence makes this precise.
Its purpose is to reconstruct the homotopy groups $\pi_*X^{-\infty}$ from the associated graded of $X$.
The non-injectivity of the transition maps on homotopy groups leads to the \emph{differentials} in a spectral sequence.
In these notes, we will further argue that instead of only reconstructing $\pi_*X^{-\infty}$ from the associated graded, the better approach is to reconstruct \emph{all} of the homotopy groups in the filtered spectrum, i.e., the bigraded homotopy groups $\oppi_{*,*}X$.

We do not include an in-depth review or motivation of the setup of a spectral sequence here.
We include a detailed but informal introduction to this in \cref{ch:informal_sseq}.

\begin{remark}
    \label{rmk:motivation_Adams_indexing_good}
    The indexing for exact couples used below is not the most common.
    Usually, Serre indexing is used in formulas and definitions, while Adams indexing is used for displaying the spectral sequence.
    We prefer to work with one indexing system throughout, and we prefer Adams indexing as it is the most practical one and arises naturally from the diagram chase of \cref{ch:informal_sseq}.
    Moreover, this indexing is more straightforward to generalise to contexts where homotopy groups have a more complicated indexing (such as filtered or synthetic spectra later in these notes).
\end{remark}

For a reference on exact couples, see, e.g., \cite[page~37 and following]{mccleary_spectral_sequences}. 

\begin{construction}
\label{constr:underlying_sseq}
Let $X$ be a filtered spectrum.
\begin{numberenum}
    \item The \defi{associated exact couple} of $X$ is the exact couple of bigraded abelian groups defined by
    \[
        \uA^{n,s}(X) \defeq \pi_{n}X^s \qquad \text{and} \qquad \uE^{n,s}(X) \defeq \pi_{n}\Gr^s X,
    \]
    with the natural maps from the long exact sequence between them, fitting into the following diagram, where each map is annotated by its $(n,s)$-bidegree.
    \[
        \begin{tikzcd}[column sep=tiny]
            \pi_n X^{s} \ar[rr,"(0{,}\,-1)"] & & \pi_n X^s \ar[dl,"(0{,}\,0)"] \\
            & \pi_n \Gr^s X \ar[ul,"(-1{,}\,1)"] & 
        \end{tikzcd}
    \]
    \item The \defi{underlying spectral sequence} of $X$, denoted by $\set{\uE_r^{n,s}(X),\ d_r}_{r\geq 1}$, is the spectral sequence in bigraded abelian groups arising from the exact couple associated with~$X$.
    By definition, this spectral sequence is of the form
    \[
        \uE_1^{n,s}(X) \defeq \pi_n \Gr^s X \implies \pi_n X^{-\infty}
    \]
    and the differential $d_r$ has bidegree $(-1,r)$ for $r\geq 1$.
    We refer to $\pi_* X^{-\infty}$ as the \defi{abutment} of the spectral sequence.
\end{numberenum}
\end{construction}

The term \emph{abutment} and the above notation are not meant to indicate any form of convergence; rather, one should think of this as the object to which the spectral sequence is trying to converge.

\begin{example}
    \label{ex:sseq_filtered_sphere}
    The associated graded of $\S^{0,0}$ is given by $\S$ in degree $0$, and zero elsewhere.
    In particular, the first page of the resulting spectral sequence is
    \[
        \uE_1^{n,s} \cong \begin{cases}
            \ \pi_n \S &\text{if }s=0,\\
            \ 0 &\text{if }s\neq 0.
        \end{cases}
    \]
    The spectral sequence of a general filtered bigraded sphere is a shift of the one for $\S^{0,0}$.
    These spectral sequences are rather uninteresting: they do not decompose $\pi_* \S$ in a new, meaningful way.
    As a result, we mostly think of the filtered bigraded spheres as useful formal objects, and not out of interest in their underlying spectral sequence.
\end{example}

To avoid confusion, let us make some indexing conventions explicit.

\begin{definition}
    \label{def:cycles_boundaries_and_infty_page}
    Let $X$ be a filtered spectrum.
    Let $r\geq 0$, and let $n,s\in\Z$.
    \begin{numberenum}
        \item Write $\uZ_r^{n,s} \subseteq \uE_1^{n,s}$ for the subgroup of the \defi{$r$-cycles}, i.e., those $x$ such that the differentials $d_1(x),\dotsc,d_r(x)$ vanish.
        This leads to a sequence of inclusions
        \[
            \dotsb \subseteq \uZ_2^{n,s} \subseteq \uZ_1^{n,s} \subseteq \uZ_0^{n,s} = \uE_{1}^{n,s}.
        \]
        We define
        \[
            \uZ_\infty^{n,s} \defeq \lim_r \uZ_r^{n,s} = \bigcap_r \uZ_r^{n,s}.
        \]
        An element of $\uZ_\infty^{n,s}$ is called a \defi{permanent cycle}.
        \item Write $\uB_r^{n,s} \subseteq \uE_1^{n,s}$ for the subgroup of the \defi{$r$-boundaries}, i.e., the image of the first~$r$ differentials.
        This leads to a sequence of inclusions
        \[
            0 = \uB_0^{n,s} \subseteq \uB_1^{n,s} \subseteq \uB_2^{n,s}\subseteq \dotsb \subseteq \uE_1^{n,s}.
        \]
        We define
        \[
            \uB_\infty^{n,s} \defeq \colim_r \uB_r^{n,s} = \bigcup_r \uB_r^{n,s}.
        \]
        Note that there is an inclusion $\uB_\infty^{n,s} \subseteq \uZ_\infty^{n,s}$.
        \item If $r \geq 1$, we define
        \[
            \uE_r^{n,s} \defeq \uZ_{r-1}^{n,s}/\uB_{r-1}^{n,s}.
        \]
        \item \label{item:infty_page}
        We define the \defi{$\infty$-term} and the \defi{derived $\infty$-term}, respectively, to be
        \begin{align*}
            \uE_\infty^{n,s} &\defeq \uZ_\infty^{n,s}/\uB_\infty^{n,s},\\
            \uR\uE_\infty^{n,s} &\defeq \derlim{1}_r\uZ_r^{n,s}.
        \end{align*}
    \end{numberenum}
\end{definition}

\begin{warning}
    In \cite{boardman_convergence}, Boardman writes $Z_{r}$ for what we would write as $Z_{r-1}$, and likewise writes $B_r$ for what we would write as $B_{r-1}$.
\end{warning}

\begin{warning}
    In general, the group $\uE_\infty^{n,s}$ is not a (co)limit of the groups $\uE_r^{n,s}$.
    In fact, as $\uE_{r+1}^{n,s}$ is a subquotient of $\uE_r^{n,s}$, there is in general no sensible map from $\uE_r^{n,s}$ to or from $\uE_{r+1}^{n,s}$.
    In certain cases, such a map does exist.
    For instance, suppose that $\uE_1^{*,s}=0$ for $s\ll 0$.
    For fixed $s$, then for $r\gg0$, we have that $\uE_{r+1}^{*,s}$ is a subgroup of $\uE_r^{*,s}$, and it follows that $\uE_\infty^{*,s}$ is the limit along the resulting sequence of inclusion maps.
\end{warning}

\begin{definition}
    \label{def:perm_cycle_and_detection}
    Let $X$ be a filtered spectrum, let $x \in \uE_1^{n,s}$ be an element, and let $\theta \in \pi_n X$.
    \begin{numberenum}
        \item Let $r\geq 1$.
        Suppose that $x$ is an $r$-cycle, so that it defines an element in $\uZ_r^{n,s}$.
        We say that $x$ \defi{survives to page $r$} if its image in $\uE_r^{n,s}$ is nonzero.
        If $x$ is a permanent cycle, then we say that $x$ \defi{survives to page $\infty$} if its image in $\uE_\infty^{n,s}$ is nonzero.
        \item Suppose that $x$ is a permanent cycle that survives to page $\infty$.
        We then say that $x$ \defi{detects} $\theta$ (or that $\theta$ is \emph{detected by} $x$) if there exists a lift $\alpha \in \pi_n X^s$ of $x$ that maps to $\theta$ under $X^s \to X^{-\infty}$.
    \end{numberenum}
\end{definition}

Beware that, in the definition of detection, the lift $\alpha$ above need not be unique if it exists.


\begin{remark}
    The underlying spectral sequence from \cref{constr:underlying_sseq} is functorial in the filtered spectrum.
    As exact couples form a $1$-category, it suffices to lift this construction to a functor from $\h \FilSp$ to exact couples.
    This is easily checked.
    Post-composing this with the functor from exact couples to spectral sequences, we obtain the desired functor
    \[
        \FilSp \to \SSeq(\grAb), \quad X \mapsto \set{\uE_r^{n,s}(X),\ d_r}.
    \]
\end{remark}

Roughly speaking, this functor forgets the homotopy of the filtered spectrum and only remembers the homotopy of the associated graded, together with the induced differentials.
The goal of the spectral sequence, then, is to reconstruct the homotopy of the filtered spectrum from this data.
Whether this is even possible (aside from extension problems) is the question of \emph{convergence}, which we discuss in \cref{ssec:convergence}.

We think of the $\infty$-category $\FilSp$ as an $\infty$-categorical enhancement of the category $\SSeq(\grAb)$.
Although not every (stable) spectral sequence arises from a filtered spectrum, in practice, they do.\footnotemark\ 
The additional homotopical structure on filtered spectra has many advantages.
For example, it allows us to talk about coherently multiplicative filtrations; see \cite[Part~II]{hedenlund_phd}.

\footnotetext{An example of a spectral sequence that does not arise from a filtration is the $p$-Bockstein spectral sequence in the way that it is set up in \cite[Section~1]{browder_torsion_Hspaces} or \cite[Chapter~10]{mccleary_spectral_sequences}.
However, this is not an essential issue: there is an alternative way to set up the $p$-Bockstein spectral sequence that does come from a filtered spectrum (namely, from the $p$-adic filtration on a spectrum or chain complex; see \cref{ex:reflection_p_adic_filtr} and \cref{ex:recover_p_BSS}).
Arguably, this latter version even has a nicer abutment, converging to $\pi_*(X_p^\wedge)$, rather than to $\F_p \otimes F$, where $F$ denotes the torsion-free quotient of $\pi_* X$.}

For completeness, we compare the above construction of the underlying spectral sequence with some other ones appearing in the literature.

\begin{remark}[Alternative constructions]
    \label{rmk:comparison_alternative_constrs_sseq}
    We used exact couples in our definition of the underlying spectral sequence, as this fits most closely with the explanation given in \cref{ch:informal_sseq}.
    There are a number of alternative approaches in the literature, including at least the following.
    \begin{itemize}
        \item Cartan--Eilenberg systems \cite[Section~XV.7]{cartan_eilenberg_homological_algebra}, which are also used by Lurie in \cite[Section~1.2.2]{HA}.
        \item Coherent cochain complexes in spectra, introduced by Ariotta \cite{ariotta_cochain_complexes}.
        \item The d\'ecalage functor on the level of filtered spectra.
        D\'ecalage was introduced by Deligne for chain complexes \cite{deligne_hodge_theory_II}, and later generalised to cosimplicial objects by Levine \cite{levine_slice_ANSS}, and generalised to filtered spectra by Antieau, as first written down by Hedenlund \cite[Part~II]{hedenlund_phd} and later by Antieau \cite{antieau_decalage}.
    \end{itemize}
    All of these approaches agree, as a consequence of the following results.
    \begin{itemize}
        \item Ariotta \cite[Theorem~3.19]{ariotta_cochain_complexes} constructs an equivalence of $\infty$-categories from $\Filhat\Sp$ to coherent cochain complexes in spectra.
        \item Let $X$ be a filtered spectrum.
        This gives rise to a Cartan--Eilenberg system via \cite[Definition~1.2.2.9]{HA}.
        Iterating the d\'ecalage functor on filtered spectra also results in a spectral sequence.
        Antieau \cite[Theorem~4.13]{antieau_decalage} showed that the spectral sequence arising from this Cartan--Eilenberg system is isomorphic to the one arising from iterating the d\'ecalage functor.
        More generally, he showed this when working with filtered objects of a stable $\infty$-category with sequential limits and colimits that is equipped with a t-structure.
        \item Let $X$ be a filtered spectrum.
        This gives rise to both a Cartan--Eilenberg system as before, as well as to an exact couple via \cref{constr:underlying_sseq}.
        The two resulting spectral sequences are isomorphic; see the discussion in \cite[page~336]{cartan_eilenberg_homological_algebra}, or, e.g., Creemers \cite[Theorem~5.1]{julie_creemers_masters} for a detailed argument.
    \end{itemize}
\end{remark}


\begin{remark}
    We use the variant of a spectral sequence that computes homotopy classes out of a compact object.
    Alternatively, as in \cite[Section~1.2.2]{HA}, one can also start with a stable $\infty$-category equipped with a t-structure, and set up a spectral sequence to compute the heart-valued homotopy groups.
    One should assume that heart-valued homotopy groups preserve sequential colimits for this to work in a reasonable generality.
\end{remark}

Finally, we make a few remarks regarding indexing and notation.

\begin{remark}[First- vs.\ second-page indexing]
    \label{rmk:E2_indexing_FilSp}
    We choose to index spectral sequences arising from filtered spectra to start on the first page; let us call this \emph{first-page indexing}.
    For various reasons (such as aesthetics, or to better fit alternative definitions of a spectral sequence), it can be useful to reindex this to start on the second (or any later) page.
    One can achieve \emph{second-page indexing} via the reindexing
    \[
        \widetilde{\uE}_{r+1}^{n,s} \defeq \uE_r^{n,\,n+s}.
    \]
    It is straightforward to check that this turns $d_r$-differentials into $\widetilde{d}_{r+1}$-differentials; in particular, this makes the spectral sequence start on the second page.
    Using this second-page indexing, the filtered bigraded spheres take the form (using the functor~$i$ from \cref{def:functor_Z_to_FilSp})
    \[
        \widetilde{\S}^{n,s} \defeq \opSigma^n i(n+s).
    \]
    Concretely, this is the filtered spectrum that is $\S^n$ in positions $n+s$ and below (with identities between them), and zero elsewhere.
    In this grading, categorical suspension takes the form $\Sigma^{1,-1}$.
\end{remark}

\begin{remark}
    \label{rmk:n_s_placement_sseq_notation}
    In accordance with our previous indexing conventions, $n$ here is indexed homologically, whereas $s$ is indexed cohomologically.
    As such, it would be more honest to write $\uE_n^s$ and $\uA_n^s$, but we do not do so, as the current notation is well established (and would leave little room for the page-index $r$).
    Depending on the context, it might be more natural to alter either of these conventions; see, e.g., \cref{rmk:convenient_indexing_left_right_concentrated}.
\end{remark}

\subsection{Convergence}
\label{ssec:convergence}

Convergence of a spectral sequence is the question whether one can reconstruct $\pi_*X^{-\infty}$ from the spectral sequence.\footnotemark\ 
More precisely, convergence concerns reconstructing $\pi_*X^{-\infty}$ by reconstructing the induced strict filtration on $\pi_*X^{-\infty}$.
We can however only hope to reconstruct the associated graded of this strict filtration, and other methods are necessary for solving the extension problems.

\footnotetext{Other times, the convergence issue is to prove that the colimit of the filtration is isomorphic to a desired spectrum. This is, of course, a question that is more specific to the situation at hand, so in this section we regard the colimit as the desired object to study.}

In most accounts of spectral sequences, convergence is additional structure on the spectral sequence.
For a spectral sequence arising from a filtered spectrum, all of this structure is supplied by the filtered spectrum, so that convergence becomes a property.

We closely follow Boardman's account \cite{boardman_convergence}.
He works with (unrolled) exact couples, but we specialise everything to the setting of filtered spectra.
We include a few detailed remarks, both for the curious reader and for use later in the more technical parts of these notes.
An alternative introduction is given by Hedenlund in \cite[Section~1.2.2]{hedenlund_phd}.

Recall from \cref{def:filtration_on_pistar} that a filtered spectrum $X$ gives rise to an \emph{induced strict filtration} on $\pi_nX^{-\infty}$, denoted by $F^s\,\pi_n X^{-\infty}$.

\begin{construction}
    \label{constr:comparison_gr_infty_page}
    Let $X$ be a filtered spectrum.
    Write $\partial_{n,s} \colon \pi_n \Gr^s X \to \pi_{n-1} X^{s+1}$ for the boundary map.
    By a diagram chase (see, e.g., \cite[Lemma~5.6]{boardman_convergence} or \cite[Lemma~2.5.10]{rognes_sseq_lecture_notes}), there is a natural isomorphism of graded abelian groups
    \[
        \frac{F^s \, \pi_nX^{-\infty}}{F^{s+1} \, \pi_nX^{-\infty}} \cong \frac{\ker \partial_{n,s}}{\uB_\infty^{n,s}}.
    \]
    Using that $\ker \partial_{n,s} / \uB_\infty^{n,s}$ naturally injects into $\uE_\infty^{n,s}$, this leads to a natural injective map
    \begin{equation}
        \Gr^s \pi_n X^{-\infty} = \frac{F^s \, \pi_nX^{-\infty}}{F^{s+1} \, \pi_nX^{-\infty}} \hookto \uE_\infty^{n,s}. \label{eq:map_assoc_graded_to_infty_page}
    \end{equation}
\end{construction}

\begin{definition}
    \label{def:strong_convergence}
    Let $X$ be a filtered spectrum.
    We say that the underlying spectral sequence \defi{converges strongly} to $\pi_*X^{-\infty}$ if
    \begin{letterenum}
        \item \label{item:strong_conv_der_complete} the induced strict filtration $\set{F^s \, \pi_* X^{-\infty}}$ is \emph{derived complete} in the sense of \cref{def:derived_complete_filtered_ab}, i.e.,
        \[
            \lim_s F^s \, \pi_* X^{-\infty} = 0 \qquad \text{and} \qquad \derlim{1}_s F^s \, \pi_* X^{-\infty} = 0;
        \]
        \item \label{item:strong_conv_Einfty_iso} the natural map \eqref{eq:map_assoc_graded_to_infty_page} is an isomorphism for all $n$ and $s$.
    \end{letterenum}
\end{definition}

\begin{warning}
    This terminology is abusive: the above definition of strong convergence is not a condition on the spectral sequence, but rather on the filtered spectrum.
    In fact, the above conditions do not make sense if we do not specify which filtered spectrum gives rise to the spectral sequence.
\end{warning}

In words, condition~\ref{item:strong_conv_der_complete} says that we can (up to extension problems) reconstruct $\pi_*X^{-\infty}$ from the induced strict filtration, and condition~\ref{item:strong_conv_Einfty_iso} says that the spectral sequence is able to recover the associated graded of this filtration.

It is more accurate to speak of convergence to the induced strict filtration on $\pi_*X^{-\infty}$, but we will usually not do this.
If the filtered spectrum is clear from the context, we may also be brief and simply say that the spectral sequence \emph{converges strongly}, which should always be understood as convergence to $\pi_*X^{-\infty}$.

\begin{remark}
    \label{rmk:can_detect_nonzeroness}
    Because the map \eqref{eq:map_assoc_graded_to_infty_page} is always injective, we learn the following, even in the absence of any of the above convergence criteria.
    If $x \in \uE_1^{n,s}$ is a permanent cycle that survives to page $\infty$, then for any lift $\alpha \in \pi_n X^s$, the image of $\alpha$ in $\pi_n X^{-\infty}$ is nonzero.
    Indeed, by injectivity of \eqref{eq:map_assoc_graded_to_infty_page}, the image of $\alpha$ defines a nonzero map in $\Gr^s \pi_n X^{-\infty}$.
    In other words, any element that $x$ detects is nonzero.
    (Without convergence hypotheses however, a lift of $x$ to $\pi_n X^s$ may not exist.)
\end{remark}

\begin{remark}
    \label{rmk:conv_implies_good_detection}
    \Cref{item:strong_conv_der_complete} in particular implies that every nonzero element of $\pi_*X^{-\infty}$ is detected by a permanent cycle that survives to page $\infty$.
    Meanwhile, \cref{item:strong_conv_Einfty_iso} implies that every permanent cycle that survives to page $\infty$ detects an element in $\pi_*X^{-\infty}$ (which is nonzero by \cref{rmk:can_detect_nonzeroness}).
\end{remark}

\begin{remark}
    Cartan--Eilenberg \cite[page~321]{cartan_eilenberg_homological_algebra} also gives names to other notions of convergence.
    In the setting of filtered spectra, these notions are the following.
    The spectral sequence is said to \emph{converge weakly} if \ref{item:strong_conv_Einfty_iso} in \cref{def:strong_convergence} holds, and said to \emph{converge} if \ref{item:strong_conv_Einfty_iso} holds and $\lim_s F^s\,\pi_*X^{-\infty}=0$.
    We will not use this terminology.
\end{remark}

As the name suggests, strong convergence is the strongest type of convergence one can hope for.
Our goal then is to find conditions that guarantee strong convergence.
In practice, we only have limited knowledge about the homotopy groups $\pi_{n,s}X = \pi_n X^s$, so we would prefer convergence criteria that involve mostly the spectral sequence rather than the filtered spectrum itself.

Boardman's notion of \emph{conditional convergence} does exactly this.
It splits the problem up into two parts.
First, one needs to establish conditional convergence, which is a structural (and often mild) condition on the filtered spectrum.
Second, once this is established, there are conditions phrased entirely in terms of the spectral sequence that guarantee strong convergence.
These conditions are more computational in nature, and need to be checked on a case-by-case basis, but are often met.
The second step is the reason for using the word `conditional': the spectral sequence converges strongly, conditional on these requirements being met.

\begin{definition}
    Let $X$ be a filtered spectrum.
    We say that the underlying spectral sequence \defi{converges conditionally} to $\pi_*X^{-\infty}$ if the filtered spectrum $X$ is complete, i.e., if the limit $X^\infty$ vanishes.
\end{definition}

\begin{warning}
    As with strong convergence, this terminology is abusive: conditional convergence is a condition on the filtration rather than the underlying spectral sequence.
\end{warning}

\begin{remark}
    \label{rmk:htpy_of_sequential_limit}
    Using the Milnor short exact sequence (see, e.g., \cite[Theorem~4.9]{boardman_convergence})
    \[
        0 \to \derlim{1}_s \pi_{n+1} X^s \to \pi_n  X^\infty \to \lim_s \pi_n X^s \to 0,
    \]
    we see that $X^\infty$ vanishes if and only if
    \[
        \lim_s \pi_* X^s = 0 \qquad \text{and} \qquad \derlim{1}_s \pi_* X^s = 0.
    \]
\end{remark}

\begin{warning}
    \label{warn:explanation_convergence_criteria}
    The previous remark may appear to suggest that the vanishing of $X^\infty$ implies that the induced strict filtration on $\pi_*X^{-\infty}$ is derived complete.
    This is \emph{not} true in general, and this is exactly what leads to Boardman's convergence criteria.
    To explain why, we introduce the following notation:
    \[
        F^s\oppi_{n,w} X = \im(\oppi_{n,\,w+s}X \to \oppi_{n,w}X).
    \]
    There are two problems.
    The first is that for any $w$, the natural map
    \[
        \pi_n X^\infty \to F^\infty \oppi_{n,w}X \defeq \lim_s F^s \oppi_{n,w}X
    \]
    need not be surjective.
    This does happen if $\uR\uE_\infty$ vanishes; see \cite[Lemma~5.9]{boardman_convergence}.
    The second problem is that the natural map
    \[
        \colim_w F^{\infty}\oppi_{n,w}X = \colim_w \lim_s F^s \oppi_{n,w}X \to \lim_s \colim_w F^s \pi_{n,w} X = F^{\infty} \,\pi_n X^{-\infty}
    \]
    also need not be surjective.
    (There is an analogous version of this map with the first derived limit in the place of the limit, but this is always surjective; see \cite[Lemma~8.11]{boardman_convergence}.)
    It is surjective in certain cases, such as when the filtered spectrum is left or right concentrated in the sense of \cref{def:left_right_concentrated} below.
    In general, Boardman's \emph{whole-plane obstruction} is the obstruction to this implication; see \cref{rmk:whole_plane_obstruction} below for a further discussion, and \cite[Lemma~8.11]{boardman_convergence} for the precise result alluded to here.
    Clearly, if all three maps above are surjective, then the vanishing of $X^\infty$ does imply derived completeness of the induced strict filtration on $\pi_*X^{-\infty}$.
\end{warning}

The `conditional' part of conditional convergence becomes easier if the filtered spectrum is of a special form.
The following terminology is nonstandard.
For the general case, see \cref{rmk:whole_plane_obstruction}.

\begin{definition}
    \label{def:left_right_concentrated}
    Let $X$ be a filtered spectrum.
    \begin{numberenum}
        \item We say that $X$ is \defi{right concentrated} if the transition maps $X^{s+1} \to X^s$ are isomorphisms for $s\gg 0$.
        \item We say that $X$ is \defi{left concentrated} if the transition maps $X^{s+1} \to X^{s}$ are isomorphisms for $s\ll 0$:
    \end{numberenum}
\end{definition}

These conditions can be checked entirely in terms of the associated graded: being right concentrated means that $\uE_1^{*,s}$ vanishes for $s\gg 0$, and being left concentrated means that it vanishes for $s \ll 0$.

\begin{remark}
    \label{rmk:convenient_indexing_left_right_concentrated}
    In practice, we usually reindex a right-concentrated filtered spectrum to be of the form
    \[
        \dotsb \congto X^{1} \to X^0 \to X^{-1} \to \dotsb,
    \]
    in which case $X^\infty = X^1$ and $\uE_1^{*,s}=0$ for $s>0$.
    Boardman calls the resulting spectral sequence a \emph{half-plane spectral sequence with exiting differentials}. 
    Likewise, we index a left-concentrated filtered spectrum to be of the form
    \[
        \dotsb \to X^{1} \to X^0 \congto X^{-1} \congto \dotsb,
    \]
    in which case $X^{-\infty}=X^0$ and $\uE_1^{*,s}=0$ for $s<0$.
    Boardman calls the resulting spectral sequence a \emph{half-plane spectral sequence with entering differentials}. 
\end{remark}

Roughly speaking, for right-concentrated filtered spectra, the `conditional' part is vacuous, while for left-concentrated filtered spectra, the only thing one has to check is the vanishing of a derived limit term.
This can be checked in terms of the spectral sequence, without requiring any further knowledge of the filtration.

\begin{theorem}[Conditional convergence, Boardman]
    \label{thm:conditional_convergence}
    Let $X$ be a filtered spectrum.
    Suppose that $X$ is right concentrated.
    Then the following are equivalent.
    \begin{enumerate}[label={\upshape(1\alph*)}]
        \item The spectral sequence underlying $X$ converges conditionally.
        \item The spectral sequence underlying $X$ converges strongly.
    \end{enumerate}
    Suppose instead that $X$ is left concentrated.
    Then any two of the following imply the third.
    \begin{enumerate}[label={\upshape(2\alph*)}]
        \item The spectral sequence underlying $X$ converges conditionally.
        \item The derived $\infty$-term $\uR\uE_\infty^{n,s}$ from \cref{def:cycles_boundaries_and_infty_page}\,\ref{item:infty_page} vanishes for all~$n$ and~$s$.
        \item The spectral sequence underlying $X$ converges strongly.
    \end{enumerate}
\end{theorem}
\begin{proof}
    This is \cite[Theorem~6.1 and Theorem~7.3]{boardman_convergence}, respectively.
    The translation between his and our notation is the following (using notation from \cref{warn:explanation_convergence_criteria}):
    \begin{align*}
        A^s &= \oppi_{*,s} X & Q^s &= F^\infty \oppi_{*,s}X \\
        A^{\infty}&=\lim_s \oppi_{*,s} X & RQ^s &= \derlim{1}_t F^t \oppi_{*,s}X\\
        RA^{\infty}&=\derlim{1}_s\pi_{*,s}X. & & \qedhere
    \end{align*}
\end{proof}

As explained in \cite[Section~7]{boardman_convergence}, it is often easy to verify that the derived $\infty$-term vanishes.
This happens if the spectral sequence collapses at a finite page (i.e., $d_r=0$ for $r>r_0$), or if for every $n$ and $s$, only finitely many differentials leaving bidegree $(n,s)$ are nonzero (but where this bound is allowed to depend on $n$ and $s$).
In general, as with any first-derived limit of abelian groups, one can use the Mittag-Leffler condition to check its vanishing.

In certain cases, the following variant will be useful as well.
It is slightly more general than working with right-concentrated filtered spectra that are complete: we only ask that every filtered abelian group $\pi_n X$ becomes zero for $s\gg 0$, but not necessarily that they become zero at the same point.

\begin{proposition}
    \label{prop:eventually_zero_implies_strong_convergence}
    Let $X$ be a filtered spectrum.
    Suppose that for every $n$, the group $\pi_n X^s$ vanishes for $s \gg 0$ \textbr{where the bound on $s$ is allowed to depend on $n$}.
    Then $X$ is complete and the spectral sequence underlying $X$ converges strongly.
\end{proposition}
\begin{proof}
    From \cref{constr:comparison_gr_infty_page}, it is clear that the map \eqref{eq:map_assoc_graded_to_infty_page} is an isomorphism if $\uZ_\infty^{n,s} = \ker \partial_{n,s}$ for all $n,s$.
    The inclusion $\ker \partial_{n,s} \subseteq \uZ_\infty^{n,s}$ always holds.
    The reverse inclusion now follows from a diagram chase involving the definition of the differential: see, e.g., \cite[Lemma~2.5.8]{rognes_sseq_lecture_notes}.
    Finally, the completeness of $X$ follows from \cref{rmk:htpy_of_sequential_limit}.
\end{proof}

\begin{remark}[Whole-plane spectral sequences]
    \label{rmk:whole_plane_obstruction}
    In the case where the filtered spectrum does not become constant in either direction, the situation becomes more difficult.
    Boardman \cite[Section~8, Equation~(8.7)]{boardman_convergence} defines a group $W$ he calls the \emph{whole-plane obstruction}.
    For any filtered spectrum $X$, Theorem~8.10 of op.\ cit.\ implies that any two of the following imply the third.
    \begin{letterenum}
        \item The spectral sequence underlying $X$ converges conditionally.
        \item The derived $\infty$-term $\uR\uE_\infty^{n,s}$ vanishes for all~$n$ and~$s$, and $W$ vanishes.
        \item The spectral sequence underlying $X$ converges strongly.
    \end{letterenum}
    He gives a criterion saying that if there is no infinite family of differentials that all cross each other in their interior, then $W=0$; see \cite[Lemma~8.1]{boardman_convergence}.
    For an alternative description of the whole-plane obstruction using Cartan--Eilenberg systems (see \cref{rmk:comparison_alternative_constrs_sseq} for more on Cartan--Eilenberg systems), see \cite{helle_rognes_whole_plane_obstruction}.
    They moreover give an alternative proof for Boardman's criterion: see \cite[Proposition~5.3]{helle_rognes_whole_plane_obstruction}.
\end{remark}

\section{Digression: reflecting}
\label{sec:reflecting}

So far, we have been thinking of a filtered spectrum (and consequently its underlying spectral sequence) as a tool to understand its colimit.
It is also possible to orient things the other way around, using a filtered spectrum to understand its limit instead.
In that case, the colimit is the object that should vanish to guarantee convergence properties.
The convergence discussion of spectral sequences does become slightly more involved in this context, because the functor $\pi_* \colon \Sp \to \grAb$ does not preserve sequential limits.
Nevertheless, for every convergence result we discussed above, Boardman \cite{boardman_convergence} also gives the limit-oriented version.

Rather than working with this limit-oriented version, we will use a duality to move back to the colimit-oriented version.
We refer to this as \emph{reflection}.
This is well known (being used, for instance, in \cite[Section~5]{bousfield_localization_spectra}), but we thought it would be helpful to make this translation and its basic properties explicit, especially its interaction with completion of filtered spectra.
It will not play a big role in the rest of these notes.
Let us also point out that this duality is a feature specific to the stable setting.

\begin{remark}
    Reflection can also be set up in the setting of exact couples.
    In that setting, when passing from a tower to a filtration, there will be a \emph{Rees system} in the sense of \cite[Section~7]{eckmann_hilton_exact_couples} witnessing that the two resulting spectral sequences are isomorphic.
    We do not detail how to specialise this argument to the setting of filtered spectra, focussing instead on the behaviour of the duality at a more categorical level.
\end{remark}

We remind the reader of the terminology we introduced in \cref{rmk:filtrations_vs_tower}: we will refer to objects of $\FilSp$ as \emph{filtrations} when working in the colimit-oriented context, and as \emph{towers} when working to the limit-oriented context.
Accordingly, if a notion of passing back and forth between towers and filtrations is to make sense, it should interchange the limit and colimit.
A natural candidate then is to reflect a filtration in its colimit, and to reflect a tower in its limit.

\begin{definition}
    \label{def:reflection_functors}
    The \defi{associated tower functor}, respectively the \defi{associated filtration functor}, are the functors defined by
    \begin{align*}
        (\blank)^\tow \colon \FilSp \to \FilSp, \quad &X \mapsto \cofib(\Sigma^{0,-1}X \shortto \Const X^{-\infty}),\\
        (\blank)^\fil \colon \FilSp \to \FilSp, \quad &X \mapsto \fib(\Const X^\infty \shortto \Sigma^{0,1}X).
    \end{align*}
\end{definition}

Concretely, if $X$ is a filtered spectrum, then
\[
    (X^\tow)^s = \cofib(X^{s+1}\to X^{-\infty}) \qquad \text{and} \qquad (X^\fil)^s = \fib(X^\infty \to X^{s-1}).
\]

We include a shift in these definitions because when working with the spectral sequence associated to a tower, one usually lets the first page consist of the fibres of the transition maps, not the cofibres.

\begin{definition}
    Let $X$ be a filtered spectrum.
    The \defi{fibre-associated graded} of $X$ is the graded spectrum $\fibGr X$ given by
    \[
        \fibGr^s X \defeq \fib(X^{s} \to X^{s-1}).
    \]
\end{definition}

In a picture, a tower $X$ together with its fibre-associated graded looks as follows:
\[
    \begin{tikzcd}
        & \fibGr^1 X \dar & \fibGr^0 X \dar & \fibGr^{-1} \dar & \\
        \dotsb \rar & X^1 \rar & X^0 \rar & X^{-1} \rar & \dotsb.
    \end{tikzcd}
\]

The following is a straightforward diagram chase.

\begin{proposition}
Let $X$ be a filtered spectrum.
Then there are natural isomorphisms of graded spectra
\[
    \fibGr (X^\tow) \cong \Gr X \qquad \text{and} \qquad \Gr (X^\fil) \cong \fibGr X.
\]
\end{proposition}

The reflection functors from \cref{def:reflection_functors} destroy some of the information contained in the original object: the associated tower only depends on the completion of the original object.
An analogous statement is true for the associated filtration, for which we introduce the following terminology.

\begin{definition}
    \label{def:cocomplete_filtered_spectrum}
    We call a filtered spectrum \defi{cocomplete} if its colimit vanishes.
    We write $\Filcheck\Sp$ for the full subcategory of $\FilSp$ on the cocomplete filtered spectra.
\end{definition}

Analogously to the case of completion, one can check that the inclusion $\Filcheck\Sp\subseteq\FilSp$ admits a right adjoint given by
\[
    \FilSp \to \Filcheck\Sp, \quad X \mapsto \widecheck{X} \defeq \fib(X \shortto X^{-\infty}).
\]
We call this functor \defi{cocompletion}, which then features in a colocalisation
\[
    \begin{tikzcd}[column sep=3.5em]
        \Filcheck\Sp \rar[shift left, hook] & \FilSp. \lar[shift left, "\widecheck{(\blank)}"]
    \end{tikzcd}
\]

We think of $\widehat{\Fil}\Sp$ as the conditionally convergent filtrations, and of $\Filcheck\Sp$ as the conditionally convergent towers.
The reflection functors of \cref{def:reflection_functors} translate between these in the following way.

\begin{proposition}
    \label{prop:properties_reflection}
    \leavevmode
    \begin{numberenum}
        \item \label{item:reflection_and_completion} We have commutative diagrams
        \[
            \begin{tikzcd}[column sep=3.5em]
                \FilSp \rar["(\blank)^\tow"] \dar["\widehat{(\blank)}"'] & \FilSp \\
                \widehat{\Fil}\Sp \rar["(\blank)^\tow"'] & \Filcheck\Sp \uar[hook]
            \end{tikzcd}
            \qquad \text{and} \qquad
            \begin{tikzcd}[column sep=3.5em]
                \FilSp \rar["(\blank)^\fil"] \dar["\widecheck{(\blank)}"'] & \FilSp \\
                \Filcheck\Sp \rar["(\blank)^\fil"'] & \widehat{\Fil}\Sp. \uar[hook]
            \end{tikzcd}
        \]
        \item \label{item:reflection_and_limits} If $X$ is a filtered spectrum, then we have natural isomorphisms
        \[
            (X^\tow)^{\infty} \cong (\widehat{X})^{-\infty} \qquad \text{and} \qquad (X^\fil)^{-\infty} \cong (\widecheck{X})^\infty.
        \]
        In particular, if $X$ is complete, then we have a natural isomorphism
        \[
            (X^\tow)^\infty \cong X^{-\infty},
        \]
        while if $X$ is cocomplete, then we have a natural isomorphism
        \[
            (X^\fil)^{-\infty} \cong X^\infty.
        \]
        \item \label{item:reflection_equivalence} The reflection functors restrict to inverse equivalences
        \[
            \begin{tikzcd}
                \Filhat\Sp \rar[shift left, "(\blank)^\tow"] &[2em]  \Filcheck\Sp. \lar[shift left, "(\blank)^\fil"]
            \end{tikzcd}
        \]
    \end{numberenum}
\end{proposition}

\begin{proof}
    For \cref{item:reflection_and_completion}, we only show the first diagram, as the argument for the second is dual to that for the first.
    Using that cofibres preserve colimits, it follows immediately from the definition that $X^\tow$ is cocomplete for all $X$.
    It remains to be verify that the map $X \to \widehat{X}$ becomes an isomorphism after taking associated towers.
    Write $T$ for $X^\tow$ and $R$ for $(\widehat{X})^\tow$.
    Then for every $s$, we have a commutative diagram where all rows and columns are cofibre sequences
    \[
        \begin{tikzcd}
            X^\infty \rar \dar[equals] & X^{s+1} \rar \dar & (\widehat{X})^{s+1} \dar\\
            X^\infty \rar \dar & X^{-\infty} \rar \dar & (\widehat{X})^{-\infty} \dar \\
            0 \rar & T^s \rar & R^s.
        \end{tikzcd}
    \]
    It follows that $T^s \to R^s$ is an isomorphism for all $s$, proving the claim.

    For \cref{item:reflection_and_limits}, we again only show the statement about associated towers.
    For a general filtered spectrum $X$, we have
    \[
        (X^\tow)^\infty = \lim_s \cofib(X^{s+1} \to X^{-\infty}) \cong \cofib(X^\infty \to X^{-\infty}).
    \]
    By \cref{item:reflection_and_completion}, it suffices to consider the case where $X$ is complete.
    In this case, the latter term is naturally isomorphic to $X^{-\infty}$, proving the claim.

    Finally, for \cref{item:reflection_equivalence}, we check that the composite $((\blank)^\tow)^\fil$ is isomorphic to the identity on complete objects; the argument for the other composite on cocomplete is analogous.
    Let $X$ be complete, and write $T=X^\tow$, and $F=T^\fil$.
    Then using that $X^\infty = 0$, we have a commutative diagram where all rows and columns are cofibre sequences
    \[
        \begin{tikzcd}
            \Omega X^s \rar \dar & 0 \rar \dar & X^s \dar \\
            0 \rar \dar & X^{-\infty} \rar[equals] \dar & X^{-\infty} \dar \\
            F^s \rar & T^\infty \rar & T^{s-1}.
        \end{tikzcd}
    \]
    This supplies a natural cofibre sequence $\Sigma^{-1,0} X \to 0 \to F$, proving the claim.
\end{proof}

\begin{example}
    \label{ex:reflection_Wh_Post}
    Let $X$ be a spectrum.
    It is a straightforward exercise to see that the reflection functors switch the Whitehead filtration and Postnikov tower of $X$: we have natural isomorphisms
    \[
        (\Wh X)^\tow \cong \Post X \qquad \text{and} \qquad (\Post X)^\fil \cong \Wh X.
    \]
    Note that this uses that the standard t-structure on spectra is complete.
    If we were to work in $\Fil(\C)$ for a stable $\infty$-category $\C$ equipped with a t-structure, then in order to interchange $\Wh X$ and $\Post X$, we would need to reflect them in $X$, rather than in their (co)limit.
\end{example}

\begin{example}
    \label{ex:reflection_p_adic_filtr}
    Let $X$ be a spectrum.
    The $p$-Bockstein filtration of $X$ is the filtered spectrum
    \[
        \begin{tikzcd}
            \dotsb \rar["p"] & X \rar["p"] & X \rar["p"] & X \rar[equals] & \dotsb,
        \end{tikzcd}
    \]
    indexed to be constant from filtration $0$ onwards.
    Observe that this filtration is complete if and only if $X$ is $p$-complete.
    Its associated tower is
    \[
        \begin{tikzcd}
            \dotsb \rar & X/p^3 \rar & X/p^2 \rar & X/p \rar & 0 \rar & \dotsb,
        \end{tikzcd}
    \]
    where $X/p^n$ appears in filtration $n-1$.
    The limit of this tower is $X_p^\wedge$.
    If we take the associated filtration of this tower, then we obtain the $p$-adic filtration on $X_p^\wedge$:
    \[
        \begin{tikzcd}
            \dotsb \rar["p"] & X_p^\wedge \rar["p"] & X_p^\wedge \rar["p"] & X_p^\wedge \rar[equals] & \dotsb.
        \end{tikzcd}\qedhere
    \]
\end{example}

\section{The Adams spectral sequence}
\label{sec:Adams_and_Tot_sseqs}

We give a brief introduction to the Adams spectral sequence.
While this is not intended as a first introduction, we take some care to explain some of the subtleties in its construction.
Our model for the Adams spectral sequence will be based on a cosimplicial object.
While this has some small downsides, it leads to an easier expression for the resulting filtered spectrum, because the cosimplicial d\'ecalage has an easy expression (\cref{def:cosimplicial_decalage}).
This limitation is not essential by any means, and can be avoided by working with filtered resolutions (\cref{rmk:cosimplicial_Adams_vs_filtered_resolution}).
Later, we will use synthetic spectra to set up this improved version; see \cref{sec:signature_synthetic_analogue}.

\subsection{The Tot spectral sequence}
\label{ssec:Tot_sseq}

We write $\Delta$ for the simplex category, and write $\Delta_{\leq n}$ for the full subcategory on the objects $[k]$ where $k\leq n$.

\begin{definition}
    \label{def:Tot_and_Tot_tower}
    Let $\C$ be a pointed $\infty$-category, and let $X^\bullet \colon \Delta \to \C$ be a cosimplicial object of $\C$.
    \begin{numberenum}
    \item The \defi{totalisation} of $X^\bullet$ is the limit
    \[
        \Tot X^\bullet = \lim_\Delta X^\bullet.
    \]
    \item For $n\geq 0$, the \defi{$n$-th partial totalisation} is the limit
    \[
        \Tot_n X^\bullet = \lim_{\Delta_{\leq n}} X^\bullet.
    \]
    \item Suppose that $\C$ admits (partial) totalisations.
    Using the filtration
    \[
        \dotsb = 0\subseteq \Delta_{\leq 0} \subseteq \Delta_{\leq 1} \subseteq \Delta_{\leq 2} \subseteq \dotsb \subseteq \Delta,
    \]
    we obtain a tower
    \[
        \Tot X^\bullet \to \dotsb \to \Tot_2 X^\bullet \to \Tot_1 X^\bullet \to \Tot_0 X^\bullet \to 0 \to \dotsb
    \]
    with limit $\Tot X^\bullet$.
    We call this tower the \defi{totalisation tower} (or \emph{Tot tower}) of $X^\bullet$.
    This is natural in $X$, resulting in a functor $\towTot\colon \C^\Delta \to \Fil(\C)$.
    \end{numberenum}
\end{definition}

We call this a \emph{tower} in accordance with \cref{rmk:filtrations_vs_tower}; below in \cref{def:Tot_filtration}, we will turn this tower into a filtration.
Note that the Tot tower is in particular a cocomplete object of $\Fil(\C)$ in the sense of \cref{def:cocomplete_filtered_spectrum}.

\begin{remark}[Cubes]
    The $n$-th partial totalisation can be computed as follows.
    We write
    \[
        \Delta_{/[n]}^\inj \subseteq \Delta_{/[n]}
    \]
    for the full subcategory of the slice on those objects given by injective maps to $[n]$.
    Note that this subcategory is a punctured $(n+1)$-cube, and in particular is a finite category.
    We have a forgetful functor
    \[
        \Delta_{/[n]}^\inj \to \Delta_{\leq n}
    \]
    and this functor is homotopy initial; see, e.g., \cite[Lemma~1.2.4.17]{HA}.
    In words, the $n$-th partial totalisation can be computed as the total fibre of a punctured $(n+1)$-cube.
\end{remark}

\begin{notation}
    Let $\Delta_\inj \subseteq \Delta$ denote the wide subcategory on those maps that are injective.
    We will write $j$ for the inclusion functor.
\end{notation}

Recall that a \defi{semicosimplicial object} of $\C$ is a functor $\Delta_\inj \to \C$.
The definition of the Tot tower can be mimicked for a semicosimplicial object, instead taking the limit over $\Delta_{\inj}$ or over $(\Delta_{\inj})_{\leq n}$.
There is a subtle difference between these two.

\begin{remark}
    \label{rmk:Tot_only_depends_semicosimplial}
    The inclusion $j\colon \Delta_\inj \subseteq \Delta$ is homotopy initial; see, e.g., \cite[Example~21.2]{dugger_htpy_colimits_primer}.
    As a result, we will also write $\Tot$ for the limit-functor
    \[
        \Sp^{\Delta_\inj} \to \Sp,
    \]
    and we therefore have a natural factorisation
    \[
        \begin{tikzcd}[row sep=small]
            \Sp^\Delta \ar[rr,"\Tot"] \drar["j^*"'] & & \Sp. \\
            & \Sp^{\Delta_\inj} \urar["\Tot"'] &
        \end{tikzcd}
    \]
    However, the same is not true for the partial totalisations: the inclusion $(\Delta_\inj)_{\leq n} \subseteq \Delta_{\leq n}$ is not homotopy initial; see, e.g., \cite[Section~21.6]{dugger_htpy_colimits_primer}.
    Limits over $(\Delta_\inj)_{\leq n}$ are bigger in some sense; see \cref{rmk:Tot_tower_for_semicosimplicial,prop:E1_of_Tot_sseq}.
\end{remark}

The semicosimplicial Tot tower is in fact an example of a cosimplicial Tot tower.

\begin{remark}
    \label{rmk:Tot_tower_for_semicosimplicial}
    Right Kan extension along $j\colon \Delta_\inj \subseteq \Delta$ is a right adjoint to the forgetful functor $\Sp^\Delta \to \Sp^{\Delta_\inj}$, and likewise for $\Ab$ in place of $\Sp$.
    If $X^\bullet$ is a semicosimplicial spectrum, then we can compute this right Kan extension to be the cosimplicial object given by
    \[
        (j_*X)^n = \prod_{[n]\twoheadto [k]} X^k
    \]
    with the product ranging over all surjective maps in $\Delta$.
    A similar formula holds for semicosimplicial abelian groups.
    We learn a number of things from this computation.
    \begin{numberenum}
        \item The functor $\pi_*$ preserves right Kan extension along $j$.
        \item The restriction of $j_* X^\bullet$ to $\Delta_{\leq n}$ is right Kan extended from the restriction of $X^\bullet$ to $(\Delta_\inj)_{\leq n}$.
        In particular, the Tot tower of $j_* X$ is the semicosimplicial Tot tower of $X^\bullet$.
    \end{numberenum}
\end{remark}


When working with cosimplicial rather than semicosimplicial objects, the Tot tower sets up a one-to-one correspondence between cosimplicial objects and certain towers.
In the following, let us write  $\Fil^{\geq 0}(\C)$ for the full subcategory of $\Fil(\C)$ on those filtered objects that vanish in negative filtration.
(This subcategory should not be confused with the connective part of the diagonal t-structure on $\Fil(\C)$.)
Then by definition, the Tot tower functor of \cref{def:Tot_and_Tot_tower} lands in $\Fil^{\geq0}(\C)$.

\begin{theorem}[Stable Dold--Kan correspondence, Lurie]
    \label{thm:stable_dold_kan}
    Let $\C$ be a stable $\infty$-category.
    Then the Tot tower functor restricts to an equivalence
    \[
        \towTot \colon \C^\Delta \simeqto \Fil^{\geq 0}(\C).
    \]
\end{theorem}
\begin{proof}
    This is a reformulation of \cite[Theorem~1.2.4.1]{HA}.
    Indeed, the notion of stability is self-dual, so the equivalence proved there dualises to an equivalence
    \[
        \towTot \colon \C^\Delta \simeqto \Fun(\Z_{\geq 0}^\op,\ \C).
    \]
    Finally, right Kan extension along the inclusion $\Z_{\geq 0}^\op \to \Z^\op$ (informally, putting zeroes in negative filtrations) results in a functor $\Fun(\Z_{\geq0}^\op,\ \C) \to \Fil(\C)$ that is fully faithful with essential image $\Fil^{\geq 0}(\C)$.
\end{proof}

The Tot tower leads to a spectral sequence.
We prefer to work with spectral sequences arising from filtrations rather than towers, so we apply the reflection duality from \cref{sec:reflecting} to land in this situation.

\begin{definition}
    \label{def:Tot_filtration}
    Let $X^\bullet$ be a cosimplicial spectrum.
    By reflecting the Tot tower of $X^\bullet$ via the functor $(\blank)^\fil$ from \cref{def:reflection_functors}, we obtain a filtered spectrum that we call the \defi{totalisation filtration} (or \emph{Tot filtration}) of $X^\bullet$.
    We denote its terms by
    \[
        \Tot^n X^\bullet \defeq \fib(\Tot X^\bullet \to \Tot_{n-1} X^\bullet).
    \]
    The resulting filtration $\filTot X^\bullet$ is of the form
    \[
        \begin{tikzcd}
            \dotsb \rar & \Tot^2 X^\bullet \rar & \Tot^1 X^\bullet \rar & \Tot X^\bullet \rar[equals] & \Tot X^\bullet \rar[equals] & \dotsb,
        \end{tikzcd}
    \]
    which is constant from filtration $0$ onwards.
    We call the spectral sequence associated to this filtered spectrum the \defi{totalisation spectral sequence} (or \emph{Tot spectral sequence}, or \emph{Bousfield--Kan spectral sequence}).
\end{definition}


Note that the Tot filtration is, by construction, a complete filtration of $\Tot X^\bullet$.
In general however, the conditional convergence of the Tot spectral sequence need not be strong.

\begin{remark}
    \label{rmk:filTot_Dold_Kan_version}
    From \cref{thm:stable_dold_kan,prop:properties_reflection}, it follows that the Tot filtration functor induces an equivalence from $\C^\Delta$ to the full subcategory of $\Filhat(\C)$ on those objects that are constant from filtration $0$ onwards.
\end{remark}

A nice feature of this spectral sequence is that there is a formula for its first page and its first differential purely in terms of the cosimplicial abelian groups $\oppi_t X^\bullet$.

\begin{notation}
    \leavevmode
    \begin{numberenum}
        \item Let $A^\bullet$ be a semicosimplicial abelian group.
        The \defi{unnnormalised cochain complex} $\uC(A^\bullet)$ of $A^\bullet$ is the cochain complex
        \[
            \dotsb \to 0 \to \uC^0 (A^\bullet) \to \uC^1 (A^\bullet) \to \uC^2 (A^\bullet) \to \dotsb
        \]
        given by, for $m\geq 0$,
        \[
            \uC^m (A^\bullet) = A^m
        \]
        and with differential $\uC^m (A^\bullet) \to \uC^{m+1} (A^\bullet)$ given by the alternating sum $\sum (-1)^i d^i$.

        \item If $A^\bullet$ is a cosimplicial group, we define its unnormalised cochain complex to be the unnormalised cochain complex of the semicosimplicial group $j^*A^\bullet$.
        
        \item Let $A^\bullet$ be a cosimplicial abelian group.
        The \defi{normalised cochain complex} $\uN(A^\bullet)$ of $A^\bullet$ is the cochain complex
        \[
            \dotsb \to 0 \to \uN^0 (A^\bullet) \to \uN^1 (A^\bullet) \to \uN^2 (A^\bullet) \to \dotsb
        \]
        given by, for $m\geq 0$,
        \[
            \uN^m (A^\bullet) = \bigcap_{i=0}^{m-1} \ker(s^i \colon A^m \to A^{m-1})
        \]
        and with differential $\uN^m (A^\bullet) \to \uN^{m+1} (A^\bullet)$ given by the alternating sum $\sum (-1)^i d^i$.
    \end{numberenum}
\end{notation}

If $A^\bullet$ is a cosimplicial abelian group, then we have an evident map of cochain complexes
\[
    \uN(A^\bullet) \to \uC(A^\bullet).
\]
This turns out to be a quasi-isomorphism, and to even admit a quasi-inverse; see the dual\footnote{To see that the definition of the normalised chain complex of a simplicial object in $\calA^\op$ used therein indeed dualises to the normalised cochain complex of a cosimplicial object in $\calA$, use \cite[Theorem~III.2.1]{goerss_jardine_simplicial}.} of \cite[Proposition~1.2.3.17]{HA}.

\begin{notation}
    Let $A^\bullet$ be a semicosimplicial abelian group, and let $s \geq 0$.
    We write $\pi^s A^\bullet$ for the $s$-th cohomology group of $\uC(A^\bullet)$.
\end{notation}


\begin{proposition}
    \label{prop:E1_of_Tot_sseq}
    \leavevmode
    \begin{numberenum}
        \item Let $X^\bullet$ be a cosimplicial spectrum, and let $\uE_r^{*,*}(X^\bullet)$ denote the resulting Tot spectral sequence.
        For all integers $n$ and $s$, there is a natural isomorphism
        \[
            \uE_1^{n,s}(X^\bullet) \cong \uN^{s} (\oppi_{n+s}X^\bullet)
        \]
        that identifies the $d_1$-differential with the differential of $\uN(\oppi_{n+s}X^\bullet)$.
        \item Let $X^\bullet$ be a semicosimplicial spectrum, and let $j_*X^\bullet$ denote the right Kan extension of it to a cosimplicial object \textbr{see \cref{rmk:Tot_tower_for_semicosimplicial}}.
        For all integers $n$ and $s$, there is a natural isomorphism
        \[
            \uE_1^{n,s}(j_*X^\bullet) \cong \uC^{s} (\oppi_{n+s}X^\bullet)
        \]
        that identifies the $d_1$-differential with the differential of $\uC(\oppi_{n+s}X^\bullet)$.
        \item \label{item:cosimp_vs_semi_second_page} Let $X^\bullet$ be a cosimplicial spectrum.
        Then under the above identifications, the unit $X^\bullet \to j_*j^*X^\bullet$ induces, on the first page of the Tot spectral sequence, the natural map
        \[
            \uN^s(\oppi_{n+s}X^\bullet) \to \uC^s(\oppi_{n+s}X^\bullet).
        \]
        In particular, this map of spectral sequences is an isomorphism from the second page onward.
    \end{numberenum}
\end{proposition}

\begin{proof}
    This follows by dualising (i.e., applying it to $\C =\Sp^\op$) the discussion of \cite[Remark~1.2.4.4, Variant~1.2.4.9]{HA}.
\end{proof}

\begin{remark}
    One can phrase \cref{item:cosimp_vs_semi_second_page} at the level of filtered spectra using the d\'ecalage functor.
    Let $\Dec \colon \FilSp \to \FilSp$ denote Antieau's d\'ecalage functor; see \cite{antieau_decalage} or \cite[Part~II]{hedenlund_phd}.
    Then if $X^\bullet$ is a cosimplicial spectrum, the above implies that the unit $X^\bullet \to j_*j^*X^\bullet$ induces an isomorphism of filtered spectra
    \[
        \Dec(\filTot(X^\bullet)) \to \Dec(\filTot(j_*j^*X^\bullet)).
    \]
    Indeed, the filtrations before applying d\'ecalage are complete, and since d\'ecalage preserves complete filtrations (\cite[Lemma~4.18]{antieau_decalage}), it is enough to check this on associated graded; see \cref{prop:pullback_tau_complete_and_invert_tau}.
    After applying d\'ecalage, the associated graded of this map is an isomorphism by \cref{item:cosimp_vs_semi_second_page} above.
\end{remark}

    

There is a variant of the Tot filtration where we have `turned a page' in the spectral sequence, and thereby starts on the second page of the Tot spectral sequence.
While there is a more general notion of such a page-turning operation (see \cref{rmk:comparison_alternative_constrs_sseq}), the cosimplicial version has a simpler expression.
These two definitions of d\'ecalage actually coincide: see \cref{prop:cosimplical_decalage_defs_agree}.

The terminology is taken from Deligne's operation of the same name for chain complexes from \cite[Definition~(1.3.3)]{deligne_hodge_theory_II}.
The generalisation to cosimplicial spectra is a result of Levine \cite[Section~6]{levine_slice_ANSS}.
Note that he also studies other (for instance, unstable) settings; we specialise his results to spectra.
Moreover, Levine compares his d\'ecalage operation to that of Deligne; see \cite[Remark~6.4]{levine_slice_ANSS}.

\begin{definition}
    \label{def:cosimplicial_decalage}
    Let $X^\bullet \colon \Delta \to \Sp$ be a semicosimplicial spectrum.
    We define the \defi{d\'ecalage} of $X^\bullet$ as the filtered spectrum given by
    \[
        \Dec^\Delta X^\bullet = \Tot(\Wh X^\bullet).
    \]
    If $X^\bullet$ is a cosimplicial spectrum, then we let $\Dec^\Delta X^\bullet$ denote the d\'ecalage of the underlying semicosimplicial object $j^*X^\bullet$; cf.\ \cref{rmk:Tot_only_depends_semicosimplial}.
\end{definition}

Concretely, the value at filtration $s$ of $\Dec^\Delta X^\bullet$ is given by $\Tot(\tau_{\geq s} X^\bullet)$, the totalisation of the (semi)cosimplicial spectrum obtained by applying $\tau_{\geq s}$ levelwise to $X^\bullet$.

\begin{remark}
    \label{rmk:decalage_lax_monoidal}
    The functor $\Dec^\Delta \colon \Sp^{\Delta_\inj} \to \FilSp$ is naturally lax symmetric monoidal.
    Indeed, it is the composite
    \[
        \begin{tikzcd}
            \Sp^{\Delta_\inj} \rar["\Wh"] & \Fil(\Sp)^{\Delta_\inj} \rar["\Tot"] & \FilSp.
        \end{tikzcd}
    \]
    The first functor is lax symmetric monoidal (for the levelwise symmetric monoidal structure on cosimplicial objects) by \cref{rmk:whitehead_lax_mon}, and the second is lax symmetric monoidal because it is the right adjoint to the constant functor $\FilSp \to \Fil(\Sp)^{\Delta_\inj}$ which is symmetric monoidal.
    The same applies for cosimplicial objects (or alternatively by noting that the forgetful functor to semicosimplicial objects is symmetric monoidal).
\end{remark}


The namesake of the d\'ecalage construction is that it has turned the page once compared to the Tot spectral sequence.
To state the comparison between these spectral sequences then, we require a reindexing; see \cref{rmk:E2_indexing_FilSp}.

\begin{theorem}[Levine]\label{decalageidentification}
    Let $X^\bullet$ be a \textbr{semi}cosimplicial spectrum, and let $\set{\mathrm{E}_r^{*,*}(X^\bullet)}_{r\geq 1}$ denote the Tot spectral sequence associated to $X^\bullet$.
    Then there is an isomorphism of spectral sequences \textbr{where $r\geq 1$}
    \[
        \uE_r^{n,s}(\Dec^\Delta X^\bullet) \cong \uE_{r+1}^{n,\, s-n}(X^\bullet).
    \]
\end{theorem}
\begin{proof}
    This follows from \cite[Proposition~6.3]{levine_slice_ANSS}, but note that Levine uses a different indexing from the Adams indexing used above.
\end{proof}

Accordingly, it often makes sense to use \emph{second-page indexing} for the spectral sequence underlying $\Dec^\Delta X^\bullet$.
Phrased like this, the above theorem says that this second-page indexed $\Dec^\Delta X^\bullet$ is isomorphic to the second page onwards of $\uE_r^{*,*}(X^\bullet)$.

For later use, we record the following more basic property of the cosimplicial d\'ecalage construction.

\begin{proposition}
    \label{lem:decalage_is_complete}
    Let $X^\bullet$ be a \textbr{semi}cosimplicial spectrum.
    Then the filtered spectrum $\Dec^\Delta X^\bullet$ is naturally a complete filtration of $\Tot X^\bullet$, meaning that its limit vanishes and its colimit is naturally isomorphic to $\Tot X^\bullet$.
\end{proposition}

\begin{proof}
    Complete filtered spectra are closed under limits, so completeness follows from the fact that Whitehead filtrations of spectra are complete.
    To compute the colimit of this filtration, note that, for every integer $s$, we have a fibre sequence
    \[
        \tau_{\ge s+1}X^\bullet\to X^\bullet\to \tau_{\le s}X^\bullet.
    \]
    Taking totalisations, one therefore has a natural fibre sequence
    \[
        \Tot (\tau_{\ge s+1}X^\bullet) \to \Tot X^\bullet \to \Tot (\tau_{\le s}X^\bullet).
    \]
    Taking colimits over $s$, one has a natural cofibre sequence
    \[
        \colim_s \Tot (\tau_{\ge s+1}X^\bullet) \to \Tot X^\bullet \to \colim_s \Tot (\tau_{\le s}X^\bullet).
    \]
    By definition, the left-hand term gives the colimit of $\Dec^\Delta(X^\bullet)$, so it suffices to show that the right-hand term vanishes. Since coconnectivity is preserved by limits, we see that for all $s$, the spectrum $\Tot\tau_{\le s}X^\bullet$ is $s$-truncated.
    As homotopy groups of spectra preserve filtered colimits, the colimit is therefore $(-\infty)$-truncated, and hence vanishes.
\end{proof}

    

Finally, we end by comparing Levine's cosimplicial d\'ecalage to Antieau's d\'ecalage.
This appears to be folklore, but for lack of a citeable reference, we give a proof here.
Note, however, that while Levine's d\'ecalage also applies to unstable settings, Antieau's d\'ecalage lives entirely in the stable world.
Let us write $\Dec \colon \FilSp \to \FilSp$ for Antieau's d\'ecalage functor.

\begin{proposition}
    \label{prop:cosimplical_decalage_defs_agree}
    Let $X^\bullet$ be a \textbr{semi}cosimplicial spectrum.
    Then there is an isomorphism, natural in $X^\bullet$, of filtered spectra
    \[
        \Dec^\Delta X^\bullet \cong \Dec(\filTot X^\bullet).
    \]
\end{proposition}

The proof uses the exact same argumentation as Antieau's argument for the Atiyah--Hirzebruch spectral sequence in \cite[Propostion~9.2]{antieau_decalage}.

\begin{proof}
    We assume familiarity with the definition of d\'ecalage via connective covers in the Be\u{\i}linson t-structure, as explained in detail in \cite{antieau_decalage} or \cite[Part~II]{hedenlund_phd}.

    Applying $\filTot$ to the diagram $\Wh X^\bullet \to X^\bullet$ in $\Sp^{\Delta_\inj}$ yields a natural diagram in $\FilSp$
    \begin{equation}
        \label{eq:Beilinson_cover_of_Totfil}
        \dotsb \to \filTot(\tau_{\geq 1}X^\bullet) \to \filTot(\tau_{\geq 0}X^\bullet) \to \dotsb \to \filTot(X^\bullet).
    \end{equation}
    We claim that this is the Whitehead filtration of $\filTot(X^\bullet)$ in the Be\u{\i}linson t-structure on $\FilSp$.
    First, we show that for every integer $n$, the filtered spectrum $\filTot(\tau_{\geq n}X^\bullet)$ is Be\u{\i}linson $n$-connective.
    For this, we ought to show that its associated graded in filtration $s$ is an $(n-s)$-connective spectrum.
    Using \cref{prop:E1_of_Tot_sseq}, we compute
    \[
        \pi_k (\Gr^s \filTot (\tau_{\geq n} X^\bullet)) = \uC^s (\oppi_{k+s}(\tau_{\geq n}X^\bullet)),
    \]
    which evidently vanishes if $k+s<n$, i.e., if $k<n-s$, so that the $s$-th associated graded is indeed $(n-s)$-connective.
    When working with cosimplicial objects and the cosimplicial Tot filtration, the same applies, using the normalised cochain complex instead.
    
    It follows that the natural map $\filTot(\tau_{\geq n}X^\bullet) \to \filTot(X^\bullet)$ factors through a map
    \[
        \filTot(\tau_{\geq n}X^\bullet) \to \optau_{\geq n}^{\mathrm{Bei}} \filTot(X^\bullet).
    \]
    We claim this is an isomorphism.
    As both filtrations are complete (the second one by \cite[Lemma~4.18]{antieau_decalage}), it is enough to show this on associated graded; see \cref{prop:pullback_tau_complete_and_invert_tau}.
    There, it follows using \cref{prop:E1_of_Tot_sseq} and \cite[Proposition~II.2.4]{hedenlund_phd}.

    We conclude that \eqref{eq:Beilinson_cover_of_Totfil} is indeed the Be\u{\i}linson Whitehead filtration of $\filTot X^\bullet$.
    From the definition of d\'ecalage, it follows that
    \[
        \Dec(\filTot X^\bullet) = \colim_s \Tot^s(\tau_{\geq *}X^\bullet) \cong \Tot(\tau_{\geq *}X^\bullet) = \Dec^\Delta X^\bullet.\qedhere
    \]
\end{proof}

Using this comparison result, the properties of d\'ecalage from \cref{decalageidentification,lem:decalage_is_complete} also follow from general properties of Antieau's d\'ecalage functor from \cite{antieau_decalage} and \cite[Part~II]{hedenlund_phd}.

\begin{remark}
    \label{rmk:only_cosimplicial_decalage_is_monoidal}
    The cosimplicial d\'ecalage functor $\Dec^\Delta$ still has an advantage over the composite functor $\Dec \circ \filTot$: the former is naturally lax symmetric monoidal (see \cref{rmk:decalage_lax_monoidal}).
    Although the functor $\Dec \colon \FilSp \to \FilSp$ is also lax symmetric monoidal, the functor $\filTot$ is not,\footnotemark\ so the lax monoidal structure only arises through the version of \cref{def:cosimplicial_decalage}.
\end{remark}

\footnotetext{If it were, then it would send $\E_\infty$-rings in $\Sp^\Delta$ (for the levelwise monoidal structure) to filtered $\E_\infty$-rings.
Instead, what we see is that it sends $\E_\infty$-rings in $\Sp^\Delta$ to filtered objects in $\E_\infty$-rings.}

\begin{remark}[Be\u{\i}linson vs.\ levelwise t-structures]
    We can deduce more than stated in the proposition above: we learn that under the stable Dold--Kan correspondence, the levelwise t-structure on cosimplicial spectra corresponds to the Be\u{\i}linson t-structure on a subcategory of $\FilSp$.
    This is remarked by Lawson in \cite[Remark~3.16]{lawson_filtered_spaces_objects}.
    In more detail: recall from \cref{rmk:filTot_Dold_Kan_version} that the functor $\filTot \colon \Sp^\Delta \to \Filhat\Sp$ defines an equivalence onto those complete filtered spectra that are constant from filtration $0$ onwards.
    In the proof above, we showed that $\filTot$ sends the levelwise Whitehead tower to the Be\u{\i}linson Whitehead tower.
    Since it is an equivalence, the result follows.
\end{remark}

\subsection{The cosimplicial Adams spectral sequence}

\label{ssec:cosimplicial_Adams_sseq}

The most general definition of the $E$-Adams spectral sequence we will use is the following, though often we will work with more specific (and more structured) spectra $E$.

\begin{definition}
    \label{def:classical_E_ASS}
    Let $E$ be a spectrum with a map $\S \to E$, and let $X$ be a spectrum.
    The map $\S \to E$ gives rise to an augmented semicosimplicial spectrum $\Delta_{\inj,+} \to \Sp$ of the form
    \[
        \begin{tikzcd}
            \S \rar & E \ar[r,shift left] \ar[r,shift right] & E \otimes E \ar[r, shift left=2] \rar \rar[shift right=2] & \dotsb.
        \end{tikzcd}
    \]
    The \defi{semicosimplicial $E$-based Adams resolution} for $X$ is the semicosimplicial spectrum
    \[
        \ASS_E^\Delta(X) \defeq E^{[\bullet] }\otimes X = \left(\begin{tikzcd}
            E \otimes X \ar[r,shift left] \ar[r,shift right] & E \otimes E \otimes X \ar[r, shift left=2] \rar \rar[shift right=2] & \dotsb
        \end{tikzcd}\right).
    \]
    Define the filtered spectrum
    \begin{align*}
        \ASS_E(X) &\defeq \filTot(\ASS_E^\Delta(X)).
    \intertext{The \defi{$E$-based Adams spectral sequence} for $X$ is the spectral sequence underlying this filtered spectrum.
    More generally, if $Y$ and $X$ are spectra, then we define the semicosimplicial spectrum}
        \ASS_E^\Delta(Y,X) &\defeq \map(Y,\, E^{[\bullet]}\otimes X)\\
    \intertext{and the filtered spectrum}
        \ASS_E(Y,X) &\defeq \filTot(\ASS_E^\Delta(Y,X)),
    \end{align*}
    and define the $E$-based Adams spectral sequence for $[Y,X]$ to be the spectral sequence underlying this filtered spectrum.
\end{definition}


We are careful to work with semicosimplicial objects rather than cosimplicial objects, because upgrading $E^{[\bullet]}$ to a cosimplicial object requires an $\E_1$-structure on $E$; see \cite[Construction~2.7]{mathew_naumann_noel_nilpotence_descent}.
Not all $E$ of interest may admit this structure, and the Adams spectral sequence should not depend on it either.
If $E$ does admit this structure, then this upgrades $\ASS_E^\Delta(Y,X)$ to a cosimplicial spectrum, which we will denote by the same notation.
We would then define $\ASS_E(Y,X)$ using the cosimplicial Tot tower, which has a more efficient first page compared to the semicosimplicial approach; see \cref{prop:E1_of_Tot_sseq}.

This brings us to a potentially confusing point about the Adams spectral sequence: usually, one is interested in it only from the \emph{second page} onward.
One should view the first page as `ill-defined' in some sense; it may admit many models, and the one from the above definition is a rather inefficient one at that.
What we call the cosimplicial Adams resolution should be regarded as only one of many potential resolutions; the one we chose is convenient as it is functorial.

Our definition above is made to align with standard conventions.
Alternatively, one can do away with the first page entirely, as follows.

\begin{variant}
    \label{var:ASS_using_cosimplicial_decalage}
    Alternatively, we could have defined the Adams spectral sequence as arising from the filtered spectrum
    \[
        \Dec^\Delta(\ASS_E^\Delta(Y,X)).
    \]
    This filtered spectrum should be viewed as the `true' incarnation of the Adams spectral sequence.
    It turns out to have better monoidality properties, though this is only possible to prove in general using more modern machinery; see \cref{sec:signature_synthetic_analogue}, particularly \cref{rmk:sigmanu_lax_monoidal,rmk:identifying_sigma_nu_monoidally}.
    To align with standard indexing, one should index the resulting spectral sequence using second-page indexing; see \cref{rmk:E2_indexing_FilSp}.
\end{variant}



\begin{remark}
    The Adams spectral sequence is not specific to spectra, but as this is the main case of interest, we will stick to it for our discussion.
    Baker--Lazarev \cite{baker_lazarev_Rmod_Adams_sseq} set up the Adams spectral sequence in modules over an $\E_\infty$-ring, and Mathew--Naumann--Noel \cite[Part~1]{mathew_naumann_noel_nilpotence_descent} set up the Adams spectral sequence in a presentably symmetric monoidal $\infty$-category.
    For an even more general setup, see the next remark.
\end{remark}

\begin{remark}
    \label{rmk:millers_Adams_sseq}
    Miller \cite{miller_adams_sseq,miller_adams_sseq_lecture_notes} showed that the $E$-based Adams spectral sequence depends on much less information than the ring spectrum $E$: it only depends on the class of morphisms that become nullhomotopic after tensoring with $E$.
    This is similar to how $E$-localisation of spectra depends on much less information than $E$.
    For a further discussion and extension of these ideas, see \cite[Sections~2 and~3]{patchkoria_pstragowski_derived_inftycats}.
\end{remark}



For the interested reader, we compare the cosimplicial approach to a filtered approach.

\begin{remark}[Cosimplicial approach and completion]
    \label{rmk:cosimplicial_Adams_vs_filtered_resolution}
    The downside to using the (semi)cosimplicial approach is that it is only able to retrieve the \emph{completion} of the filtered spectrum giving rise to the Adams spectral sequence.
    Let us explain this by comparing it with the other approach.
    If $E$ is a spectrum with a map $\S \to E$, let $\overline{E}$ denote the fibre of this map.
    The \emph{filtered $E$-Adams resolution}\footnote{Various people refer to this as the \emph{$E$-Adams tower}, but our use of the words \emph{tower} and \emph{filtration} (see \cref{rmk:filtrations_vs_tower}) prevents us from using that terminology here.} of the sphere is the filtered spectrum given by
    \[
        \begin{tikzcd}
            \dotsb \rar & \overline{E} \otimes \overline{E} \rar & \overline{E} \rar & \S \rar[equals] & \S \rar[equals] & \dotsb
        \end{tikzcd}
    \]
    indexed to be constant from filtration $0$ onwards.
    If $X$ is a spectrum, then the $E$-Adams filtration of $X$ is by definition obtained by tensoring this with $X$ levelwise.
    Clearly the resulting filtered spectrum has colimit $X$, i.e., it is a filtration of $X$.
    It need not be complete however, and this results in a convergence problem.
    If $E$ is an $\E_1$-ring, then the semicosimplicial object $E^{[\bullet]}$ naturally upgrades to a cosimplicial diagram; see \cite[Construction~2.7]{mathew_naumann_noel_nilpotence_descent}.
    Mathew--Naumann--Noel \cite[Proposition~2.14]{mathew_naumann_noel_nilpotence_descent} show that, under the stable Dold--Kan correspondence of \cref{thm:stable_dold_kan}, the Tot tower of the cosimplicial spectrum $E^{[\bullet]} \otimes X$ is matched up with the associated tower (in the sense of \cref{def:reflection_functors}) of the filtered $E$-Adams resolution for $X$.
    Because the associated tower functor factors through completion (\cref{prop:properties_reflection}), this shows that the cosimplicial approach recovers only the completion of the approach based on filtered Adams resolutions in the above sense.
\end{remark}

In the generality of \cref{def:classical_E_ASS}, it is very hard to say much about the second page of the resulting spectral sequence.
Things improve if we impose conditions on~$E$.
The following is a rather restrictive one, but luckily covers some of the main cases of interest.

\begin{definition}[\cite{adams_blue_book}, Condition~III.13.3]
    \label{def:Adams_type}
    Let $E$ be a homotopy associative ring spectrum.
    \begin{numberenum}
        \item A finite spectrum $P$ is called \defi{finite $E$-projective} if $E_*P$ is a projective $E_*$-module.\footnotemark
        \item We say that $E$ is of \defi{Adams type} if it can be written as a filtered colimit of finite $E$-projective spectra $E_\alpha$ such that for every $\alpha$, the natural map
        \begin{equation}
            \label{eq:fp_spectrum_dualising_condition}
            E^*E_\alpha \to \Hom_{E_*}(E_*E_\alpha,\, E_*)
        \end{equation}
        is an isomorphism.
    \end{numberenum}
\end{definition}

\footnotetext{This should not be confused with what one might call an \emph{$E$-finite projective} spectrum, meaning a spectrum $P$ such that $E_*P$ is a finite projective $E_*$-module. This notion need not imply that the spectrum $P$ is itself finite, but this is a condition we very intentionally require on $P$.}

\begin{remark}
    If $E$ is an $\E_1$-ring, then the condition on a finite-projective $E_\alpha$ that the map \eqref{eq:fp_spectrum_dualising_condition} is an isomorphism is automatic.
    Indeed, if $E$ is $\E_1$, then we have a good $\infty$-category $\Mod_E(\Sp)$ of $E$-modules, which we can use to set up an Ext spectral sequence as in \cite[Chapter~IV]{EKMM}.
    The fact that $E_*E_\alpha$ is projective implies that the resulting K\"unneth spectral sequence computing $E^*E_\alpha$ is concentrated in filtration~$0$, implying that \eqref{eq:fp_spectrum_dualising_condition} is an isomorphism.
    (In fact, the definition of Adams type is engineered to be able to set up a K\"unneth and universal coefficient spectral sequence; see \cite[Chapter~III.13]{adams_blue_book}.)
\end{remark}

\begin{remark}
    \label{rmk:Adams_implies_flat}
    If $E$ is of Adams type, then $E_*E$ is in particular a flat (left and right) $E_*$-module.
    Indeed, homotopy groups preserve filtered colimits, and projective modules are flat.
\end{remark}

\begin{example}
    \label{ex:adams_type}
    \leavevmode
    \begin{numberenum}
        \item The sphere is of Adams type: it is the colimit of the one-point diagram $\set{\S}$.
        
        \item If $E$ is $\F_p$, or more generally if $\pi_*E$ is a graded field (e.g., if $E$ is a Morava K-theory), then every finite spectrum is finite $E$-projective.
        Since every spectrum is a filtered colimit of finite spectra, it follows that $E$ is of Adams type if $\pi_*E$ is a graded field.

        The case $E=\F_p$ is the one originally considered by Adams in \cite{Adams_Adams_sseq}, and is often simply referred to as the \emph{Adams spectral sequence}.
        (Strictly speaking, the original version of \cite{Adams_Adams_sseq} is based on \emph{cohomology} rather than homology, but for $E=\F_p$, this difference is less material.)
        
        \item The ring spectrum $\MU$ is of Adams type, being witnessed by it being the colimit of Thom spectra of Grassmannians.
        Moreover, every Landweber-exact homotopy-associative ring spectrum is of Adams type; see \cite[Proposition~1.3]{devinatz_morava_brown_comenetz}.
        In particular, Morava E-theories are of Adams type.

        The cases $E=\MU$ and $E=\BP$ are both referred to as the \emph{Adams--Novikov spectral sequence} (ANSS).
        We refer to \cite{ravenel_novice_guide_ANSS} for an introduction to the Adams--Novikov spectral sequence and its interplay with the $\F_p$-Adams spectral sequence.
        
        \item A non-example is $\Z$: this follows since $\pi_*(\Z\otimes \Z)$ contains $p$-torsion for every prime $p$.
        In particular, it is not flat over the integers, so \cref{rmk:Adams_implies_flat} shows it cannot be of Adams type.
        Likewise, for every prime $p$, the ring $\Z_{(p)}$ is not of Adams type.
        
        However, one can modify the notion of Adams type and work with $\Z$ and $\Z_{(p)}$ as if they were of Adams type; see \cite{burklund_pstragowski_quivers}.
        
        \item Another non-example is real K-theory, both in its connective and periodic variants.
        Mahowald \cite{mahowald_bo_resolutions} nevertheless computes the second page of the $\ko$-based Adams spectral sequence, leading to a proof of the telescope conjecture at height $1$ and at the prime $2$; see \cite{mahowald_image_J_EHP}.\qedhere
    \end{numberenum}
\end{example}

The reason for imposing these restrictions is to obtain an abelian category that computes the second page of the Adams spectral sequence.
This abelian category is defined for any Hopf algebroid, not just ones arising from a ring spectrum.
For a further and more detailed treatment of (the category of) comodules over a Hopf algebroid, we refer to \cite[Section~1]{hovey_htpythy_comodules} or \cite[Appendix~A.1]{ravenel_green_book}.

\begin{definition}
    Let $(A,\Gamma)$ be a graded Hopf algebroid.
    We write $\grComod_{(A,\Gamma)}$ for the category of comodules over $(A,\Gamma)$ in graded abelian groups.
    If $n$ is an integer, then we write $[n]$ for the $n$-fold shift operator,
    \[
        (M[n])_i = M_{i-n}.
    \]
    The category $\grComod_{(A,\Gamma)}$ is naturally symmetric monoidal, with the underlying $A$-module being given by the tensor product over $A$.
\end{definition}

To ensure that we can do homological algebra in this setting, we need an assumption on $\Gamma$.
We say that a (graded) Hopf algebroid $(A,\Gamma)$ is \defi{flat} if $\Gamma$ is flat as an $A$-module via either the left or right unit; note that it does not matter which unit we choose, as they differ by an automorphism of $\Gamma$.

\begin{proposition}
    Let $(A,\Gamma)$ be a graded Hopf algebroid.
    Suppose that $(A,\Gamma)$ is flat.
    Then the category $\grComod_{(A,\Gamma)}$ is a Grothendieck abelian category, and the forgetful functor to $\grAb$ preserves small colimits and finite limits \textbr{in particular, it is exact}.
    In particular, $\grComod_{(A,\Gamma)}$ has enough injectives.
\end{proposition}
\begin{proof}
    Graded comodules over $(A,\Gamma)$ are equivalent to comodules over the comonad
    \[
        \grMod_A \to \grMod_A, \quad M \mapsto \Gamma \otimes_A M.
    \]
    If $(A,\Gamma)$ is flat, then this comonad is exact, from which it follows that the category of comodules is abelian and that the forgetful functor to $\grMod_A$ is exact; see, e.g., \cite[Proposition~A.2]{patchkoria_pstragowski_derived_inftycats}.
    Consequently, the forgetful functor also detects exactness (since it is conservative).
    Combining this with the fact that $\grMod_A$ is a Grothendieck abelian category and that the comonad $\Gamma\otimes_A \blank$ preserves colimits, one can check directly that $\grComod_{(A,\Gamma)}$ is Grothendieck abelian also.
\end{proof}

The forgetful functor $\grComod_{(A,\Gamma)} \to \grMod_A$ has a right adjoint, given on underlying modules by $M \mapsto \Gamma \otimes_A M$.
It follows immediately that the image of an injective $A$-module under this right adjoint is an injective comodule.

In general, if $\Gamma$ is not flat over $A$, then $\grComod_{(A,\Gamma)}$ is merely an additive category.

\begin{notation}
    Let $(A,\Gamma)$ be a graded Hopf algebroid, and let $M$ and $N$ be comodules over it.
    If $s$ and $t$ are integers, then we write
    \[
        \Ext_{(A,\Gamma)}^{s,t}(M,N) \defeq \Ext_{(A,\Gamma)}^s(M[t],\, N).
    \]
    Note that this is a \emph{homological} indexing convention for $t$, and a \emph{cohomological} indexing convention for $s$.
\end{notation}

\begin{construction}
    If $E$ is a homotopy-associative ring spectrum, then the pair $(E_*,\, E_*E)$ is naturally a Hopf algebroid.
    If $E$ is moreover of Adams type, then $E_*E$ is in particular flat over $E_*$; see \cref{rmk:Adams_implies_flat}.
    This gives us a good symmetric monoidal abelian category of comodules over $(E_*,\, E_*E)$.
    We will abbreviate this category by
    \[
        \grComod_{E_*E}.
    \]
    The fact that $E$ is of Adams type tells us more about this category: it tells us that it is generated under colimits by dualisable objects; see \cite[Section~3.1]{pstragowski_synthetic}.
\end{construction}

\begin{remark}
    \label{rmk:flat_vs_Adams_type}
    Let $E$ be a homotopy-associative ring spectrum.
    Then $E$ being of Adams type implies that $E_*E$ is flat; see \cref{rmk:Adams_implies_flat}.
    It is not known whether the converse is true; see \cite{burklund_levy_pstragowski_compositions_Adams_type} for a discussion (specifically Question~(1) following Remark~7 in op.\ cit.).
    As we will need the stronger assumption of Adams type, and because of the lack of a known counterexample to flatness implying Adams type, we content ourselves with assuming this potentially stronger condition.
\end{remark}


One potentially confusing point of the Adams spectral sequence is that there are two variants of what one might guess the induced strict filtration on $[Y,X]_*$ could be.
One of them is more accessible algebraically, while the other is what comes out of the definition of \cref{def:classical_E_ASS}.
In some cases, such as when $Y$ is a sphere, these agree; see \cref{prop:different_Adams_filtrations_agree} below.
We use the following (nonstandard) names to distinguish between them.

\begin{definition}
    \label{def:algebraic_topological_Adams_filtrations}
    Let $X$, $Y$ and $E$ be spectra.
    \begin{numberenum}
        \item The \defi{algebraic $E$-Adams filtration} on $[Y,X]$ is the strict filtration where, for every $s \geq 1$, a map $f$ is in $F^s [Y,X]$ when it can be written as a composite
        \[
            f = f_1 \circ \dotsb \circ f_s
        \]
        where each map $f_i$ induces the zero map on $E_*$-homology.
        We put $F^0[Y,X]$ equal to $[Y,X]$ by definition.
        \item The \defi{topological $E$-Adams filtration} on $[Y,X]$ is the strict filtration where, for every $s \geq 1$, a map $f$ is in $F^s [Y,X]$ when it can be written as a composite
        \[
            f = f_1 \circ \dotsb \circ f_s
        \]
        where each map $f_i$ becomes a nullhomotopic of spectra after tensoring with $E$.
        We put $F^0[Y,X]$ equal to $[Y,X]$ by definition.
    \end{numberenum}
\end{definition}

Taking $Y=\S^n$ results in filtrations on $\pi_n X$, which is the case we will be interested in most of the time.

In general, the filtration captured by the spectral sequence is the topological one.

\begin{proposition}
    \label{prop:ASS_induced_filtration_is_Adams_filtration}
    Let $E$ be a spectrum with a homotopy-unital multiplication, and let~$X$ and~$Y$ be spectra.
    The induced strict filtration on $[Y,X]_*$ by the filtered spectrum $\ASS_E(Y,X)$ is the topological $E$-Adams filtration.
\end{proposition}

\begin{proof}
    We freely use the notation introduced in \cref{rmk:cosimplicial_Adams_vs_filtered_resolution}.
    By loc.\ cit., we may compute the induced strict filtration on $[Y,X]_*$ using the Adams resolution
    \begin{equation}
        \label{eq:Adams_res_X}
        \dotsb \to \overline{E} \otimes \overline{E} \otimes X \to \overline{E} \otimes X \to X.
    \end{equation}
    Note that the map $\overline{E} \to \S$ becomes nullhomotopic after tensoring with $E$: indeed, the cofibre sequence defining $\overline{E}$ becomes, after tensoring with $E$, the cofibre sequence
    \[
        E\otimes \overline{E} \to E \to E \otimes E.
    \]
    The unital multiplication on $E$ provides a retraction of the second map, so indeed the first map is nullhomotopic.
    
    We see therefore that each map in the Adams resolution \eqref{eq:Adams_res_X} becomes null after tensoring with $E$.
    As a result, a map $\Sigma^n Y \to X$ that lifts to $\overline{E}{}^{\otimes s} \otimes X$ for $s \geq 1$ is of topological Adams filtration $\geq s$.
    For the converse, it is enough to consider the case $s=1$ and $n=0$, i.e., we have a map $f\colon Y \to X$ that becomes nullhomotopic after tensoring with $E$.
    Then also the composite $Y \to E\otimes Y \to E\otimes X$ is nullhomotopic, so it follows that there exists a dashed factorisation
    \[
        \begin{tikzcd}
            & Y \dar["f"] \dlar[dashed] & \\
            \overline{E} \otimes X \rar & X \rar & E\otimes X,
        \end{tikzcd}
    \]
    so that $f$ lifts to filtration $1$ in the induced strict filtration.
\end{proof}

Meanwhile, the algebraic filtration is preferable computationally.
Clearly, the topological Adams filtration of a map is a lower bound for its algebraic Adams filtration.
In general, the converse need not be true.
If $E$ is of Adams type, then it is true so long as $E_*Y$ is projective over $E_*$.

\begin{proposition}
    \label{prop:different_Adams_filtrations_agree}
    Let $E$ be a homotopy-associative ring spectrum of Adams type, and let~$X$ and~$Y$ be spectra.
    Suppose that $E_*Y$ is projective over $E_*$ \textbr{e.g., if $\pi_*E$ is a graded field, or if $Y = \S^n$}.
    Then the following hold.
    \begin{numberenum}
        \item There is an isomorphism
        \[
            \uE_2^{n,s}(Y,X) \cong \Ext_{E_*E}^{s,\,n+s}(E_*Y,\, E_*X),
        \]
        where the left-hand side denotes the second page of the $E$-Adams spectral sequence.
        \item The algebraic and topological $E$-Adams filtrations on $[Y,X]_*$ coincide.
    \end{numberenum}
    In particular, the above are true if $Y=\S$.
\end{proposition}

\begin{proof}
    For the first part, see \cite[Chapter III.15]{adams_blue_book}.
    The second part then follows from \cite[Warning~3.21]{patchkoria_pstragowski_derived_inftycats}.
\end{proof}

\begin{remark}
    \label{rmk:Estar_based_ASS}
    The idea of the proof is that one can construct a modified version of the $E$-Adams spectral sequence whose second page is \emph{always} given by Ext groups of comodules, and whose filtration is \emph{always} the algebraic one.
    One might refer to this as the \emph{$E_*$-based Adams spectral sequence}.
    There is a comparison map from $E$-based to the $E_*$-based version, and if $Y$ has projective homology, then this is an isomorphism on second pages (and hence also on all later pages).
    For a further explanation of this point, we refer to the work of Patchkoria--Pstr\k{a}gowski \cite{patchkoria_pstragowski_derived_inftycats}; see Section~3 of op.\ cit., particularly Warning~3.21 therein.
\end{remark}

\subsection{Convergence}

Next, we deal with the issue of convergence.
In this cosimplicial formulation, the object $\ASS_E(X)$ is a complete filtration by \cref{rmk:filTot_Dold_Kan_version}, so it might appear as if there are no convergence problems.
The issue, however, is whether the underlying spectrum of $\ASS_E(X)$ is isomorphic to $X$.
Only in this case are we willing to speak of (conditional) convergence of the Adams spectral sequence.

We closely follow Bousfield's seminal discussion on convergence \cite[Sections~5 and~6]{bousfield_localization_spectra}.
Although those results are specific to the setting of spectra, some results have been generalised by Mantovani \cite[Section~7]{mantovani_localisations_stable}.

We begin with some terminology.

\begin{definition}
    Let $X$ be a spectrum, and let $E$ be spectrum with a map $\S \to E$.
    The \defi{$E$-nilpotent completion} of $X$ is the spectrum
    \[
        X_E^\wedge \defeq \Tot (E^{[\bullet]} \otimes X).
    \]
    The map $\S \to E^{[\bullet]}$ induces a natural map $X \to X_E^\wedge$.
    We say that $X$ is \defi{$E$-nilpotent complete} if this map is an isomorphism.
\end{definition}
Bousfield shows that one may alternatively compute the $E$-nilpotent completion by a choice of what he calls an \emph{$E$-nilpotent resolution}; see \cite[Proposition~5.8]{bousfield_localization_spectra}.

Essentially by definition, the $E$-Adams spectral sequence converges conditionally to $\map(Y,\,X_E^\wedge)$.
If $Y$ is the sphere, then this is simply $X_E^\wedge$.
In both cases then, the question of convergence is whether the map $X \to X_E^\wedge$ is an isomorphism.

One should think of the $E$-nilpotent completion as a more computable approximation to the $E$-localisation of $X$.
As it is a limit of $E$-local objects, the $E$-nilpotent completion is $E$-local, so that the natural map $X \to X_E^\wedge$ factors through a map
\begin{equation}
    \label{eq:LE_to_Enilp_completion}
    L_E X \to X_E^\wedge.
\end{equation}
Unfortunately, this map can fail to be an isomorphism.
This failure is connected to some unexpected behaviour of nilpotent completion: it can fail to be idempotent.

\begin{proposition}[Bousfield]
    Let $E$ and $X$ be as above.
    Then the map \eqref{eq:LE_to_Enilp_completion} is an isomorphism if and only if the natural map $X_E^\wedge \to (X_E^\wedge)_E^\wedge$ is an isomorphism.
\end{proposition}
\begin{proof}
    See the discussion following the proof of Proposition~5.5 in \cite{bousfield_localization_spectra}.
\end{proof}


In the case where both $E$ and $X$ are bounded below, the $E$-nilpotent completion and $E$-localisation of $X$ coincide and are of an arithmetic nature.
For the following, recall that if $R$ is an ordinary ring, then its \emph{core} is its subring
\[
    \set{x \in R \mid x \otimes 1 = 1 \otimes x \text{ holds in }R \otimes R}.
\]

\begin{theorem}[Bousfield]
    Let $E$ be a bounded-below homotopy-associative ring spectrum.
    \begin{numberenum}
        \item Suppose that the core of the ring $\pi_0 E$ is $\Z[J^{-1}]$ for some set of primes $J$.
        Then for every bounded-below spectrum~$X$, the natural map $L_E X \to X_E^\wedge$ is an isomorphism, and moreover
        \[
            L_E X \cong X [J^{-1}].
        \]
        \item Suppose that the core of the ring $\pi_0 E$ is $\Z/n$ for some nonzero integer $n$.
        Then for every bounded-below spectrum~$X$, the natural map $L_E X \to X_E^\wedge$ is an isomorphism, and moreover
        \[
            L_E X \cong X_n^\wedge.
        \]
    \end{numberenum}
\end{theorem}
\begin{proof}
    See \cite[Theorems~6.5 and~6.6]{bousfield_localization_spectra}.
    Note that Bousfield uses the term \emph{connective} to mean \emph{bounded below}.
\end{proof}

\begin{remark}
    Bousfield computes all rings that can appear as the core of a ring.
    These fall into four types; the above two are the first two types.
    Bousfield shows \cite[Theorem~6.7]{bousfield_localization_spectra} that if the core of $\pi_0 E$ is of one of the other two types, then the map $L_E\Z \to \Z_E^\wedge$ is not an isomorphism.
\end{remark}

\begin{example}
    \label{ex:MU_BP_Fp_nilpotent_completion}
    The ring $\MU$ satisfies the first condition for $J = \varnothing$, so that for all spectra $X$ that are bounded below, we have
    \[
        L_\MU X \cong X_\MU^\wedge \cong X.
    \]
    For a fixed prime $p$, the ring $\BP$ satisfies the first condition for $J$ the set of primes different from $p$, so that for all bounded below spectra $X$, we have
    \[
        L_\BP X \cong X_\BP^\wedge \cong X_{(p)}.
    \]
    Lastly, the ring $\F_p$ satisfies the second condition for $n=p$, so that for all bounded below spectra $X$, we have
    \[
        L_{\F_p} X \cong X_{\F_p}^\wedge \cong X_p^\wedge.
    \]
    Beware that for a general $X$ not bounded below, none of these isomorphisms need to hold.
\end{example}

%% file: fil/tau_abridged.tex
So far, we have introduced spectral sequences in the way that they are normally introduced.
We will now rephrase them in terms of the \emph{$\tau$-formalism}.
In \cref{sec:tau_in_fil_ab}, we introduce $\tau$ in filtered abelian groups.
Like in the previous chapter, this is both as a warm-up and to describe the structure present on the homotopy groups of filtered spectra.
In this context, $\tau$ is a helpful variable for keeping track of filtrations and hidden extensions; we showcase some examples in \cref{ssec:eg_filtered_ab_tau_mod}.
Next, in \cref{sec:tau_in_fil_sp} we lift this story to filtered spectra, and begin to rephrase spectral sequences in the language of $\tau$.
We compare various operations in the $\tau$-formalism in spectra and abelian groups; for example, the issue of convergence can also be phrased as the difference between $\tau$-completion in filtered spectra and in filtered abelian groups (\cref{warn:tau_complete_filtered_vs_abelian}).
This formalism becomes particularly convenient when discussing \emph{total differentials}, which we introduce in \cref{sec:total_differentials}, and where we discuss strengthened versions of the ordinary Leibniz rule.

After introducing this more foundational setup, our goal becomes to prove the \emph{Omnibus Theorem}, showing that if $X$ is a filtered spectrum, then the $\Z[\tau]$-module $\oppi_{*,*}X$ captures all of the structure of the spectral sequence underlying $X$.
To help us prove this, we introduce the \emph{$\tau$-Bockstein spectral sequence} in \cref{sec:tau_Bockstein_sseq}.
Although it is not strictly necessary to use this, it is a helpful organisational tool, and is useful to have available in general.
With this in hand, we prove the Omnibus Theorem in \cref{sec:filtered_omnibus}.

Finally, we would like to export all of the aforementioned results and tools to other settings.
We discuss the notion of a \emph{deformation} in \cref{sec:deformations}.
For deformations arising in a special way, which we refer to as \emph{monoidal deformations}, we show that indeed everything discussed above exports directly.
Our main example of a deformation is that of \emph{synthetic spectra} and is deferred to the next chapter, but we include a few additional examples at the end of this chapter.

Many of the results in this chapter are well known to experts.
We draw from \cite{barkan_monoidal_algebraicity}, \cite[Appendices~A--C]{burklund_hahn_senger_Rmotivic}, \cite{lurie_rotation_invariance}, and \cite{pstragowski_perfect_even_filtration}.
Our proof of the Omnibus Theorem is heavily inspired by \cite[Appendix~A]{burklund_hahn_senger_manifolds}.

\section{Filtered abelian groups}
\label{sec:tau_in_fil_ab}

Recall that a module over the polynomial ring $\Z[x]$ is the same as an abelian group together with an endomorphism.
Using this, we can give a different description of filtered abelian groups, as follows.
We reserve the letter $\tau$ as a formal variable for the polynomial ring $\Z[\tau]$.
We turn this into a graded ring by giving $\tau$ degree~$-1$.
By forgetting the transition maps, a filtered abelian group has an underlying graded abelian group, and the transition maps can be viewed as a graded $\Z[\tau]$-module structure on this graded abelian group.
The forgetful functor thus factors through a functor $\FilAb \to \Mod_{\Z[\tau]}(\grAb)$.

\begin{proposition}
    The functor
    \[
        \FilAb \simeqto \Mod_{\Z[\tau]}(\grAb)
    \]
    is a symmetric monoidal equivalence, where we regard $\grAb$ as having the symmetric monoidal structure \emph{without} any signs in the swap maps.
\end{proposition}

Mathematically, there is nothing deep about this statement.
The value is in the human aspect: it can be less mentally taxing to think in terms of algebraic equations involving $\tau$, than it is to picture the diagram that is a filtered abelian group.
Even for strict filtrations this is very helpful, particularly when recording filtration jumps and hidden relations.
We give a few examples in \cref{ssec:eg_filtered_ab_tau_mod}.

Various properties of filtered abelian groups can be rephrased using $\tau$.
The ones we focus on are the following.

\begin{numberenum}
    \item A filtered abelian group is \emph{strict} if and only if the corresponding $\Z[\tau]$-module is \emph{$\tau$-power torsion free}.
    \item The \emph{associated graded} of a filtered abelian group $A$ is given by the graded abelian group $A/\tau$.
    \item The underlying object $A^{-\infty}$ of a filtered abelian group $A$ can be identified with the \emph{$\tau$-inversion} of the corresponding $\Z[\tau]$-module.
    \item A filtered abelian group is \emph{derived complete} (\cref{def:derived_complete_filtered_ab}) if and only if the corresponding $\Z[\tau]$-module is \emph{$\tau$-adically complete}.
\end{numberenum}

The first two claims are obvious; let us elaborate on the other two.

\begin{definition}
    The \defi{constant filtration functor} $\Const \colon \Ab \to \FilAb$ is given by sending an abelian group $A$ to the constant filtration on $A$, given by
    \[
        \begin{tikzcd}
            \dotsb \rar[equals] & A \rar[equals] & A \rar[equals] & A \rar[equals] & \dotsb.
        \end{tikzcd}
    \]
\end{definition}

It is easy to check that $\Const$ restricts to an equivalence from $\Ab$ onto the filtered abelian groups whose transition maps are all isomorphisms; let us say that such a filtered abelian group is \defi{constant}.
In terms of $\tau$, this says that $\tau$ acts invertibly on it.
For this reason, we may also refer to such filtered abelian groups as being \defi{$\tau$-invertible}, and we write $\FilAb[\tau^{-1}]$ for the full subcategory on these objects.

\begin{proposition}
    \label{prop:tau_invertible_fil_ab}
    The inclusion $\FilAb[\tau^{-1}]\subseteq \FilAb$ admits both a left and a right adjoint.
    Under the equivalence $\Const \colon \Ab \simeq \FilAb[\tau^{-1}]$, the adjunctions
    \[
        \begin{tikzcd}[column sep=5em]
            \FilAb \ar[r, bend left, "\tau^{-1}"] \ar[r, bend right] & \FilAb[\tau^{-1}] \ar[l, hook']
        \end{tikzcd}
        \qquad \text{become} \qquad
        \begin{tikzcd}[column sep=6em]
            \FilAb \ar[r, bend left, "A \, \mapsto \, A^{-\infty}"] \ar[r, bend right, "A\, \mapsto \, A^\infty"'] & \Ab. \ar[l, "\Const"', hook']
        \end{tikzcd}
    \]
\end{proposition}
\begin{proof}
    By definition of the limit and colimit, the functor $\Const$ admits both a left and right adjoint, being given by the colimit and limit functor, respectively.
    The claims then follow immediately from the fact that $\Const$ restricts to an equivalence $\Ab \simeq \FilAb[\tau^{-1}]$.
\end{proof}





To allow for a distinction to be made between an abelian group and its corresponding constant filtered abelian group, we introduce the following notation.

\begin{notation}
    \label{not:tau_inv_vs_tau_is_one}
    We write $\tau^{-1} \colon \FilAb \to \FilAb[\tau^{-1}]$ for the functor sending $A$ to
    \[
        A[\tau^{-1}] \defeq \colim(\begin{tikzcd}
            A \rar["\tau"] & A \rar["\tau"] & \dotsb
        \end{tikzcd}).
    \]
    We write $(\blank)^{\tau=1}\colon \FilAb \to \Ab$ for the composite
    \[
        \begin{tikzcd}
            \FilAb \rar["\tau^{-1}"] & \FilAb[\tau^{-1}] \rar["\simeq"] & \Ab.
        \end{tikzcd}
    \]
    Both of these functors are naturally symmetric monoidal.
\end{notation}

In terms of $\Z[\tau]$-modules, these functors take the following form.

\begin{variant}   
    The subcategory of $\tau$-invertible filtered abelian groups is (symmetric monoidally) equivalent to the category of graded $\Z[\tau^\pm]$-modules.
    Rephrased like this, the functors $\tau^{-1}$ and $(\blank)^{\tau=1}$ are given, respectively, by sending a graded $\Z[\tau]$-module $M$ to
    \[
        \Z[\tau^\pm ]\otimes_{\Z[\tau]} M, \qquad \text{respectively} \qquad \Z \otimes_{\Z[\tau]} M,
    \]
    where in the latter we let $\tau$ act on $\Z$ by the identity.\footnotemark\ 
    Note that these two functors are related by the (symmetric monoidal) equivalence
    \[
        \Mod_{\Z[\tau^{\pm}]}(\grAb) \simeqto \Ab
    \]
    given by evaluation at degree zero.
\end{variant}

\footnotetext{Note that in the second of these cases, we lose a grading, because letting $\tau$ act by the identity on $\Z$ does not turn $\Z$ into a \emph{graded} $\Z[\tau]$-module.
Said differently, $(\blank)^{\tau=1}$ is given by taking the quotient by $\tau-1$, which is not a homogeneous element, thus resulting in a loss of grading.}

This explains what we mean by the colimit of a filtered abelian group being the same as the $\tau$-inversion of the corresponding $\Z[\tau]$-module.
The claim about completion follows from the following.

\begin{definition}
    \label{def:adic_filtration_module}
    Let $R$ be a (commutative) ring, let $x \in R$, and let $M$ be an $R$-module.
    The \defi{$x$-adic filtration} on $M$ is the filtered $R$-module
    \[
        \begin{tikzcd}
            \dotsb \ar[r,"x"] & M \ar[r,"x"] & M \ar[r,"x"] & M \ar[r,equals] & \dotsb,
        \end{tikzcd}
    \]
    which we index to be constant from degree $0$ onwards.
    If $M$ is a commutative $R$-algebra, then this is naturally a filtered commutative $R$-algebra.
\end{definition}

\begin{proposition}
    \label{prop:when_is_adic_filtration_derived_complete}
    Let $R$ be a ring, let $x \in R$, and let $M$ be an $R$-module.
    Then the $x$-adic filtration on $M$ is derived complete if and only if $M$ is $x$-complete as an $R$-module, i.e., if the natural map
    \[
        M \to M_x^\wedge \defeq \lim_k M/x^k
    \]
    is an isomorphism.
\end{proposition}

\begin{warning}
    The above should not be confused with the notion of \emph{derived $x$-completion} (such as derived $p$-completion, a.k.a.\ L-completion; see \cite[Appendix~A]{hovey_strickland_Ktheories_localisation}) in the sense of \cite{greenlees_may_derived_completion}.
    In fact, for filtered abelian groups, derived $\tau$-completion is in general different from $\tau$-completion; see \cref{warn:tau_complete_filtered_vs_abelian} below.
\end{warning}

\begin{proof}[Proof of \cref{prop:when_is_adic_filtration_derived_complete}]
    Consider the diagram
    \[
        \begin{tikzcd}
            \vdots \ar[d] & \vdots \ar[d] & \vdots \ar[d]\\
            M \ar[d, "x"'] \ar[r,"x^3"] & M \ar[r] \ar[d,equals] & M/x^3 \ar[d] \\
            M \ar[d, "x"'] \ar[r,"x^2"] & M \ar[r] \ar[d,equals] & M/x^2 \ar[d] \\
            M \ar[r, "x"] & M \ar[r] & M/x.
        \end{tikzcd}
    \]
    Taking limits in the vertical direction, we get an exact sequence
    \[
        0 \to L \to M \to M_x^\wedge \to K \to 0
    \]
    where
    \[
        L = \lim(\begin{tikzcd} \dotsb \rar["x"] & M \rar["x"] & M \end{tikzcd}) \quad \text{and} \quad K = \derlim{1} (\begin{tikzcd} \dotsb \rar["x"] & M \rar["x"] & M \end{tikzcd} ).
    \]
    These are precisely the limit and first-derived limit of the $x$-adic filtration on $M$.
    In other words, we see that $M \to M_x^\wedge$ is an isomorphism if and only if the $x$-adic filtration on $M$ is derived complete.
\end{proof}

\begin{corollary}
    \label{cor:filtered_ab_derived_complete_iff_tau_complete}
    A filtered abelian group is derived complete \textbr{\cref{def:derived_complete_filtered_ab}} if and only if its corresponding graded $\Z[\tau]$-module is $\tau$-adically complete.
\end{corollary}
\begin{proof}
    Let $A$ be a filtered abelian group.
    Applying the previous to the ring $R=\Z[\tau]$ and the element $x = \tau$, the limit and the first-derived limit of the $\tau$-adic filtration on $A$ are the constant filtered abelian groups on the limit and first-derived limit of $A$, respectively.
\end{proof}

\begin{variant}
    \label{var:filtered_graded_abelian_groups}
    There is an obvious variant of all of the above for graded abelian groups, and all the analogous versions of the previous results hold true.
    A \defi{filtered graded abelian group} is a functor $\Z^\op \to \grAb$.
    Note that this by definition means that the transition maps preserve degrees.
    We write $\FilgrAb \defeq \Fun(\Z^\op,\grAb)$ for the category of filtered graded abelian groups.
    
    We give $\grAb$ the symmetric monoidal structure with the Koszul sign rule.
    We then regard $\FilgrAb$ as a symmetric monoidal category under Day convolution.

    We turn $\Z[\tau]$ into a bigraded ring by giving $\tau$ bidegree $(0,-1)$.
    We give $\bigrAb$ the symmetric monoidal structure with the Koszul sign rule according to the \emph{first} grading.
    This results in a symmetric monoidal equivalence
    \[
        \FilgrAb \simeqto \Mod_{\Z[\tau]}(\bigrAb)
    \]
    where the first grading is the internal grading, and the second grading is the grading arising from the filtration.
    There is a sign rule for swapping elements according to their \emph{internal grading}; the filtration does not play a role in these signs.
    This indexing is designed to fit with the indexing conventions for spectral sequences from \cref{sec:spectral_sequences}.
\end{variant}

\subsection{Examples}
\label{ssec:eg_filtered_ab_tau_mod}


Consider the $p$-adic filtration (see \cref{def:adic_filtration_module}) on an abelian group $A$.
The induced strict filtration on $A$ is the \defi{strict $p$-adic filtration}, given by
\[
    F^s = p^s A \subseteq A.
\]
Note, however, that the maps in the $p$-adic filtration itself need not be injective, as $A$ might contain $p$-torsion.
The additional information in the $p$-adic filtration is that it remembers all possible choices of $p$-divisions of elements.

\begin{example}
    \label{ex:p_adic_filtration_Z}
    The abelian group $\Z$ is $p$-torsion free for all $p$.
    The graded $\Z[\tau]$-algebra corresponding to the $p$-adic filtration on $\Z$ is
    \[
        \Z[\tau,\ptil]/(\tau \cdot \ptil = p) \qquad \text{where }\abs{\ptil} = 1.
    \]
    We think of $\ptil$ as a refinement of $p\in \Z$ that records the fact that $p$ has filtration $1$.
    
    The fact that the filtration is constant from filtration $0$ onward translates to the fact that in the $\Z[\tau]$-module, multiplication by $\tau$ is an isomorphism in degrees zero and below.
    The element $1$ is not $\tau$-divisible however, reflecting the fact that the transition map from filtration 1 to filtration 0 is not surjective.
    As with all $\Z[\tau]$-modules, the filtration of an element corresponds to the $\tau$-divisibility.

    This filtration is not $\tau$-complete however; for this, we would have to pass to $\Z_p$.
\end{example}

In \cref{rmk:subadditive_filtration}, we remarked that in a strictly filtered ring, the filtration of elements is subadditive.
The corresponding $\Z[\tau]$-algebra records this very elegantly.

\begin{example}
    \label{ex:subadditive_filtration}
    Consider the ring
    \[
        A = \Z[\eta,\nu]/ (2\eta,\ 8\nu,\ 4\nu = \eta^3).
    \]
    We give $A$ a strict filtration by letting both $\eta$ and $\nu$ be of filtration $1$, and all of $\Z$ be of filtration $0$.
    The relation $4\cdot \nu = \eta^3$ is then a jump in filtration: the product $4\cdot \nu$ is in $F^1$, but happens to land in the smaller subgroup $F^3$.
    In particular, we do not see this relation on the associated graded.

    The corresponding $\Z[\tau]$-algebra keeps track of this more clearly.
    Saying that $4\nu$ lands in filtration $3$ means that $4\nu$ is the $\tau^2$-multiple of an element in filtration $3$.
    Due to the lack of $\tau$-torsion, in this case there is a unique such element, namely $\eta^3$.
    The resulting graded $\Z[\tau]$-algebra is given by
    \[
        \Z[\tau,\eta,\nu]/ (2\eta,\ 8\nu, \ 4\nu = \tau^2 \eta^3) \qquad \text{where } \abs{\eta} = \abs{\nu} = 1.
    \]
    In an informal sense, we were forced to insert a $\tau^2$-term in the last relation: unlike filtered rings, graded rings do not allow for a grading-jump under multiplication.
    Since $\tau$ has degree $-1$, the relation $4\nu = \tau^2 \eta^3$ now respects this rule.
    Observe that if we put $\tau=1$, then this indeed recovers the original ring $A$.
    
    Note that, unlike in \cref{ex:p_adic_filtration_Z}, we do not write $\xtilde{\eta}$ or $\xtilde{\nu}$, but instead use the symbols $\eta$ and $\nu$ to directly record the filtration of the elements in the ring $A$.
    We do this because, once we fix a filtration on $A$, we think of the filtration of an element as an intrinsic property, not something to be witnessed by another element.
    We cannot use this type of notation in \cref{ex:p_adic_filtration_Z} however, because the symbol $p$ is usually reserved for $1+\dotsb+1$, and it is a bad idea to break this convention.
    (For instance, using the symbol $2$ to denote anything other than $1+1$ is not advisable.)
\end{example}

Elements that are in the kernel of a transition map now translate to elements that are $\tau$-torsion.
This will become especially important when dealing with spectral sequences: there, $\tau$-power torsion will encode the presence of \emph{differentials}.

\section{Filtered spectra}
\label{sec:tau_in_fil_sp}

By the Yoneda lemma, the natural transition map $\pi_{n,s} \to \pi_{n,\, s-1}$ is induced by a map $\S^{n,\,s-1}\to\S^{n,s}$.

\begin{definition}
    \label{def:fil_tau}
    \leavevmode
    \begin{numberenum}
        \item The map $\tau \colon \S^{0,-1} \to \S$ is the image of the morphism $-1\to 0$ in $\Z$ under the functor $i\colon \Z \to \FilSp$ from \cref{def:functor_Z_to_FilSp}.
        \item If $X$ is a filtered spectrum, then tensoring $\tau \colon \S^{0,-1} \to \S$ with $X$ results in a map $\Sigma^{0,-1}X \to X$.
        We will denote this map by $\tau_X$, or simply by $\tau$ if there is little risk of confusion.
    \end{numberenum}
\end{definition}

In a diagram, writing $\S^{0,-1}$ in the top row and $\S$ in the bottom row, the map $\tau$ looks like
\[
    \begin{tikzcd}
        \dotsb \ar[r] & 0 \ar[r] \ar[d] & 0 \ar[r] \ar[d] & \S \ar[r,equals] \ar[d,equals] & \dotsb\\
        \dotsb \ar[r] & 0 \ar[r] & \S \ar[r,equals] & \S \ar[r,equals] & \dotsb.
    \end{tikzcd}
\]

If $X$ is a filtered spectrum, then the map $\tau_X$ looks like
\[
    \begin{tikzcd}
        \dotsb \ar[r] & X^2 \ar[r,"f_1"] \ar[d,"f_1"] & X^1 \ar[r,"f_0"] \ar[d,"f_0"] & X^0 \ar[r] \ar[d,"f_{-1}"] & \dotsb\\
        \dotsb \ar[r] & X^1 \ar[r,"f_0"] & X^0 \ar[r,"f_{-1}"] & X^{-1} \ar[r] & \dotsb.
    \end{tikzcd}
\]
In words: the components of the map $\tau_X$ are the transition maps of $X$.

\begin{remark}
    \label{rmk:interpretation_of_i_as_tau_tower}
    Using the above notation, the functor $i \colon \Z \to \FilSp$ from \cref{def:functor_Z_to_FilSp} can be depicted as the diagram in $\FilSp$ given by
    \[
        \begin{tikzcd}
            \dotsb \ar[r,"\tau"] & \S^{0,-1} \ar[r,"\tau"] & \S \ar[r,"\tau"] & \S^{0,1} \ar[r,"\tau"] & \dotsb.
        \end{tikzcd}
    \]
\end{remark}

Previously, we considered the homotopy groups of a filtered spectrum as a filtered graded abelian group.
Now, we rephrase this in terms of $\Z[\tau]$-modules.

\begin{definition}
    \label{def:filtered_htpy_as_tau_module}
    We define a functor
    \[
        \pi_{*,*} \colon \FilSp \to \Mod_{\Z[\tau]}(\bigrAb), \quad X \mapsto \oppi_{*,*}X
    \]
    where the $\Z[\tau]$-module structure is given by the map $\tau$ of \cref{def:fil_tau}.
    This is naturally a lax symmetric monoidal functor, where $\bigrAb$ is given the Koszul sign rule according to the first grading.
\end{definition}

Our ultimate goal is to make precise that the bigraded $\Z[\tau]$-module $\oppi_{*,*}X$ captures the data of the spectral sequence underlying $X$.
Before we can do this, we start by reformulating the basic building blocks of spectral sequences in terms of $\tau$.
We also compare these notions to the analogous notions for $\Z[\tau]$-modules from the previous sections.
Note, however, that these do not always coincide: modding out by $\tau$ and completing at $\tau$ are different in the stable than in the abelian setting.

\subsection{Inverting \texorpdfstring{$\tau$}{tau}}
\label{ssec:filtered_inverting_tau}

\begin{definition}
    A filtered spectrum $X$ is called \defi{$\tau$-invertible} if the map $\tau \colon \opSigma^{0,-1}X \to X$ is an isomorphism.
    We write $\FilSp[\tau^{-1}]$ for the full subcategory of $\FilSp$ on the $\tau$-invertible filtered spectra.
    Write $\tau^{-1}\colon \FilSp \to \FilSp$ for the functor sending $X$ to the colimit
    \[
        X[\tau^{-1}] \defeq \colim(\begin{tikzcd} X \ar[r,"\tau"] & \Sigma^{0,1} X \ar[r,"\tau"] & \dotsb\end{tikzcd}).
    \]
\end{definition}

It is easy to check that $\tau$ acts invertibly on $X[\tau^{-1}]$, so that $\tau$-inversion lands in $\tau$-invertible filtered spectra.
Moreover, it participates in an adjunction
\[
    \begin{tikzcd}[column sep=4em]
        \FilSp \ar[r, "\tau^{-1}", shift left=1.25] & \FilSp[\tau^{-1}]. \ar[l, shift left=1.25, hook']
    \end{tikzcd}
\]

Inverting $\tau$ is a particularly good kind of localisation: it is a \emph{smashing localisation}.
We refer to \cite[Section~3]{gepner_groth_nikolaus_multiplicative} for an introduction to such localisations.
The practical upshot is that $\tau$-invertible objects are closed under limits, colimits and tensor products, and that $\tau$-local objects get an essentially unique structure of a $\S[\tau^{-1}]$-module.

\begin{proposition}
    The functor of $\tau$-localisation is a smashing localisation, i.e., it is given by tensoring with the idempotent object $\S[\tau^{-1}]$.
    In particular, the inclusion functor $\FilSp[\tau^{-1}] \subseteq \FilSp$ preserves colimits and has a further right adjoint.
\end{proposition}
\begin{proof}
    The tensor product of filtered spectra preserves colimits.
    It follows that
    \begin{align*}
        X[\tau^{-1}] &= \colim(\begin{tikzcd}[ampersand replacement=\&] X \ar[r,"\tau"] \& \Sigma^{0,1} X \ar[r,"\tau"] \& \dotsb\end{tikzcd})\\
        &\cong \colim(\begin{tikzcd}[ampersand replacement=\&] \S \ar[r,"\tau"] \& \S^{0,1} \ar[r,"\tau"] \& \dotsb\end{tikzcd}) \otimes X\\
        &= \S[\tau^{-1}]\otimes X.\qedhere
    \end{align*}
\end{proof}

The notion of a $\tau$-invertible filtered spectrum is not new: it is the same as a filtered spectrum whose transition maps are invertible.
More precisely, we have the following identifications.

\begin{proposition}
    \label{prop:identification_tau_invertible_FilSp}
    The symmetric monoidal functor $\Const \colon \Sp \to \FilSp$ restricts to an equivalence onto the $\tau$-invertible filtered spectra:
    \[
        \Const \colon \Sp \simeqto \FilSp[\tau^{-1}].
    \]
    Under this equivalence, the adjunctions
    \[
        \begin{tikzcd}[column sep=5em]
            \FilSp \ar[r, bend left, "\tau^{-1}"] \ar[r, bend right] & \FilSp[\tau^{-1}] \ar[l, hook']
        \end{tikzcd} \qquad \text{become} \qquad
        \begin{tikzcd}[column sep=6em]
            \FilSp \ar[r, bend left, "X \, \mapsto \, X^{-\infty}"] \ar[r, bend right, "X\, \mapsto \, X^\infty"'] & \Sp. \ar[l, "\Const"', hook']
        \end{tikzcd}
    \]
\end{proposition}

Paralleling \cref{not:tau_inv_vs_tau_is_one}, we will sometimes use the following notation to distinguish between the two equivalent $\infty$-categories $\Sp$ and $\FilSp[\tau^{-1}]$.
In practice however, we may refer to both functors as ``$\tau$-inversion''.

\begin{notation}
    We write $(\blank)^{\tau=1}$ for the composite
    \[
        \begin{tikzcd}
            \FilSp \ar[r,"\tau^{-1}"] & \FilSp[\tau^{-1}] \rar["\simeq"] & \Sp.
        \end{tikzcd}
    \]
\end{notation}

\begin{remark}
    The operation of $\tau$-inversion of filtered spectra is compatible with $\tau$-inversion of filtered (graded) abelian groups from \cref{sec:tau_in_fil_ab}.
    More specifically, if $X$ is a filtered spectrum, then the natural map provides an isomorphism
    \[
        (\oppi_{n,*}X)[\tau^{-1}] \congto \oppi_{n,*}(X[\tau^{-1}]),
    \]
    due to the fact that bigraded homotopy groups preserve filtered colimits (as the filtered spheres are compact).
    As a result, this also induces an isomorphism
    \[
        (\pi_{n,*}X)^{\tau=1} \cong \oppi_n (X^{\tau=1}).
    \]
\end{remark}

\begin{remark}[Detection and $\tau$-divisibility]
    Let $X$ be a filtered spectrum, and let $\theta \in \oppi_n X^{\tau=1}$ be nonzero.
    The statement that $\theta$ is detected (see \cref{def:perm_cycle_and_detection}) in the spectral sequence underlying $X$ in filtration $s$ translates to the statement that $\theta$ lifts to a non-$\tau$-divisible element $\alpha$ in $\oppi_{n,s}X$.
    Indeed, saying that $\alpha$ is not $\tau$-divisible is another way of saying that $\alpha$ does not lift to $\oppi_{n,\,s+1}X$, which is equivalent to its image in $\uE_1^{n,s}$ not being zero, which is implied by the definition of detection from \cref{def:perm_cycle_and_detection}.
\end{remark}

\begin{example}
    \label{ex:invert_tau_on_filtered_sphere}
    Recall the definition $\S^{n,s}=\opSigma^{n} i(s)$ from \cref{def:filtered_bigraded_spheres}.
    It is immediate from the definition of $i$ that we have a natural isomorphism $i(\blank)[\tau^{-1}] \cong \Const\S$.
    As $\tau$-inversion and $\Const$ are exact functors, they preserve suspensions, so we find that
    \[
        \S^{n,s}[\tau^{-1}] = \opSigma^{n}i(s)[\tau^{-1}] \cong \Sigma^{n} \Const \S \cong \Const \S^n.
    \]
    In other words, $(\S^{n,s})^{\tau=1} \cong \S^n$.
    We can think of this as saying that inverting $\tau$ forgets the filtration.
\end{example}

\subsection{Modding out by \texorpdfstring{$\tau$}{tau}}

\begin{notation}
    Let $k\geq 1$.
    \begin{itemize}
        \item We write $C\tau^k$ for the cofibre of the map $\tau^k \colon \S^{0,\, -k}\to\S$.
        \item If $X$ is a filtered spectrum, then we write $X/\tau^k$ for $C\tau^k \otimes X$.
    \end{itemize}
\end{notation}

Concretely, $C\tau^k$ is the filtered spectrum
\[
    \begin{tikzcd}
        \dotsb \rar & 0 \rar & \S \rar[equals] & \dotsb \rar[equals] & \S \rar & 0 \rar & \dotsb,
    \end{tikzcd}
\]
where the nonzero terms are in filtrations $0,-1,\dotsc,-k+1$.
If $X$ is a filtered spectrum, then $X/\tau^k$ is in filtration $s$ given by the cofibre of $X^{s+k} \to X^s$.

Unlike in the case of filtered abelian groups, in higher algebra, taking quotients in a monoidal way is a treacherous matter.
In this case, it turns out we can do this in the best possible way.

\begin{theorem}[Lurie]
    \label{thm:Ctau_Einfty}
    For every $k\geq 0$, the filtered spectrum $C\tau^k$ admits \textbr{uniquely up to contractible choice} the structure of a filtered $\E_\infty$-ring such that its unit map $\S \to C\tau^k$ is an isomorphism in filtrations $0,-1,\dotsc,-k+1$.
\end{theorem}
\begin{proof}
    See \cite[Proposition~3.2.5]{lurie_rotation_invariance}, bearing in mind that Lurie writes $\mathrm{Rep}(\Z)$ for $\FilSp$, and writes $\mathbb{A}$ for $C\tau$.
\end{proof}

Using this ring structure on $C\tau$, we can make precise the way in which tensoring with $C\tau$ recovers the associated graded.
We first define a functor
\[
    d\colon \grSp \to \FilSp
\]
given by left Kan extension along the functor $\Z^\discr \to \Z^\op$, with $\Z^\discr$ denoting the discrete category with objects $\Z$.
Informally, this functor is given by sending a graded spectrum $(X_n)_n$ to the filtered spectrum
\[
    \dotsb \to \bigoplus_{n\geq 1} X_n \to \bigoplus_{n\geq 0} X_n \to \bigoplus_{n\geq -1} X_n \to \dotsb,
\]
with maps the natural inclusions.
Being defined as the left Kan extension along a symmetric monoidal functor, this is naturally a symmetric monoidal functor.

\begin{theorem}[Lurie]
    \label{thm:mod_tau_is_associated_graded}
    The composite
    \[
        \begin{tikzcd}
            \grSp \ar[r,"d"] & \FilSp \ar[r,"C\tau\otimes \blank"] &[1.5em] \Mod_{C\tau}(\FilSp)
        \end{tikzcd}
    \]
    is a symmetric monoidal equivalence.
    Moreover, this equivalence fits into a commutative diagram
    \[
        \begin{tikzcd}
            \FilSp \ar[r,"\Gr"] \ar[dr,"C\tau\otimes\blank"'] & \grSp \ar[d,"\simeq"] \\
            & \Mod_{C\tau}(\FilSp).
        \end{tikzcd}
    \]
\end{theorem}
\begin{proof}
    See \cite[Proposition~3.2.7]{lurie_rotation_invariance}, bearing in mind that Lurie writes $\mathrm{Rep}(\Z)$ for $\FilSp$, writes $\mathrm{Rep}(\Z^{\mathrm{ds}})$ for $\grSp$, writes $\mathbb{A}$ for $C\tau$, and writes $I$ for $d$.
\end{proof}

\begin{remark}
    \label{rmk:Ctau_gives_Gr_symm_mon_str}
    The above in particular puts a symmetric monoidal structure on the associated graded functor $\Gr \colon \FilSp\to\grSp$, because the functor $C\tau \otimes \blank \colon \FilSp \to \Mod_{C\tau}(\FilSp)$ is canonically symmetric monoidal.
    One could have also done this more directly: see \cite[Section~II.1.3]{hedenlund_phd}.
\end{remark}

\begin{warning}
    Being a module over $C\tau$ is not a property, but additional structure.
    One can see this by observing that $C\tau$ is not an idempotent, i.e., that $C\tau \otimes C\tau$ is not isomorphic to $C\tau$: we instead have
    \[
        C\tau \otimes C\tau \cong C\tau \oplus \opSigma^{1,\,-1} C\tau.
    \]
    Indeed, the map $\tau$ on $C\tau$ is nullhomotopic (because $C\tau$ is a ring), so the cofibre sequence defining $C\tau$ splits after tensoring with $C\tau$.
    Alternatively, note that the associated graded of $C\tau$ is concentrated in filtrations $0$ and $-1$, where it is $\S$ and $\S^1$, respectively.
    Via the equivalence of \cref{thm:mod_tau_is_associated_graded}, this precisely gives us the above splitting.
\end{warning}

\begin{warning}
The notion of modding out by $\tau$ in the spectral setting is decidedly different from modding out by $\tau$ in the setting of filtered abelian groups.
Indeed, this is exactly the difference between a cofibre and a cokernel: the former also sees the kernel.
\end{warning}

\subsection{Completing at \texorpdfstring{$\tau$}{tau}}
\label{ssec:filtered_tau_completion}

\begin{definition}
    The functor of \defi{$\tau$-adic completion} \textbr{or \emph{$\tau$-completion} for short} is $C\tau$-localisation of $\FilSp$, i.e., inverting those maps that become an isomorphism after tensoring with $C\tau$.
\end{definition}

We write $\FilSp_\tau^{\wedge}$ for the full subcategory of $\FilSp$ on the $\tau$-adically complete filtered spectra.
This results in an adjunction
\[
    \begin{tikzcd}[column sep=4em]
        \FilSp \ar[r,"(\blank)_{\tau}^{\wedge}", shift left=1.25] & \FilSp_{\tau}^{\wedge}. \ar[l, shift left=1.25, hook']
    \end{tikzcd}
\]

This notion is, in fact, nothing but the notion of completeness previously introduced in \cref{def:complete_filtered_spectrum}.

\begin{proposition}
    \label{prop:tau_completion_filsp_concretely}
    \leavevmode
    \begin{numberenum}
        \item A filtered spectrum $X$ is $\tau$-adically complete if and only if its limit $X^\infty$ vanishes.
        \item A map of filtered spectra is a $C\tau$-equivalence if and only if it induces an isomorphism on associated graded.
        \item For every filtered spectrum $X$, the natural map $X \to \cofib(X^\infty \to X)$ is the $\tau$-completion of $X$.
    \end{numberenum}
\end{proposition}
\begin{proof}
    By definition, a filtered spectrum $X$ is $\tau$-complete if and only if $\Map(Y,X)$ is contractible for all $\tau$-invertible $Y$.
    The first claim therefore follows immediately from \cref{prop:identification_tau_invertible_FilSp}.
    The second item is immediate from \cref{thm:mod_tau_is_associated_graded}.
    For the final claim, we have to show that $\cofib(X^\infty \to X)$ is $\tau$-complete and that the map from $X$ into it is a $C\tau$-equivalence.
    The former is clear, and the latter follows from the tetrahedral axiom.
\end{proof}

The operations of inverting and completing at $\tau$ are related in the following way.

\begin{proposition}
    \label{prop:pullback_tau_complete_and_invert_tau}
    For a filtered spectrum $X$, there is a natural pullback square of lax symmetric monoidal functors
    \[
        \begin{tikzcd}
            X \ar[r] \ar[d] \pullback & X_\tau^{\wedge} \ar[d] \\
            X[\tau^{-1}] \ar[r] & (X_\tau^{\wedge})[\tau^{-1}].
        \end{tikzcd}
    \]
    In particular, a map of filtered spectra $X \to Y$ is an isomorphism if and only if the maps
    \[
        X[\tau^{-1}] \to Y[\tau^{-1}] \qquad \text{and} \qquad C\tau \otimes X \to C\tau \otimes Y
    \]
    are both an isomorphism.
\end{proposition}
\begin{proof}
    This is a standard result; see \cite[Proposition~4.1.1]{mantovani_localisations_stable}.
\end{proof}

Concretely, this says that a map of filtered spectra is an isomorphism if and only if it is an isomorphism on the colimit and on the associated graded.

Finally, we compare $\tau$-completeness of filtered spectra and filtered abelian groups.

\begin{warning}
    \label{warn:tau_complete_filtered_vs_abelian}
    If $X$ is a filtered spectrum, then $\tau$-completeness of $X$ need not imply $\tau$-completeness of the $\Z[\tau]$-module $\oppi_{*,*}X$.
    Indeed, this is part of the discussion of convergence.
    In detail: using \cref{cor:filtered_ab_derived_complete_iff_tau_complete}, we see that part~\ref{item:strong_conv_der_complete} of the definition of strong convergence in \cref{def:strong_convergence} is asking $\pi_{*,*}X$ to be $\tau$-complete.
    As explained in \cref{warn:explanation_convergence_criteria}, the $\tau$-completeness of $X$ (i.e., the vanishing of $X^\infty$) need not imply this.
    The conditional convergence criteria of Boardman from \cref{thm:conditional_convergence,rmk:whole_plane_obstruction} give the further conditions needed to go from $\tau$-completeness to $\tau$-completeness of $\pi_{*,*}X$.
\end{warning}

\section{Total differentials}
\label{sec:total_differentials}

Our next topic concerns the differentials in the spectral sequence.
Although the previous concepts involving $\tau$ are reformulations of ones from \cref{ch:filtered_spectra_sseqs}, the \emph{total differentials} to be introduced in this section did not make an appearance there.
However, this is not because one cannot phrase these without $\tau$, but because we find that~$\tau$ provides for an easier notational setup to introduce these concepts.

Many authors have used total differentials in both synthetic and motivic spectra; see \cite{burklund_hahn_senger_manifolds}, \cite{chua_E3_ASS} and \cite{isaksen_etal_motivic_ANSS_tmf}, for example.
They naturally arise in the filtered setting as well, and they are one of the big benefits of working at the filtered level: they lead to strengthened versions of the usual Leibniz rule.
Further, knowledge of total differentials allows one to deduce hidden extensions from non-hidden extensions.
We explain these applications at the end of this section.

\begin{notation}
    \label{not:tau_Bocksteins}
    Let $X$ be a filtered spectrum, and let $n \geq 1$.
    Write $\partial_n^\infty$ for the boundary map in the cofibre sequence
    \[
        \begin{tikzcd}
            \opSigma^{0,\, -n} X \ar[r,"\tau^n"] & X \ar[r] & X/\tau^n \ar[r,"\partial_n^\infty"] & \opSigma^{1,\, -n} X.
        \end{tikzcd}
    \]
    For $N \geq n$, write $\partial_n^N$ for the boundary map in
    \[
        \begin{tikzcd}
            \opSigma^{0,\, -n} X/\tau^{N-n} \ar[r,"\tau^n"] & X/\tau^{N} \ar[r] & X/\tau^n \ar[r,"\partial_n^N"] & \opSigma^{1,\,-n} X/\tau^{N-n}.
        \end{tikzcd}
    \]
    We call $\partial_1^\infty$ the \defi{total differential} on $X$, and call $\partial_1^N$ the \defi{$N$-truncated total differential} on $X$.
\end{notation}




The map $\partial_n^N$ captures information about the $d_n,\dotsc,d_{N-1}$-differentials in $X$.
Decreasing $N$ results in a loss of information, but increasing the lower index $n$ should be thought of as an increase of information: we will see that, roughly speaking, $\partial_n^N$ is only defined on elements on which the differentials $d_1,\dotsc,d_{n-1}$ vanish.


We will now make these ideas precise.
Before we begin, let us point out what these maps are concretely.
As colimits in $\FilSp$ are computed levelwise, if we evaluate the cofibre sequence defining $\partial_1^\infty$ at filtration $s$, we obtain the cofibre sequence
\[
    X^{s+1} \to X^s \to \Gr^s X \to \Sigma X^{s+1},
\]
where the first map is the transition map.
Recall from \cref{ch:informal_sseq} that the boundary map $\Gr^s \to \Sigma X^{s+1}$ is precisely the map used to define all the differentials in the underlying spectral sequence.
The other total differentials are variants on this map, and a similar diagram chase will make precise the intuition for $\partial_n^N$ we gave above.

With this in mind as our intuition, we proceed to the formal proofs.

\begin{proposition}
    \label{prop:relation_different_tau_Bocksteins}
    Let $X$ be a filtered spectrum, let $n \geq k \geq 1$, and let $\infty \geq N \geq n$.
    Then we have commutative diagrams
    \[
        \begin{tikzcd}
            X/\tau^n \ar[r,"\partial_n^\infty"] \ar[dr, "\partial_n^N"'] & \opSigma^{1,\,-n} X \ar[d]\\
            & \opSigma^{1,\,-n} X/\tau^{N-n}
        \end{tikzcd}
        \qquad \text{and} \qquad
        \begin{tikzcd}
            X/\tau^n \ar[d] \ar[r,"\partial_n^\infty"] & \opSigma^{1,\,-n} X \ar[d,"\tau^{n-k}"] \\
            X/\tau^k \ar[r,"\partial_k^\infty"] & \opSigma^{1,\,-k} X,
        \end{tikzcd}
    \]
    where the unlabelled maps are the reduction maps.
    In words: $\partial_n^N$ is the mod $\tau^{N-n}$ reduction of $\partial_n^\infty$, and $\tau^{n-k} \cdot \partial_n^\infty = \partial_k^\infty$.
\end{proposition}

\begin{proof}
    For readability, we omit the bigraded suspensions in this proof.
    We start with the commutative diagram
    \[
        \begin{tikzcd}
            X \ar[d,"\tau^{N-n}"'] \ar[r,equals] & X \ar[d,"\tau^N"] \\
            X \ar[r,"\tau^n"] & X.
        \end{tikzcd}
    \]
    Taking pushouts in the vertical direction once, and repeatedly in the horizontal direction, we arrive at a commutative diagram
    \[
        \begin{tikzcd}
            X \ar[d] \ar[r,"\tau^n"] & X \ar[d] \ar[r] & X/\tau^n \ar[d,equals] \ar[r,"\partial_n^\infty"] & X \ar[d] \\
            X/\tau^{N-n} \ar[r,"\tau^n"] & X/\tau^N \ar[r] & X/\tau^n \ar[r,"\partial_n^N"] & X/\tau^{N-n}.
        \end{tikzcd}
    \]
    The right-most square is the first claimed diagram.
    
    The second diagram comes from the commutative diagram
    \[
        \begin{tikzcd}
            X \ar[d,"\tau^{n-k}"'] \ar[r, "\tau^n"] & X \ar[d,equals] \ar[r] & X/\tau^n \ar[d] \ar[r,"\partial_n^\infty"] & X \ar[d,"\tau^{n-k}"] \\
            X \ar[r,"\tau^k"] & X \ar[r] & X/\tau^k \ar[r,"\partial_k^\infty"] & X
        \end{tikzcd}
    \]
    obtained by taking horizontal pushouts of the left-most square.
\end{proof}

Next, we relate the $\tau$-divisibility of the (truncated) total differentials to the vanishing of differentials in the underlying spectral sequence.

For the case of truncated total differentials, we run into the subtlety that there are two kinds of $\tau$-multiples in $\oppi_{*,*}X/\tau^k$.
Namely, we can either speak of the $\tau$-multiple of an element from $X/\tau^k$, or of the $\tau$-multiple of an element in $X/\tau^{k-1}$ regarded as an element of $X/\tau^k$.
The latter of these is more general.
This difference will come up a number of times, particularly when describing the structure of $\oppi_{*,*}X/\tau^k$ later in \cref{ssec:truncated_filtered_omnibus}.
We use the following notation to distinguish between these.

\begin{notation}
    \label{not:two_different_meanings_multiplication_by_tau}
    Let $X$ be a filtered spectrum, and let $k\geq m \geq 0$ be integers.
    \begin{itemize}
        \item For $\theta\in\oppi_{*,*}X/\tau^k$, we write $\tau^m \cdot \theta$ for the $\tau^m$-multiple of $\theta$ in the $\Z[\tau]$-module $\oppi_{*,*}(X/\tau^k)$.
        \item For $\theta \in \oppi_{*,*}X/\tau^{m-k}$, we write $\tau^m(\theta)$ for the image of $\theta$ under the map
        \[
            \tau^m \colon \opSigma^{0,-m}X/\tau^{k-m} \to X/\tau^k.
        \]
    \end{itemize}
\end{notation}

Both versions have their advantages and disadvantages.
The former of the two is, by definition, captured by the $\Z[\tau]$-module $\oppi_{*,*}(X/\tau^k)$, while the map $\tau^m$ in the latter of the two participates in a cofibre sequence
\[
    \opSigma^{0,-m}X/\tau^{k-m} \to X/\tau^k \to X/\tau^m,
\]
and as a result is closely tied to the truncated total differential $\partial_m^k$.
They are, however, related in the following way.

\begin{remark}
    \label{rmk:relation_two_different_tau_multiplications}
    There is a commutative diagram
    \[
        \begin{tikzcd}
            \opSigma^{0,\,-m} X/\tau^k \rar["\tau^m"] \dar & X/\tau^k \\
            \opSigma^{0,\,-m} X/\tau^{k-m} \ar[ur, "\tau^m"'] & 
        \end{tikzcd}
    \]
    where the left vertical map is the reduction.
    It follows that for $\theta \in \oppi_{*,*}X/\tau^k$, the elements $\tau^m \cdot \theta$ and $\tau^m(\theta)$ (where in the latter, we implicitly reduce $\theta$ mod $\tau^{k-m}$) are the same.
    The difference between the notations of \cref{not:two_different_meanings_multiplication_by_tau}, then, is that there could be more elements of the form $\tau^m(\alpha)$ (where $\alpha \in \oppi_{*,*}X/\tau^{k-m}$) than elements of the form $\tau^m \cdot \theta$ (where $\theta \in \oppi_{*,*}X/\tau^k$).
\end{remark}

\begin{proposition}
    \label{prop:divisibility_total_diff_and_ordinary_diff}
    Let $X$ be a filtered spectrum, and let $x \in \uE_1^{n,s}$.
    View $x$ as an element of $\oppi_{n,s}X/\tau$.
    Let $r\geq 0$.
    \begin{enumerate}[label={\upshape(1\alph*)}]
        \item \label{item:total_diff_divisible} If $\partial_1^\infty(x)$ is $\tau^r$-divisible, then $d_1(x),\dotsc,d_r(x)$ vanish.

        \item More generally, if $k > r$ and $\partial_1^k(x)$ is of the form $\tau^r(\alpha)$ for some $\alpha \in \oppi_{*,*}X/\tau^{k-r-1}$, then $d_1(x),\dotsc,d_r(x)$ vanish.
    \end{enumerate}
    \begin{enumerate}[label={\upshape(2\alph*)}]
        \item Suppose that $\alpha \in \oppi_{n-1,\, s+r+1}X$ is an element such that $\tau^r \cdot \alpha = \partial_1^\infty(x)$.
        Then the mod $\tau$ reduction of $\alpha$ is a representative for $d_{r+1}(x)$ in $\uE_1^{n-1,\, s+r+1}$.
        
        \item More generally, if $k > r+1$ and $\alpha \in \oppi_{*,*}X/\tau^{k-r-1}$ is such that $\tau^r(\alpha) = \partial_1^k(x)$, then the mod $\tau$ reduction of $\alpha$ is a representative for $d_{r+1}(x)$.
    \end{enumerate}
\end{proposition}
\begin{proof}
    We begin with the claims regarding $\partial_1^\infty$.
    These are a rephrasing of the definition of the differentials in the spectral sequence associated to $X$.
    Indeed, evaluating the cofibre sequence of filtered spectra
    \[
        \begin{tikzcd}
            \opSigma^{0,\, -r} X \ar[r,"\tau^r"] & X \ar[r] & X/\tau^r \ar[r,"\partial_r^\infty"] & \opSigma^{1,\, -r} X
        \end{tikzcd}
    \]
    at filtration $s$ is exactly the cofibre sequence of spectra
    \[
        \begin{tikzcd}
            X^{s+r} \ar[r] & X^s \ar[r] & \Gr^s X \ar[r] & \opSigma X^{s+r}.
        \end{tikzcd}
    \]
    Finally, the analogous claims for $\partial_1^k$ follow by combining the ones for $\partial_1^\infty$ with \cref{prop:relation_different_tau_Bocksteins}.
\end{proof}

Later, when we have more of an understanding of how $\oppi_{*,*}X/\tau^k$ relates to the spectral sequence, we will be able to rephrase the above result in a convenient way; see \cref{constr:map_from_mod_tau_r_to_Er,cor:diagram_differential_and_total_differential}.

\begin{warning}
    \label{warn:diff_not_imply_tau_torsion}
    The converse of either item~(1a) or~(1b) is not true in general.
    The reason is that the $r$-th differential is only well defined up to the images of shorter differentials.
    As a result, one cannot in general use $d_r(x) = 0$ to deduce that $\partial_1^\infty(x)$ is $\tau^r$-divisible if $r>1$.
    This can be done, of course, if previous differentials vanish in the appropriate range: more specifically, if bidegree $(n-1,\, s+r)$ receives no differentials of length shorter than $r$.
\end{warning}

\subsection{Applications}
\label{ssec:total_diff_applications}

The reason for going to the trouble of using total differentials is that they allow for more sophisticated differential-deduction techniques.
Particularly, it allows one to deduce longer differentials from shorter ones.

A simple example of this is to use linearity of the total differential.

\begin{proposition}[Linearity of the total differential]
    \label{prop:linearity_total_differential}
    Let $n\geq 1$, and let $n\leq N\leq \infty$.
    Let $X$ be a \textbr{left} homotopy-module over a homotopy-associative ring $R$ in $\FilSp$.
    Then the map $\oppi_{*,*}\partial_n^N$ on $X$ is $\oppi_{*,*}R$-linear.
    If $N<\infty$, then the map $\oppi_{*,*}\partial_n^N$ is also $\oppi_{*,*}R/\tau^N$-linear.
\end{proposition}
\begin{proof}
    This is immediate.
\end{proof}

A more powerful version of this is the following result, due to Burklund \cite[Chapter~3]{burklund_cookware_draft}.
It is also referred to as the \emph{synthetic Leibniz rule}, but as we will see, the synthetic version follows directly from the filtered version.

\begin{theorem}[Filtered Leibniz rule, Burklund]\label{syntheticleibnizrule}
    Let $R$ be a homotopy-associative ring in $\FilSp$.
    For any $n\geq 1$, the map
    \[
        \partial_n^{2n} \colon \oppi_{*,*}(R/\tau^n) \to \oppi_{*-1,\, *+n} (R/\tau^n)
    \]
    is a derivation.
    In particular, for any two classes $x,y \in \oppi_{\ast,\ast}(R/\tau^n)$, we have the relation
    \[
        \partial_n^{2n}(xy) = \partial_n^{2n}(x)\cdot y + (-1)^{\abs{x}}x\cdot  \partial_{n}^{2n}(y).
    \]
\end{theorem}

\begin{proof}
    See \cite[Theorem~2.34]{CDvN_part2}.
\end{proof}

Because mod $\tau$ reduction is a ring map, it follows from this that $\partial_n^{n+1}$ is also a derivation.
This recovers the ordinary Leibniz rule for $d_n$, or more precisely, for a first-page representative of $d_n$.


Finally, we discuss how knowledge of total differentials is, to some extent, the same as hidden extensions.
This was used in, e.g., \cite[Proof of Proposition~A.20\,(14)]{burklund_hahn_senger_manifolds} and \cite[Remark~4.1.3]{marek_nu_tmf}, and is explained in detail in \cite[Method~2.17, Example~2.18, Proposition~4.5]{isaksen_etal_motivic_ANSS_tmf}.

\begin{remark}[Total differentials and hidden extensions]
    \label{rmk:total_diffs_and_hidden_extensions}
    Fix a filtered homotopy-ring spectrum $X$.
    Let $x,y\in \uE_1^{*,*}$ be elements, and write $\alpha = \partial_1^\infty(x)$ and $\beta = \partial_1^\infty(y)$.
    In other words, we have a differential on $x$ hitting (the mod $\tau$ reduction of) $\alpha$, and likewise from $y$ to (the mod $\tau$ reduction of) $\beta$, possibly of different lengths.
    We describe a technique for deducing (possibly hidden) extensions between $\alpha$ and $\beta$ from a (non-hidden) extension between $x$ and $y$ on the $\uE_1$-page.
    
    Let $\theta \in \oppi_{*,*}X$ be another element, write $t \in \uE_1^{*,*}$ for its mod $\tau$ reduction, and suppose that we have a multiplicative relation in $\uE_1^{*,*}$ given by
    \[
        t\cdot x = y.
    \]
    Using linearity of $\partial_1^\infty$ over $\oppi_{*,*}X$, we find that
    \[
        \beta = \partial_1^\infty(y) = \partial_1^\infty(t\cdot x) = \partial_1^\infty(\theta \cdot x) = \theta \cdot \partial_1^\infty(x) = \theta \cdot \alpha.
    \]
    In words: the extension $t\cdot x = y$ allows us to deduce the extension $\theta \cdot \alpha = \beta$, provided that we know that $\alpha$ and $\beta$ are the total differentials on $x$ and $y$, respectively.
    In fact, reading the argument backwards, we see that knowing the relation $\theta\cdot \alpha = \beta$ and the total differential $\partial_1^\infty(x) = \alpha$ allows us to deduce the total differential on $y$.

    The resulting extension $\theta\cdot \alpha = \beta$ corresponds to a hidden extension when the differential on $y$ is longer than the differential on $x$.
    Indeed, set up properly, the length of the differential corresponds to the $\tau$-divisibility of the total differential.
    If $\alpha = \tau^r \alpha'$ and $\beta = \tau^s \beta'$, where $\alpha'$ and $\beta'$ are not $\tau$-divisible, and $0<r < s$, then we learn that
    \[
        \theta \cdot \alpha' = \tau^{s-r} \beta' \qquad \text{up to $\tau^r$-divisible elements.}
    \]
    Modulo $\tau$, this relation is hidden: the right-hand side reduces to zero modulo~$\tau$.
    Yet, the relation $t \cdot x = y$ is \emph{not} a hidden extension, so that the knowledge of total differentials allows us to turn non-hidden extensions into hidden extensions between the targets of the differentials.

    The only downside of this approach is that it requires $\alpha$ and $\beta$ to be images of a total differential, or equivalently, that they have to be $\tau$-power torsion (in other words, they come from classes that are hit by differentials).
    Such classes define the zero element in $\pi_*X^{-\infty}$, making it seem like this method has limited use for deducing extensions between elements in $\pi_*X^{-\infty}$.
    In practice, one can get around this using the following trick.
    If we start with non-$\tau$-power torsion elements $\alpha$ and $\beta$ we would like to find a relation between, we may be able to find another element $\gamma$ such that $\alpha \cdot \gamma$ and $\beta \cdot \gamma$ are $\tau$-power torsion.
    By exactness this must mean they are in the image of a total differential.
    If the sources of these total differentials are related by a multiplicative extension on the $\uE_1$-page, then (up to dividing by $\gamma$) we learn about a (possibly hidden) extension between $\alpha$ and $\beta$.
    
    This method can be modified with other total differentials in the place of $\partial_1^\infty$.
    If we work with $\partial_1^N$, then we can learn about hidden extensions of length at most $N-1$.
    These truncated total differentials may be easier to obtain, giving it a practical upshot at the cost of a less powerful outcome.
    On the other hand, we can also work with $\partial_n^\infty$ (or even $\partial_n^N$).
    The key difference is that there, we can take multiplicative relations $X/\tau^n$ as input, which may themselves correspond to hidden extensions of length smaller than $n$.
    In this way, this technique becomes a method to turn hidden extensions into hidden extensions of a possibly greater length.
\end{remark}

\begin{example}
    This is a slightly simplified example of a hidden extension in $\Tmf$ that is shown to hold in \cite{CDvN_part2}.
    In $\Tmf$, there is a $2$-extension in the $110$-stem:
    \begin{equation}
        \label{eq:example_hidden_ext}
        2 \cdot \kappa_4 = \eta_1^2 \kappabar^3.
    \end{equation}
    In the DSS for $\Tmf$, this is a hidden extension, with filtration jump $12$: the element $\kappa_4$ has filtration $2$, the element $\eta_1$ has filtration $1$, and the element $\kappabar$ has filtration~$4$.
    Working in the filtered spectrum giving rise to this DSS, we have canonical lifts of these elements to their respective filtrations, which we will denote by the same name.
    Although neither $\kappa_4$ nor $\eta_1^2\kappabar^3$ is $\tau$-power torsion, they become so after multiplying with $\kappabar^3$: we even have the total differentials
    \[
        \partial_4^{28}(2\nu\Delta^7) = \tau^8 \kappa_4 \kappabar^3 \qquad \text{and} \qquad \partial_4^{28}(\tau^2 \eta^3 \Delta^7) = \tau^{20} \, \eta_1^2\, \kappabar^6.
    \]
    The hidden extension $4\nu = \tau^2\eta^3$, which is of length $2$, now stretches to a length $12$ extension:
    \[
        2 \cdot \tau^8 \,\kappa_4 \, \kappabar^3 = \partial_4^{28}(4 \nu \Delta^7) = \partial_4^{28}(\tau^2 \eta^3 \Delta^7) = \tau^{20}\, \eta_1^2 \, \kappabar^6.
    \]
    From this, the desired extension \eqref{eq:example_hidden_ext} follows.
\end{example}

%% file: fil/taubss.tex
\section{Digression: the \texorpdfstring{$\tau$}{tau}-Bockstein spectral sequence}
\label{sec:tau_Bockstein_sseq}

To prove the Omnibus Theorem, we will use the \emph{$\tau$-Bockstein spectral sequence} of a filtered spectrum.
This is a spectral sequence that compute the bigraded homotopy of a filtered spectrum $X$.
Its usefulness is due to the fact that the $\tau$-Bockstein differentials can be identified with the differentials in the spectral sequence underlying~$X$.

This is one of the more technical sections in these notes.
In this text, we only need this spectral sequence for the proof of the Omnibus Theorem in the next section, so the reader who is willing to take the statements of the Omnibus Theorem on faith does not need to read this section.
Although one could have instead phrased that proof with a more hands-on argument, the $\tau$-Bockstein spectral sequence is a useful device in and of itself, so we thought it worthwhile to give an account of this spectral sequence.


\begin{remark}
    The setup of the $\tau$-Bockstein spectral sequence is, in some sense, the most general setup of a (stable) Bockstein spectral sequence.
    Through the use of deformations of \cref{sec:deformations} (specifically \cref{thm:universal_property_FilSp}), one can recover other Bockstein spectral sequences; see \cref{ex:recover_p_BSS} for the case of the $p$-Bockstein spectral sequence of spectra, for instance.
\end{remark}

\begin{definition}
    Let $X$ be a filtered spectrum.
    The \defi{$\tau$-Bockstein filtration} on $X$, denoted $\BF_\tau X$, is the bifiltered spectrum $\Z^\op \to \FilSp$ given by
    \[
        \begin{tikzcd}
            \dotsb \ar[r,"\tau"] & \opSigma^{0,-2} X \ar[r,"\tau"] & \opSigma^{0,-1}X \ar[r,"\tau"] & X \ar[r,equals] & \dotsb
        \end{tikzcd}
    \]
    indexed to be constant from filtration $0$ onwards.
\end{definition}

\begin{remark}
    It is by design that the filtration becomes constant from filtration $0$ onwards.
    This way, its underlying object is $X$, and consequently this filtrations helps us understand $\oppi_{*,*}X$.
    If we had continued the pattern of placing $\tau$'s going all the way to the right, then the underlying object would be $X[\tau^{-1}]$, and in fact, we would learn no more than if we had used the underlying spectral sequence of $X$ directly.
\end{remark}

This leads to a spectral sequence in the same way as for singly filtered spectra, except that the indexing is more involved.
To help us index it, we use the philosophy from \cref{rmk:motivation_Adams_indexing_good}.
First of all, this means that the filtration variable $s$ records the location in the diagram as depicted above.
Second, homotopy groups in $\FilSp$ are naturally bigraded, so we would like to understand $\pi_{n,w}$ of the colimit.
Accordingly, we apply $\pi_{n,w}$ to the above diagram, and study the resulting behaviour using the long exact sequences involving the associated graded.

To avoid potential confusion regarding indexing, we discuss this spectral sequence in the case of a general bifiltered spectrum.
In this text, we will only apply this in the case of (variants of) the $\tau$-Bockstein filtration.

\begin{construction}
    \label{constr:trigraded_sseq_bifiltered_spectrum}
    Let $X \colon \Z^\op \to \FilSp$ be a bifiltered spectrum.
    We define a trigraded exact couple
    \[
        \uA^{n,w,s} = \pi_{n,w}(X^s) \qquad \text{and} \qquad \uE^{n,w,s} = \pi_{n,w}(\Gr^s X),
    \]
    fitting into the following diagram, where each map is annotated by its tridegree.
    \[
        \begin{tikzcd}
            \uA^{n,w,s} \ar[rr,"(0{,}\, 0{,}\, -1)"] & & \uA^{n,w,s} \ar[dl,"(0{,}\,0{,}\,0)"]\\
             & \uE^{n,w,s} \ar[ul,"(-1{,}\,0{,}\,1)"]
        \end{tikzcd}
    \]
    In this indexing, the differential $d_r$ has tridegree $(-1,0,r)$: indeed, a $d_r$-differential is given by applying a boundary map of degree $(-1,0,1)$ once, lifting against a map of degree $(0,0,-1)$ a total of $r-1$ times, and then projecting down by a map of degree zero.
    The resulting spectral sequence is of signature
    \[
        \uE_1^{n,w,s} = \oppi_{n,w}(\Gr^s X) \implies \oppi_{n,w}X^{-\infty}.
    \]
    The induced strict filtration on $\oppi_{n,w} X^{-\infty}$ is
    \[
        F^s\oppi_{n,w} X^{-\infty} = \im(\oppi_{n,w}X^s \to \oppi_{n,w}X^{-\infty}).
    \]
\end{construction}

\begin{example}
    Let $X$ be a filtered spectrum.
    We refer to the trigraded spectral sequence underlying $\BF_\tau X$ as the \defi{$\tau$-Bockstein spectral sequence} ($\tau$-BSS) of $X$.
    Plugging in the definitions, we see that
    \[
        \pi_{n,w}(\BF_\tau^s X) = \pi_{n,w}(\opSigma^{0,-s}X) = \oppi_{n,\, w+s}X \qquad \text{and} \qquad \Gr^s (\BF_\tau X) = \opSigma^{0,-s} X/\tau.
    \]
    In particular, the first page is of the form
    \begin{equation}
        \label{eq:identification_E1_filtered_tau_BSS}
        \uE_1^{n,w,s} = \begin{cases}
            \ \pi_{n,\, w+s} (X/\tau) & \text{if }s\geq 0,\\
            \ 0 & \text{else.}
        \end{cases}
    \end{equation}
    The induced strict filtration on $\oppi_{*,*}X$ is
    \[
        F^s \oppi_{n,w}X = \im(\tau^s \colon \oppi_{n,\, w+s}X \to \oppi_{n,w}X).\qedhere
    \]
\end{example}

It will be useful to think of an element in filtration $s$ in the $\tau$-BSS as a formal $\tau^s$-multiple.

\begin{notation}
    \label{not:tau_tilde}
    We define a formal variable $\taubar$ to have tridegree $(0,-1,1)$.
    Let $X$ be a filtered spectrum, and $\set{\uE_r^{*,*,*}}$ denote its $\tau$-BSS.
    We put a $\Z[\taubar]$-module structure on $\uE_1^{*,*,*}$ in such a way that the isomorphism \eqref{eq:identification_E1_filtered_tau_BSS} becomes an isomorphism trigraded $\Z[\taubar]$-modules
    \[
        \uE_1^{*,*,*} \cong \Z[\taubar] \otimes \pi_{*,*}(X/\tau),
    \]
    where $\oppi_{n,w}(X/\tau)$ is placed in tridegree $(n,w,0)$.
\end{notation}

\begin{remark}
    Even though the $\tau$-BSS is a trigraded spectral sequence, one can to a certain extent depict it as a bigraded one, as follows.
    Since the $d_r^\tau$-differential has tridegree $(-1,0,s)$, by fixing a constant value for $w$, we get a bigraded spectral sequence trying to converge to $\oppi_{*,w}X$.
    An element in filtration $s$ will be a formal $\tau^s$-multiple, the only catch being that the class it is a formal multiple of lives in a spectral sequence for a different $w$-value (namely $w+s$).
\end{remark}

The behaviour of the $\tau$-Bockstein spectral sequence is the same as any Bockstein spectral sequence: differentials capture $\tau$-power torsion in $\oppi_{*,*}X$, where a differentials of length $r$ corresponds to $\tau^r$-torsion.
What is special to the $\tau$-BSS of $X$ is that the differentials are exactly the differentials in the spectral sequence underlying~$X$.
More specifically, a differential $d_r(x)=y$ in the spectral sequence underlying $X$ corresponds to a $\tau$-Bockstein differential
\[
    d_r^\tau(x) = \taubar^r \cdot y.
\]
The insertion of this $\taubar$-power means that, instead of killing elements directly, the differentials in the spectral sequence underlying $X$ introduce $\tau$-power torsion in $\oppi_{*,*}X$.
Together with $\taubar$-linearity of the $d_r^\tau$-differentials, this determines all differentials in  the $\tau$-BSS.

Stating this in a precise way requires some care.
For the sake of completeness, we prove this here in detail.
We base the statement of this result on \cite{palmieri_bockstein_sseq} and \cite[Appendix~A]{burklund_hahn_senger_manifolds}.
The full statement about $\tau$-power torsion will be given by the Omnibus Theorem of the next section.



In the statement of the following theorem, all claims about differentials should be interpreted as up to boundaries of lower differentials.
That is, we allow ourselves to speak of a differential $d_r(x)=y$ where $x$ and $y$ are elements of $\uE_1$ that are $d_{\leq r-1}$-cycles, implicitly considering these elements as defining classes in $\uE_r$.

\begin{theorem}
    \label{thm:tau_BSS_structure}
    Let $X$ be a filtered spectrum.
    Let $\set{\uE_r^{*,*},d_r}$ denote the spectral sequence underlying $X$, and let $\set{\uE_r^{*,*,*},d_r^\tau}$ denote the $\tau$-Bockstein spectral sequence of $X$.
    \begin{numberenum}
        \item \label{item:E1_tau_BSS} There is a natural isomorphism of trigraded $\Z[\taubar]$-modules
        \[
            \uE_1^{*,*,*} \cong \Z[\taubar] \otimes \pi_{*,*}(X/\tau) = \Z[\taubar] \otimes \uE_1^{*,*},
        \]
        where $\pi_{n,w}(X/\tau)$ is placed in tridegree $(n,w,0)$ and $\taubar$ has tridegree $(0,-1,1)$.
        \item \label{item:differentials_tau_linear} The differentials are $\taubar$-linear and $\taubar$-divisible: for $x \in \uE_1^{n,w,s}$ and $y\in \uE_1^{n-1,\, w,\,s+r}$, there is a differential
        \[
            d_r^\tau(x) = y,
        \]
        if and only if for any \textbr{hence all} $m\geq 0$, there is a differential
        \[
            d_r^\tau(\taubar^m x) = \taubar^m y.
        \]
        In particular, the $\Z[\taubar]$-module structure on the first page induces a $\Z[\taubar]$-module structure on later pages.
        \item \label{item:differentials_hit_tau_multiples} The target of a $d_r^\tau$-differential is a $\taubar^r$-multiple: for every $x \in \uE_r^{n,w,s}$, there is an element $y \in \uE_r^{n-1,\, w+r,\,s}$ such that
        \[
            d_r^\tau(x) = \taubar^r \cdot y,
        \]
        where the multiplication denotes the $\Z[\taubar]$-module structure on $\uE_r^{*,*,*}$ from \ref{item:differentials_tau_linear}.
        \item \label{item:tau_BSS_diff_are_underlying_diff}
        If $x\in \uE_1^{n,w}$ and $y\in\uE_{r}^{n-1,\, w+r}$, then there is a differential
        \begin{align*}
            d_r(x) &= y\\
        \shortintertext{if and only if there is a differential}
            d_r^\tau(x) &= \taubar^r \cdot y.
        \end{align*}
        \item \label{item:tau_BSS_detection} Suppose that $x \in \uE_1^{n,w,s}$ detects an element $\theta \in \oppi_{n,w}X$.
        Then $\theta$ is $\tau^s$-divisible.
        Moreover, for every $m\geq 0$, the class $\taubar^m \cdot x$ detects $\tau^m \cdot \theta$.
        
        \item \label{item:tau_BSS_convergence} The $\tau$-Bockstein spectral sequence for $X$ converges conditionally to $\oppi_{*,*}X$ if and only if the spectral sequence underlying $X$ converges conditionally to $\pi_*X^{-\infty}$.
        If this is the case, then the $\tau$-Bockstein spectral sequence converges strongly if and only if $\uR\uE_\infty^{n,s}$ \textbr{the derived $\infty$-term for the spectral sequence underlying $X$} vanishes for all $n,s$.
    \end{numberenum}
\end{theorem}

Here, \cref{item:differentials_tau_linear,item:tau_BSS_diff_are_underlying_diff} should be read inductively, in the following way.
\begin{itemize}
    \item \Cref{item:E1_tau_BSS} provides a $\Z[\taubar]$-module structure on the first page, so that the statement of $\taubar$-linearity of the $d_1^\tau$-differential is well defined.
    As the second page is the homology of the first page, this endows the second page with a natural $\Z[\taubar]$-module structure, so that we can talk about $\taubar$-linearity of $d_2^\tau$-differentials, and so forth.
    \item Using \cref{item:E1_tau_BSS}, the comparison of $d_1$ with $d_1^\tau$ makes sense.
    Once this is established, it follows that the isomorphism from \cref{item:E1_tau_BSS} induces an isomorphism
    \[
        \uE_2^{n,w,s}\cong \uE_2^{n,\,w+s}
    \]
    for all $s\geq 1$.
    Since the $d_2^\tau$-differential only hits filtrations $s\geq 2$, we can use these isomorphisms to make sense of the comparison between $d_2$ and $d_2^\tau$, and so forth.
\end{itemize}

\begin{proof}
    \Cref{item:E1_tau_BSS} is true by definition of $\taubar$ in \cref{not:tau_tilde}.
    
    \Cref{item:differentials_tau_linear} is clear after unwrapping definitions.
    Let us do this in the case $m=1$, which is sufficient to prove the entire statement.
    If $x$ is an element of tridegree $(n,w,s)$, this means that $x$ is an element of
    \[
        \oppi_{n,w}(\Gr^s \BF_\tau X) = \oppi_{n,\,w+s}X/\tau.
    \]
    The meaning of a differential $d_r^\tau(x) = y$ is as follows: there is an element $\alpha$ in $\oppi_{n,w}\BF_\tau^{s+r} X=\oppi_{n,\,w+s+r}X$ such that
    \[
        \tau^{r-1}\cdot \alpha = \partial_1^\infty(x) \quad \text{in } \oppi_{n,w} \BF_\tau^{s+1} X = \oppi_{n,\,w+s+1}X,
    \]
    and such that the mod $\tau$ reduction of $\alpha$ is equal to $y$ (up to boundaries of shorter differentials).
    The element $\taubar\cdot x$ is given by considering $x$ as an element in
    \[
        \oppi_{n,w-1}(\Gr^{s+1}\BF_\tau X) \cong \oppi_{n,\, w-1+s+1} X/\tau = \oppi_{n,w}X/\tau,
    \]
    and the differential on it is calculated in exactly the same way.
    This proves the `only if' statement.
    To also deduce the `if' statement, it suffices to observe that by induction on $r$, the shorter differentials entering tridegree $(n,w,s+r)$ are $\taubar$-power multiples of differentials originating from filtration $0$.
    As a result, both the differential $d_r^\tau(x)=y$ and $d_r^\tau(\taubar x) = \taubar y$ have the same boundaries as their indeterminacy, so the claim follows.
    
    
    \Cref{item:differentials_hit_tau_multiples} will follow from \cref{item:differentials_tau_linear}, combined with the following claim: for every $r$, $n$, $w$, and $s$, multiplication by $\taubar$ induces a surjection
    \[
        \taubar \colon \uE_r^{n,w,s} \twoheadto \uE_{r}^{n,\, w-1,\, s+1}.
    \]
    This claim, in turn, we prove by induction on $r$.
    For $r=1$ it is clear from \cref{item:E1_tau_BSS}.
    Assume that for some $r\geq 1$, multiplication by $\taubar$ induces a surjection $\uE_r^{n,w,s} \twoheadto \uE_r^{n,\, w-1,\, s+1}$.
    Write $(\ker d_r)^{n,w,s}$ for the kernel of the $d_r$-differential out of $\uE_r^{n,w,s}$, and write $(\im d_r)^{n,w,s}$ for the image of $d_r$ into $\uE_r^{n,w,s}$.
    Then we have
    \[
        \uE_{r+1}^{n,w,s} = \frac{(\ker d_r)^{n,w,s}}{(\im d_r)^{n,w,s}} \qquad \text{and} \qquad \uE_{r+1}^{n,\, w-1,\, s+1} = \frac{(\ker d_r)^{n,\, w-1,\,s+1}}{(\im d_r)^{n,\, w-1,\, s+1}}.
    \]
    The $\taubar$-divisibility of the differentials from~\ref{item:differentials_tau_linear} implies that multiplication by $\taubar$ restricts to a surjection $(\ker d_r)^{n,w,s} \to (\ker d_r)^{n,\, w-1,\, s+1}$, which implies the desired statement.

    To prove item~\ref{item:tau_BSS_diff_are_underlying_diff}, it suffices to prove this for one fixed value of $w=w_0$ at a time.
    The exact couple computing the differentials going out of filtration $w$ is also computed by the exact couple associated with the truncated version of $X$:
    \[
        \begin{tikzcd}
            \dotsb \rar &  X^{w_0+2} \rar & X^{w_0+1} \rar & X^{w_0} \rar[equals] & X^{w_0} \rar[equals] & \dotsb,
        \end{tikzcd}
    \]
    where $X^{w_0}$ is placed in filtration $w_0$.
    Let us denote this filtered spectrum by $Y$.
    (In a picture, the spectral sequence underlying $Y$ is obtained from the one associated to $X$ by removing the elements in filtration strictly below $w_0$.)
    Next, we observe that by levelwise evaluating $\BF_\tau X\colon \Z^\op \to \FilSp$ at degree $w_0$, we obtain $\opSigma^{0,-w_0}Y$.
    As a result, the functor of levelwise evaluating at degree $w_0$ induces an isomorphism of exact couples
    \[
        \uA^{n,\,w_0,\,s}(\BF_\tau X) \congto \uA^{n,\,w_0+s}(Y) \qquad \uE^{n,\,w_0,\,s}(\BF_\tau X) \congto \uE^{n,\,w_0+s}(Y)
    \]
    where $s\geq 0$.
    This identifies the $d_r^\tau$-differentials with the $d_r$-differentials going out of in filtration $w_0$.

    
    
    For \cref{item:tau_BSS_detection}, recall from \cref{def:perm_cycle_and_detection} that $x$ detecting $\theta$ means that there is a lift $\alpha \in \oppi_{n,w}\BF_\tau^s X = \oppi_{n,\, w+s}X$ of $x$ that maps to $\theta$ under $\oppi_{n,w}\BF_\tau^s X \to \oppi_{n,w}\BF_\tau^0 X$.
    This transition map is given by multiplication by $\tau^s$, so this is saying that $\tau^s \cdot \alpha = \theta$.
    This proves the first clause.
    We prove the second clause for $m=1$, from which the general case follows by iterating.
    Using the $\E_\infty$-structure on $C\tau$ from \cref{thm:Ctau_Einfty}, it follows that multiplication by $\tau$ is zero on $\oppi_{*,*}(X/\tau)$, so that $\tau\cdot \theta$ is detected in filtration at least $s+1$.
    Next, the element $\alpha$, considered as an element of $\oppi_{n,\,w-1}\BF_\tau^{s+1}X = \oppi_{n,\,w+s}X$, is a lift of $\taubar \cdot x \in \uE_1^{n,\,w-1,\,s+1}$.
    Clearly, the element $\alpha$ in $\oppi_{n,\,w-1}\BF_\tau^{s+1}X$ maps to $\tau^{s+1}\cdot \alpha = \tau\cdot \theta$ under $\tau^{s+1} \colon \oppi_{n,\,w+s}X\to \oppi_{n,w-1}X$, so indeed $\taubar\cdot x$ detects $\tau\cdot \theta$.
    
    Finally, we address the claims regarding convergence.
    Conditional convergence of the $\tau$-Bockstein spectral sequence is the claim that the limit of $\BF_\tau X$ vanishes.
    By \cref{prop:identification_tau_invertible_FilSp}, this limit is isomorphic the constant spectrum on $X^\infty$, proving the claim about conditional convergence.
    Suppose now that both spectral sequences converge conditionally.
    By the trigraded analogue\footnotemark\ of \cref{thm:conditional_convergence}, the $\tau$-Bockstein spectral sequence converges strongly if and only if $\uR\uE_\infty^{n,w,s}$ vanishes for all $n,w,s$.
    Using the isomorphism from \ref{item:tau_BSS_diff_are_underlying_diff}, we find that (for $s\geq 0$)
    \[
        \uR\uE_\infty^{n,w,s} = \derlim{1}_r \uZ_r^{n,w,s} \cong \derlim{1}_r \uZ_r^{n,\,w+s} = \uR\uE_\infty^{n,\,w+s}.\qedhere
    \]
    \footnotetext{As explained in \cite[Remark~1.2.4]{hedenlund_phd}, one needs certain requirements on an abelian category $\calA$ to ensure that Boardman's arguments apply to spectral sequences valued in $\calA$.
    Here, we are working with spectral sequences valued in bigraded abelian groups, where these conditions are certainly met, and Boardman's arguments go through without any change (merely tagging on an additional grading).}
\end{proof}

\begin{remark}
    \label{rmk:cycles_boundaries_tau_BSS}
    A different way of stating (part of) \cref{item:differentials_hit_tau_multiples,item:tau_BSS_diff_are_underlying_diff} above is to say that the isomorphism of \cref{item:E1_tau_BSS} induces isomorphisms (for $s\geq 0$)
    \[
        \uZ_r^{n,w,s} \cong \uZ_r^{n,\,w+s} \qquad \uB_r^{n,w,s} \cong \uB_{\min(r,s)}^{n,\, w+s} \qquad \uE_{r+1}^{n,w,s} \cong \uZ_r^{n,\,w+s}/\uB_{\min(r,s)}^{n,\,w+s}.
    \]
    In particular, we have isomorphisms
    \[
        \uZ_\infty^{n,w,s} \cong \uZ_\infty^{n,\,w+s} \qquad \uB_\infty^{n,w,s} \cong \uB_{s}^{n,\, w+s} \qquad \uE_\infty^{n,w,s} \cong \uZ_\infty^{n,\,w+s}/\uB_{s}^{n,\,w+s}.
    \]
    The appearance of the $s$-boundaries instead of the $\infty$-boundaries in the last expression is the reason that a $d_r$-differential in the spectral sequence underlying $X$ leads to $\tau^r$-power torsion in $\oppi_{*,*}X$.
\end{remark}

\begin{remark}[Convergence of the $\tau$-BSS]
    Notably, in part~\ref{item:tau_BSS_convergence}, we do not need to assume strong convergence of the spectral sequence underlying $X$ to ensure strong convergence of the $\tau$-BSS for $X$.
    This distinction is relevant in the case where $X$ is neither left nor right concentrated, in which case additionally Boardman's whole-plane obstruction $W$ from \cref{rmk:whole_plane_obstruction} needs to vanish to ensure strong convergence of the spectral sequence underlying $X$.
    The reason this does not appear for the $\tau$-BSS is explained by \cref{warn:explanation_convergence_criteria}.
    Indeed, the group $W$ is the obstruction to a limit commuting with a colimit, but the abutment of the $\tau$-BSS is $\oppi_{n,w}X$, where we have not yet taken the colimit over $w$ (which would result in $\pi_n X^{-\infty}$).
    Of course, to be able to study $\pi_*X^{-\infty}$ using $\oppi_{*,*}X$, one would then need $W$ to vanish.
\end{remark}

\begin{remark}[Second-page indexing]
    \label{rmk:second_page_indexing_tau_BSS}
    As per usual up till this point, we have used first-page indexing for the $\tau$-Bockstein spectral sequence.
    Even in situations where one indexes filtered spectra using second-page indexing, there is something to be said for using first-page indexing for the $\tau$-BSS: only in this indexing does filtration correspond to $\tau$-power divisibility.
    In that case, a $d_r^\tau$-differential would correspond to a $d_{r+1}$-differential.
    (See \cref{var:tau_BSS_synthetic} for an example.)
    If desired however, second-page indexing for the $\tau$-BSS can be achieved via
    \[
        \widetilde{\uE}_{r+1}^{n,w,s} \defeq \uE_r^{n,\,w,\, s+w}.
    \]
    In this indexing, $\taubar$ has tridegree $(0,-1,0)$, and the non-$\taubar$-divisible groups are located in tridegrees of the form $(n,w,w)$.
\end{remark}

\begin{warning}
    The $\taubar$-divisions appearing in \cref{item:differentials_tau_linear,item:differentials_hit_tau_multiples} are not necessarily unique, due to the $\taubar$-torsion caused by shorter differentials.
\end{warning}

\subsection{The truncated \texorpdfstring{$\tau$}{tau}-Bockstein spectral sequence}
\label{ssec:truncated_tau_Bockstein_sseq}

There is a truncated version of the $\tau$-Bockstein spectral sequence, which instead computes the bigraded homotopy groups of $X/\tau^k$.
One could, of course, apply the $\tau$-BSS directly to $X/\tau^k$, but it is more efficient to use the following modification of the $\tau$-BSS.

\begin{definition}
    \label{def:truncated_tau_BSS}
    Let $X$ be a filtered spectrum and let $k\geq 1$.
    The \defi{$k$-truncated $\tau$-Bockstein filtration} on $X$ is the bifiltered spectrum $\tr_k\BF_\tau X$ given in nonnegative filtrations by
    \[
        \begin{tikzcd}
            \dotsb \ar[r] & 0 \ar[r] & \opSigma^{0,\,-k+1} X/\tau \ar[r,"\tau"] & \opSigma^{0,\,-k+2} X/\tau^2 \ar[r,"\tau"] & \dotsb \ar[r,"\tau"] & X/\tau^k
        \end{tikzcd}
    \]
    and indexed to be constant from filtration $0$ onwards.
    The resulting spectral sequence we call the \defi{$k$-truncated $\tau$-Bockstein spectral sequence} for $X$.
\end{definition}

\begin{remark}
    Upon evaluation at a fixed filtration, we obtain the filtered spectrum described by \cite[Construction~3.17]{antieau_decalage}.
\end{remark}

Rather than having to re-prove the analogous version of \cref{thm:tau_BSS_structure} from the ground up, we can instead deduce this from \cref{thm:tau_BSS_structure} by means of the following map.


\begin{construction}
    \label{constr:maps_between_tauBSS}
    Let $X$ be a filtered spectrum.
    Then for every $k\geq 1$, we have a morphism of bifiltered spectra $\BF_\tau X \to \tr_k \BF_\tau X$ of the form
    \[
        \begin{tikzcd}
            \dotsb \ar[r,"\tau"] & \opSigma^{0,-k} X \ar[r,"\tau"] \ar[d] & \opSigma^{0,-k+1}X \ar[r,"\tau"] \ar[d] & \opSigma^{0,-k+2}X \ar[r,"\tau"] \ar[d] & \dotsb \ar[r,"\tau"] & X \ar[d] \\
            \dotsb \ar[r] & 0 \ar[r] & \opSigma^{0,-k+1} X/\tau \ar[r,"\tau"] & \opSigma^{0,-k+2} X/\tau^2 \ar[r,"\tau"] & \dotsb \ar[r,"\tau"] & X/\tau^k
        \end{tikzcd}
    \]
    constructed by letting each nontrivial square be a pushout.
    Formally, it is left Kan extended from the subdiagram
    \[
        \begin{tikzcd}
            \dotsb \ar[r,"\tau"] & \opSigma^{0,-k-1} X \ar[r,"\tau"] \ar[d] & \opSigma^{0,-k}X \ar[r,"\tau"] \dar & \dotsb \ar[r,"\tau"] & X. \\
            \dotsb \ar[r] & 0 \rar & 0 & &
        \end{tikzcd}
    \]
    Likewise, for $m\geq k\geq 1$, we have a morphism $\tr_m \BF_\tau X \to \tr_k \BF_\tau X$, and these fit into a tower
    \[
        \BF_\tau X \to \dotsb \to \tr_k \BF_\tau X \to \tr_{k-1} \BF_\tau X \to\dotsb \to \tr_1\BF_\tau X.
    \]
    On colimits (equivalently, on filtration $0$), this tower is the \emph{$\tau$-adic tower} of $X$:
    \[
        X \to \dotsb \to X/\tau^k \to X/\tau^{k-1} \to \dotsb\to X/\tau.
    \]
\end{construction}



The only subtlety with the truncated version is the distinction between the two kinds of $\tau$-multiples, as discussed in \cref{not:two_different_meanings_multiplication_by_tau}.

\begin{theorem}
    \label{thm:truncated_tau_BSS_structure}
    Let $X$ be a filtered spectrum and let $k\geq 1$.
    Let $\set{\uE_r^{*,*},d_r}$ denote the spectral sequence underlying $X$, and let $\set{\uE_r^{*,*,*},d_r^\tau}$ denote the $k$-truncated $\tau$-Bockstein spectral sequence of $X$.
    \begin{numberenum}
        \item \label{item:morph_tauBSS_to_truncated} The map on spectral sequences induced by $\BF_\tau X \to \tr_k \BF_\tau X$ from \cref{constr:maps_between_tauBSS} is, on first pages, given by the quotient
        \[
            \Z[\taubar] \otimes \oppi_{*,*}(X/\tau) \twoheadto \Z[\taubar]/\taubar^k \otimes \oppi_{*,*}(X/\tau).
        \]
        For $m\leq k$, the map on spectral sequences induced by $\tr_k \BF_\tau X \to \tr_m \BF_\tau X$ from \cref{constr:maps_between_tauBSS} is, on first pages, given by the quotient
        \[
            \Z[\taubar]/\tau^k \otimes \oppi_{*,*}(X/\tau) \twoheadto \Z[\taubar]/\taubar^m \otimes \oppi_{*,*}(X/\tau).
        \]
        \item \label{item:trunc_diffs_tau_linear} The differentials are $\taubar$-linear: for $x \in \uE_1^{n,w,s}$ and $y\in \uE_1^{n-1,\, w,\,s+r}$, if there is a differential
        \[
            d_r^\tau(x) = y,
        \]
        then for all $m\geq 0$, there is a differential
        \[
            d_r^\tau(\taubar^m x) = \taubar^m y.
        \]
        In particular, the $\Z[\taubar]$-module structure on the first page induces a $\Z[\taubar]$-module structure on later pages.
        \item \label{item:trunc_diffs_hit_tau_multiples} The target of a $d_r^\tau$-differential is a $\taubar^r$-multiple: for every $x \in \uE_r^{n,w,s}$, there is an element $y \in \uE_r^{n-1,\, w+r,\,s}$ such that
        \[
            d_r^\tau(x) = \taubar^r \cdot y,
        \]
        where the multiplication denotes the $\Z[\taubar]$-module structure on $\uE_r^{*,*,*}$ from \ref{item:trunc_diffs_tau_linear}.
        \item \label{item:trunc_diffs_are_underlying_diff}
        Let $r \leq k-1$.
        If $x\in \uE_1^{n,w}$ and $y\in\uE_{r}^{n-1,\, w+r}$, then there is a differential
        \begin{align*}
            d_r(x) &= y\\
        \shortintertext{if and only if there is a differential}
            d_r^\tau(x) &= \taubar^r \cdot y.
        \end{align*}
        \item \label{item:trunc_tauBSS_detection} Suppose that $x \in \uE_1^{n,w,s}$ detects an element $\theta \in \oppi_{n,w}X/\tau^k$.
        Then $\theta$ is in the image of the map
        \[
            \tau^{s} \colon \opSigma^{0,\,-s} X/\tau^{k-s} \to X/\tau^k,
        \]
        i.e., it is of the form $\tau^s(\alpha)$ for some $\alpha \in \oppi_{n,\,w+s}X/\tau^{k-s}$.
        
        Moreover, for every $m$ such that $s+m \leq k-1$, the class $\taubar^m \cdot x$ detects $\tau^m \cdot \theta$.
        For every $m$ such that $s+m\geq k$, we instead have $\tau^m \cdot \theta = 0$.
        \item \label{item:trunc_tauBSS_convergence} The $k$-truncated $\tau$-Bockstein spectral sequence converges strongly to $\oppi_{*,*}X/\tau^k$.
    \end{numberenum}
\end{theorem}

Take particular note that in \cref{item:trunc_diffs_tau_linear}, the truncated differentials are no longer $\taubar$-divisible in general, and that \cref{item:trunc_diffs_are_underlying_diff} only applies to differentials of length $r\leq k-1$ (longer $d_r^\tau$-differentials vanish for degree reasons).

\begin{proof}
    The definition of the map $\BF_\tau X \to \tr_k \BF_\tau X$ from \cref{constr:maps_between_tauBSS} as a left Kan extension shows that in filtrations $0 \leq s \leq k-1$, the map induces an isomorphism on associated graded.
    This proves the first part of \cref{item:morph_tauBSS_to_truncated}, and the second is analogous.
    
    \Cref{item:trunc_diffs_tau_linear,item:trunc_diffs_hit_tau_multiples,item:trunc_diffs_are_underlying_diff} follow immediately from \cref{item:differentials_tau_linear,item:differentials_hit_tau_multiples,item:tau_BSS_diff_are_underlying_diff} of \cref{thm:tau_BSS_structure}, using the map from \cref{item:morph_tauBSS_to_truncated}.
    For \cref{item:trunc_tauBSS_convergence}, it is enough to note that $\tr_k\BF_\tau X$ vanishes in filtrations $s \geq k$, so that (the trigraded analogue of) \cref{prop:eventually_zero_implies_strong_convergence} implies the strong convergence.

    It remains to prove \cref{item:trunc_tauBSS_detection}.
    The statement that $x$ detects $\theta$ means that there is a lift $\alpha \in \oppi_{n,\,w+s}X/\tau^{k-s}$ of $x$ that maps to $\theta$ under the transition map
    \[
        \tr_k \BF_\tau^s X = \opSigma^{0,\,-s}X/\tau^{k-s} \to X/\tau^k = \tr_k \BF_\tau^0 X,
    \]
    in other words, that $\tau^s(\alpha)=\theta$.
    If $s+m\leq k-1$, then
    \[
        \oppi_{n,\,w-m}(\Gr^{s+m}\tr_k\BF_\tau X) = \oppi_{n,w}X/\tau,
    \]
    and $\taubar^m\cdot x$ is given by the element $x$ considered as an element of this group via this identification.
    Write $\beta$ for the mod $\tau^{k-s-m}$ reduction of $\alpha$.
    Then $\beta$ is a lift of $\taubar^m \cdot x$ to $\oppi_{n,\,w-m}(\tr_k\BF_\tau^{s+m}X)=\oppi_{n,\,w+s}X/\tau^{k-s-m}$.
    Evidently, the image of this element in $\tr_k\BF_\tau^0 X$ is given by
    \[
        \tau^{s+m}(\beta) = \tau^s(\tau^m(\beta)) = \tau^s (\tau^m \cdot \alpha) = \tau^m \cdot \tau^s(\alpha) = \tau^m \cdot \theta.
    \]
    Lastly, if instead $s+m \geq k$, then this means $m \geq k-s$.
    Since $\tau^m \cdot \tau^s(\alpha) = \tau^s(\tau^m \cdot \alpha)$ and $\alpha$ lives in $\oppi_{*,*}X/\tau^{k-s}$, it follows that $\tau^m \cdot \alpha = 0$ (using the ring structure on $C\tau^{k-s}$), proving the final claim.
\end{proof}


\begin{warning}
    Unlike in the integral case, if an element in the truncated $\tau$-BSS is detected in filtration $s$, this does not imply that it is a $\tau^s$-multiple in $\oppi_{*,*}X/\tau^k$.
    Instead, in general this only implies that it is of the form $\tau^s(\theta)$ for some $\theta$ in $\oppi_{*,*}X/\tau^{k-s}$.
    (See \cref{not:two_different_meanings_multiplication_by_tau} for the distinction between these two.)
\end{warning}

%% file: fil/omnibus.tex
\section{The Omnibus Theorem}
\label{sec:filtered_omnibus}

We can summarise \cref{sec:tau_in_fil_sp} as follows: if $X$ is a filtered spectrum, then its underlying spectral sequence is of the form
\[
    \uE_1^{n,s} = \pi_{n,s}(C\tau\otimes X) \implies \pi_n(X^{\tau=1}).
\]
In \cref{ch:informal_sseq}, we argued that the bigraded homotopy groups $\oppi_{*,*}X$ as a $\Z[\tau]$-module should capture the differentials in this spectral sequence.
The \emph{Omnibus Theorem} makes this precise, by describing the structure of $\oppi_{*,*}X$ in terms of the underlying spectral sequence of $X$.
We prove it here in the context of filtered spectra, which should be regarded as the most general (stable) setting for it: later in \cref{sec:deformations}, we will show how it extends to any \emph{monoidal deformation}.
In \cref{ch:synthetic_spectra}, we will see how this recovers and extends the synthetic Omnibus Theorem of Burklund--Hahn--Senger, and compare our proof to theirs; see \cref{sec:synthetic_omnibus} in particular.


\begin{theorem}[Omnibus]
    \label{thm:filtered_omnibus}
    Let $X$ be a complete filtered spectrum, and assume that in the spectral sequence underlying $X$, we have $\uR\uE_\infty^{*,*}=0$ \textbr{for instance, this happens if the spectral sequence converges strongly}.
    Let $x \in \uE_1^{n,s} = \oppi_{n,s}(X/\tau)$ be a nonzero class.
    Then the following are equivalent.
    \begin{enumerate}[label={\upshape(1\alph*)}]
        \item \label{item:fil_ombs_perm_cycle} The element $x$ is a permanent cycle.
        \item \label{item:fil_ombs_lifts} The element $x \in \oppi_{n,s}(X/\tau)$ lifts to an element of $\oppi_{n,s}X$.
    \end{enumerate}
    For any such lift $\alpha$ to $\oppi_{n,s}X$, the following are true.
    \begin{enumerate}[label={\upshape(2\alph*)}]
        \item \label{item:fil_omnbs_not_tau_torsion} If $x$ survives to page $r$, then $\tau^{r-1}\cdot \alpha$ is nonzero.
        \item \label{item:fil_ombs_detection} If $x$ survives to page $\infty$, then $\alpha$ maps to a nonzero element in $\pi_n X^{\tau=1} = \pi_n X^{-\infty}$, and this element is detected by~$x$.
    \end{enumerate}
    Moreover, if $x$ lifts to $X$, then there exists a lift $\alpha$ with either one of the following additional properties.
    \begin{enumerate}[label={\upshape(3\alph*)}]
        \item \label{item:fil_omnbs_torsion_lift} If $x$ is the target of a $d_r$-differential, then $\tau^r \cdot \alpha = 0$.
        \item \label{item:fil_omnbs_detection_lift} If $\theta \in \pi_n X^{\tau=1}$ is detected by $x$, then $\alpha$ is sent to $\theta$ under $\oppi_{n,s}X \to \pi_n X^{\tau=1}$.
    \end{enumerate}
    Finally, we have the following generation statement.
    \begin{enumerate}[label={\upshape(4)}]
        \item \label{item:fil_omnbs_generation} Let $\set{\alpha_i}$ be a collection of elements of $\oppi_{n,*}X$ such that their mod $\tau$ reductions generate the permanent cycles in stem $n$.
        Then the $\tau$-adic completion of the $\Z[\tau]$-submodule of $\oppi_{n,*}X$ generated by the $\set{\alpha_i}$ is equal to $\oppi_{n,*}X$.
    \end{enumerate}
\end{theorem}
\begin{proof}
    By \cref{thm:tau_BSS_structure}\,\ref{item:tau_BSS_convergence}, the hypotheses on $X$ are equivalent to strong convergence of the $\tau$-Bockstein spectral sequence for $X$.

    We begin with (1).
    The map $\BF_\tau X \to \tr_1\BF_\tau X$ from \cref{constr:maps_between_tauBSS} is, on underlying objects, equal to the reduction map $X\to X/\tau$.
    Meanwhile the induced map on the first page of the resulting spectral sequences is given by quotienting by $\taubar$:
    \[
        \Z[\taubar] \otimes \oppi_{*,*}(X/\tau) \twoheadto \oppi_{*,*}(X/\tau),
    \]
    so in particular it is an isomorphism in filtration $0$.
    By strong convergence, every permanent cycle in the $\tau$-BSS lifts to $\oppi_{*,*}X$.
    Accordingly, it follows that $x \in \oppi_{n,s}(X/\tau)$ lifts to $\oppi_{n,s}X$ if and only if the element $x$ considered in tridegree $(n,s,0)$ of the $\tau$-BSS is a permanent cycle.
    The equivalence of \ref{item:fil_ombs_perm_cycle} and \ref{item:fil_ombs_lifts} therefore follows from the identification of the differentials in the $\tau$-BSS from \cref{thm:tau_BSS_structure}\,\ref{item:tau_BSS_diff_are_underlying_diff}.

    Next, we check \ref{item:fil_omnbs_not_tau_torsion}.
    Suppose that $\alpha \in \oppi_{n,s}X$ is a lift of $x$.
    Then since $\BF_\tau X$ is constant from filtration $0$ and onwards, this means that the class $x$ detects $\alpha$.
    By \cref{thm:tau_BSS_structure}\,\ref{item:tau_BSS_detection}, the element $\taubar^{r-1}\cdot x$ detects $\tau^{r-1}\cdot \alpha$.
    To show that $\tau^{r-1}\cdot \alpha$ is nonzero, by (the trigraded analogue of) \cref{rmk:can_detect_nonzeroness} we need to show that $\taubar^{r-1}\cdot x$ survives to page $\infty$.
    Again using \cref{thm:tau_BSS_structure}\,\ref{item:tau_BSS_diff_are_underlying_diff}, we see that $\taubar^{r-1}\cdot x$ survives to page~$r$ if and only if $x$ in the spectral sequence underlying $X$ survives to page~$r$.
    If this happens, then for degree reasons $\taubar^{r-1}\cdot x$ also survives to page $\infty$.

    
    Properties \ref{item:fil_ombs_detection} and \ref{item:fil_omnbs_detection_lift} are a restatement of the definition of detection from \cref{def:perm_cycle_and_detection}, with \ref{item:fil_ombs_detection} specifically following from \cref{rmk:can_detect_nonzeroness}.

    We prove \ref{item:fil_omnbs_torsion_lift} by induction on $r$.
    Suppose that $x$ is the target of a $d_r$-differential, say $d_r(y) = x$ for $y \in \oppi_{n+1,\, s-r}(X/\tau)$.
    This implies that $d_r^\tau(y) = \taubar^r \cdot x$.
    Unrolling the definition of the $d_r^\tau$-differential, this means the following.
    There exists an element $\alpha \in \oppi_{n,s}X$ such that
    \[
        \tau^{r-1} \cdot \alpha = \partial_1^\infty(y) \quad \text{in } \oppi_{n-1,\,s-r+1}X,
    \]
    and such that the mod $\tau$ reduction of $\alpha$ agrees with $x$ up to the images of differentials of length shorter than $r$.
    \[
        \begin{tikzcd}
            \opSigma^{0,\,-r} X \ar[r,"\tau^{r-1}"] \ar[d] & \Sigma^{0,\, -1} X \ar[r,"\tau"] & X \ar[d]\\
            \opSigma^{0,\,-r} X/\tau & &  X/\tau \ar[ul, dashed, "\partial_1^\infty"]
        \end{tikzcd}
    \]
    By induction, we may assume that the images of shorter differentials have $\tau^{r-1}$-torsion lifts, so that without loss of generality, we may assume $\alpha$ is a lift of $x$.
    By exactness, we have $\tau \cdot \partial_1^\infty(y) = 0$.
    But then $\tau^r\cdot \alpha = \tau \cdot \partial_1^\infty(y)=0$, so that $\alpha$ is a lift satisfying \ref{item:fil_omnbs_torsion_lift}.

    For the final claim, let $M$ denote the $\Z[\tau]$-submodule of $\oppi_{n,*}X$ generated by the $\alpha_i$.
    By strong convergence of the $\tau$-BSS, the $\Z[\tau]$-module $\oppi_{n,*}X$ is $\tau$-complete; see \cref{warn:tau_complete_filtered_vs_abelian}.
    It follows that $M_\tau^\wedge$ is naturally a submodule of $\oppi_{n,*}X$.
    To show that this inclusion is an equality, it suffices to show that it becomes surjective after quotienting by $\tau$, as both modules are $\tau$-complete.
    After quotienting by $\tau$, the inclusion becomes
    \[
        M/\tau \to (\oppi_{n,*}X)/\tau = F^0_\tau \oppi_{n,*}X,
    \]
    where the last identification is \cref{thm:tau_BSS_structure}\,\ref{item:tau_BSS_detection}.
    By the identifications of the differentials in the $\tau$-BSS of \cref{thm:tau_BSS_structure}\,\ref{item:tau_BSS_diff_are_underlying_diff}, the assumption on the $\alpha_i$ translates to this map being a surjection.
    This finishes the proof.
\end{proof}

\begin{warning}
    Suppose that $x \in \uE_1^{n,s}$ is the target of a $d_r$-differential, and suppose that $\alpha \in \oppi_{n,s}X$ is a lift of $x$.
    Then \cref{item:fil_omnbs_torsion_lift} does \emph{not} necessarily imply that $\alpha$ is $\tau^r$-torsion (in fact, in general $\alpha$ need not even be $\tau$-power torsion): the theorem only guarantees that there exists \emph{some} lift of $x$ that is $\tau^r$-torsion.
\end{warning}

We can draw a number of simpler conclusions from this result.

\begin{corollary}
    \label{cor:tau_torsion_free_iff_no_differentials}
    Let $X$ be a filtered spectrum satisfying the conditions of \cref{thm:filtered_omnibus}, and let $n$ be an integer.
    Then the $\Z[\tau]$-module $\oppi_{n,*}X$ is $\tau$-power torsion free if and only if the $n$-stem in the spectral sequence underlying $X$ does not receive any nonzero differentials.
\end{corollary}
\begin{proof}
    The lack of incoming differentials implies that every permanent cycle survives to page $\infty$, so one direction follows from \cref{item:fil_omnbs_not_tau_torsion}.
    Conversely, if the $n$-stem does receive a nonzero differential, then by \cref{item:fil_omnbs_torsion_lift} there exists a $\tau$-power torsion element in $\oppi_{n,*}X$.
\end{proof}

\begin{corollary}
    Let $X$ be a filtered spectrum satisfying the conditions of \cref{thm:filtered_omnibus}, and let $n,s$ be integers.
    If $\oppi_{n,\,s+d}X/\tau$ vanishes for all $d\geq 0$, then $\oppi_{n,s}X$ vanishes also.
\end{corollary}
\begin{proof}
    This follows directly from \cref{item:fil_omnbs_generation}.
\end{proof}

\begin{remark}
    By inspecting the proof of \cref{thm:filtered_omnibus}, we see that the convergence conditions on $X$ are only needed for items~(1) and~\ref{item:fil_omnbs_generation}.
    Without these assumptions, item~(2) still holds for any lift, but such a lift is no longer guaranteed to exist.
    In item~(3) meanwhile, the assumptions on $x$ in both \ref{item:fil_omnbs_torsion_lift} and \ref{item:fil_omnbs_detection_lift} imply that a lift exists (see \cref{def:perm_cycle_and_detection}).
\end{remark}

\subsection{The truncated Omnibus Theorem}
\label{ssec:truncated_filtered_omnibus}

There is also a version of the Omnibus Theorem that describes the structure of $\oppi_{*,*}(X/\tau^k)$ in terms of the spectral sequence underlying $X$.
In this case, we no longer need any convergence conditions, but the generation statement is more involved to state.
Let us begin therefore with the other parts of the Omnibus Theorem.

\begin{theorem}[Truncated Omnibus, part 1]
    \label{thm:truncated_filtered_omnibus}
    Let $X$ be a filtered spectrum, and let $k\geq 1$.
    Let $x \in \uE_1^{n,s} = \oppi_{n,s}(X/\tau)$ be a class.
    Then the following are equivalent.
    \begin{enumerate}[label={\upshape(1\alph*)}]
        \item \label{item:tr_fil_omnbs_cycle} The differentials $d_1(x),\dotsc,d_{k-1}(x)$ vanish.
        \item \label{item:tr_fil_omnbs_lifts} The element $x \in \oppi_{n,s}(X/\tau)$ lifts to an element of $\oppi_{n,s}(X/\tau^k)$.
    \end{enumerate}
    For any such lift $\alpha$ to $\oppi_{n,s} (X/\tau^k)$, the following are true.
    \begin{enumerate}[label={\upshape(2\alph*)}]
        \item \label{item:tr_fil_omnbs_not_tau_torsion} If $x$ survives to page $r$ for $r \leq k$, then $\tau^{r-1}\cdot \alpha$ is nonzero.
        \item \label{item:tr_fil_omnbs_differential} The image of $\alpha$ under $\partial_{k}^{k+1} \colon \oppi_{n,s} (X/\tau^k) \to \oppi_{n-1,\, s+r} (X/\tau)$ is a representative for $d_k(x)$.
    \end{enumerate}
    Moreover, if $x$ lifts to $X/\tau^k$, then there exists a lift $\alpha$ with the following additional property.
    \begin{enumerate}[label={\upshape(3)}]
        \item \label{item:tr_fil_omnbs_torsion_lift} If $x$ is the target of a $d_r$-differential for $r<k$, then $\tau^r \cdot \alpha = 0$.
    \end{enumerate}
\end{theorem}

\begin{proof}
    The $k$-truncated $\tau$-Bockstein spectral sequence for $X$ converges strongly to $\oppi_{*,*}X/\tau^k$ by \cref{thm:truncated_tau_BSS_structure}\,\ref{item:trunc_tauBSS_convergence}.

    The morphism $\tr_k\BF_\tau X \to \tr_1\BF_\tau X$ from \cref{constr:maps_between_tauBSS} is, on underlying objects, the reduction map $X/\tau^k \to X/\tau$.
    On first pages of the underlying spectral sequences on the other hand, the map takes the form of quotienting $\taubar$:
    \[
        \Z[\taubar]/\taubar^k\otimes \oppi_{*,*}X/\tau \twoheadto \oppi_{*,*}X/\tau.
    \]
    In the same way as in the proof of \cref{thm:filtered_omnibus}, the equivalence between \ref{item:tr_fil_omnbs_cycle} and \ref{item:tr_fil_omnbs_lifts} follows by strong convergence and from the identification of the differentials in the truncated $\tau$-BSS arising from \cref{thm:truncated_tau_BSS_structure}\,\ref{item:morph_tauBSS_to_truncated}.

    Next, we prove item~\ref{item:tr_fil_omnbs_not_tau_torsion}.
    Let $\alpha \in \oppi_{n,s}X/\tau^k$ be a lift of $x$.
    By \cref{thm:truncated_tau_BSS_structure}\,\ref{item:trunc_tauBSS_detection}, the element $\tau^{r-1}\cdot \alpha$ is detected by $\taubar^{r-1}\cdot x$.
    The assumption that $x$ survives to page $r$ implies that $\taubar^{r-1}\cdot x$ is a permanent cycle that survives to page $r$, and for degree reasons it then also survives to page $\infty$.
    Using \cref{rmk:can_detect_nonzeroness}, it follows that $\tau^{r-1}\cdot \alpha$ is nonzero.
    
    Next, we prove item~\ref{item:tr_fil_omnbs_differential}.
    Let $\alpha \in \oppi_{n,s}X/\tau^k$ be a lift of $x$.
    As a consequence of \cref{thm:truncated_tau_BSS_structure}\,\ref{item:morph_tauBSS_to_truncated}, it is enough to show that $\partial_{k}^{k+1}(\alpha)$ is a representative for $d_{k}^\tau(x)$ in the non-truncated $\tau$-BSS, where we regard $x$ in tridegree $(n,s,0)$.
    Recall how the $d_{k}$-differential on $x$ in the non-truncated $\tau$-BSS is computed: we apply the boundary map $\partial_1^\infty(x)$, choose a $\tau^{k-1}$-division of this element, and reduce this mod $\tau$.
    \[
        \begin{tikzcd}
            \opSigma^{0,\,-k} X \ar[r,"\tau^{k-1}"] \ar[d] & \Sigma^{0,\, -1}X \ar[r,"\tau"] & X \ar[d]\\
            \opSigma^{0,\,-k} X/\tau & &  X/\tau \ar[ul, dashed, "\partial_1^\infty"]
        \end{tikzcd}
    \]
    We have $\tau^{k-1} \cdot \partial_{k}^\infty (\alpha) = \partial_1^\infty(x)$.
    In other words, $\partial_{k}^\infty(\alpha)$ is a valid choice of $\tau^{k-1}$-division of $\partial_1^\infty(a)$, so that its projection to $\oppi_{n,\, s+k}X/\tau$ is a representative for $d_{k}^\tau(x)$.
    But the mod~$\tau$ reduction of $\partial_{k}^\infty$ is $\partial_{k}^{k+1}$, proving item~\ref{item:tr_fil_omnbs_differential}.

    Next, we prove item~\ref{item:tr_fil_omnbs_torsion_lift}.
    Suppose that $x$ is the target of a $d_r$-differential for $r<k$.
    We consider the element $\taubar^{r-1}\cdot x$ in tridegree $(n,\ s-r+1,\ r-1)$ of the $k$-truncated $\tau$-Bockstein spectral sequence.
    Then the $d_r$-differential hitting $x$ translates to a $d_{r}^\tau$-differential
    \[
        d_r^\tau(y) = \taubar^r x.
    \]
    Unrolling what this means, we learn the following.
    There exists an element $\beta \in \oppi_{n,s}(X/\tau^{k-r})$ such that $\tau^{r-1}(\beta) = \partial_1^k(y)$, and such that $\beta$ reduces to $x$ modulo~$\tau$ up to the images of shorter differentials.
    \[
        \begin{tikzcd}
            \opSigma^{0,\,-r} X/\tau^{k-r} \ar[r,"\tau^{r-1}"] \ar[d] & \Sigma^{0,\, -1} X/\tau^{k-1} \ar[r,"\tau"] & X/\tau^k \ar[d]\\
            \opSigma^{0,\,-r} X/\tau & &  X/\tau \ar[ul, dashed, "\partial_1^k"]
        \end{tikzcd}
    \]
    Like in the proof of \cref{thm:filtered_omnibus}, by induction we may assume without loss of generality that $\beta$ reduces to $x$. 
    By exactness, the element $\partial_1^{k}(y)$ satisfies $\tau(\partial_1^k(y))=0$, so any choice of $\beta$ satisfies $\tau^r(\beta)=0$.
    As a result, it suffices to show that there is a choice of $\beta \in \oppi_{n,s}X/\tau^{k-r}$ that lifts to an element in $\oppi_{n,s} X/\tau^k$.
    Indeed, if such a lift $\alpha$ exists, then by \cref{rmk:relation_two_different_tau_multiplications}, we have $\tau^r \cdot \alpha = \tau^r(\beta) = 0$, which would mean that $\alpha$ is the lift proving \cref{item:tr_fil_omnbs_torsion_lift}.
    
    To produce such an $\alpha$, we first show that $y \in \oppi_{n+1,\,s-r} X/\tau$ lifts to $X/\tau^{r}$.
    This follows from (1) because $y$ is a $d_{\leq r-1}$-cycle.
    Choose a lift $\widetilde{y}$.
    It then follows that $\partial_{r}^{k}(\widetilde{y})$ is a valid choice for $\beta$ as above.
    To show that this $\beta$ lifts to $X/\tau^k$, we need to show that $\partial_{k-r}^{k}(\beta)=0$.
    Note that $\partial_{k-r}^{k}\circ \partial_{r}^{k}$ can be written as (omitting shifts for readability)
    \[
        \begin{tikzcd}
            X/\tau^{k-1} \rar["\partial_{k}^\infty"] & X \rar & X/\tau^{k-r} \rar["\partial_{k-r}^\infty"] &[1em] X \rar & X/\tau^{r}.
        \end{tikzcd}
    \]
    The middle two maps are part of a cofibre sequence, so in particular their composition is zero.
    This means that indeed $\partial_{k-r}^k(\beta) = 0$, showing that $\beta$ lifts to the desired~$\alpha$, thus proving item~\ref{item:tr_fil_omnbs_torsion_lift}.
\end{proof}

\begin{construction}
    \label{constr:map_from_mod_tau_r_to_Er}
    Let $X$ be a filtered spectrum, and let $n,s$ be integers.
    \begin{numberenum}
        \item \label{item:map_X_to_Einfty} By \cref{thm:filtered_omnibus}\,(1), the reduction map $\oppi_{n,s}X\to \oppi_{n,s}X/\tau = \uE_1^{n,s}$ has image given by the subgroup $\uZ_\infty^{n,s}$ of permanent cycles.
        Postcomposing with the quotients, we in particular obtain, for every $1\leq r \leq \infty$, a map
        \[
            \oppi_{n,s}X \to \uE_{r}^{n,s}
        \]
        which is surjective in the case $r=\infty$.
        \item Let $r\geq 1$.
        By \cref{thm:truncated_filtered_omnibus}\,(1), the reduction map $\oppi_{n,s}X/\tau^r \to \oppi_{n,s}X/\tau = \uE_1^{n,s}$ has image given by the subgroup $\uZ_{r-1}^{n,s}$ of $(r-1)$-cycles.
        Postcomposing with the quotient, we in particular obtain a surjective map
        \[
            \oppi_{n,s}X/\tau^r \twoheadto \uE_{r}^{n,s}.
        \]
    \end{numberenum}
    These maps are compatible with each other in the obvious way.
\end{construction}

Using these maps, we can now reformulate \cref{prop:divisibility_total_diff_and_ordinary_diff}.

\begin{corollary}
    \label{cor:diagram_differential_and_total_differential}
    Let $X$ be a filtered spectrum, let $r\geq 1$, and let $n$ and $s$ be integers.
    Then we have a commutative diagram
    \[
        \begin{tikzcd}
            \oppi_{n,s}X/\tau^{r} \dar \rar["\partial_{r}^\infty"] & \oppi_{n-1,\, s+r} X \dar \\
            \uE_r^{n,s} \rar["d_{r}"] & \uE_r^{n-1,\, s+r}.
        \end{tikzcd}
    \]
    More generally, for $R \geq r$, we have a commutative diagram
    \[
        \begin{tikzcd}
            \oppi_{n,s}X/\tau^{r} \ar[r,"\partial_{r}^R"] \ar[d] & \oppi_{n-1,\, s+r} X/\tau^{R-r} \ar[d] \\
            \uE_r^{n,s} \ar[r,"d_r"] & \uZ_{m}^{n-1,\, s+r}/\uB_{r-1}^{n-1,\, s+r}
        \end{tikzcd}
    \]
    where $m$ is the minimum of $r-1$ and $R-r-1$.
\end{corollary}

\begin{proof}
    This follows directly from \cref{prop:divisibility_total_diff_and_ordinary_diff}.
\end{proof}

Next, we turn to the question of finding generators for $\oppi_{*,*}X/\tau^k$.
In the non-truncated version (\cref{thm:filtered_omnibus}\,\ref{item:fil_omnbs_generation}), we could start with lifts of permanent cycles, and take the ($\tau$-complete) $\Z[\tau]$-module that they generate to reconstruct all of $\oppi_{*,*}X$.
In the truncated case, taking $\tau$-multiples is a more subtle notion.
In order to generate $\oppi_{*,*}X/\tau^k$, we need to take the $\tau$-multiples of elements from all lower truncations $X/\tau^{i}$ for $i\leq k$ into account.
Unfortunately, stating this precisely makes the indexing get a little out of hand.

For applications, we also need a relative version describing the kernel of $X/\tau^k \to X/\tau^m$.
To compensate for the more intricate phrasing of this result, we give a simplified, more coarse description in \cref{cor:easy_lifting_higher_power_tau} below.

\begin{theorem}[Truncated Omnibus, part 2]
    \label{thm:truncated_filtered_omnibus_generation}
    Let $X$ be a filtered spectrum, let $n,s\in\Z$, and let $r\geq 1$.
    \begin{numberenum}
        \item Let $k\geq 1$ be fixed.
        Suppose that for every $1 \leq i \leq k$, we have a collection of elements
        \[
            \set{\beta^i_j}_j \subseteq \oppi_{n,\, s+k-i}X/\tau^i
        \]
        whose mod $\tau$ reductions generate the abelian group
        \begin{equation}\label{eq:cycles_mod_boundaries}
            \uZ_{i-1}^{n,\,s+k-i}/\uB_{k-i}^{n,\, s+k-i}.
        \end{equation}
        \textbr{Note that by \cref{thm:truncated_filtered_omnibus}\,\textup{(1)}, such a collection exists for every $i$.}
        Write
        \[
            \alpha^i_j \defeq \tau^{k-i}(\beta^i_j) \in \oppi_{n,s}X/\tau^k.
        \]
        Then $\set{\alpha^i_j}_{i,j}$ is a set of generators for the abelian group $\oppi_{n,s}X/\tau^k$.
        
        \item Let $1 \leq m \leq k$ be fixed.
        Suppose that for every $1 \leq i \leq k-m$, we have a collection of elements
        \[
            \set{\beta^i_j}_j \subseteq \oppi_{n,\, s+k-i}X/\tau^i
        \]
        whose mod $\tau$ reductions generate the abelian group \eqref{eq:cycles_mod_boundaries}.
        Write $\alpha^i_j \defeq \tau^{k-i}(\beta^i_j)$.
        Then $\set{\alpha^i_j}_{i,j}$ is a set of generators for the abelian group
        \[
            \ker (\oppi_{n,s}X/\tau^k \to \oppi_{n,s}X/\tau^m).
        \]
    \end{numberenum}
\end{theorem}

\begin{proof}
    The $k$-truncated Bockstein spectral sequence converges strongly to $\oppi_{*,*}X/\tau^k$ by \cref{thm:truncated_tau_BSS_structure}\,\ref{item:trunc_tauBSS_convergence}.
    The first result therefore follows from \cref{thm:truncated_tau_BSS_structure}\,\ref{item:trunc_diffs_are_underlying_diff}; see also \cref{rmk:cycles_boundaries_tau_BSS}.
    The second follows analogously by considering the natural map from the $k$-truncated $\tau$-BSS to the $m$-truncated $\tau$-BSS for~$X$ induced by the map from \cref{constr:maps_between_tauBSS}.
\end{proof}

Sometimes, the following simplified criterion is sufficient.

\begin{corollary}
    \label{cor:easy_lifting_higher_power_tau}
    Let $X$ be a filtered spectrum, let $n,s\in\Z$, and let $k \geq m \geq 1$.
    \begin{numberenum}
        \item If $\oppi_{n,\, s+d} X/\tau$ vanishes for $0 \leq d \leq k-1$, then $\oppi_{n,s}X/\tau^k$ vanishes also.
        \item If $\oppi_{n,\, s+d} X/\tau$ vanishes for $m \leq d \leq k-1$, then the reduction map $\oppi_{n,s}X/\tau^k \to \oppi_{n,s}X/\tau^m$ is injective.
    \end{numberenum}
\end{corollary}
\begin{proof}
    In the notation of \cref{thm:truncated_filtered_omnibus_generation}, we have $\uE_1^{n,s} = \oppi_{n,s}X/\tau$, and $\uZ_r^{n,s}$ is a subgroup of this.
    It follows that the relevant groups in \eqref{eq:cycles_mod_boundaries} vanish, so the claims follow.
\end{proof}

%% file: fil/deform2.tex
\section{Deformations}
\label{sec:deformations}

So far, we have seen that the $\infty$-category $\FilSp$ is the natural home for (stable) spectral sequences.
For specific purposes however, this category might be a bit unwieldy, and it might be helpful to find a modification of $\FilSp$ that is more suited to the problem at hand.
The main example of such a modification in these notes is that of \emph{synthetic spectra}.
However, much of the setup of synthetic spectra holds much more generally, and leads one to a broad theory of `modifications' of $\FilSp$.
These have become known as \emph{deformations}.
Readers only interested in synthetic spectra may move on to the next chapter, referring back to this section as needed.

This section is concerned with the general properties of deformations.
Of particular interest is the case where this deformation structure arises from a (symmetric) monoidal left adjoint out of $\FilSp$; we call these \emph{\textbr{symmetric} monoidal deformations}, which are the subject of \cref{ssec:monoidal_deformations}.
For monoidal deformations, we can prove much more: we prove all that we need in order to import all results about filtered spectra into a monoidal deformation; see \cref{prop:deformation_sigma_tauinv}.
In particular, the Omnibus Theorem holds in any monoidal deformation, where the underlying spectral sequence is replaced by what we call the \emph{signature spectral sequence}.
We will make particular use of this in the case of synthetic spectra in the next chapter.

Later, in \cref{ch:variants_Syn}, we will continue our study of deformations and discuss \emph{cellularity}, \emph{filtered models}, and \emph{evenness}.
For now, our goals are more modest, and our main aim is to show how to deduce results in a deformation from results in filtered spectra.

Much of the material in this section is based on the treatment of deformations given by Barkan \cite[Section~2]{barkan_monoidal_algebraicity} and Burklund--Hahn--Senger \cite[Appendices~A--C]{burklund_hahn_senger_Rmotivic}.
While these sources also discuss constructing new deformations out of old ones, we will focus on studying phenomena within a fixed deformation.

\begin{definition}[\cite{barkan_monoidal_algebraicity}, Definition~2.2]
    \label{def:deformation}
    A \defi{(1-parameter, stable) deformation} is a left module over $\FilSp$ in $\PrL_\st$.
\end{definition}

A module over $\FilSp$ in $\PrL_\st$ is also called an \emph{$\FilSp$-linear $\infty$-category}.
We refer to \cite[Section~3]{naumann_pol_ramzi_monoidal_fracture} for further background on $\D$-linear $\infty$-categories where $\D$ is a presentably symmetric monoidal $\infty$-category.

A deformation $\C$ is in particular left tensored over $\FilSp$.
As a result, we will use the same notation of bigraded shifts $\Sigma^{n,s}$ on $\C$ to mean tensoring with the filtered spectrum $\S^{n,s}$.
Likewise, if $X \in \C$, we will generally denote by $\tau \colon \Sigma^{0,-1}X\to X$ the map given by tensoring the map $\tau \colon \S^{0,-1}\to \S$ of filtered spectra with $X$.
Moreover, if $A$ is a filtered ring spectrum, then we can speak of modules over $A$ in $\C$.
Two cases of $A$ deserve a special name.

\begin{notation}
    \label{not:special_generic_fibres_deformation}
    Let $\C$ be a deformation.
    The \defi{generic fibre} of $\C$ is defined by
    \[
        \C[\tau^{-1}] \defeq \Mod_{\S[\tau^{-1}]}(\C),
    \]
    and we refer to the functor
    \[
        \C \to \C[\tau^{-1}], \quad X \mapsto X[\tau^{-1}]\defeq \S[\tau^{-1}] \otimes X
    \]
    as the \defi{$\tau$-inversion functor}.
    Further, the \defi{special fibre} of $\C$ is defined by
    \[
        \Mod_{C\tau}(\C),
    \]
    and for $X \in \C$, we write
    \[
        X/\tau \defeq C\tau \otimes X.
    \]
    Finally, we write $\C_\tau^\wedge$ for the full subcategory of $\C$ on the $C\tau$-local objects, which we call \defi{$\tau$-complete}.
    This results in a localisation
    \[
        \begin{tikzcd}[column sep=4em]
            \C \rar["(\blank)_\tau^\wedge", shift left] & \C_{\tau}^\wedge. \lar[shift left, hook']
        \end{tikzcd}
    \]
\end{notation}

\begin{remark}
    \label{rmk:expressions_inv_tau_and_modtau}
    Using that the tensoring over $\FilSp$ preserves colimits in each variable, we find that
    \[
        X[\tau^{-1}] = \colim(\begin{tikzcd}
            X \rar["\tau"] & \Sigma^{0,1} X \rar["\tau"] & \dotsb
        \end{tikzcd})
    \]
    and
    \[
        X/\tau = \cofib(\tau \colon \opSigma^{0,-1}X \to X).
    \]
    Moreover, since $\S[\tau^{-1}]$ is an idempotent in $\FilSp$, it follows that $\C\to\C[\tau^{-1}]$ is a smashing localisation, and the forgetful functor $\C[\tau^{-1}] \to \C$ is fully faithful with essential image consisting of those $X$ on which $\tau$ is an isomorphism.
    In most examples, the forgetful functor $\Mod_{C\tau}(\C) \to \C$ is not fully faithful.
\end{remark}

In a moment, we will see a way to obtain examples of deformations.
For now, let us mention the following.

\begin{example}
    The universal example of a deformation is that of $\FilSp$ itself.
    In this case, \cref{prop:identification_tau_invertible_FilSp} identifies the generic fibre with $\Sp$, and \cref{thm:mod_tau_is_associated_graded} identifies the special fibre with $\grSp$.
\end{example}

\begin{example}[\cite{burklund_hahn_senger_Rmotivic}, Appendix~C.1]
    Let $R$ be a filtered $\E_1$-ring.
    Then the $\infty$-category $\Mod_R(\FilSp)$ is naturally a deformation.
    Its generic and special fibres are equivalent to, respectively,
    \[
        \Mod_{R^{\tau=1}}(\Sp) \qquad \text{and} \qquad \Mod_{R/\tau}(\grSp).
    \]
    If $R$ is a filtered $\E_\infty$-ring, then these equivalences are naturally symmetric monoidal.
    In fact, in this case the deformation $\Mod_R(\FilSp)$ is a \emph{symmetric monoidal deformation}, a concept to be introduced in \cref{ssec:monoidal_deformations} below.
\end{example}

\begin{remark}
    The use of the terms \emph{deformation}, \emph{special fibre} and \emph{generic fibre} is inspired by algebraic geometry.
    We think of $\tau$ as the deformation parameter, and the `deformation' is one of the special fibre to the generic fibre.
    Geometrically, $\FilSp$ plays the role of $\bA^1/\bG_m$.
    One can in fact make this comparison precise; see \cite{moulinos_geometry_filtrations}.
\end{remark}

The analogous statement to \cref{prop:pullback_tau_complete_and_invert_tau} holds in any deformation.
In this sense, the parameter~$\tau$ governs the large-scale structure of a deformation.

\begin{proposition}
    \label{prop:deformation_pullback_tau_complete_and_invert_tau}
    Let $\C$ be a deformation.
    For $X \in \C$, there is a natural pullback square
    \[
        \begin{tikzcd}
            X \ar[r] \ar[d] \pullback & X_\tau^{\wedge} \ar[d] \\
            X[\tau^{-1}] \ar[r] & (X_\tau^{\wedge})[\tau^{-1}].
        \end{tikzcd}
    \]
    In particular, a map $X \to Y$ in $\C$ is an isomorphism if and only if the maps
    \[
        X[\tau^{-1}] \to Y[\tau^{-1}] \qquad \text{and} \qquad C\tau \otimes X \to C\tau \otimes Y
    \]
    are both an isomorphism.
\end{proposition}
\begin{proof}
    This follows by tensoring the pullback square from \cref{prop:pullback_tau_complete_and_invert_tau} for the unit in $\FilSp$ with the object $X \in \C$.
    As tensoring preserves colimits in each variable, and we are in the stable setting, the resulting square is again a pullback square.
\end{proof}

\begin{remark}
    Because the inclusion $\C[\tau^{-1}]\subseteq \C$ admits both a left and a right adjoint, it follows from \cite[Proposition~A.8.20]{HA} that this equips $\C$ with the structure of a stable recollement
    \[
        \begin{tikzcd}
            \C[\tau^{-1}] \rar[hook] &[2em] \C \rar["(\blank)_\tau^\wedge"] \lar[bend left=50] \lar[bend right=50,"\tau^{-1}"'] &[2em] \C_\tau^\wedge \lar[bend left=50, hook'] \lar[bend right=50, hook']
        \end{tikzcd}
    \]
    and the pullback square of \cref{prop:deformation_pullback_tau_complete_and_invert_tau} is the one corresponding to this recollement. 
\end{remark}

A deformation is naturally enriched in filtered spectra; more loosely speaking, it is enriched in spectral sequences.

\begin{construction}
    \label{constr:deformation_enriched_FilSp}
    Let $\C$ be a deformation.
    Let $X,Y \in \C$.
    Recall the functor $i\colon \Z \to \FilSp$ from \cref{def:functor_Z_to_FilSp}.
    Define the \defi{filtered mapping spectrum} $\filmap_\C(X,Y)$ from $X$ to $Y$ as the filtered spectrum
    \[
        \map_\C(i(\blank)\otimes X,\, Y) \colon \Z^\op \to \Sp,
    \]
    where $\map_\C(\blank,\blank)$ denotes the mapping spectrum of $\C$.
    This is naturally functorial in $X$ and $Y$, leading to a functor $\filmap \colon \C^\op \times \C \to \FilSp$.
    Concretely, the $\filmap_\C(X,Y)$is given by
    \[
        \dotsb \to \map_\C(\opSigma^{0,1}X,\,Y) \to \map_\C(X,Y) \to \map_\C(\opSigma^{0,-1}X,\,Y) \to \dotsb
    \]
    with transition maps induced by $\tau$.
\end{construction}

\subsection{Monoidal deformations}
\label{ssec:monoidal_deformations}

One way to obtain the structure of a deformation on $\C$ is to give it the additional structure of an \emph{algebra} over $\FilSp$.
The universal property of $\FilSp$ provides a way to construct this.
In the following, we are careful with the distinction between monoidal and symmetric monoidal functors, because of the existence of important examples where these functors are not symmetric.
Nevertheless, in the case of synthetic spectra, these issues do not arise, so the distinction will not appear much in the later text.

Recall the symmetric monoidal functor $i\colon \Z \to\FilSp$ from \cref{def:functor_Z_to_FilSp} which, by \cref{rmk:interpretation_of_i_as_tau_tower}, is of the form
\[
    \begin{tikzcd}
        \dotsb \ar[r,"\tau"] & \S^{0,-1} \ar[r,"\tau"] & \S \ar[r,"\tau"] & \S^{0,1} \ar[r,"\tau"] & \dotsb.
    \end{tikzcd}
\]
The universal property of $\FilSp$ says that it is the universal presentable stable $\infty$-category on a diagram of this form.
In more loose terms, this says it is the universal category on the endomorphism $\tau$ of the unit.

\begin{proposition}[Universal property of filtered spectra]
    \label{thm:universal_property_FilSp}
    \leavevmode
    \begin{numberenum}
        \item \label{item:univ_prop_FilSp} Let $\C$ be a presentable stable $\infty$-category.
        Then there is an equivalence
        \[
            \Fun(\Z,\C) \simeq \LFun(\FilSp,\, \C)
        \]
        such that, if $f \colon \Z\to \C$ corresponds to $F \colon \FilSp\to \C$, then we have a natural isomorphism
        \[
            F\circ i \cong f.
        \]
        \item Let $\C$ be a presentably \textbr{symmetric} monoidal stable $\infty$-category.
        Then there is an equivalence
        \[
            \Fun^\otimes(\Z,\C) \simeq \LFun^\otimes(\FilSp,\, \C),
        \]
        where $\Fun^\otimes$ denotes the $\infty$-category of \textbr{symmetric} monoidal functors.
        Moreover, when forgetting the \textbr{symmetric} monoidal structure on the functors, this equivalence coincides with the equivalence from item~\ref{item:univ_prop_FilSp}, and the natural isomorphism $F \circ i \cong f$ is naturally \textbr{symmetric} monoidal.
    \end{numberenum}
\end{proposition}
\begin{proof}
    Note that we have a symmetric monoidal equivalence $\FilSp \simeq \Sp(\PSh(\Z))$, where $\PSh(\Z)$ is also equipped with the Day convolution symmetric monoidal structure.
    The first universal property is therefore the combination of the universal property of presheaves from \cite[Tag~\href{https://kerodon.net/tag/03W9}{03W9}]{kerodon} and the presentable universal property of stabilisation from \cite[Corollary~1.4.4.5]{HA}.
    The universal property of Day convolution in the monoidal (respectively, symmetric monoidal) case from \cite[Example~2.2.6.10]{HA} (respectively, Example~2.2.6.9 of op.\ cit.) then upgrades this to the second claimed equivalence.
\end{proof}

\begin{notation}
    \label{not:rho_sigma_general_deformation}
    Let $\C$ be a presentably monoidal $\infty$-category, and let $f \colon \Z \to \C$ be a monoidal functor.
    \begin{itemize}
        \item We typically reserve the letter $\rho$ for for the monoidal functor $\FilSp \to \C$ induced by~$f$.
        This functor in particular turns $\C$ into a (left) $\FilSp$-module, which informally is given by (where $A \in \FilSp$ and $X \in \C$)
        \[
            A \otimes X \defeq \rho(A) \otimes X.
        \]
        We summarise this by saying that $\rho$ gives $\C$ the structure of a \defi{monoidal deformation}.\footnotemark

        \item We typically reserve the letter $\sigma$ for the (lax monoidal) right adjoint $\C \to \FilSp$ to $\rho$.
        We call $\sigma$ the \defi{signature functor}; if $X \in \C$, then we refer to $\sigma X$ as the \defi{signature} of $X$.
        \item If $\C$ is symmetric monoidal and $f$ is a symmetric monoidal functor, then $\rho$ is naturally symmetric monoidal, and $\sigma$ is naturally lax symmetric monoidal.
        In this case, we say $\rho$ gives $\C$ the structure of a \defi{symmetric monoidal deformation}.
    \end{itemize}
\end{notation}

\footnotetext{This is a bad choice of terminology.
If $\calO$ is an $\infty$-operad, one should define an \emph{$\calO$-monoidal deformation} as an $\calO$-algebra in $\PrL_\st$.
A monoidal functor $f\colon \Z \to \C$ gives rise to an $\E_1$-monoidal functor $\FilSp \to \C$, but this does not equip $\C$ with the structure of an $\E_1$-algebra in $\PrL_\st$.
(This is the usual difference between an $\E_1$-algebra in $\Mod_A(\C)$ and an $\E_1$-map $A \to R$ in $\C$.)
However, in these notes, we are mostly concerned with the symmetric monoidal (i.e., $\E_\infty$) case where this distinction goes away, so this abuse of terminology does not have many ramifications for the rest of these notes.}

Roughly speaking, the functor $\rho$ is characterised by preserving colimits and by sending $\tau$ in $\FilSp$ to the map $f(-1 \to 0)$.
Thus, we can think of $f$ as a `prescription' for what the map $\tau$ in $\C$ ought to be.

\begin{example}
    Let $R$ be a filtered $\E_2$-ring.
    Then $\Mod_R(\FilSp)$ is naturally an $\E_1$-monoidal $\infty$-category, and the functor $R\otimes \blank \colon \FilSp \to \Mod_R(\FilSp)$ is a monoidal functor, turning $\Mod_R(\FilSp)$ into a monoidal deformation.
    If $R$ is a filtered $\E_\infty$-ring, then this turns it into a symmetric monoidal deformation.

    In a precise sense, this example includes many deformations: all symmetric monoidal deformations with compact unit that are generated by the spheres are of this form; see \cref{sec:filtered_models_syn}.    
\end{example}

Using that $\rho$ is monoidal, we can import structure from $\FilSp$ into $\C$.
For example, the $\E_\infty$-structure on $C\tau$ induces an $\E_1$-structure on $\rho(C\tau)$.
If $\rho$ is symmetric monoidal, then $\rho(C\tau)$ also acquires an $\E_\infty$-structure.

\begin{remark}
    \label{rmk:expression_for_sigma}
    Using that $\sigma$ is right adjoint to $\rho$, it follows that for $X \in \C$, the filtered spectrum $\sigma X$ is given by
    \[
        \dotsb \to \map_\C(f(1), X) \to \map_\C(f(0), X) \to \map(f(-1),X) \to \dotsb,
    \]
    where $\map_\C(\blank,\blank)$ denotes the mapping spectrum of $\C$, and where the transition maps are induced by $f$.
    We see that, after forgetting the lax monoidal structure on~$\sigma$, it is naturally isomorphic to $\filmap(\1_\C,\blank)$ from \cref{constr:deformation_enriched_FilSp}.
\end{remark}

\begin{remark}
    It is difficult in general to obtain symmetric monoidal functors out of $\Z$, as it is not free as a symmetric monoidal $\infty$-category.
    We learned from Shaul Barkan that t-structures are a source of such functors, using the Whitehead filtration; see \cite[Section~2]{barkan_vN_cellular_chromatic}.
    We will use this in \cref{sec:Syn_as_deformation} to give synthetic spectra the structure of a symmetric monoidal deformation.
\end{remark}

So far, we have used a deformation structure on an $\infty$-category $\C$ to define similar-looking operations in $\C$ as we have in $\FilSp$, such as inverting and modding out by $\tau$.
In the monoidal case, the functor $\sigma$ shows that these operations translate back to the respective operations in $\FilSp$.
This both gives a more concrete interpretation of these operations, and also ties the story of deformations into the study of spectral sequences.

To prove the desired properties of $\sigma$, we use the following lemma.
We refer to, e.g., \cite[Definition~3.5]{naumann_pol_ramzi_monoidal_fracture} for a definition of the projection map, and Section~3 of op.\ cit.\ for an introduction to these ideas.

\begin{lemma}
    \label{lem:projection_formula_deformation}
    Let $\C$ be a monoidal deformation.
    Let $A \in \FilSp$ and $X \in \C$.
    Consider the natural \emph{projection map}
    \[
        A \otimes \sigma(X) \to \sigma(\rho(A) \otimes X).
    \]
    \begin{numberenum}
        \item \label{item:proj_formula_dualisable} If $A$ is dualisable, then the projection map is an isomorphism for all $X \in \C$.
        \item \label{item:proj_formula_all} If $\sigma$ preserves colimits, then the projection map is an isomorphism for all $A \in \FilSp$ and all $X \in \C$.
    \end{numberenum}
\end{lemma}
\begin{proof}
    The first part follows from \cite[Lemma~3.8\,(b)]{naumann_pol_ramzi_monoidal_fracture}.
    The second part follows from the fact that, in this case, both sides preserve colimits in $A$ and that $\FilSp$ is generated under colimits by dualisable objects.
\end{proof}

In particular, the functor $\sigma$ commutes with bigraded suspensions.


\begin{theorem}
    \label{prop:deformation_sigma_tauinv}
    Let $\C$ be a \textbr{symmetric} monoidal deformation which arises from a \textbr{symmetric} monoidal functor $f\colon \Z \to \C$.
    \begin{numberenum}
        \item \label{item:sigma_preserves_tau} For every $X \in \C$, we have a natural isomorphism
        \[
            \tau_{\sigma X} \cong \sigma(\tau_X).
        \]
        In particular, $\sigma$ sends $\tau$-invertible objects in $\C$ to $\tau$-invertible \textbr{a.k.a.\ constant} filtered spectra.
        \item \label{item:sigma_pres_modtau_taucompl} There are commutative diagrams of lax \textbr{symmetric} monoidal functors
        \[
            \begin{tikzcd}
                \C \rar["\sigma"] \dar["(\blank)/\tau"'] & \FilSp \dar["(\blank)/\tau"] \\
                \Mod_{C\tau}(\C) \rar["\sigma"] & \Mod_{C\tau}(\FilSp).
            \end{tikzcd}
            \qquad \text{and} \qquad
            \begin{tikzcd}
                \C \rar["\sigma"] \dar["(\blank)_\tau^\wedge"'] & \FilSp \dar["(\blank)_\tau^\wedge"] \\
                \C_\tau^\wedge \rar["\sigma"] & \FilSp_\tau^\wedge.
            \end{tikzcd}
        \]
        In particular, for $X \in \C$, we have natural isomorphisms
        \[
            \sigma(C\tau \otimes X) \cong C\tau \otimes \sigma(X) \qquad \text{and} \qquad \sigma(X_\tau^\wedge) \cong \sigma(X)_\tau^\wedge.
        \]
        \item \label{item:colims_sigma_condition} The functor $\sigma$ preserves small colimits if and only if the monoidal unit of $\C$ is compact.

        If this happens, then there is a commutative diagram of lax \textbr{symmetric} monoidal functors
        \[
            \begin{tikzcd}
                \C \rar["\sigma"] \dar["\tau^{-1}"'] & \FilSp \dar["\tau^{-1}"] \\
                \C[\tau^{-1}] \rar["\sigma"] & \FilSp[\tau^{-1}].
            \end{tikzcd}
        \]
        In particular, if $\sigma$ preserves colimits, then for $X \in \C$, we have a natural isomorphism
        \[
            \sigma(X[\tau^{-1}]) \cong \sigma(X)[\tau^{-1}].
        \]
        \item \label{item:conservative_sigma_condition} The functor $\sigma$ is conservative if and only if the image of $f$ generates $\C$ as a stable $\infty$-category under colimits.
    \end{numberenum}
\end{theorem}
\begin{proof}
    \Cref{item:sigma_preserves_tau} and the $C\tau$-part of \cref{item:sigma_pres_modtau_taucompl} follow from \cref{lem:projection_formula_deformation}\,\ref{item:proj_formula_dualisable}.
    To prove that $\sigma$ also preserves $\tau$-completion, we need to check that $\sigma$ preserves $C\tau$-equivalences and preserves $C\tau$-local objects.
    The first again follows from the projection formula, and the second follows because $\rho$ preserves $C\tau$-acyclics (being an $\FilSp$-linear functor).
    

    For \cref{item:colims_sigma_condition}, note that since $\sigma$ is exact, it preserves colimits if and only if it preserves filtered colimits.
    The latter is equivalent to its left adjoint $\rho$ preserving compact objects.
    For this, it is equivalent to check that $\rho$ sends a collection of compact generators of $\FilSp$ to compact objects of $\C$.
    As the filtered spheres $\S^{0,s}$ for $s\in \Z$ form stable generators for $\FilSp$, and $\rho(\S^{0,s}) \cong f(s)$, we see that $\sigma$ preserves colimits if and only if $f$ lands in compact objects of $\C$.
    Since $f$ is monoidal and hence sends the unit to the unit, the latter condition implies that the unit of $\C$ is compact.
    Conversely, if the unit of $\C$ is compact, then all dualisable objects of $\C$ are compact.
    Because $f$ is monoidal and all objects of $\Z$ are dualisable, it then follows that all values of $f$ are compact.
    We conclude that $\sigma$ preserves colimits if and only if the unit of $\C$ is compact.

    If $\sigma$ preserves colimits, then \cref{lem:projection_formula_deformation}\,\ref{item:proj_formula_all} implies that $\sigma$ preserves $\tau$-inversion.
    
    Finally, for \cref{item:conservative_sigma_condition}, we use again that the filtered spheres are generators, so that $\sigma X$ is zero if and only if for all $s\in \Z$, the mapping spectrum
    \[
        \map_\FilSp(\S^{0,s},\, \sigma X) \cong \map_\C(\rho(\S^{0,s}),\, X) \cong \map_\C(f(s),\, X)
    \]
    vanishes.
    This shows the final claim.
\end{proof}

\begin{remark}
    \label{rmk:sigma_internal_right_adjoint}
    As explained in \cite[Definition~3.9, Remark~3.10]{naumann_pol_ramzi_monoidal_fracture}, it follows from \cref{prop:deformation_sigma_tauinv}\,\ref{item:colims_sigma_condition} that $\rho$ is an \emph{internal left adjoint} in $\FilSp$-linear $\infty$-categories (i.e., its right adjoint $\sigma$ is itself an $\FilSp$-linear functor) if and only if the unit of $\C$ is compact.
\end{remark}

We now introduce the notion of the signature spectral sequence.
We will be more brief, as we will develop this in detail in the next chapter in the case of synthetic spectra (which is our main case of interest).

\begin{definition}
    Let $\C$ be a monoidal deformation.
    If the unit of $\C$ is compact, then the functor $\pi_{n,s}\colon \C \to \Ab$ defined by
    \[
        \pi_{n,s}(\blank) \defeq [\opSigma^{n,s}\1_\C,\ \blank]
    \]
    preserves filtered colimits.
    It follows that for $X \in \C$, the spectral sequence underlying $\sigma X$ is of the form
    \[
        \uE_1^{n,s} = \oppi_{n,s}(X/\tau) \implies [\opSigma^n \1_\C,\ X[\tau^{-1}]] \cong \pi_{n,*}(X)[\tau^{-1}].
    \]
    Accordingly, we refer to this as the \defi{signature spectral sequence} of $X$.
\end{definition}

\Cref{prop:deformation_sigma_tauinv} tells us that we can understand this spectral sequence $\sigma X$ through computations in $\C$.
For instance, \cref{item:sigma_pres_modtau_taucompl} tells us that $\sigma X$ is conditionally convergent whenever $X$ is $\tau$-complete, and that the associated graded is given by applying~$\sigma$ to $X/\tau$.
This is particularly useful in cases where $\Mod_{C\tau}(\C)$ is simpler than the original case of filtered spectra, where $\Mod_{C\tau}(\FilSp)\simeq\grSp$ is topological in nature.

The above also tells us that the Omnibus Theorems from \cref{sec:filtered_omnibus} carry over to the deformation $\C$, with all occurrences of the underlying spectral sequence replaced by the signature spectral sequence.
We have to be a little careful if the unit of $\C$ is not compact, as then $\sigma$ need not preserve $\tau$-inversion.
Even if the unit of $\C$ is not compact, the truncated Omnibus Theorem of \cref{thm:truncated_filtered_omnibus,thm:truncated_filtered_omnibus_generation} apply in $\C$ without change.
In the non-truncated case of \cref{thm:filtered_omnibus}, all results that do not compare $\oppi_{*,*}X$ with $\pi_* X^{\tau=1}$ apply in $\C$ as well.
If the unit of $\C$ is compact, then all of the Omnibus Theorems apply in their entirety.
In summary then, we learn that the $\Z[\tau]$-module $\oppi_{*,*}X$ captures the signature spectral sequence of $X$.


Not just the Omnibus Theorem carries over to a deformation, but also the $\tau$-Bockstein spectral sequence we used to prove it.

\begin{variant}
    \label{var:tau_BSS_deformation}
    Let $\C$ be a monoidal deformation, and $X \in \C$.
    The \defi{$\tau$-adic filtration} on $X$ is the filtered object $\Z^\op \to \C$ given by
    \[
        \begin{tikzcd}
            \dotsb \ar[r,"\tau"] & \opSigma^{0,-2} X \ar[r,"\tau"] & \opSigma^{0,-1}X \ar[r,"\tau"] & X \ar[r,equals] & \dotsb.
        \end{tikzcd}
    \]
    Analogously to \cref{constr:trigraded_sseq_bifiltered_spectrum}, this leads to a trigraded spectral sequence which we call the \defi{$\tau$-Bockstein spectral sequence} of $X$, which is of the form
    \[
        \uE_1^{n,w,s} \cong \oppi_{n,\, w+s}(X/\tau) \implies \oppi_{n,w} X.
    \]
    From \cref{prop:deformation_sigma_tauinv}, it follows that when we apply $\sigma$ to the $\tau$-adic filtration on $X$, we obtain the $\tau$-adic filtration on $\sigma X$.
    It follows that the $\tau$-BSS of $X$ is the $\tau$-BSS of $\sigma X$.
    We may therefore freely use the results of \cref{sec:tau_Bockstein_sseq} for this spectral sequence; in particular, it captures the signature spectral sequence of $X$.
\end{variant}

Two structural properties of deformations, namely \emph{cellularity} and \emph{evenness}, will be discussed later in \cref{ssec:cellular,ssec:even_synthetic_spectra}.
For now, we finish this chapter with a few examples of deformations, which are of a different flavour than the one we will meet in the next chapter.

We learned the following example from Christian Carrick and Lennart Meier.
See also \cite[Examples~A.8 and~A.9]{burklund_hahn_senger_Rmotivic}.

\begin{example}
    Let $\Sp_{C_2}$ denote the $\infty$-category of genuine $C_2$-spectra.
    The Euler class $a_\sigma$ gives $\Sp_{C_2}$ the structure of a deformation.
    (To avoid notational confusion, we will avoid the names $\rho$ and $\sigma$ for the deformation functors.)
    More specifically, let $\sigma$ denote the sign representation of $C_2$.
    Then the inclusion of fixed points results in a map $S^0 \to S^\sigma$, which stably results in a map called the \emph{Euler class}
    \[
        a_\sigma \colon \S^{-\sigma} \to \S.
    \]
    We expect, but do not check in detail, that this assembles to a monoidal functor
    \[
        \begin{tikzcd}
            \dotsb \rar["a_\sigma"] & \S^{-\sigma} \rar["a_\sigma"] & \S \rar["a_\sigma"] & \S^{\sigma} \rar["a_\sigma"] & \dotsb.
        \end{tikzcd}
    \]
    (Note, however, that this cannot be made symmetric monoidal, due to a nontrivial switch map for $\S^\sigma\otimes\S^\sigma$.)
    This monoidal functor induces a monoidal left adjoint $\FilSp\to\Sp_{C_2}$, resulting in the promised monoidal deformation structure.
    The generic fibre is equivalent to $\Sp$.
    Interestingly, its special fibre is also spectra: it is given by modules over $\S/a_\sigma \cong \Sigma^\infty_+ C_2$, which by \cite[Theorem~1.1]{balmer_dellambrogio_sanders} is equivalent to $\Sp$ (via a lift of the restriction-coinduction adjunction).

    Under these identifications, $\tau$-inversion is identified with geometric fixed points, and $\tau$-completion is identified with Borel completion (i.e., inverting those maps that induce an isomorphism on homotopy fixed points).
    As a result, the pullback square of \cref{prop:deformation_pullback_tau_complete_and_invert_tau} becomes, after taking (genuine) fixed points, the \emph{Tate square}
    \[
        \begin{tikzcd}
            X^{C_2} \rar \dar \pullback & X^{hC_2} \dar \\
            X^{\Phi C_2} \rar & X^{tC_2}.
        \end{tikzcd}
    \]
    The resulting $\tau$-BSS is identified with the $a_\sigma$-BSS, which is the homotopy fixed-point spectral sequence.

    One could mimic these constructions for $\Sp_G$ for a general finite group $G$ in the place of $C_2$.
    However, the resulting deformation would be rather contrived: for general $G$, the structure of $\Sp_G$ is better captured by taking all subgroups of $G$ into account.
    Only if $G=C_2$ is the resulting structure exactly that of a deformation.
\end{example}

Finally, we end with an example that is not strictly speaking a deformation, but which is close enough in that it also allows theorems from filtered spectra to be imported over.

\begin{example}[Recovering the $p$-Bockstein spectral sequence]
    \label{ex:recover_p_BSS}
    Fix a prime $p$, and consider the functor $f\colon \Z \to \Sp$ given by
    \[
        \begin{tikzcd}
            \dotsb \rar["p"] & \S \rar["p"] & \S \rar["p"] & \S \rar["p"] & \dotsb.
        \end{tikzcd}
    \]
    This functor induces an adjunction
    \[
        \begin{tikzcd}[column sep=3.5em]
            \FilSp \rar[shift left,"\rho"] & \Sp. \lar[shift left, "\sigma"]
        \end{tikzcd}
    \]
    However, neither $f$ nor $\rho$ can be made monoidal.
    Indeed, the functor $\rho$ sends $\tau$ to~$p$, and hence sends $C\tau$ to $\S/p$.
    The latter does not admit an $\E_1$-structure for any $p$, preventing $\rho$ (and hence $f$) from being monoidal.
    
    While strictly speaking not a deformation, we can still use this adjunction to import information from the $\tau$-Bockstein spectral sequence.
    Namely, we can check by hand that $\sigma$ sends $\tau_X$ to $\tau_{\sigma X}$.
    In $\Sp$, the map $\tau$ is given by $p$, which $\sigma$ sends to $p$ because it is additive.
    Further, by \cref{rmk:expression_sigma}, the functor $\sigma$ sends a spectrum $X$ to
    \[
        \begin{tikzcd}
            \dotsb \rar["p"] & X \rar["p"] & X \rar["p"] & X \rar["p"] & \dotsb
        \end{tikzcd}
    \]
    so indeed $\sigma$ sends $\tau_X$ to $\tau_{\sigma X}$.
    It follows that $\sigma$ sends the $p$-Bockstein filtration on~$X$ to the $\tau$-Bockstein filtration on $\sigma X$.
    The resulting trigraded spectral sequence we compute to be the $p$-BSS for $X$ with an additional filtration tagged on.
    More specifically, at every level of this new filtration, it is the $p$-BSS for $X$, and the transition maps for this new filtration are all given by multiplication by $p$.
    Using this, we can deduce the analogous version of \cref{thm:tau_BSS_structure} for the $p$-BSS, at least those parts that do not depend on any monoidality properties of $\tau$ (such as \cref{thm:tau_BSS_structure}\,\ref{item:tau_BSS_detection}).
    

    For a monoidal version of this example, consider instead the functor $\Z \to \D(\Ab)$ given by
    \[
        \begin{tikzcd}
            \dotsb \rar["p"] & \Z \rar["p"] & \Z \rar["p"] & \Z \rar["p"] & \dotsb.
        \end{tikzcd}
    \]
    This functor is symmetric monoidal: it lands in the heart of $\D(\Ab)$, so equivalently is given by a functor $\Z\to\Ab$.
    It is easy to check that this functor describes a strict, multiplicative filtration on the commutative ring $\Z[\tfrac{1}{p}]$, and as a result is canonically a symmetric monoidal functor.
    It follows that this gives $\D(\Ab)$ the structure of a symmetric monoidal deformation.
    We can therefore directly deduce the analogue of \cref{thm:tau_BSS_structure} by using \cref{prop:deformation_sigma_tauinv}, thereby recovering, e.g., \cite[Theorem~3.8]{palmieri_bockstein_sseq}.

    Informally, we may summarise both situations by saying that we ``put $\tau$ equal to~$p$'' and thereby recover the $p$-Bockstein spectral sequence.
    However, this slogan should be taken with a grain of salt, as it ignores the monoidality issues raised above, which depend on the specific category one is working with.
    In particular, only in the monoidal case will we be able to import multiplicative properties of the $\tau$-BSS.
\end{example}

%% file: syn/syn.tex
Previously in \cref{sec:Adams_and_Tot_sseqs}, we discussed the classical definition of the Adams spectral sequence.
The more modern way to interact with Adams spectral sequences is to use \emph{synthetic spectra}.
This chapter is intended both as a first introduction to and as a manual for working with synthetic spectra.
Compared to most of the existing literature, our distinctive focus is to understand these through the lens of filtered spectra.

We review the main categorical features of synthetic spectra in \cref{sec:Syn_categorical_basics,sec:homological_t_structure}.
There, among other things, we encounter the synthetic map $\tau$.
In \cref{sec:Syn_as_deformation}, we show that this gives the $\infty$-category of synthetic spectra the structure of a (symmetric monoidal) deformation in the sense of \cref{ssec:monoidal_deformations}, and we study the structure of this deformation.
In particular, this results in a functor
\[
    \sigma \colon \Syn_E \to \FilSp
\]
that preserves limits and colimits and sends $\tau$ to the filtered map $\tau$.
We refer to $\sigma X$ as the \emph{signature} of~$X$.
Consequently, we obtain a synthetic Omnibus Theorem, describing the homotopy groups of a synthetic spectrum $X$ in terms of the spectral sequence underlying $\sigma X$.

The question of which spectral sequence this is requires a computation.
In \cref{sec:signature_synthetic_analogue}, our goal is to compute this for so-called \emph{synthetic analogues}; this is both an important foundational result, and also showcases how to work with synthetic spectra.
There is a lax symmetric monoidal functor
\[
    \nu \colon \Sp \to \Syn_E
\]
called the \emph{synthetic analogue functor}.
We show that the signature of the synthetic analogue of a spectrum is its $E$-Adams spectral sequence; see \cref{sec:signature_synthetic_analogue}, particularly \cref{thm:signature_nu}.
Finally, in \cref{sec:synthetic_omnibus} we collect some implications of this result, and briefly discuss some notational conventions.

\begin{remark}
    In addition to providing a new interface for interacting with Adams spectral sequences, synthetic spectra have also been used in setting up obstruction theories; see \cite{barkan_monoidal_algebraicity,Hopkins_Lurie_Brauer_E,pstragowski_vankoughnett_obstruction_theory}.
    We do not (yet) discuss these obstruction theories in this version of the notes.
\end{remark}

For the most part, this chapter consists of an overview of results from \cite{pstragowski_synthetic}.
We learned much of the relationship with filtered spectra from \cite{barkan_monoidal_algebraicity}, \cite[Appendix~C]{burklund_hahn_senger_Rmotivic}, \cite[Section~1]{christian_jack_synthetic_j}, and \cite{pstragowski_perfect_even_filtration}.
Some results in this chapter are modifications of results appearing in \cite[Section~1]{CDvN_part1}.

\section{Categorical properties}
\label{sec:Syn_categorical_basics}

Before we can do more serious work, we need to know the basic properties and structure of the category we are dealing with.
We do not review the constructions given by Pstr\k{a}gowski, but content ourselves with summarising the main properties, preferring to work with synthetic spectra in a `model-independent way'.\footnote{We mean this in a loose way: we specifically choose the construction from \cite{pstragowski_synthetic} over the other available ones. However, phrasing our arguments and computations in this way should make them more robust and more easily adapted to other synthetic contexts.}

\begin{construction}
    Let $E$ be a homotopy associative ring spectrum of Adams type.
    In \cite{pstragowski_synthetic}, Pstr\k{a}gowski constructs a symmetric monoidal $\infty$-category $\Syn_E$ of \defi{$E$-based synthetic spectra}, together with a unital lax symmetric monoidal functor $\nu \colon \Sp \to \Syn_E$.
    We call $\nu$ the \defi{synthetic analogue functor}.
\end{construction}

We may refer to $E$-based synthetic spectra as \emph{$E$-synthetic spectra}, or even simply by \emph{synthetic spectra} if $E$ is clear from the context.
On the opposite end, when we want to vary the variable $E$, we would write $\nu_E$ for $\nu$, emphasising it as the \emph{$E$-synthetic analogue}.
It would be more principled to write $\Syn_E(\Sp)$ instead of $\Syn_E$, indicating that it is a modification of the $\infty$-category $\Sp$.
For simplicity, we will stick to the shorter name in this chapter.

\begin{notation}
    An \defi{$E$-synthetic $\E_\infty$-ring} is an $\E_\infty$-algebra object in $\Syn_E$.
    If $E$ is clear from the context, then we may also refer to such an object as a \emph{synthetic $\E_\infty$-ring}.
\end{notation}

\begin{remark}
    Just as with the Adams spectral sequence, the role of the $\infty$-category of spectra here is not of fundamental importance: there should be a similar modification of any nice enough stable $\infty$-category with a type of Adams spectral sequence.
    We stick to the case of spectra here to more conveniently cite Pstr\k{a}gowski's construction.
    A more general theory can be found in \cite{patchkoria_pstragowski_derived_inftycats}, but this construction differs from the one in \cite{pstragowski_synthetic} in various ways; see \cite[Section~6.5]{patchkoria_pstragowski_derived_inftycats}.
\end{remark}

\begin{remark}
    \label{rmk:SynE_depends_on_less_than_E}
    Although the notation seems to suggest otherwise, the symmetric monoidal $\infty$-category $\Syn_E$ depends on much less data than the ring spectrum~$E$.
    This is because the Adams spectral sequence depends on less data than the ring spectrum~$E$; see \cref{rmk:millers_Adams_sseq}.
    In particular, $\Syn_E$ is not sensitive to a potential coherent multiplicative structure on $E$, nor does it require it for its construction as a symmetric monoidal $\infty$-category.
\end{remark}

\begin{warning}
    \label{warn:E_should_be_Adams_type_SynE}
    The $\infty$-category $\Syn_E$ is functorial in $E$, and in fact the construction does not require the assumption that $E$ is of Adams type.
    However, it only gives the correct answer if $E$ is of Adams type.
    For instance, as observed by Pstr\k{a}gowski and explained by Sch\"appi in \cite[Theorem~2.3.7]{schappi_flat_replacements}, taking $E = \MU$ or $E=\Z$ results in the \emph{same} symmetric monoidal $\infty$-category, even though the $\Z$-Adams spectral sequence is wildly different from the $\MU$-Adams spectral sequence (see \cref{ex:adams_type}).
\end{warning}

Because of this warning, throughout this chapter we stick to the following assumption.

\begin{notation}
    For the remainder of this chapter, $E$ denotes a fixed choice of a homotopy-associative ring spectrum of Adams type.
\end{notation}

We begin by studying some categorical properties of synthetic spectra.
Recall the notion of a \emph{finite $E$-projective spectrum} from \cref{def:Adams_type}.

\begin{proposition}
    \label{prop:properties_of_Syn}
    \leavevmode
    \begin{numberenum}
        \item The $\infty$-category $\Syn_E$ is stable.\label{item:Syn_stable}
        \item The $\infty$-category $\Syn_E$ is presentable, and the symmetric monoidal structure preserves colimits in each variable separately; that is to say, $\Syn_E$ is presentably symmetric monoidal.\label{item:Syn_pres_mon}
        \item \label{item:fp_spectra_are_compact_dualisable_generators} If $P$ is a finite $E$-projective spectrum, then $\nu P$ is a compact and dualisable object in $\Syn_E$, with dual $\nu(P^\vee)$.
        In particular, the monoidal unit is compact.
        \item \label{item:Syn_generated_by_spheres} As a stable $\infty$-category, $\Syn_E$ is compactly generated under colimits by the synthetic analogues of finite $E$-projectives.
        That is, the collection of $\opSigma^k \nu P$, for $P$ finite $E$-projective and $k\in \Z$, forms a set of compact dualisable generators.
        In particular, $\Syn_E$ is compactly generated by dualisables.
        \item The monoidal $\infty$-category is \emph{rigid} in the sense that an object is compact if and only if it is dualisable.\label{item:Syn_ridigly_compactly_gen}
    \end{numberenum}
\end{proposition}
\begin{proof}
    \Cref{item:Syn_stable,item:Syn_pres_mon} are \cite[Proposition~4.2]{pstragowski_synthetic}, and \cref{item:fp_spectra_are_compact_dualisable_generators,item:Syn_generated_by_spheres} are \cite[Remark~4.14]{pstragowski_synthetic}.
    \Cref{item:Syn_ridigly_compactly_gen} then follows formally from the fact that the unit is compact and that it has a set of compact dualisable generators; see, e.g., (the footnote to) \cite[Terminology~4.8]{naumann_pol_ramzi_monoidal_fracture}.
\end{proof}

Next, we turn to properties of the functor $\nu$.

\begin{proposition}
    \label{prop:properties_of_nu}
    \leavevmode
    \begin{numberenum}
        \item \label{item:categorical_properties_nu} The functor $\nu \colon \Sp \to \Syn_E$ is fully faithful, additive, and preserves filtered colimits.
        In particular, $\nu$ preserves arbitrary coproducts.
        \item \label{item:nu_and_cofibre_seqs} Consider a cofibre sequence of spectra
        \[
            \begin{tikzcd}
                X \ar[r,"f"] & Y \ar[r,"g"] & Z.
            \end{tikzcd}
        \]
        Then the induced sequence
        \[
            \begin{tikzcd}
                \nu X \ar[r,"\nu f"] & \nu Y \ar[r,"\nu g"] & \nu Z
            \end{tikzcd}
        \]
        is a cofibre sequence of synthetic spectra if and only if
        \[
            \begin{tikzcd}
                0 \ar[r] & E_*X \ar[r,"f_*"] & E_*Y \ar[r,"g_*"] & E_*Z \ar[r] & 0
            \end{tikzcd}
        \]
        is short exact, or in other words, if the boundary map $Z \to \Sigma X$ is zero on $E_*$-homology.
        \item \label{item:strong_monoidality_nu} The comparison map $\nu X \otimes \nu Y \to \nu(X\otimes Y)$ coming from the lax monoidal structure on $\nu$ is an isomorphism whenever $X$ or $Y$ is a filtered colimit of finite $E$-projective spectra.
        
        More generally, if the $E_*$-homology of $X$ or $Y$ is flat as an $E_*$-module, then the map $\nu X\otimes \nu Y \to \nu(X\otimes Y)$ is a $\nu E$-equivalence.
    \end{numberenum}
\end{proposition}
\begin{proof}
    \Cref{item:categorical_properties_nu} follows from \cite[Lemma~4.4 and Corollary~4.38]{pstragowski_synthetic}, \cref{item:nu_and_cofibre_seqs} is \cite[Lemma~4.23]{pstragowski_synthetic}, and \cref{item:strong_monoidality_nu} is \cite[Lemma~4.24]{pstragowski_synthetic}.
\end{proof}

Both conditions of \cref{prop:properties_of_nu}\,\ref{item:strong_monoidality_nu} are a type of flatness condition.
This is obvious for the second one.
For the first, compare this with the algebraic result that a module over a ring is flat if and only if it can be written as a filtered colimit of finite free modules; see \cite[\href{https://stacks.math.columbia.edu/tag/058G}{Tag 058G}]{stacks_project}.

\begin{example}
    \label{ex:symm_mon_if_one_is_E}
    The definition of Adams type directly implies that $\nu E \otimes X \to \nu(E \otimes X)$ is an isomorphism for all spectra $X$.
\end{example}

\begin{example}
    Suppose $E = \F_p$, or more generally that $E$ is a homotopy-associative ring spectrum such that $\pi_*E$ is a graded field.
    Then every finite spectrum is finite $E$-projective.
    The smallest subcategory of $\Sp$ that contains all finite spectra and is closed under filtered colimits is equal to all of $\Sp$.
    We therefore learn from \cref{prop:properties_of_nu}\,\ref{item:strong_monoidality_nu} that $\nu$ is actually a strong symmetric monoidal functor if $\pi_*E$ is a graded field.
\end{example}

\begin{remark}
    The reader may gain intuition for the above properties by thinking of $\nu$ as being similar to the Whitehead filtration functor $\Wh \colon \Sp \to \FilSp$.
    In fact, this is more than a formal analogy: in the case $E=\S$, the $\infty$-category $\Syn_\S$ is equivalent to $\Mod_{\Wh \S}(\FilSp)$, and this equivalence identifies $\nu$ with the Whitehead filtration functor; see \cref{cor:sphere_synthetic_filtered_model}.
\end{remark}

Take particular note that $\nu$ is \emph{not} an exact functor, even though it is a functor between stable $\infty$-categories.
For example, \cref{prop:properties_of_nu}\,\ref{item:nu_and_cofibre_seqs} implies that $\Sigma(\nu X) \cong \nu(\Sigma X)$ if and only if $E_*X=0$.
In terms of the $E$-Adams spectral sequences, having vanishing $E$-homology is a very degenerate case, so the functor $\nu$ practically never preserves suspensions.

The difference between suspending in spectra and in synthetic spectra has a conceptual meaning as well: the former has the effect of shifting its Adams spectral sequence one to the right, while suspending its synthetic analogue also shifts it down by one filtration.
This is made precise by the following definition of the \emph{synthetic bigraded spheres}.
The indexing convention we use here turns out to be the most practical; we defer a more detailed explanation to \cref{rmk:why_second_page_indexing_Syn}.

\begin{definition}[Synthetic bigraded spheres]
    \label{def:synthetic_bigraded_spheres}
    Let $n$ and $s$ be integers.
    \begin{numberenum}
        \item The \defi{synthetic $(n,s)$-sphere} is
        \[
            \S^{n,s} \defeq \Sigma^{-s}\opnu(\S^{n+s}).
        \]
        We refer to $n$ as the \defi{stem}, and to $s$ as the \defi{filtration}.
        \item We write $\Sigma^{n,s} \colon \Syn_E \to \Syn_E$ for the functor given by tensoring with $\S^{n,s}$ on the left.
        \item We write $\pi_{n,s}\colon \Syn_E \to \Ab$ for the functor
        \[
            \pi_{n,s}(\blank) := [\S^{n,s},\ \blank].
        \]
        \item \label{item:def_syn_tau} The map $\tau \colon \S^{0,-1} \to \S^{0,0}$ is the colimit-comparison map
        \[
            \tau \colon \S^{0,-1} = \Sigma(\nu \S^{-1}) \to \nu\S =\S^{0,0}.
        \]
        If $X$ is a synthetic spectrum, then tensoring it with the map $\tau\colon \S^{0,-1}\to \S$ results in a map $\Sigma^{0,-1}X \to X$, which we denote by $\tau_X$, or by $\tau$ when there is no risk of confusion.
    \end{numberenum}
\end{definition}

Note that $\pi_{*,*}$ naturally lifts to a functor $\Syn_E \to \Mod_{\Z[\tau]}(\bigrAb)$.

\begin{remark}
    If $X$ is a spectrum, then we also have the natural colimit-comparison map $\Sigma \nu (\Sigma^{-1}X) \to \nu X$.
    This coincides with the map $\tau \otimes \nu X$ by \cite[Proposition~4.28]{pstragowski_synthetic}.
\end{remark}

\begin{remark}[Cellularity]
    \label{rmk:cellularity_first_advertisement}
    It is not necessarily true that bigraded homotopy groups detect isomorphisms of synthetic spectra.
    If this is the case, we say that $\Syn_E$ is \emph{cellular}.
    For many $E$, this is the case.
    We discuss this issue more in \cref{ssec:cellular}.
    For applications to spectral sequences, one can equally well work with the cellularisation of $\Syn_E$, so we view this as a technicality.
\end{remark}

Remembering the precise definition of $\S^{n,s}$ is not the most important; it is enough to remember the following key facts.

\begin{example}
    \label{ex:synthetic_bigraded_spheres}
    \leavevmode
    \begin{numberenum}
        \item For every $n$, the synthetic spectrum $\nu(\S^n)$ is the bigraded sphere $\S^{n,0}$.
        This is the first instance where we see that $\nu$ places everything in \emph{Adams filtration zero}; see \cref{sec:synthetic_omnibus} below for a further discussion.
        Later we will also see that these are the only synthetic spheres that are in the essential image of~$\nu$: see \cref{ex:connectivity_synthetic_analogue_explained}.
        
        We will abuse notation and abbreviate $\S^{0,0}$ simply by $\S$, and refer to it as the \defi{synthetic sphere}.
        Most of the time, the context will allow one to infer whether the sphere spectrum or the synthetic sphere is meant by this notation.

        \item More generally, if $X$ is a spectrum, then we have a natural isomorphism
        \[
            \opSigma^{n,0} \nu X \cong \nu(\Sigma^n X).
        \]
        Indeed, this follows since $\nu$ is strong symmetric monoidal when one factor is a sphere; see \cref{prop:properties_of_nu}\,\ref{item:strong_monoidality_nu}.
        
        \item \label{item:coordinates_synthetic_suspension} Categorical suspension is given by the bigraded suspension $\Sigma^{1,-1}$.\qedhere
    \end{numberenum}
\end{example}

\begin{remark}[Koszul sign rule]
    It is possible to set up matters so that if $A$ is a homotopy-commutative algebra in $\Syn_E$, then $\oppi_{*,*}A$ becomes a bigraded ring with a Koszul sign rule according to the first variable (i.e., the stem).
    Doing this involves choices, as explained by Dugger \cite{dugger_coherence_invertible_objects}, see also \cite{motivic_sign_rules_explained}.
    The choice described by Pstr\k{a}gowski in \cite[Remark~4.10]{pstragowski_synthetic}, and explained in detail by Chua in \cite[Section~6]{chua_E3_ASS}, results in this sign on homotopy groups.
\end{remark}

As with any symmetric monoidal stable $\infty$-category, synthetic spectra have an internal notion of homology.
We regard it as a bigraded object.

\begin{notation}
    Let $A$ and $X$ be synthetic spectra, and let $n$ and $s$ be integers.
    We write $A_{n,s}(X)$ for $\pi_{n,s}(A\otimes X)$.
\end{notation}

Finally, we make a few comments about notation and indexing compared to the literature.

\begin{remark}
    \label{rmk:lambda_instead_of_tau}
    In the specific case of $\F_p$-synthetic spectra, it is becoming more and more common to use the letter $\lambda$ to denote the map otherwise denoted by $\tau$; see, e.g., \cite{burklund_isaksen_xu_F2synthetic} (particularly Section~1.1 therein).
    This is done to allow for computations that involve both $\BP$-synthetic and $\F_p$-synthetic arguments at the same time.
    Because these notes are not aimed at these computations, we will still use the letter $\tau$ even in the $\F_p$-synthetic case.
\end{remark}

Although the grading convention of \cref{def:synthetic_bigraded_spheres} has become more standard, it is not the only one in the literature.
We refer to the indexing of \cref{def:synthetic_bigraded_spheres} as \defi{Adams grading} of synthetic spectra.
This is not the convention used in \cite{pstragowski_synthetic}, which instead follows the \emph{motivic grading}.
Unless explicitly said otherwise, we will not use motivic grading in these notes.

\begin{remark}[Motivic grading]
    \label{rmk:motivic_grading_Syn}
    The \defi{motivic grading} on synthetic spectra is to define
    \[
        \S^{t,w} := \opSigma^{t-w}\nu(\S^w).
    \]
    Conversion from Adams to motivic grading, and vice versa, is given respectively by
    \[
        (n,s) \mapsto (n,\, n+s) \qquad \text{and} \qquad (t,w) \mapsto (t,\, w-t).
    \]
    The only cases in which motivic grading agrees with Adams grading are those where the stem is $0$.
    (In particular, $\tau$ has bidegree $(0,-1)$ in both conventions.)
    In \cite{pstragowski_synthetic}, the degree $w$ is called the \defi{weight}, and the difference $t-w$ is called the \defi{Chow degree}.
    In Adams grading, the weight of $\S^{n,s}$ is given by~$n+s$, while the Chow degree is given by~$-s$.
    The motivic grading is designed to match with the indexing conventions of motivic homotopy theory; see \cref{ssec:motivic} for more information.
\end{remark}

\section{The homological t-structure}
\label{sec:homological_t_structure}

One of the most important features of $\Syn_E$ that distinguishes it from the category of $\FilSp$ is the existence of a particular t-structure.
Following Burklund--Hahn--Senger \cite{burklund_hahn_senger_manifolds}, we refer to this t-structure as the \emph{homological t-structure}; it is referred to as the \emph{natural t-structure} in \cite{pstragowski_synthetic}.
Although this is indeed the default t-structure on synthetic spectra for our purposes, we prefer this more descriptive name.

As the name suggests, the defining feature of this t-structure is that it looks at \emph{$\nu E$-homology} of $E$-synthetic spectra to measure (co)connectivity, \emph{not} at the bigraded homotopy groups.
This has both its upsides and downsides.
On the one hand, this homology tends to be a lot simpler than the homotopy, due to the special role that $E$ plays for $E$-synthetic spectra.
On the other hand, this means that taking truncations or connective covers can have very unpredictable effects on bigraded homotopy groups.

\begin{remark}
    Carrick and Davies \cite[Section~2]{christian_jack_synthetic_j} define t-structures on $\Syn_E$ based on the bigraded homotopy groups.
\end{remark}

\begin{definition}
    The \defi{homological t-structure} on $\Syn_E$ is the t-structure where a synthetic spectrum $X$ is connective if and only if
    \[
        \nu E_{n,s}(X) = 0 \qquad \text{whenever } s > 0.
    \]
    We will write $\tau_{\geq n}$ and $\tau_{\leq n}$ for the $n$-connective cover and $n$-truncation functors with respect to this t-structure, respectively.
\end{definition}

\begin{theorem}
    \label{thm:homological_t_structure}
    The homological t-structure is an accessible t-structure on $\Syn_E$ that satisfies the following.
    \begin{letterenum}
        \item \label{item:syn_connectivity_property} A synthetic spectrum $X$ is connective if and only if
        \[
            \nu E_{n,s}(X) = 0 \qquad \text{whenever } s > 0.
        \]
        \item \label{item:syn_truncatedness_property} A synthetic spectrum $X$ is $0$-truncated if and only if $X$ is $\nu E$-local and
        \[
            \nu E_{n,s}(X) = 0 \qquad \text{whenever }s<0.
        \]
        \item \label{item:syn_conn_cover_iso_homology} Let $X$ be a synthetic spectrum.
        The connective cover $\tau_{\geq0} X \to X$ induces an isomorphism
        \[
            \nu E_{n,s}(\tau_{\geq 0} X) \congto \nu E_{n,s}(X) \qquad \text{whenever }s\leq 0.
        \]
        Likewise, the $0$-truncation $X \to \tau_{\leq 0} X$ induces an isomorphism
        \[
            \nu E_{n,s}(X) \congto \nu E_{n,s}(\tau_{\leq 0} X) \qquad \text{whenever }s\geq 0.
        \]
        \item \label{item:synthetic_analogue_connective} For every spectrum $X$, the synthetic spectrum $\nu X$ is connective.
        \item \label{item:heart_SynE} There exists a monoidal equivalence of categories
        \[
            \Syn_E^\heart \simeq \grComod_{E_*E}
        \]
        under which $\opSigma^{n,0}$ becomes the functor $[n]$, and which fits into a commutative diagram of lax monoidal functors
        \[
            \begin{tikzcd}
                \Sp \ar[r,"\tau_{\leq0} \nu"] \ar[dr, "E_*(\blank)"'] & \Syn_E^\heart \ar[d,"\simeq"] \\
                & \grComod_{E_*E}.
            \end{tikzcd}
        \]
        Moreover, if $E$ is homotopy commutative, then this equivalence and the above diagram are naturally symmetric monoidal.
        \item The t-structure is right complete \textbr{but for general $E$, not even left separated}.\label{item:syn_right_complete}
        \item The t-structure is compatible with filtered colimits.\label{item:syn_is_Grothendieck}
        \item The t-structure is compatible with the monoidal structure.
    \end{letterenum}
\end{theorem}

\begin{proof}
    Property~\ref{item:syn_connectivity_property} alone determines this t-structure uniquely.
    By Theorem~4.18 of \cite{pstragowski_synthetic}, this definition of the t-structure agrees with the definition from Proposition~2.16 of op.\ cit., which has the remaining desired properties by Propositions~4.16, 4.18, and~4.21 of op.\ cit.
\end{proof}

In fact, as we will re-prove in \cref{ex:connectivity_synthetic_analogue_explained} below, the $\nu E$-homology of $\nu X$ has a very simple form: there is an isomorphism of bigraded $\Z[\tau]$-modules
\[
    \nu E_{*,*}(\nu X) \cong E_*X[\tau]
\]
where $E_nX$ is placed in bidegree $(n,0)$.
This property has useful consequences for working with $\nu$; see \cref{rmk:nu_is_right_adjoint_when_landing_in_connective} for instance.
It is also one of the first instances where we see that  $\Syn_E$ is more suited for working with Adams spectral sequences than $\FilSp$: as we will see in \cref{warn:filtered_vs_synthetic_sphere}, the filtered spectrum underlying $\nu X$ is rarely connective in the diagonal t-structure on $\FilSp$, making it harder to work with that filtered spectrum directly.

\begin{remark}
    It follows from \cref{item:syn_is_Grothendieck} that $(\Syn_E)_{\geq 0}$ is a \emph{Grothendieck prestable $\infty$-category}; see \cite[Proposition~C.1.4.1]{SAG}.
\end{remark}

\begin{remark}
    There exist examples for which the monoidal equivalence $\Syn_E^\heart \simeq \grComod_{E_*E}$ cannot be made \emph{symmetric} monoidal.
    One can think of this as saying that the `correct' braiding on $\grComod_{E_*E}$ is not the usual algebraic one, but rather a more exotic `topological' one.
    For an example of this phenomenon, see \cite[Section~6]{Hopkins_Lurie_Brauer_E} and \cite[Section~4]{barthel_pstragowski_Morava_Ktheory}, where $E = \K(n)$ and where $\Sp$ is replaced by ($\K(n)$-local) modules over Morava E-theory.
\end{remark}

Finally, let us make a few comments regarding notation and indexing.

\begin{remark}
    \label{rmk:cohom_vs_hom_indexing_t_structure_Syn}
    The description of the connective objects is somewhat confusing, in that an object is connective when certain groups in a \emph{positive} degree vanish.
    This clash is because the filtration in Adams spectral sequences is indexed cohomologically, while (at least in homotopy theory) we usually index t-structures homologically.
    Arguably, it would be less confusing to index this t-structure cohomologically instead (as is more common in algebraic geometry), writing $\tau^{\leq n}$ for what we normally denote by $\tau_{\geq -n}$, and $\tau^{\geq n}$ for $\tau_{\leq -n}$.
    However, to prevent confusion with the standard convention in homotopy theory, we will refrain from doing this.
\end{remark}

\begin{warning}
    \label{warn:syn_heart_htpy_groups}
    Often with t-structures, one writes $\pi_n^\heart$ for the functor $\opSigma^{-n}\tau_{\leq n}\tau_{\geq n}$ considered as landing in the heart of the t-structure.
    Because this t-structure is measured by homology instead of homotopy, this notation can get confusing: the functor $\pi_n^\heart$ is \emph{not} given by bigraded homotopy groups.
    Instead, for $X \in \Syn_E$, by (the comment following) \cite[Theorem~4.18]{pstragowski_synthetic}, we have an isomorphism of graded $E_*E$-comodules
    \[
        \pi_n^\heartsuit(X) \cong \nu E_{*+n,\,-n}(X) = \nu E_{*,\,-n}(X) [-n].
    \]
    Note also the minus sign in the filtration on the right-hand side; this is again due to the difference between homological and cohomological grading (cf.~\cref{rmk:cohom_vs_hom_indexing_t_structure_Syn}).
    To avoid the potential confusion with the bigraded homotopy groups, we will generally avoid the notation $\pi_n^\heartsuit$.
\end{warning}

\section{Synthetic spectra as a deformation}
\label{sec:Syn_as_deformation}

With the foundational properties and structure in hand, we can relate synthetic spectra to spectral sequences.
We begin by defining a deformation structure.

\begin{lemma}
    \label{lem:symm_mon_tau_tower_on_synthetic_sphere_first_mention}
    There is a natural symmetric monoidal structure on the functor $\Z \to \Syn_E$ given by the multiplication-by-$\tau$ tower on the unit:
    \[
        \begin{tikzcd}
            \dotsb \ar[r,"\tau"] & \S^{0,-1} \ar[r,"\tau"] & \S \ar[r,"\tau"] & \S^{0,1} \ar[r,"\tau"] & \dotsb.
        \end{tikzcd}
    \]
\end{lemma}
\begin{proof}
    We follow the argument given in the proof of \cite[Corollary~6.1]{lawson_cellular}.
    In Pstr\k{a}gowski's model \cite{pstragowski_synthetic}, the synthetic sphere $\S^{0,s}$ is defined as the sheafification of the presheaf $\tau_{\geq -s}\map(\blank,\S)$, with $\tau$ induced by the suspension-comparison map.
    The Whitehead filtration functor is lax symmetric monoidal (\cref{rmk:whitehead_lax_mon}), as is sheafification, so the $\E_\infty$-structure on the sphere spectrum induces a symmetric monoidal structure on the multiplication-by-$\tau$ tower.
\end{proof}

We will use the notation and terminology introduced in \cref{not:rho_sigma_general_deformation}; let us repeat it here for convenience.

\begin{notation}
    \label{not:rho_sigma_Syn}
    By the universal property of $\FilSp$ from \cref{thm:universal_property_FilSp}, the symmetric monoidal functor $\Z \to \Syn_E$ from \cref{lem:symm_mon_tau_tower_on_synthetic_sphere_first_mention} induces an adjunction
    \[
        \begin{tikzcd}[column sep=3.5em]
            \FilSp \ar[r, shift left, "\rho"] & \Syn_E \ar[l, shift left, "\sigma"]
        \end{tikzcd}
    \]
    where the left adjoint $\rho$ is a symmetric monoidal functor.
    As a result, the functor $\sigma$ is naturally lax symmetric monoidal.
    If $X$ is a synthetic spectrum, then we refer to the filtered spectrum $\sigma X$ as its \defi{signature}.
\end{notation}

\begin{remark}[History]
This functor has appeared before in \cite[Appendix~C]{burklund_hahn_senger_Rmotivic} under the name $i_*$.
The name \emph{signature} was introduced in \cite{christian_jack_synthetic_j}, with the letter $\sigma$ starting to be used in \cite{CDvN_part1,CDvN_part2} (and in later revisions of \cite{christian_jack_synthetic_j}).
\end{remark}

The functor $\rho$ lets us import important structure from $\FilSp$.
For instance, for every $k\geq 1$, the synthetic spectrum $C\tau^k$ inherits an $\E_\infty$-structure from the filtered spectrum $C\tau^k$.\footnote{Alternatively, one can use \cref{prop:nu_mod_tau} and the monoidality of the homological t-structure to give $C\tau$ an $\E_\infty$-structure; this is how it is done in \cite[Corollary~4.30]{pstragowski_synthetic}.}

As explained in \cref{sec:deformations}, the functor $\sigma$ can be thought of as an `underlying spectral sequence' functor.
More precisely, as a consequence of \cref{prop:deformation_sigma_tauinv}, it sends the synthetic map $\tau_X$ to the transition map of $\sigma X$, and preserves modding out by $\tau$.
We now check that it also preserves colimits, and in particular preserves $\tau$-inversion.

\begin{proposition}
    \leavevmode
    \label{cor:properties_syn_sigma}
    \begin{numberenum}
        \item The functor $\sigma$ preserves colimits, and is even $\FilSp$-linear.
        In particular, $\sigma$ preserves $\tau$-inversion, and $\rho$ is an internal left adjoint in $\FilSp$-linear $\infty$-categories.
        \item The functor $\sigma$ is conservative if and only if $\Syn_E$ is cellular.
    \end{numberenum}
\end{proposition}
\begin{proof}
    This follows from \cref{prop:deformation_sigma_tauinv}\,\ref{item:colims_sigma_condition} and~\ref{item:conservative_sigma_condition}, and \cref{rmk:sigma_internal_right_adjoint}, using that the synthetic sphere is compact, and that $\S^{0,s}$ for $s\in\Z$ form stable generators if and only if $\Syn_E$ is cellular.
\end{proof}

\begin{remark}
    \label{rmk:expression_sigma}
    As a special case of \cref{rmk:expression_for_sigma}, the functor $\sigma$ can be described as follows.
    Write $\map(\blank,\blank)$ for the mapping spectrum functor of the stable $\infty$-category $\Syn_E$.
    Then $\sigma$ is given by levelwise applying $\map(\S,\blank)$ to the multiplication-by-$\tau$ tower functor.
    In diagrams: for $X \in \Syn_E$, the filtered spectrum $\sigma X$ is given by
    \[
        \begin{tikzcd}
            \dotsb \ar[r,"\tau"] & \map(\S,\,  \Sigma^{0,-1} X) \ar[r,"\tau"] & \map(\S,\, X) \ar[r,"\tau"] & \map(\S,\, \Sigma^{0,1} X) \ar[r,"\tau"] & \dotsb.
        \end{tikzcd}
    \]
\end{remark}

\begin{remark}
    The adjunction $\rho \dashv \sigma$ is very close to a monadic adjunction.
    More precisely, it is a monadic adjunction if and only if $\Syn_E$ is cellular.
    What requires more assumptions is to then identify the monad on filtered spectra without making reference to the synthetic category.
    We discuss these things more in \cref{sec:filtered_models_syn}.
\end{remark}

\subsection{The signature spectral sequence}
\label{ssec:synthetic_htpy_signature}
The deformation picture tells us how to understand the bigraded homotopy groups of a synthetic spectrum.
Namely, the synthetic bigraded spheres are in the image of the left adjoint $\rho$; as a result, understanding synthetic homotopy groups is equivalent to understanding the filtered homotopy groups of $\sigma$ applied to the synthetic spectrum.
The latter, as explained by the Omnibus Theorem, captures a spectral sequence.

The only subtlety in this story is that there is a reindexing taking place when passing between filtered and synthetic spectra.
To avoid confusion, we will for the moment distinguish the filtered and synthetic settings by writing
\[
    \S_\fil^{n,s} \quad \text{and} \quad \S_\syn^{n,s}
\]
for the filtered and synthetic spheres, respectively, and similarly $\pi^\fil_{*,*}$ and $\pi^\syn_{*,*}$ for the homotopy groups.
\begin{proposition}
    \label{prop:synthetic_vs_filtered_htpy_groups}
    \leavevmode
    \begin{numberenum}
        \item \label{item:fil_sphere_vs_syn_sphere} For all $n$ and $s$, we have an isomorphism
        \[
            \rho(\S^{n,s}_\fil) \cong \S^{n,\, s-n}_\syn.
        \]
        \item \label{item:pi_fil_vs_pi_syn} For all $n$, we have a natural isomorphism of graded $\Z[\tau]$-modules \textbr{where $X \in \Syn_E$}
        \[
            \pi_{n,*}^\syn(X) \cong \pi_{n,\,*+n}^\fil(\sigma X).
        \]
        \item The functor $\rho$ is right t-exact \textbr{with respect to the diagonal t-structure on $\FilSp$ and the  homological t-structure on $\Syn_E$}; equivalently, the functor $\sigma$ is left t-exact.
    \end{numberenum}
\end{proposition}
\begin{proof}
    The functor $\rho$ is characterised by preserving colimits and sending $\S^{0,s}_\fil$ to $\S^{0,s}_\syn$ for all $s$.
    In particular, $\rho$ is exact, so it preserves arbitrary suspensions.
    Using the identifications
    \[
        \S^{n,s}_\fil \cong \opSigma^n \S^{0,s}_\fil \qquad \text{and} \qquad \Sigma_\syn(\blank) \cong \S^{1,\,-1}_\syn \otimes \blank,
    \]
    the first isomorphism follows.
    Using that $\rho$ is left adjoint to $\sigma$, this implies that for every $n$ and $s$, we have a natural isomorphism of abelian groups (where $X \in \Syn_E$)
    \[
        \pi_{n,s}^\syn(X) \cong \pi_{n,\,s+n}^\fil(\sigma X).
    \]
    By \cref{prop:deformation_sigma_tauinv}, the functor $\sigma$ sends $\tau_X$ to $\tau_{\sigma X}$, so this assembles to the claimed isomorphism of graded $\Z[\tau]$-modules.

    For the final claim, recall that $(\FilSp)_{\geq0}$ is the smallest subcategory generated under colimits by the objects $\S^{n,s}_\fil$ for $n-s\geq 0$.
    Since synthetic analogues are connective, it follows that $\S^{k,u}_\syn = \opSigma^{-u}\nu(\S^{k+u})$ is $(-u)$-connective in the homological t-structure.
    It follows that $\rho$ sends $\S^{n,s}_\fil$ for $n-s\geq 0$ to a connective synthetic spectrum.
    Because $\rho$ preserves colimits, it follows that it restricts to a functor $(\FilSp)_{\geq0} \to (\Syn_E)_{\geq 0}$, proving the claim.
\end{proof}

\begin{warning}
    \label{warn:filtered_vs_synthetic_sphere}
    Even though $\rho(\S^{n,s}_\fil)$ is a synthetic sphere, the filtered spectrum $\sigma(\S^{n,s}_\syn)$ is \emph{very} different from a filtered sphere.
    Indeed, the spectral sequence associated to a filtered sphere is uninteresting (see \cref{ex:sseq_filtered_sphere}), while the spectral sequence underlying $\sigma(\S^{n,s}_\syn)$ is (a shift of) the $E$-Adams spectral sequence for the sphere spectrum (see \cref{thm:signature_nu} below), which is very interesting and highly nontrivial for many $E$.
    In particular, $\rho$ is very far from preserving bigraded homotopy groups.
    This also implies $\sigma$ is not right t-exact.
\end{warning}

\Cref{prop:synthetic_vs_filtered_htpy_groups} tells us that the synthetic homotopy groups capture a spectral sequence; we give it a special name.

\begin{definition}
    Let $X$ be a synthetic spectrum.
    The \defi{signature spectral sequence} of~$X$ is the spectral sequence underlying the filtered spectrum $\sigma X$.
\end{definition}

We will use second-page indexing for this spectral sequence, so that it is of the form
\[
    \uE_2^{n,s} = \oppi_{n,s}(C\tau \otimes X) \implies \pi_n(X[\tau^{-1}]).
\]
The Omnibus Theorem makes precise the way in which the synthetic homotopy groups capture this spectral sequence.
One has to be slightly careful in that we reindexed the above to start on the second page, which substracts one from powers of $\tau$ in the Omnibus Theorem.
For example, $d_r$-differentials in the signature spectral sequence of $X$ introduce $\tau^{r-1}$-torsion in $\oppi_{*,*}X$.
We spell this out in the non-truncated case in \cref{sec:synthetic_omnibus} below.

\begin{remark}
    \label{rmk:why_second_page_indexing_Syn}
    We use second-page indexing because, for $X$ a synthetic analogue, the filtered spectrum $\sigma X$ is the \emph{d\'ecalage} of an Adams spectral sequence.
    As a result, it makes most sense to index this spectral sequence to agree with the usual indexing for Adams spectral sequences.
    The definition of synthetic spheres from \cref{def:synthetic_bigraded_spheres} was chosen exactly to fit with second-page indexing.
    Because we use first-page indexing on filtered spectra, this has the unfortunate side effect of causing the reindexing as in \cref{prop:synthetic_vs_filtered_htpy_groups}\,\ref{item:fil_sphere_vs_syn_sphere}.
    The reindexing $(n,s) \mapsto (n,\, s-n)$ of \cref{prop:synthetic_vs_filtered_htpy_groups} is precisely the reindexing of \cref{rmk:E2_indexing_FilSp}.
\end{remark}

\begin{remark}
The distinction between the similar, but different, terms \emph{signature} and \emph{signature spectral sequence} is intentional.
The former is a filtered spectrum, and as a result is able to capture more intricate structures (e.g., $\E_n$-structures), while the latter is only an algebraic object.
However, we will not need to make this distinction very often.
\end{remark}

The bare formalism only takes us so far: it does not tell us which spectral sequences arise in this way, nor does it tell us what the structure of $\tau$-inverted synthetic spectra or $C\tau$-modules in synthetic spectra are.
We will investigate the second question first, leaving the computation of signatures for \cref{sec:signature_synthetic_analogue}.
As advertised in the introduction, we will show that $C\tau$-modules are of an algebraic nature, forming a type of derived $\infty$-category of an abelian category.
This gives signature spectral sequences a structural advantage over the one coming from a bare filtered spectrum: the starting page is, in a sense, entirely algebraic.

Before we begin, let us briefly record the definition (and reindexing) of the $\tau$-BSS in the synthetic setting.

\begin{variant}
    \label{var:tau_BSS_synthetic}
    If $X$ is a synthetic spectrum, then its \defi{$\tau$-adic filtration} is the filtered synthetic spectrum $\Z^\op \to \Syn_E$ given by
    \[
        \begin{tikzcd}
            \dotsb \ar[r,"\tau"] & \opSigma^{0,-2} X \ar[r,"\tau"] & \opSigma^{0,-1}X \ar[r,"\tau"] & X \ar[r,equals] & \dotsb.
        \end{tikzcd}
    \]
    Analogously to \cref{constr:trigraded_sseq_bifiltered_spectrum}, this leads to a trigraded spectral sequence that we call the \defi{$\tau$-Bockstein spectral sequence} of $X$, which is of the form
    \[
        \uE_1^{n,w,s} \cong \begin{cases}
        \ \oppi_{n,\, w+s}(X/\tau) &\text{if }s \geq 0\\
        \ 0 &\text{else}\end{cases}
        \quad \implies \oppi_{n,w} X.
    \]
    Its differential $d_r^\tau$ is of tridegree $(-1,1,r)$ for $r\geq 1$.
    Note that when we apply $\sigma$ to the $\tau$-adic filtration on $X$, we obtain the $\tau$-adic filtration on $\sigma X$.
    Using \cref{prop:synthetic_vs_filtered_htpy_groups}, it follows that the $\tau$-BSS of $X$ is merely a reindexing of the $\tau$-BSS of $\sigma X$; the reindexing is given by
    \[
        \uE_1^{n,w,s}(X) \cong \uE_1^{n,\,w+n,\,s}(\sigma X).
    \]
    We may therefore freely use the results of \cref{sec:tau_Bockstein_sseq} for this spectral sequence; in particular, it captures the signature spectral sequence of $X$.
    Note that this is a situation where we use second-page indexing for ordinary spectral sequences (coming to us from the conventions for synthetic spectra), but nevertheless use first-page indexing for the corresponding $\tau$-BSS; see \cref{rmk:second_page_indexing_tau_BSS} for a further discussion of this.
    If we use second page indexing for both (see the previously cited remark), then we obtain the indexing of the $\tau$-BSS used in, e.g., \cite[Theorem~A.8]{burklund_hahn_senger_manifolds}.
\end{variant}

\subsection{Inverting \texorpdfstring{$\tau$}{tau}}
\label{ssec:synthetic_invert_tau}

Recall that the functor $\nu \colon \Sp \to \Syn_E$ is fully faithful.
However, since it is not an exact functor, we should not think too strongly of the image of $\nu$ as an embedding of spectra into synthetic spectra.
If we make $\nu$ exact in a universal way, then this does result in a good embedding of spectra into synthetic spectra, and these happen to be exactly the $\tau$-invertible synthetic spectra.

As in \cref{not:special_generic_fibres_deformation}, we write $\Syn_E[\tau^{-1}]$ for the generic fibre of the deformation $\Syn_E$.
By \cref{rmk:expressions_inv_tau_and_modtau}, this is the full subcategory of $\Syn_E$ on the $\tau$-invertible synthetic spectra, and is moreover a smashing localisation of $\Syn_E$.

\begin{definition}
    Write $\yo \colon \Sp \to \Syn_E$ for the functor $\nu(\blank)[\tau^{-1}]$.
\end{definition}

The functor $\yo$ is also referred to as the \emph{spectral Yoneda embedding}.
By definition, $\yo$ lands in $\tau$-invertible synthetic spectra.

\begin{theorem}[\cite{pstragowski_synthetic}, Theorem~4.37]
    \label{thm:tau_inv_Syn}
    The functor $\yo$ is fully faithful, exact, and symmetric monoidal, and restricts to a symmetric monoidal equivalence
    \[
        \yo \colon \Sp \simeqto \Syn_E[\tau^{-1}].
    \]
\end{theorem}

\begin{notation}
We write $(\blank)^{\tau=1}$ for the composite
\[
    \begin{tikzcd}
        \Syn_E \ar[r,"\tau^{-1}"] & \Syn_E[\tau^{-1}] \simeq \Sp.
    \end{tikzcd}
\]  
\end{notation}

The following example is the analogous one to \cref{ex:invert_tau_on_filtered_sphere}.

\begin{example}
    \label{ex:invert_tau_on_bigraded_syn_sphere}
    Recall the definition $\S^{n,s}=\opSigma^{-s}\nu(\S^{n+s})$ from \cref{def:synthetic_bigraded_spheres}.
    As $\tau$-inversion is an exact functor on synthetic spectra, it preserves suspensions, so we find that
    \[
        \S^{n,s}[\tau^{-1}] = \opSigma^{-s}\nu(\S^{n+s})[\tau^{-1}] = \Sigma^{-s}\yo (\S^{n+s}) \cong \yo(\S^n).
    \]
    In other words, $(\S^{n,s})^{\tau=1} \cong \S^n$.
    We can think of this as saying that inverting $\tau$ forgets the Adams filtration.
\end{example}

This identification is compatible with our earlier identification of $\FilSp[\tau^{-1}]$ from \cref{ssec:filtered_inverting_tau}.

\begin{proposition}
    \label{prop:Syn_sigma_special_fibres}
    The adjunction $\rho \dashv \sigma$ restricts to an adjoint equivalence between $\tau$-invertible objects.
    Moreover, this equivalence fits into a commutative diagram of symmetric monoidal equivalences
    \[
        \begin{tikzcd}
            & \Sp \ar[dl,"\Const"'] \ar[dr,"\yo"] & \\
            \FilSp[\tau^{-1}] \ar[rr,"\rho",shift left] & & \Syn_E[\tau^{-1}]. \ar[ll,"\sigma",shift left]
        \end{tikzcd}
    \]
    In particular, for every spectrum $X$, the colimit of the filtered spectrum $\sigma(\nu X)$ is naturally isomorphic to $X$.
\end{proposition}
\begin{proof}
    The functors $\rho \Const$ and $\yo$ are symmetric monoidal colimit-preserving functors $\Sp \to \Syn_E$.
    By the universal property of $\Sp$, it follows that they are naturally isomorphic as symmetric monoidal functors; see \cite[Corollary~4.8.2.19]{HA}.
    By two-out-of-three, the functor $\rho$ restricts to an equivalence between $\tau$-invertible objects.
    Because $\sigma$ is right adjoint to $\rho$, it follows that $\sigma$ restricts to an inverse for it.

    The final claim follows from the isomorphism $\sigma(\nu(X)[\tau^{-1}])\cong \Const X$ and the fact that $\sigma$ preserves $\tau$-inversion by \cref{cor:properties_syn_sigma}.
\end{proof}

So far, we have focussed on $\nu$ and defined $\yo$ in terms of it.
It is also possible to go in the other direction and characterise $\nu$ in terms of $\yo$, using the homological t-structure.

\begin{proposition}
    \label{prop:nu_vs_yo}
    Let $X$ be a spectrum.
    \begin{numberenum}
        \item The $\tau$-inversion map
        \[
            \nu X \to \yo(X)
        \]
        is a connective cover with respect to the homological t-structure.
        \item There is a natural isomorphism of functors $\Z^\op \to \Syn_E$ between the Whitehead filtration of $\yo(X)$,
        \[
            \begin{tikzcd}
                \dotsb \rar & \tau_{\geq 1} \yo(X) \rar & \tau_{\geq 0} \yo(X) \rar & \tau_{\geq -1} \yo(X) \rar & \dotsb
            \end{tikzcd}
        \]
        and the multiplication-by-$\tau$ tower on $\nu X$,
        \[
            \begin{tikzcd}
                \dotsb \ar[r,"\tau"] & \opSigma^{0,-1} \nu X \ar[r,"\tau"] & \nu X \ar[r,"\tau"] & \opSigma^{0,1}\nu X \ar[r,"\tau"] & \dotsb.
            \end{tikzcd}
        \]
    \end{numberenum}
\end{proposition}
\begin{proof}
    The first is \cite[Proposition~4.36]{pstragowski_synthetic}, and the second follows from the same argument as in \cref{lem:symm_mon_tau_tower_on_synthetic_sphere_first_mention}.
\end{proof}

\begin{remark}
    \label{rmk:nu_is_right_adjoint_when_landing_in_connective}
    The functor $\nu \colon \Sp \to \Syn_E$ is neither a left nor right adjoint, as it is not even an exact functor.
    When considered as landing in connective synthetic spectra however, its categorical properties improve: it is then right adjoint to inverting~$\tau$.
    This follows from the isomorphism $\nu \cong \tau_{\geq 0} \circ \yo$ and by pasting adjunctions: the horizontal composites in
    \[
        \begin{tikzcd}[column sep=4em]
            \Sp \ar[r, shift right=1.25,"\yo"'] & \Syn_E[\tau^{-1}] \ar[r, shift right=1.25, hook] \ar[l, shift right=1.25,"(\blank)^{\tau=1}"'] & \Syn_E \ar[r, shift right=1.25,"\tau_{\geq0}"'] \ar[l, shift right=1.25,"\tau^{-1}"'] & (\Syn_E)_{\geq 0} \ar[l, shift right=1.25, hook']
        \end{tikzcd}
    \]
    form the adjunction
    \[
        \begin{tikzcd}[column sep=4em]
            \Sp  \ar[r, shift right=1.25,"\nu"'] & (\Syn_E)_{\geq 0}. \ar[l, shift right=1.25, "(\blank)^{\tau=1}"']
        \end{tikzcd}
    \]
\end{remark}

\subsection{Modding out by \texorpdfstring{$\tau$}{tau}}

Although $\tau$-invertible synthetic spectra are equivalent to $\tau$-invertible filtered spectra, modules over the cofibre of $\tau$ are very different in synthetic spectra compared to filtered spectra.
This distinction is controlled by the homological t-structure.

\begin{proposition}
    \label{prop:nu_mod_tau}
    For every spectrum $X$, the natural map $\nu X \to C\tau\otimes \nu X$ exhibits the target as the $0$-truncation of the source.
    In particular, we have a natural isomorphism
    \[
        C\tau \otimes \nu X \cong E_*(X)
    \]
    where we regard the right-hand side as an element of $\grComod_{E_*E} \simeq \Syn_E^\heart$.
\end{proposition}
\begin{proof}
    The first statement is \cite[Lemma~4.29]{pstragowski_synthetic}.
    The final claim follows by combining this with \cref{thm:homological_t_structure}\,\ref{item:heart_SynE}.
\end{proof}

Morally, $C\tau$-modules in $\Syn_E$ are equivalent to the derived $\infty$-category of graded $E_*E$-comodules.
This is not entirely true, as for most $E$, the unit of $\D(\grComod_{E_*E})$ is not a compact object.\footnote{Note that this is a phenomenon that does not occur in $\D(\Mod_A)$ or $\D(\mathrm{grMod}_A)$ for a (graded) commutative ring $A$, and is specific to the comodule setting.}
Accordingly, to make such a statement true, we need a modification of this derived $\infty$-category to have a compact unit.
It turns out that this modification is actually an enlargement.

The idea behind this construction is that instead of inverting homology-isomorphisms of chain complexes of comodules, we should invert homotopy-isomorphisms.
Hovey \cite{hovey_htpythy_comodules} constructs a model category doing this; see also the introduction to op.\ cit.\ for a further motivation for this construction.
Barthel--Heard--Valenzuela \cite{barthel_heard_valenzuela_comodules} give the following description of the underlying $\infty$-category of Hovey's model category; see also \cite[Section~3.2]{pstragowski_synthetic} for a summary.

\begin{recall}
    Let $(A,\Gamma)$ be a graded Hopf algebroid.
    Write $\Perf_\Gamma$ for the thick subcategory of $\D(\grComod_{(A,\Gamma)})$ generated by the dualisable (ordinary, non-derived) comodules over $(A,\Gamma)$ (considered as objects in the heart the derived $\infty$-category).
    Define the \defi{stable comodule $\infty$-category} as
    \[
        \Stable_{(A,\Gamma)} = \Ind(\Perf_\Gamma).
    \]
    The inclusion functor $\Perf_\Gamma \to \D(\grComod_{(A,\Gamma)})$ induces a functor $\Stable_{(A,\Gamma)} \to \D(\grComod_{(A,\Gamma)})$, and this turns out to have a fully faithful right adjoint:
    \begin{equation}
        \label{eq:Stable_derived_adjunction}
        \begin{tikzcd}
            \Stable_{(A,\Gamma)} \rar[shift left] & \D(\grComod_{(A,\Gamma)}). \lar[shift left, hook']
        \end{tikzcd}
    \end{equation}
    As $\Perf_\Gamma$ is closed under tensor products, the $\infty$-category $\Stable_{(A,\Gamma)}$ is naturally a symmetric monoidal functor, and the localisation functor to the derived is symmetric monoidal.
    By definition, the unit of $\Stable_{(A,\Gamma)}$ is compact.
    Moreover, the localisation \eqref{eq:Stable_derived_adjunction} is precisely given by $\Gamma$-localisation, i.e., inverting those maps that become isomorphisms after tensoring with $\Gamma$.
\end{recall}

If $E$ is a homotopy-associative ring spectrum, then we also write $\Stable_{E_*E}$ for the stable comodule $\infty$-category of the Hopf algebroid $(E_*,E_*E)$.

\begin{theorem}
    \label{thm:Syn_mod_tau_is_Stable}
    There is a right t-exact fully faithful left adjoint of monoidal $\infty$-categories
    \[
        \Mod_{C\tau}(\Syn_E) \hookto \Stable_{E_*E}
    \]
    with the following properties.
    \begin{numberenum}
        \item \label{item:nu_mod_tau_into_Stable} This functor sits in a commutative diagram of lax monoidal functors
        \[
            \begin{tikzcd}
                \Sp \rar["\nu(\blank)/\tau"] \dar["E_*(\blank)"'] & \Mod_{C\tau}(\Syn_E) \dar[hook] \\
                \grComod_{E_*E} \rar[hook] & \Stable_{E_*E}.
            \end{tikzcd}
        \]
        \item \label{item:Stable_is_Syn_Landweber_exact} If $E$ is Landweber exact or is the sphere spectrum, then this functor is an equivalence.
        \item \label{item:Stable_symm_mon} If $E$ is homotopy commutative, then this functor is naturally symmetric monoidal, and the diagram of \ref{item:nu_mod_tau_into_Stable} is naturally one of lax symmetric monoidal functors.
    \end{numberenum}
\end{theorem}

\begin{proof}
    The main result and \cref{item:Stable_symm_mon} are \cite[Theorem~4.46]{pstragowski_synthetic}.
    \Cref{item:nu_mod_tau_into_Stable} follows by combining this with \cref{thm:homological_t_structure}\,\ref{item:heart_SynE}, and \cref{item:Stable_is_Syn_Landweber_exact} is \cite[Proposition~4.53]{pstragowski_synthetic}.
\end{proof}

In particular, the special fibre of synthetic spectra is entirely algebraic (except for possibly the braiding if $E$ is not homotopy commutative).
The failure of the above functor to be essentially surjective is the problem of $\Stable_{E_*E}$ not being generated by the objects $E_*P$ where $P$ ranges over the finite $E$-projective spectra.
We regard this as a minor technical issue.

\begin{remark}
    We now have two different generalisations of $E$-homology for synthetic spectra, namely $\nu E$-homology and $C\tau$-homology.
    They both differ from $E$-homology for ordinary spectra, but in different ways.
    \begin{itemize}
        \item For $\nu E$-homology, we obtain a second grading, and even a $\Z[\tau]$-module structure.
        Although the $\nu E$-homology of a synthetic analogue is very simple (having no $\tau$-torsion for instance, see \cref{ex:connectivity_synthetic_analogue_explained} below), the $\nu E$-homology of a general synthetic spectrum can be a highly nontrivial $\Z[\tau]$-module.
        \item For $C\tau$-homology, this takes values in (a modification of) the derived $\infty$-category of $E_*E$-comodules.
        For a synthetic analogue, this lands in the heart, and by \cref{prop:nu_mod_tau} is identified with $E$-homology in the ordinary sense.
        For a general synthetic spectrum however, the resulting object will rarely be an honest $E_*E$-comodule, but will generally be a derived or stable one.
    \end{itemize}
\end{remark}

\begin{remark}
    Strengthening \cref{prop:nu_mod_tau}, the essential image of $\nu$ in fact consists precisely of those synthetic spectra $X$ that are connective and for which $C\tau \otimes X$ is discrete in the homological t-structure; see \cite[Proposition~2.16]{pstragowski_vankoughnett_obstruction_theory}.
    Using \cref{thm:Syn_mod_tau_is_Stable}, we can equivalently state this condition as asking $C\tau\otimes \nu X$ to be an honest comodule, rather than a derived one.
\end{remark}

Recall that the derived $\infty$-category is obtained from the stable comodule $\infty$-category by localising at $E_*E$-equivalences.
Translated into synthetic terms, $E_*E$ corresponds to $\nu E/\tau$, leading to the following.
Moreover, note that the generation issues go away after $\nu E$-localisation, and we obtain an actual equivalence.

\begin{theorem}[\cite{pstragowski_synthetic}, Theorem~4.54]
    \label{thm:hypercomplete_Syn_mod_tau}
    The functor from \cref{thm:Syn_mod_tau_is_Stable} restricts to a monoidal equivalence
    \[
        L_{\nu E} \, \Mod_{C\tau}(\Syn_E) \simeq \D(\grComod_{E_*E})
    \]
    which is naturally symmetric monoidal if $E$ is homotopy commutative.
\end{theorem}

\begin{notation}
    \label{not:Synhat}
    Following \cite{pstragowski_synthetic}, we will also write $\Synhat_E$ for $L_{\nu E}\, \Syn_E$.
    In op.\ cit., objects of $\Synhat_E$ are called \emph{hypercomplete}, stemming from their definition as sheaves; we will not use this name, and instead refer to them simply as $\nu E$-local objects.
    We warn the reader that, while this notation is convenient when working with a fixed $E$, it could lead to confusion when working with various $E$ at once.
    In these notes, we will always work with a fixed $E$, so this confusion should not arise.
\end{notation}

\begin{warning}
    For general $E$, the unit in $\Synhat_E$ is not compact for general $E$, because the unit of $\D(\grComod_{E_*E})$ is usually not compact.
    As a result, although $\Synhat_E$ is a symmetric monoidal deformation in its own right (obtained by $\nu E$-localising the functor $\rho$), the resulting right adjoint $\Synhat_E \to \FilSp$ does \emph{not} preserve colimits; see \cref{prop:deformation_sigma_tauinv}\,\ref{item:colims_sigma_condition}.
\end{warning}

At this point, we can begin to see the Adams spectral sequence appearing, at least its second page.

\begin{example}
    \label{ex:mod_tau_maps_between_nu}
    Let $X$ and $Y$ be spectra.
    Combining \cref{prop:nu_mod_tau} and \cref{thm:hypercomplete_Syn_mod_tau}, we learn that $C\tau$-linear maps between synthetic analogues are computed by maps of comodules:
    \[
        [\nu Y/\tau,\ \nu X/\tau]_{C\tau} \cong \Hom_{E_*E}(E_*Y,\, E_*X).
    \]
    For an integer $k$, let us denote the $k$-fold grading-shift functor on $\D(\grComod_{E_*E})$ by $[k]$, and let us write $\Sigma^k$ for the $k$-fold $\infty$-categorical suspension as usual.
    By \cref{prop:nu_mod_tau}, the synthetic spectra $\nu Y /\tau$ and $\nu X/\tau$ are $0$-truncated, which by \cref{thm:homological_t_structure}\,\ref{item:syn_truncatedness_property} in particular means they are $\nu E$-local.
    Using \cref{thm:hypercomplete_Syn_mod_tau} and the definition $\S^{n,s} = \opSigma^{-s}\nu(\S^{n+s})$, it therefore follows that
    \[
        [\opSigma^{n,s}\nu Y/\tau,\ \nu X/\tau]_{C\tau} \cong [(E_*Y)[n+s],\, \Sigma^s \, E_*X]_{\D(\grComod_{E_*E})} = \Ext^{s,\,n+s}_{E_*E}(E_*Y,\, E_*X).
    \]
    In particular, we have
    \[
        \oppi_{n,s}(\nu X/\tau) = [\S^{n,s},\, \nu X/\tau] \cong [C\tau\otimes \S^{n,s},\, \nu X/\tau]_{C\tau} \cong \Ext^{s,\,n+s}_{E_*E}(E_*,\,E_*X).\qedhere
    \]
\end{example}

\subsection{Synthetic lifts}
\label{ssec:synthetic_lifts}

Previously, we argued that the $\tau$-inversion of a synthetic spectrum can be thought of as an `underlying spectrum'.
We can also turn this question around, fixing a spectrum and asking how many synthetic spectra have this as their underlying spectrum.
It is useful to introduce some terminology for this.

\begin{definition}
    Let $X$ be a spectrum.
    A \defi{synthetic lift} of $X$ is a synthetic spectrum $S$ such that $S^{\tau=1} \cong X$.
\end{definition}

\begin{example}
\Cref{thm:tau_inv_Syn} says that $\nu$ provides a functorial synthetic lift.
\end{example}

We think of a synthetic lift of a spectrum $X$ as encoding a \emph{modified Adams spectral sequence} for $X$.
The synthetic analogue of $X$, from this perspective, is the standard synthetic lift; as we will see in \cref{thm:signature_nu}, it encodes the ordinary Adams spectral sequence for $X$.
For a further discussion of these ideas, see \cite{christian_jack_synthetic_j}.

We end this section by discussing how to construct synthetic lifts out of old ones.
As $\tau$-inversion preserves colimits, taking colimits results in a synthetic lift of the colimit of the underlying spectra.
More subtle is the use of limits, since $\tau$-inversion does not preserve all limits.
For instance, if $S$ is a synthetic spectrum, then every term in its $\tau$-adic tower
\[
    \dotsb \to S/\tau^3 \to S/\tau^2 \to S/\tau
\]
vanishes upon $\tau$-inversion; meanwhile, the $\tau$-inversion of the limit is $(S_\tau^\wedge)[\tau^{-1}]$, which is nontrivial for many $S$ (e.g., if $S = \nu X$ for $X$ an $E$-nilpotent complete spectrum, by \cref{thm:signature_nu}).
Nevertheless, if the diagram is of a special form, then $\tau$-inversion does preserve the limit.

\begin{proposition}
    \label{prop:limits_synthetic_analogues}
    Let $X\colon I \to \Sp$ be a diagram of spectra.
    Then we have an isomorphism
    \[
        \nu\br*{\lim X} \cong \tau_{\geq 0} \br*{\lim \nu (X)},
    \]
    and the limit-comparison map
    \[
        \br*{\lim \nu(X)}^{\tau=1} \to \lim X
    \]
    is an isomorphism of spectra.
    In particular, $\lim  \nu(X)$ is a synthetic lift of $\lim X$.
\end{proposition}

Because $\nu$ is fully faithful, this in fact says that $\tau$-inversion preserves the limit of any diagram that takes values in synthetic analogues.

The key input for the proof is the homological t-structure, particularly \cref{rmk:nu_is_right_adjoint_when_landing_in_connective} and the following lemma.

\begin{lemma}[\cite{pstragowski_synthetic}, Lemma~4.35]
    \label{lem:bounded_above_zero_tauinv}
    If $S$ is a bounded above synthetic spectrum, then $S[\tau^{-1}]$ is zero.
    In particular, if $S$ is any synthetic spectrum, then for every $n$, the map $\tau_{\geq n} S \to S$ becomes an isomorphism upon $\tau$-inversion.
\end{lemma}
\begin{proof}
    For every integer $s$, the suspension $\opSigma^{0,s}$ decreases coconnectivity by $s$; this follows directly from the coconnectivity criterion of \cref{thm:homological_t_structure}\,\ref{item:syn_truncatedness_property}.
    Since $\nu E$-homology preserves filtered colimits, it follows that if $S$ is bounded above, then the colimit
    \[
        S[\tau^{-1}] = \colim(\begin{tikzcd}
            S \rar["\tau"] & \opSigma^{0,1} S \rar["\tau"]& \opSigma^{0,2} S \rar["\tau"] & \dotsb 
        \end{tikzcd}),
    \]
    is $(-\infty)$-coconnective.
    Since the t-structure on $\Syn_E$ is right complete by \cref{thm:homological_t_structure}\,\ref{item:syn_right_complete}, it follows that $S[\tau^{-1}]$ vanishes.

    The final claim follows from the fact that the cofibre of $\tau_{\geq n}S \to S$ is $\tau_{\leq n-1}S$, which in particular is bounded above, and that $\tau$-inversion is an exact functor.
\end{proof}

\begin{proof}[Proof of \cref{prop:limits_synthetic_analogues}]
    Since $\nu$ takes values in connective synthetic spectra, we may consider $\nu\circ X$ as landing in $(\Syn_E)_{\geq 0}$.
    By \cref{rmk:nu_is_right_adjoint_when_landing_in_connective}, the functor $\nu \colon \Sp \to (\Syn_E)_{\geq 0}$ is right adjoint to $\tau$-inversion, implying that $\nu$ sends limits of spectra to limits in $(\Syn_E)_{\geq 0}$.
    A limit in the connective subcategory is computed as the connective cover of the limit in $\Syn_E$, so we find that
    \[
        \nu\br*{\lim X} = \tau_{\geq 0} \br*{\lim\nu(X)}.
    \]
    As $\tau$-inversion is left inverse to $\nu$, it follows that the right-hand side $\tau$-inverts to $\lim X$.
    The claim now follows from \cref{lem:bounded_above_zero_tauinv}.
\end{proof}

\section{The signature of a synthetic analogue}
\label{sec:signature_synthetic_analogue}

Previously in \cref{ssec:cosimplicial_Adams_sseq}, we defined a (cosimplicial) model for the Adams spectral sequence, resulting in a functor $\Sp \to \FilSp$.
The better and more modern definition is the following.

\begin{definition}
    \label{def:synthetic_version_of_Adams_filtration}
    The \defi{$E$-based Adams filtration} is the functor $\sigma \circ \nu_E \colon \Sp \to \FilSp$.
\end{definition}

The goal of this section is to give a justification for this name: in \cref{thm:signature_nu}, we show that this spectral sequence agrees with the d\'ecalage of the $E$-based Adams spectral sequence as defined in \cref{def:classical_E_ASS}.
(For an explanation why the d\'ecalage appears, see \cref{var:ASS_using_cosimplicial_decalage} and the discussion preceding it.)
However, as we pointed out before, the point of this is not to let go of the synthetic origins of this functor, but rather to demonstrate that this recovers the correct notion.

One concrete reason for preferring this definition over the old one is the following.

\begin{remark}[Lax monoidality; \cite{patchkoria_pstragowski_derived_inftycats}, Section~5.5]
    \label{rmk:sigmanu_lax_monoidal}
    The functor $\sigma \circ \nu_E$ is a composite of two lax symmetric monoidal functors, making it a lax symmetric monoidal functor.
    This only requires $E$ to be homotopy-associative.
    By contrast, to turn the classical definition of the $E$-Adams filtration into a lax symmetric monoidal functor, one would need an $\E_\infty$-structure on $E$.
    Such a structure does not always exist in cases of interest (e.g., $\BP$ or Morava K-theories), and the Adams spectral sequence does not depend on it, so it is not desirable to require these structures.
    For a further discussion, and a way to construct this for $E$ having only a left-unital multiplication, see \cite[Section~5.5]{patchkoria_pstragowski_derived_inftycats}.
\end{remark}

Our proof strategy is to first show this comparison on resolution objects, which for the $E$-Adams spectral sequence are the homotopy $E$-modules.
This relies on a computation of the synthetic homotopy groups of $E$-modules.
The general case follows from this by descending from a resolution by such objects.

\begin{remark}[History]
    \Cref{thm:signature_nu} is not new and is well-known to experts, but has not been written down in this specific form.
    In \cite[Proposition~5.56 and Theorem~5.60]{patchkoria_pstragowski_derived_inftycats}, \citeauthor{patchkoria_pstragowski_derived_inftycats} prove this result for the synthetic-like categories they construct therein; while formally these categories are different, the proofs follow the same ideas.
    Another closely-related result is \cite[Theorem~6.26]{pstragowski_perfect_even_filtration}, whose proof we follow closely in this section.
    Similar results in the nilpotent-complete case can be found in \cite[Remark~4.64]{pstragowski_synthetic} and \cite[Appendix~A.1]{burklund_hahn_senger_manifolds}.
    A proof of the nilpotent-complete case also appeared in \cite[Section~1.4]{CDvN_part1}, which this section is an adaptation of.
\end{remark}

The following holds for homotopy classes of maps $Y \to X$ between two spectra, but for simplicity we record it only for homotopy groups.
In words, it says that the spectral sequence underlying $\nu X$ is concentrated in nonnegative filtrations, and that the underlying spectrum of the filtration $\sigma (\nu X)$ is given by $X$.

\begin{proposition}[\cite{pstragowski_synthetic}, Theorem~4.58]
    \label{lem:htpy_synthetic_analogue_in_negative_filtration}
    Let $X$ be a spectrum.
    Then for all $s \leq 0$ and all $n$, inverting $\tau$ induces a natural isomorphism
    \[
        \pi_{n,s}(\nu X) \congto \pi_{n}X.
    \]
    Phrased differently: inverting $\tau$ induces a natural isomorphism of bigraded $\Z[\tau]$-modules
    \[
        \pi_{*,\, \leq 0}(\nu X) \cong \pi_*(X)[\tau],
    \]
    where $\pi_nX$ is placed in bidegree $(n,0)$.
\end{proposition}

For the didactic value, we include Pstr\k{a}gowski's proof.

\begin{proof}
    The cofibre sequence
    \[
        \begin{tikzcd}
            \nu X \rar["\tau"] & \opSigma^{0,1} \nu X \rar & \opSigma^{0,1} \nu X/\tau
        \end{tikzcd}
    \]
    gives rise to a long exact sequence on bigraded homotopy groups; by \cref{ex:mod_tau_maps_between_nu}, part of this reads
    \[
        \begin{tikzcd}[column sep=1.5em]
            \Ext^{s-2,\, n+s-1}_{E_*E}(E_*,\, E_*X) \ar[r] & \oppi_{n, s} (\nu X) \ar[r,"\tau"] & \oppi_{n,\,s-1}(\nu X) \ar[r] & \Ext_{E_*E}^{s-1,\, n+s-1}(E_*, E_*X).
        \end{tikzcd}
    \]
    Note that $\Ext^{s,t}_{E_*E} (E_*,\, E_*X) = 0$ whenever $s < 0$.
    Therefore if $s\leq 0$, we see that the two outer terms vanish, so that the map in the middle is an isomorphism.
    As a result, to prove the claim, we only have to show that $\tau$-inversion induces an isomorphism
    \[
        \oppi_{n,0}(\nu X) = [\nu \S^n,\, \nu X ] \to [\S^n,\, X] = \pi_n X.
    \]
    This follows from the fact that $\tau$-inversion is a left inverse to $\nu$ (see \cref{thm:tau_inv_Syn}) and that $\nu$ is fully faithful.
\end{proof}

For a particularly nice class of spectra, this computes the entirety of their synthetic homotopy groups.

\begin{proposition}[\cite{pstragowski_synthetic}, Proposition~4.60]
    \label{lem:synthetic_homotopy_of_module}
    Let $M$ be a spectrum admitting a homotopy $E$-module structure.
    Then inverting $\tau$ induces a natural isomorphism of bigraded $\Z[\tau]$-modules
    \[
        \oppi_{*,*}(\nu M) \congto \pi_*(M) [\tau],
    \]
    where $\pi_n M$ is placed in bidegree $(n,0)$.
\end{proposition}

In words, this says that the signature spectral sequence for $\nu M$ is concentrated in filtration zero, and as a result collapses without any differentials.
Again we include Pstr\k{a}gowski's proof.

\begin{proof}
    Using the previous result, we only have to show that $\pi_{n,s}(\nu M)$ vanishes when $s\geq 1$.
    We first show this for $s=1$.
    Since $M$ is a homotopy $E$-module, the Hurewicz homomorphism
    \[
        E_*(\blank) \colon \pi_n M \to \Hom_{E_*E}(E_*[n], \, E_*M)
    \]
    is an isomorphism; see \cite[Remark~3.18]{pstragowski_synthetic}.
    Under the isomorphism $\pi_n M \cong \pi_{n,0}(\nu M)$, the Hurewicz homomorphism is the right-most map in the exact sequence
    \[
        \begin{tikzcd}
            \Ext^{-1,\, n}_{E_*E}(E_*,\, E_*M) \ar[r] & \oppi_{n, 1} (\nu M) \ar[r,"\tau"] & \oppi_{n,0}(\nu M) \ar[r] & \Ext_{E_*E}^{0,\, n}(E_*, E_*M).
        \end{tikzcd}
    \]
    As the Ext group on the left vanishes, we learn that $\pi_{n,1}(\nu M)=0$ for all $n$.

    Next, we consider the case $s>1$.
    Since $M$ is a homotopy $E$-module, it follows from \cite[Remark~3.18]{pstragowski_synthetic} that we have an isomorphism of graded comodules
    \[
        E_*M \cong E_*E \otimes_{E_*} M_*,
    \]
    implying that
    \[
        \Ext^{s,t}_{E_*E}(E_*,\, E_*M) \cong \Ext_{E_*}^{s,t}(E_*,\, M_*).
    \]
    In particular, we see that these Ext groups vanish whenever $s\geq 1$.
    By the long exact sequence, this means that multiplication by $\tau$ induces an isomorphism
    \[
        \tau \colon \oppi_{n,\, s+1}(\nu X) \congto \oppi_{n,s}(\nu X)
    \]
    for all $s\geq 1$.
    We previously showed that $\pi_{*,1}(\nu M)=0$, so we are done.
\end{proof}

\begin{example}
    \label{ex:connectivity_synthetic_analogue_explained}
    Recall from \cref{ex:symm_mon_if_one_is_E} that for all spectra $X$, we have an isomorphism $\nu E \otimes \nu X \cong \nu(E\otimes X)$.
    Because $E \otimes X$ is a homotopy $E$-module, we learn from \cref{lem:synthetic_homotopy_of_module} that
    \[
        \nu E_{*,*}(\nu X) = \pi_{*,*}(\nu E\otimes \nu X) \cong \pi_{*,*}(\nu(E\otimes X)) \cong E_*(X)[\tau].
    \]
    In particular, by \cref{thm:homological_t_structure}\,\ref{item:syn_connectivity_property}, this shows that $\nu X$ is connective in the homological t-structure.
    We now see why this is independent of the connectivity of the spectrum~$X$: the connectivity of $\nu X$ is about the collapse of the $E$-Adams spectral sequence for $\pi_*(E\otimes X)$.
    
    We learn a number of things from this computation.
    \begin{numberenum}
        \item The shift $\opSigma^{0,s}\nu X$ for $s \neq 0$ is not in the essential image of $\nu$ (unless $E_*(X)$ vanishes).
    
        \item The $\Z[\tau]$-module $\nu E_{*,*}(\nu X)$ is $\tau$-torsion free.
        As a result, we learn that
        \[
            \nu E_{*,*}(C\tau \otimes \nu X) \cong (\nu E_{*,*}(\nu X))/\tau \cong E_*(X),
        \]
        where we mean the quotient by $\tau$ in the (non-derived) algebraic sense.
        This explains (apart from the $\nu E$-locality) why $C\tau \otimes \nu X$ is $0$-truncated in the homological t-structure; cf.\ \cref{thm:homological_t_structure}\,\ref{item:syn_truncatedness_property}.
    
        \item Inverting $\tau$ on $\nu E_{*,*}(\nu X)$ yields
        \[
            \nu E_{*,*}(\yo X) = \nu E_{*,*}(\nu X[\tau^{-1}]) \cong E_*(X)[\tau^\pm].
        \]
        This gives an indication of why $\nu X \to \yo X$ is a connective cover, and more generally, why the Whitehead tower of $\yo X$
        \[
            \dotsb \to \tau_{\geq 1} (\yo X) \to \tau_{\geq 0} (\yo X) \to \tau_{\geq -1} (\yo X) \to \dotsb
        \]
        can be identified with the multiplication-by-$\tau$ tower on $\nu X$
        \[
            \begin{tikzcd}
                \dotsb \rar["\tau"] & \opSigma^{0,-1}\nu X \rar["\tau"] & \nu X \rar["\tau"] & \opSigma^{0,1}\nu X \rar["\tau"] & \dotsb.
            \end{tikzcd}\qedhere
        \]
    \end{numberenum}
\end{example}

We can restate \cref{lem:synthetic_homotopy_of_module} in terms of the signature of $\nu M$.

\begin{corollary}
    \label{cor:signature_htpy_module_is_whitehead}
    Let $M$ be a spectrum admitting a homotopy $E$-module structure.
    Then there is a natural isomorphism of filtered spectra
    \[
        \sigma(\nu M)\cong \Wh M
    \]
    which is naturally a symmetric monoidal natural transformation in $M$.
\end{corollary}

\begin{proof}
    Because $\sigma$ preserves $\tau$-inversion by \cref{cor:properties_syn_sigma}, we find that for any spectrum $X$, applying $\sigma$ to the $\tau$-inversion map $\nu X \to \nu X [\tau^{-1}]$ results in a natural (symmetric monoidal) transformation
    \begin{equation}
        \label{eq:sigma_on_tauinv}
        \sigma(\nu X) \to \sigma (\nu X)[\tau^{-1}] \cong \Const X,
    \end{equation}
    where we use the identification from \cref{prop:Syn_sigma_special_fibres}.
    If now $M$ is a homotopy $E$-module spectrum, then \cref{lem:synthetic_homotopy_of_module} implies that $\sigma(\nu M)$ is connective in the diagonal t-structure on filtered spectra.
    Indeed, combining \cref{lem:synthetic_homotopy_of_module} with \cref{prop:synthetic_vs_filtered_htpy_groups}\,\ref{item:pi_fil_vs_pi_syn}, we see that the group
    \[
        \pi_{n,s} (\opsigma(\nu M)) \cong \oppi_{n,\, s-n}(\nu M)
    \]
    vanishes whenever $s-n > 0$, that is, whenever $n < s$.
    As a result, the natural map \eqref{eq:sigma_on_tauinv} in the case $X = M$ factors through a natural map
    \[
        \sigma(\nu M) \to \tau^\diag_{\geq 0} (\Const M) = \Wh M.
    \]
    Moreover, this factorisation is through a symmetric monoidal transformation, because the diagonal t-structure on filtered spectra is monoidal (\cref{prop:diagonal_t_structure_filsp}\,\ref{item:diag_t_str_is_monoidal}).
    To establish that it is an isomorphism, it suffices to show that each component $\sigma(\nu M)^s \to \optau_{\geq s}M$ is an isomorphism for all $s$.
    As this map is induced by $\tau$-inversion, this is the other part of \cref{lem:synthetic_homotopy_of_module}.
\end{proof}

Although it follows directly from \cref{lem:synthetic_homotopy_of_module} that the signature spectral sequence of $\nu M$ converges to $M$, a more refined argument even shows that $\nu M$ is $\tau$-complete in $\Syn_E$ itself.
(If $\Syn_E$ is not cellular, then this is not automatic from completeness of $\sigma(\nu M)$.)

\begin{corollary}[\cite{burklund_hahn_senger_manifolds}, Lemma~A.15]
    \label{prop:htpy_module_tau_complete}
    Let $M$ be a spectrum admitting a homotopy $E$-module structure.
    Then $\nu M$ is $\tau$-complete.
\end{corollary}

Again for didactic value, we include the proof given by Burklund--Hahn--Senger.

\begin{proof}
    We have to show that the limit (as $s \to \infty$) of
    \[
        \begin{tikzcd}
            \dotsb \rar["\tau"] & \opSigma^{0,-s}\nu M \rar["\tau"] & \dotsb \rar["\tau"] & \opSigma^{0,-1} \nu M \rar["\tau"] & \nu M
        \end{tikzcd}
    \]
    vanishes.
    It is enough to check this on mapping spectra out of $\nu P$ for all finite $E$-projective spectra $P$.
    As $\nu P$ is dualisable with dual $\nu (P^\vee)$ by \cref{prop:properties_of_Syn}\,\ref{item:fp_spectra_are_compact_dualisable_generators}, we find that
    \[
        \map(\nu P, \lim_s \opSigma^{0,-s} \nu M) \cong \lim_s \map(\S^{0,s},\, \nu (P^\vee) \otimes \nu M) \cong \lim_s \map(\S^{0,s},\, \nu(P^\vee \otimes M)),
    \]
    where for the latter isomorphism we use \cref{prop:properties_of_nu}\,\ref{item:strong_monoidality_nu}.
    The $n$-th homotopy group of the mapping spectrum at stage $s$ on the right-hand side is given by
    \[
        \pi_{n,\,s-n}(\nu (P^\vee \otimes M)).
    \]
    Because $P^\vee \otimes M$ admits a homotopy $E$-module structure (since $M$ does), these groups vanish when $n<s$ by \cref{lem:synthetic_homotopy_of_module}.
    It follows that this mapping spectrum is $s$-connective, so that the limit over $s$ is $\infty$-connective, and therefore vanishes.
\end{proof}

We are now ready to show that the spectral sequence underlying $\sigma (\nu X)$ is the $E$-based Adams spectral sequence for general $X$, or more precisely, that it captures its d\'ecalage.
To show this, we may choose any preferred $E$-resolution of $X$ to compute the $E$-based Adams spectral sequence; we use the one from \cref{def:classical_E_ASS}.

\begin{construction}
    \label{constr:comparison_map_sigmanu_EASS}
    Since $E$ is a homotopy ring spectrum, its unit map $\S \to E$ gives rise to a semicosimplicial spectrum $E^{[\bullet]}\colon \Delta_\inj \to \Sp$ which receives a map from $\S$.
    Tensoring this resulting diagram with a spectrum $X$, we obtain a diagram of spectra
    \[
        X \to E^{[\bullet]} \otimes X.
    \]
    Applying $\sigma \circ \nu$ to this diagram, we obtain a diagram of filtered spectra
    \[
        \sigma (\nu X) \to \sigma (\nu(E^{[\bullet]}\otimes X)).
    \]
    Note that for every $n\geq 1$, the spectrum $E^{\otimes n}\otimes X$ has the structure of a homotopy $E$-module.
    By \cref{cor:signature_htpy_module_is_whitehead}, it follows that $\sigma(\nu(E^{\otimes n}\otimes X)) \cong \Wh(E^{\otimes n}\otimes X)$.
    As such, the above diagram induces a map
    \[
        \sigma (\nu X) \to \Tot (\Wh(E^{[\bullet]}\otimes X)) = \Dec^\Delta(E^{[\bullet]}\otimes X).
    \]
    The target of this map is $\Dec(\ASS_E(X))$, as follows by combining \cref{prop:cosimplical_decalage_defs_agree,def:classical_E_ASS}.
\end{construction}

\begin{theorem}
    \label{thm:signature_nu}
    Let $X$ be any spectrum.
    \begin{numberenum}
        \item \label{item:nuX_complete_iff_nilp_complete} The spectrum $X$ is $E$-nilpotent complete if and only if $\nu X$ is $\tau$-complete.
        \item \label{item:comp_map_is_tau_completion} The natural comparison map from \cref{constr:comparison_map_sigmanu_EASS}
        \[
            \sigma (\nu X) \to \Dec(\ASS_E(X))
        \]
        is completion \textbr{a.k.a.\ $\tau$-completion} of filtered spectra, i.e., it is an isomorphism on associated graded and the target is complete \textbr{a.k.a.\ $\tau$-complete}.
        \item \label{item:converse_comp_map} The above comparison map is an isomorphism of filtered spectra if $X$ is $E$-nilpotent complete.
        If $\Syn_E$ is cellular, then the converse is true, i.e., it is an isomorphism if and only if $X$ is $E$-nilpotent complete.
    \end{numberenum}
\end{theorem}
\begin{proof}
    Because $\sigma$ preserves limits, the map from \cref{constr:comparison_map_sigmanu_EASS} is obtained by applying $\sigma$ to the natural map
    \begin{equation}
        \label{eq:comp_map_before_sigma}
        \nu X \to \Tot(\nu (E^{[\bullet]}\otimes X)).
    \end{equation}
    We begin by showing that this map is $\tau$-completion, i.e., that it is an isomorphism after tensoring with $C\tau$ and that the target is $\tau$-complete.
    By \cref{prop:htpy_module_tau_complete}, the target is a limit of $\tau$-complete objects, and is therefore $\tau$-complete.
    Because $C\tau$ is finite, tensoring with it preserves limits, so we find that tensoring \eqref{eq:comp_map_before_sigma} with $C\tau$ results in a map in $\calD(\grComod_{E_*E}) \subseteq \Mod_{C\tau}(\Syn_E)$ of the form
    \[
        E_*X \to \Tot((E_*E)^{\otimes [\bullet]}\otimes E_*X),
    \]
    where the tensor products are over $E_*$.
    As $E_*E$ is flat, we may take these tensor products to be underived.
    This map is an isomorphism in $\calD(\grComod_{E_*E})$, because for any $M \in \grComod_{E_*E}$, the cobar complex
    \[
        (E_*E)^{\otimes [\bullet]} \otimes M
    \]
    constitutes a cosimplicial resolution of $M$ by relative injectives, so that the map from $M$ into its totalisation is a quasi-isomorphism.
    We conclude that \eqref{eq:comp_map_before_sigma} is indeed $\tau$-completion in $\Syn_E$.

    Next, if we invert $\tau$ on \eqref{eq:comp_map_before_sigma}, then by \cref{prop:limits_synthetic_analogues} we obtain the map of spectra
    \[
        X \to \Tot(E^{[\bullet]} \otimes X).
    \]
    This map is, by definition, an isomorphism if and only if $X$ is $E$-nilpotent complete.
    Because any map of synthetic spectra is an isomorphism if and only if it so after inverting $\tau$ and after quotienting by $\tau$ (see \cref{prop:deformation_pullback_tau_complete_and_invert_tau}), we find that \eqref{eq:comp_map_before_sigma} is an isomorphism if and only if $X$ is $E$-nilpotent complete.
    As \eqref{eq:comp_map_before_sigma} is $\tau$-completion, this proves \cref{item:nuX_complete_iff_nilp_complete}.

    \Cref{item:comp_map_is_tau_completion} follows immediately from the fact that $\sigma$ preserves $\tau$-completion; see \cref{prop:deformation_sigma_tauinv}\,\ref{item:sigma_pres_modtau_taucompl}.
    If $X$ is $E$-nilpotent complete, then \eqref{eq:comp_map_before_sigma} is an isomorphism, and hence so is the map from \cref{constr:comparison_map_sigmanu_EASS}.
    Finally, cellularity of $\Syn_E$ is equivalent to $\sigma$ being conservative (\cref{cor:properties_syn_sigma}), thereby showing the final claim of \cref{item:converse_comp_map}.
\end{proof}

\begin{remark}
    \label{rmk:sigma_nu_and_completion}
    The reason why completion appears in the comparison between $\sigma(\nu X)$ and the Adams filtration is due to our use of cosimplicial objects in defining the latter.
    Note that the colimit of the filtration $\sigma (\nu X)$ is always naturally isomorphic to~$X$, but that this filtration may not be complete.
    On the other hand, the cosimplicial definition of $\ASS_E(X)$ is always a complete filtration, but its colimit may not be~$X$.
    These convergence problems are in fact the same, since the natural map from the former to the (d\'ecalage of the) latter is completion.
    We expect that $\sigma(\nu X)$ is isomorphic, as a filtered spectrum, to the d\'ecalage of a filtered Adams resolution of $X$ (see \cref{rmk:cosimplicial_Adams_vs_filtered_resolution}), so that \cref{rmk:cosimplicial_Adams_vs_filtered_resolution} would also explain this difference in packaging of the convergence problem.
    (On the other hand, incorporating monoidal structures as in the next remark would require working from the second page onward, for reasons explained by \cref{rmk:only_cosimplicial_decalage_is_monoidal}.)
\end{remark}

\begin{remark}[Monoidal version]
    \label{rmk:identifying_sigma_nu_monoidally}
    The left-hand side of the comparison map of \cref{constr:comparison_map_sigmanu_EASS} is naturally a lax symmetric monoidal functor in $X$; see \cref{rmk:sigmanu_lax_monoidal}.
    Recall from \cref{rmk:decalage_lax_monoidal} that $\Dec^\Delta$ is naturally lax symmetric monoidal.
    If $E$ carries considerably more structure, then we can give the functor $X\mapsto E^{[\bullet]} \otimes X$ a lax symmetric monoidal structure as well, in which case the comparison map
    \[
        \sigma(\nu X) \to \Dec^\Delta(E^{[\bullet]}\otimes X)
    \]
    of \cref{constr:comparison_map_sigmanu_EASS} matches up these monoidal structures, as follows.
    \begin{itemize}
        \item If $E$ carries an $\E_1$-structure, then the semicosimplicial spectrum $E^{[\bullet]}$ naturally extends to a cosimplicial spectrum $\Delta \to \Sp$; see \cite[Construction~2.7]{mathew_naumann_noel_nilpotence_descent}.
        If $E$ carries an $\E_n$-structure for $1\leq n \leq \infty$, then using Dunn additivity, this construction turns $E^{[\bullet]}$ into a cosimplicial $\E_{n-1}$-ring.
        \item Consequently, if $E$ is $\E_n$ for $1\leq n\leq \infty$, then the functor
        \[
            \Sp \to \Sp^\Delta, \quad X \mapsto E^{[\bullet]}\otimes X
        \]
        is naturally a lax $\E_{n-1}$-monoidal functor.
        Postcomposing with the lax symmetric monoidal functor $\Dec^\Delta$, we obtain a lax $\E_{n-1}$-monoidal functor
        \[
            \Sp \to \FilSp, \quad X \mapsto \Dec^\Delta(E^{[\bullet]}\otimes X).
        \]
        In this case, the comparison map of \cref{constr:comparison_map_sigmanu_EASS} is an $\E_{n-1}$-monoidal natural transformation.
    \end{itemize}
    Said differently, by completing the filtered spectrum $\sigma(\opnu_E (\blank))$, through \cref{thm:signature_nu} we obtain a lax symmetric monoidal structure structure on $\Dec(\ASS_E(\blank))$, even if $E$ is merely homotopy-associative.
    Only when $E$ is $\E_\infty$ can we identify this structure concretely in terms of $E$; in general, we have to work at the level of synthetic spectra.
\end{remark}

\begin{remark}
    A slightly different way of showing that $\nu X$ is $\tau$-complete if and only if $X$ is $E$-nilpotent complete is given by \cite[Proposition~A.13]{burklund_hahn_senger_manifolds}.
    They also show that this is equivalent to $\nu X$ being $\nu E$-nilpotent complete.
    Note that this also follows from our proof: using \cref{ex:symm_mon_if_one_is_E}, we see that \eqref{eq:comp_map_before_sigma} is precisely the $\nu E$-nilpotent completion of $\nu X$.
\end{remark}

\begin{remark}
    \label{rmk:general_signature_maps_Y_X}
    One could try to run the same argument for the Adams spectral sequence for maps $[Y,X]_*$ for a general spectrum $Y$.
    Due to the nature of the definition of $\Syn_E$ as in \cite{pstragowski_synthetic}, this is only sensible when $Y$ is a filtered colimit of finite $E$-projective spectra.
    (The above proof breaks down in general because when taking \eqref{eq:comp_map_before_sigma} mod $\tau$, we obtain a resolution only by \emph{relative injectives}.
    This is also related to $\Syn_E$ being generated by the synthetic analogues of finite $E$-projectives; see \cref{prop:properties_of_Syn}\,\ref{item:Syn_generated_by_spheres}.)
    If $\pi_*E$ is a graded field, then this condition on $Y$ is vacuous, but not in general.
    To obtain the Adams spectral sequence for general $E$ and $Y$, one has to work with a different version of synthetic spectra, as in the work by Patchkoria--Pstr\k{a}gowski \cite{patchkoria_pstragowski_derived_inftycats}.
    In Theorem~5.60 of op.\ cit., when working in this different version, they identify the spectral sequence underlying the filtered mapping spectrum (see \cref{constr:deformation_enriched_FilSp}) from $\nu Y$ to $\nu X$ with the Adams spectral sequence for $[Y,X]_*$.
\end{remark}

\section{The synthetic Omnibus Theorem}
\label{sec:synthetic_omnibus}

Previously in \cref{sec:Syn_as_deformation}, we explained how the filtered Omnibus \cref{thm:filtered_omnibus}, as well as its truncated variants of \cref{thm:truncated_filtered_omnibus,thm:truncated_filtered_omnibus_generation}, directly imply a synthetic version, up to a reindexing.
For the convenience of the reader, we state this synthetic version here.
We leave the analogous re-indexing of the truncated versions to the reader, or refer to \cite[Theorems~2.21 and~2.28]{CDvN_part2} for a recorded version.

\begin{theorem}[Synthetic Omnibus]
    \label{thm:synthetic_omnibus}
    Let $X$ be a $\tau$-complete synthetic spectrum, and assume that in its signature spectral sequence, we have $\uR\uE_\infty^{*,*}=0$ \textbr{for instance, this happens if the spectral sequence converges strongly}.
    Let $x \in \uE_2^{n,s} = \oppi_{n,s}(X/\tau)$ be a nonzero class.
    Then the following are equivalent.
    \begin{enumerate}[label={\upshape(1\alph*)}]
        \item \label{item:syn_ombs_perm_cycle} The element $x$ is a permanent cycle.
        \item \label{item:syn_ombs_lifts} The element $x \in \oppi_{n,s}(X/\tau)$ lifts to an element of $\oppi_{n,s}X$.
    \end{enumerate}
    For any such lift $\alpha$ to $\oppi_{n,s}X$, the following are true.
    \begin{enumerate}[label={\upshape(2\alph*)}]
        \item \label{item:syn_omnbs_not_tau_torsion} If $x$ survives to page $r$, then $\tau^{r-2}\cdot \alpha$ is nonzero.
        \item \label{item:syn_ombs_detection} If $x$ survives to page $\infty$, then $\alpha$ maps to a nonzero element in $\pi_n X^{\tau=1}$, and this element is detected by~$x$.
    \end{enumerate}
    Moreover, if $x$ lifts to $X$, then there exists a lift $\alpha$ with either of the following additional properties.
    \begin{enumerate}[label={\upshape(3\alph*)}]
        \item \label{item:syn_omnbs_torsion_lift} If $x$ is the target of a $d_r$-differential, then $\tau^{r-1} \cdot \alpha = 0$.
        \item \label{item:syn_omnbs_detection_lift} If $\theta \in \pi_n X^{\tau=1}$ is detected by $x$, then $\alpha$ is sent to $\theta$ under $\oppi_{n,s}X \to \pi_nX^{\tau=1}$.
    \end{enumerate}
    Finally, we have the following generation statement.
    \begin{enumerate}[label={\upshape(4)}]
        \item \label{item:syn_omnbs_generation} Let $\set{\alpha_i}$ be a collection of elements of $\oppi_{n,*}X$ such that their mod $\tau$ reductions generate the permanent cycles in stem $n$.
        Then the $\tau$-adic completion of the $\Z[\tau]$-submodule of $\oppi_{n,*}X$ generated by the $\set{\alpha_i}$ is equal to $\oppi_{n,*}X$.
    \end{enumerate}
\end{theorem}
\begin{proof}
    This follows by combining \cref{thm:filtered_omnibus} with \cref{prop:synthetic_vs_filtered_htpy_groups} and \cref{prop:deformation_sigma_tauinv}.
\end{proof}

If the synthetic spectrum $X$ in \cref{thm:synthetic_omnibus} is the synthetic analogue of a spectrum~$Y$, then \cref{thm:signature_nu} tells us that the signature of $\nu Y$ is the $E$-Adams spectral sequence for $Y$.
Moreover, in this case the convergence condition is precisely asking for the Adams spectral sequence for $Y$ to be strongly convergent.
(Indeed, because this filtration is left-concentrated in the sense of \cref{def:left_right_concentrated}, this follows from \cref{thm:conditional_convergence}.)
In this way, the filtered Omnibus Theorem together with the computation of $\sigma \circ \nu$ recovers the Omnibus Theorem for synthetic analogues of Burklund--Hahn--Senger \cite[Theorem~9.19]{burklund_hahn_senger_manifolds}.

This combination of \cref{thm:signature_nu,thm:synthetic_omnibus} lets information flow both ways.
On the one hand, we now see that if we use synthetic spectra to compute the homotopy of a synthetic analogue, this gives us new information about the underlying Adams spectral sequence.
On the other hand, we can also use this to import existing knowledge about Adams spectral sequence into synthetic spectra, thereby giving us a starting point for new computations.

\begin{remark}[Comparison of proofs]
    \label{rmk:comparison_omnibus_proofs}
    Our proof of the Omnibus Theorem is inspired by the one of Burklund--Hahn--Senger in \cite[Appendix~A]{burklund_hahn_senger_manifolds}.
    They identify the $\nu E$-Adams spectral sequence of $\nu X$ with (to use our terminology) the $\tau$-BSS for the $E$-ASS of $X$; see Theorem~A.8 of op.\ cit.\ for the precise meaning of this.\footnote{Note that they use second-page indexing, which can be obtained from ours using \cref{rmk:second_page_indexing_tau_BSS}. See also \cref{var:tau_BSS_synthetic}.}
    Generalising this to an arbitrary synthetic spectrum requires finding suitable replacements for these three spectral sequences.
    In our approach, we view the signature spectral sequence as the appropriate replacement for the $E$-ASS, and we do away with the the $\nu E$-ASS, going straight to the $\tau$-BSS of the signature spectral sequence.
    As the signature spectral sequence is defined using filtered spectra, the proof naturally takes place there, so that the synthetic version is a special case of the filtered version of \cref{sec:filtered_omnibus}.
    The computation of $\sigma \circ \nu$ is then the result that gives this a concrete meaning.
\end{remark}

As a consequence of \cref{thm:signature_nu}, we learn that the strict filtration on $\oppi_n X$ induced by $\oppi_{n,*}\nu X$ coincides with the $E$-Adams filtration, recovering \cite[Corollary~9.21]{burklund_hahn_senger_manifolds}.
See \cref{def:algebraic_topological_Adams_filtrations} for the definition of the (algebraic) Adams filtration.

\begin{corollary}[Geometric Adams filtration]
    \label{thm:geometric_Adams_filtration}
    Let $X$ be a spectrum, and let $n$ be an integer.
    Let $f \colon \S^n \to X$ be a map, and let $s\geq 0$.
    Then $f$ has \textbr{algebraic} $E$-Adams filtration \textbr{see \cref{def:algebraic_topological_Adams_filtrations}} at least $s$ if and only if there exists a factorisation
    \[
        \begin{tikzcd}
            \nu \S^n = \S^{n,0} \rar["\nu f"] \dar["\tau^s"'] & \nu X.\\
            \opSigma^{0,s} \nu \S^n = \S^{n,s} \ar[ur,dashed]
        \end{tikzcd}
    \]
    In other words, $f$ is of \textbr{algebraic} $E$-Adams filtration at least $s$ if and only if $\nu f$ is divisible by $\tau^s$.
    
    In particular, the \textbr{algebraic} $E$-Adams filtration on $\pi_nX$ coincides with
    \[
        F^s \, \pi_nX = \im\left(\begin{tikzcd}[ampersand replacement=\&]
    {[}\S^{n,s}{,}\ \nu X{]} \rar["\tau=1"] \& {[}\S^n{,}\ X{]}
    \end{tikzcd}\right).
    \]
\end{corollary}
\begin{proof}
    Combine \cref{prop:ASS_induced_filtration_is_Adams_filtration,prop:different_Adams_filtrations_agree} with \cref{thm:signature_nu}.
\end{proof}

Phrased differently, the functor $\nu$ sends a map $f\colon \S^n \to X$ of spectra to the map $\nu f$, which we think of as $f$ placed in Adams filtration $0$.
Then a $\tau^s$-division of $\nu f$ (should it exist) is a witness that $f$ has Adams filtration at least $s$.

More generally, the analogous result holds when $\S^n$ is replaced by a spectrum $Y$ that can be written as a filtered colimit of finite $E$-projective spectra.
For a general spectrum $Y$, one has to work with the different construction of synthetic spectra from \cite{patchkoria_pstragowski_derived_inftycats}, as explained in \cref{rmk:general_signature_maps_Y_X}; the general result in this case is \cite[Theorem~5.60\,(2)]{patchkoria_pstragowski_derived_inftycats}.

\begin{notation}
    \label{not:labelling_synthetic_htpy_elements}
    If $Y$ is a spectrum and $\alpha \in \oppi_n Y$ is an element, then there are two conventions one can take for naming elements of $\oppi_{*,*}\nu Y$.
    \begin{numberenum}
        \item One could denote $\nu(\alpha) \in \oppi_{n,0}\nu Y$ by the same symbol $\alpha$ again, and write $\widetilde{\alpha}$ for (a choice of) the maximal $\tau$-division thereof.
        This is the convention used by Burklund--Hahn--Senger in \cite{burklund_hahn_senger_manifolds}.
        \item One could denote (a choice of) the maximal $\tau$-division of $\nu(\alpha)$ by the same symbol $\alpha$.
        This is used, e.g., by Burklund in \cite{burklund_extension54}.
    \end{numberenum}
    The second convention is the one we will use.
    Conceptually, it emphasises the Adams filtration of an element in $\pi_n Y$ as an intrinsic invariant, so that we should be thinking of its maximal $\tau$-division as its `true origin'.
\end{notation}

\begin{remark}[Uniqueness]
    Note that the maximal $\tau$-division of $\nu(\alpha)$ may not be uniquely defined, as $\oppi_{n,*}\nu Y$ might contain $\tau$-power torsion.
    Specifically, if $\alpha$ has filtration $s$, then the maximal $\tau$-division of $\nu(\alpha)$ is defined up to $\tau^s$-power torsion in $\oppi_{n,s}\nu Y$.
    Via the Omnibus Theorem (if the convergence conditions are met), this is more or less saying that it is defined up to permanent cycles in filtrations $k\geq s$ that are hit by a differential of length $\leq k-1$.
\end{remark}

\begin{example}
    Consider the $1$-stem of the $\MU$-synthetic sphere.
    In the ordinary sphere, the element $\eta \in \oppi_1 \S$ has Adams--Novikov filtration $1$.
    Translated in synthetic terms, we have $\oppi_{1,s}C\tau = 0$ for $s \neq 1$, and $\oppi_{1,1}C\tau$ is isomorphic to $\Z/2$ with generator $h_0$.
    Since the $1$-stem does not receive or support any differentials (for degree reasons), it follows from \cref{cor:tau_torsion_free_iff_no_differentials} that $\oppi_{1,*}\S$ is $\tau$-power torsion free.
    It follows that $\tau$-inversion is an injection, and as $\oppi_1\S\cong \Z/2\cdot \eta$, we obtain an isomorphism
    \[
        \oppi_{1,*} \S \cong \Z/2[\tau] \cdot x,
    \]
    where $x$ is in bidegree $(1,1)$, where $x$ maps to $h_0$ under mod $\tau$ reduction, and where $\tau \cdot x$ is the synthetic analogue of the map $\eta \in \oppi_1\S$.
    
    The second convention above would be to use the symbol $\eta$ to denote the generator~$x$.
\end{example}

As an another example, in \cref{ex:subadditive_filtration}, we studied essentially $\oppi_{n,s}\S$ of the ($2$-local) $\MU$-synthetic sphere for $n\leq 3$, and also followed the second notational convention mentioned above.

%% file: syn/syn_models.tex
In this chapter, we discuss certain variants and modifications of synthetic spectra.
The first main result is that, under a mild condition, the $\infty$-category of synthetic spectra is equivalent to modules in $\FilSp$ over a filtered ring spectrum; see \cref{sec:filtered_models_syn}.
The second is that $\MU$-synthetic spectra are equivalent to (a small subcategory of) $\bC$-motivic spectra; see \cref{ssec:motivic}.
The first of these results requires a discussion of cellularity, and the second requires a discussion of evenness; these are discussed in \cref{ssec:cellular} and \cref{ssec:even_synthetic_spectra}, respectively.

Most of the results in this section are not new, but are presented slightly differently compared to the literature.
In particular, our discussion of filtered models is once again focussed on the signature adjunction.
One minor new result is that our discussion of evenness shows that even synthetic spectra are closed under limits (\cref{cor:even_synthetic_closed_limits}), but our approach has the downside of only applying to the cellular case.
Our discussion of the relation to motivic spectra is no more than a quick review, and will not be used heavily.

\section{Cellularity}
\label{ssec:cellular}

\begin{definition}
    Let $\C$ be a monoidal deformation.
    The \defi{cellular subcategory} $\C^\cell$ of $\C$ is the smallest subcategory closed under colimits that contains $\rho(\S^{n,s})$ for all $n$ and $s$.
    We say that $\C$ is \defi{cellular} if $\C^\cell = \C$.
\end{definition}

\begin{remark}
    A different way of stating the definition of cellularity is that the object $\1_\C$ of $\C$ generates $\C$ as a $\FilSp$-linear category under colimits.
    Phrased in this way, the definition can be extended to a general deformation, where one has to provide a choice of an object to play the role of $\1_\C$.
    In this case, one might speak of the deformation being a \emph{monogenic} $\FilSp$-linear $\infty$-category.
\end{remark}

\begin{proposition}
    Let $\C$ be a monoidal deformation.
    \begin{numberenum}
        \item The $\infty$-category $\C^\cell$ is a presentable stable $\infty$-category.
        \item The subcategory $\C^\cell \subseteq \C$ is closed under colimits and tensor products.
    \end{numberenum}
    Consequently, $\C^\cell$ inherits the structure of a presentably monoidal $\infty$-category from $\C$, and we have a colocalisation adjunction
    \[
        \begin{tikzcd}[column sep=3.5em]
            \C^\cell \rar[shift left, hook] & \C \lar[shift left, "(\blank)^\cell"]
        \end{tikzcd}
    \]
    with a monoidal left adjoint.
    If $\C$ is a symmetric monoidal deformation, then $\C^\cell$ is a symmetric monoidal subcategory and the inclusion is a symmetric monoidal functor.
\end{proposition}
\begin{proof}
    Presentability follows because it is generated under small colimits by a small set of objects.
    Closure under colimits is clear, and closure under tensor products follows because the tensor product of $\C$ preserves colimits in each variable separately.
    The Adjoint Functor Theorem then yields the desired adjunction.
\end{proof}

We refer to the right adjoint $\C \to \C^\cell$ as the \defi{cellularisation functor}, which is canonically lax monoidal (symmetric if $\C$ is a symmetric deformation).
It follows from this adjunction that for every $X \in \C$, the counit $X^\cell \to X$ induces an isomorphism
\[
    \oppi_{*,*}X \cong \oppi_{*,*}(X^\cell).
\]
In this sense, one may think of the passage from $\C$ to $\C^\cell$ as similar to the passage from topological spaces up to homotopy equivalence, to topological spaces up to weak homotopy equivalence.

This is also indicative of a larger phenomenon: when using a deformation to study the resulting signature spectral sequences, we may as well work with its cellularisation.

\begin{proposition}
    Let $\C$ be a \textbr{symmetric} monoidal deformation.
    \begin{numberenum}
        \item The functor $\sigma$ is conservative if and only if $\C$ is cellular.
        \item We have natural commutative diagrams of lax \textbr{symmetric} monoidal functors
        \[
            \begin{tikzcd}
                & \C^\cell \ar[dr,hook]  & \\
                \FilSp \ar[rr,"\rho"'] \ar[ur,"\rho"] & & \C 
            \end{tikzcd} \qquad \text{and} \qquad \begin{tikzcd}
                \C \ar[rr,"\sigma"] \drar["(\blank)^\cell"'] &  & \FilSp.\\
                & \C^\cell \urar["\sigma"'] & 
            \end{tikzcd}
        \]
        \item The adjunction $\rho \dashv \sigma$ restricts to an adjunction between $\FilSp$ and $\C^\cell$, which we denote by
        \[
            \begin{tikzcd}[column sep=3.5em]
                \FilSp \rar[shift left, "\rho^\cell"] & \C^\cell. \lar[shift left, "\sigma^\cell"]
            \end{tikzcd}
        \]
    \end{numberenum}
\end{proposition}
\begin{proof}
    The first claim is a restatement of \cref{prop:deformation_sigma_tauinv}\,\ref{item:conservative_sigma_condition}.
    The functor $\rho$ lands in $\C^\cell$ because it preserves colimits and because $\FilSp$ is generated under colimits by the bigraded spheres.
    To see that $\sigma$ factors over cellularisation, we have to show that $\sigma$ sends the counit $X^\cell \to X$ to an isomorphism for every $X \in \C$.
    By adjunction, this follows from $\rho$ landing in the cellular subcategory.
    The final claim follows immediately from this.
\end{proof}

\begin{proposition}
    Let $\C$ be a monoidal deformation.
    Then $\C$ is cellular if and only if the generic fibre of $\C$ is generated under colimits by $\rho(\S^{n,0}[\tau^{-1}])$ for $n\in\Z$ and the special fibre is generated under colimits by $\rho(C\tau\otimes \S^{n,s})$ for $n,s\in\Z$.
\end{proposition}
\begin{proof}
    The deformation $\C$ is cellular if and only if the counit $X^\cell \to X$ is an isomorphism for all $X\in\C$.
    We may check whether this map is an isomorphism after inverting and after modding out by $\tau$; see \cref{prop:deformation_pullback_tau_complete_and_invert_tau}.
    One checks that the cellularisation adjunction reduces to the analogous adjunction for generation by the $\tau$-invertible or the mod $\tau$ bigraded spheres, respectively.
    Since the filtered spheres $\S^{n,s}$ and $\S^{n,t}$ become isomorphic upon $\tau$-inversion, the case of the generic fibre reduces to the stated claim.
\end{proof}

Next, we turn to the case of synthetic spectra.
To the author's knowledge, it is not known if for every $E$ of Adams type, the deformation $\Syn_E$ is cellular.
It is known to hold in many cases.
In the case $E=\F_p$, it is relatively easy to see that it is cellular; see \cite[Lemma~2.3]{christian_jack_synthetic_j}.
Pstr\k{a}gowski proved in his foundational work that $\MU$-synthetic spectra are cellular; see \cite[Theorem~6.2]{pstragowski_synthetic}.
Later, Lawson generalised this to arbitrary connective $E$.

\begin{theorem}[Lawson \cite{lawson_cellular}]
    \label{thm:E_connective_implies_cellular_SynE}
    Let $E$ be a homotopy-associative ring spectrum of Adams type.
    If $E$ is connective, then $\Syn_E$ is cellular.
\end{theorem}

In joint work with Barkan, we prove cellularity in an important nonconnective case.

\begin{theorem}[\cite{barkan_vN_cellular_chromatic}]
    \label{thm:Morava_E_theory_cellular}
    Let $E$ be a Morava E-theory at an arbitrary height and prime.
    Then $\Syn_E$ is cellular.
\end{theorem}

Note that \cref{thm:E_connective_implies_cellular_SynE} does not imply \cref{thm:Morava_E_theory_cellular}.
Indeed, if $E$ denotes Morava E-theory, then (except at height $0$) we have that $\F_p$ is $E$-acyclic for all $p$, while if $R$ is a connective homotopy-ring spectrum, then there is at least one $p$ for which $\F_p$ is not $R$-acyclic.
It follows that $\Syn_E$ is not equivalent to $\Syn_R$ for any connective homotopy-ring spectrum $R$.

\section{Filtered models}
\label{sec:filtered_models_syn}

Under certain conditions, there is a description of $\Syn_E$ as modules in $\FilSp$ over a certain filtered ring spectrum.
These results have appeared before \cite[Appendix~C]{burklund_hahn_senger_Rmotivic} \cite{lawson_cellular} \cite[Sections~3.3 and~3.5]{pstragowski_perfect_even_filtration}.
We include a discussion here to highlights its connection to the signature adjunction.
We learned this approach from Shaul Barkan.
See also \cite[Appendix~A]{burklund_hahn_senger_Rmotivic}.

\begin{theorem}
    \label{thm:deformation_filtered_model}
    Let $\C$ be a monoidal deformation.
    Suppose that $\C$ satisfies the following conditions.
    \begin{letterenum}
        \item The monoidal unit $\1_\C$ of $\C$ is compact.
        \item The deformation $\C$ is cellular.
    \end{letterenum}
    Then the adjunction
    \[
        \begin{tikzcd}[column sep=3.5em]
            \FilSp \ar[r, shift left, "\rho"] & \C \ar[l, shift left, "\sigma"]
        \end{tikzcd}
    \]
    is a monadic adjunction that identifies the $\infty$-category underlying $\C$ with $\Mod_{\sigma(\1_C)}(\FilSp)$.
    If $\C$ is a symmetric monoidal deformation, then $\sigma(\1_\C)$ is a filtered $\E_\infty$-ring, and the identification of $\C$ with $\Mod_{\sigma(\1_\C)}(\FilSp)$ is naturally one of symmetric monoidal $\infty$-categories.
\end{theorem}


\begin{proof}
    It is enough to check the three conditions of \cite[Proposition~5.29]{mathew_naumann_noel_nilpotence_descent}.
    The second, that $\sigma$ preserves colimits, follows from \cref{prop:deformation_sigma_tauinv}\,\ref{item:colims_sigma_condition}.
    The first, that the adjunction satisfies the projection formula, follows from $\sigma$ preserving colimits and \cref{lem:projection_formula_deformation}.
    Finally, the third condition, that $\sigma$ is conservative, follows from \cref{prop:deformation_sigma_tauinv}\,\ref{item:conservative_sigma_condition}.
\end{proof}

\begin{remark}
    Alternatively, as explained by Pstr\k{a}gowski, the hypotheses on $\C$ can be interpreted as saying that it has a single compact filtered generator, so that a filtered version of Schwede--Shipley applies; see \cite[Proposition~3.16]{pstragowski_synthetic}.
\end{remark}

Specialising to the synthetic setting, we immediately obtain the following.

\begin{corollary}
    \label{thm:filtered_model_SynE}
    Let $E$ be a homotopy-associative ring spectrum of Adams type.
    Then the adjunction
    \[
        \begin{tikzcd}[column sep=3.5em]
            \FilSp \ar[r, shift left, "\rho^\cell"] & \Syn_E^\cell \ar[l, shift left, "\sigma^\cell"]
        \end{tikzcd}
    \]
    is a monadic adjunction that identifies $\Syn_E^\cell$ with $\Mod_{\sigma(\nu\S)}(\FilSp)$ as a symmetric monoidal $\infty$-category.
    Under this equivalence, the $\tau$-completion of $(\nu X)^\cell$ \textbr{where $X$ is a spectrum} is sent to $\Tot(\Wh(E^{[\bullet]}\otimes X))$.
\end{corollary}
\begin{proof}
    This follows by combining \cref{thm:deformation_filtered_model} and \cref{thm:signature_nu}.
\end{proof}

Note, however, that this does not yet give an alternative description of $\Syn_E$ (even after forgetting the monoidal structure): we have not yet identified the filtered spectrum $\sigma(\nu\S)$ as a filtered $\E_1$-ring (let alone as a filtered $\E_\infty$-ring) without reference to $\Syn_E$.
We can do this using \cref{thm:signature_nu}, but only when $E$ is sufficiently structured; see \cref{rmk:identifying_sigma_nu_monoidally}.
However, this result only identifies the \emph{completion} of $\sigma(\nu\S)$; by the same theorem, $\sigma(\nu\S)$ is complete if and only if $\S$ is $E$-nilpotent complete.



The prime example of such a setting is $E = \MU$, which is an $\E_\infty$-ring, and the sphere is $\MU$-nilpotent complete (see \cref{ex:MU_BP_Fp_nilpotent_completion}).
As mentioned before, $\Syn_\MU$ is cellular, resulting in the following.

\begin{corollary}
    \label{cor:filtered_model_Syn_MU}
    There is a symmetric monoidal equivalence
    \[
        \Syn_\MU \simeq \Mod_{\Tot(\Wh(\MU^{[\bullet]}))}(\FilSp)
    \]
    under which $\sigma$ becomes the forgetful functor \textbr{as a lax symmetric monoidal functor}.
    Moreover, on $\MU$-nilpotent complete spectra, the functor $\nu$ is identified with the functor sending $X$ to $\Tot(\Wh(\MU^{[\bullet]}\otimes X))$ \textbr{as a lax symmetric monoidal functor}.
\end{corollary}


\begin{remark}[\cite{burklund_hahn_senger_Rmotivic}, Appendix~C]
    Because we can identify the $\tau$-completion of $\sigma(\nu \S)$ by \cref{thm:signature_nu}, we can ignore the nilpotent-completeness issues by passing to a slight modification of synthetic spectra.
    Namely, if $E$ is an $\E_1$-ring of Adams type, then the functor~$\sigma$ also induces an equivalence
    \[
        \Mod_{(\nu \S)_\tau^\wedge}(\Syn_E^\cell) \simeqto \Mod_{\Tot(\Wh(E^{[\bullet]}))}(\FilSp),
    \]
    which is symmetric monoidal if $E$ is an $\E_\infty$-ring.
    This follows by applying \cref{thm:deformation_filtered_model} to the symmetric monoidal deformation $\Mod_{(\nu\S)_\tau^\wedge}(\Syn_E)$ and using that $\sigma$ preserves $\tau$-completion (\cref{prop:deformation_sigma_tauinv}\,\ref{item:sigma_pres_modtau_taucompl}).
    This is the formulation of the result given by Burklund--Hahn--Senger in \cite[Proposition~C.22]{burklund_hahn_senger_Rmotivic}.
    (Of course, the same reasoning holds in the context of \cref{thm:deformation_filtered_model}, yielding a filtered model for modules over the $\tau$-completion of the unit of $\C$.)
\end{remark}

A more elementary but curious-looking example is the case where $E=\S$.
In this case, all nilpotent completeness conditions become vacuous, and we learn the following.

\begin{corollary}
    \label{cor:sphere_synthetic_filtered_model}
    There is a symmetric monoidal equivalence
    \[
        \Syn_\S \simeq \Mod_{\Wh\S}(\FilSp)
    \]
    under which $\sigma$ becomes the forgetful functor and under which $\nu$ becomes the Whitehead filtration functor \textbr{and these identifications are as lax symmetric monoidal functors}, and under which the homological t-structure is identified with the diagonal t-structure.
\end{corollary}

One can check that for every Adams type $E$, the functor $\sigma \colon \Syn_E \to \FilSp$ factors through modules over $\Wh \S$, so that this is in a sense the minimal structure present on $E$-synthetic spectra as $E$ varies.

\section{Evenness}
\label{ssec:even_synthetic_spectra}

The goal of this section is to discuss \emph{even synthetic spectra} introduced in \cite[Section~5.2]{pstragowski_synthetic}, also studied in a slightly different form in \cite{pstragowski_perfect_even_filtration}.
For spectral sequences, \emph{evenness} is nothing more than the condition that the starting page is concentrated in certain even degrees (where the specific meaning of this depends on the indexing system being used).
On the synthetic side, even objects turn out to have geometric interpretations: the special fibre of even $\MU$-synthetic is the derived $\infty$-category of quasi-coherent sheaves over $\M_\fg$, see \cref{ex:even_MU_synthetic}.
Moreover, in the next section, we will see that even $\MU$-synthetic spectra can be viewed as \emph{motivic spectra}.

We will reinterpret the notion of evenness defined by Pstr\k{a}gowski through the language of deformations.
Accordingly, we begin by discussing this in the setting of filtered spectra.

\subsection{Even filtered spectra}

We remind the reader that we use first-page indexing on filtered spectra.

\begin{lemma}
    \label{lem:fil_characterisation_even}
    Let $X$ be a filtered spectrum.
    Then the following conditions are equivalent.
    \begin{letterenum}
        \item \label{item:even_transition_map_fil} For every $s$, the transition map
        \[
            X^{2s} \to X^{2s-1}
        \]
        is an isomorphism.
        In other words, for every $n$ and $s$, the map
        \[
            \tau \colon \oppi_{n,\,2s}X \to \oppi_{n,\,2s-1}X
        \]
        is an isomorphism.
        \item For every odd integer $s$, the associated graded spectrum $\Gr^s X$ vanishes.
        In other words, the homotopy groups $\oppi_{n,s}(X/\tau)$ vanish whenever $s$ is odd.\label{item:even_assoc_gr_fil}
        \item The filtered spectrum $X$ belongs to the smallest subcategory of $\FilSp$ generated under colimits by the filtered spheres $\S^{n,\,2s}$ for all $n$ and $s$.\label{item:X_is_even_filtered}
    \end{letterenum}
\end{lemma}
\begin{proof}
    By exactness, it follows that \ref{item:even_transition_map_fil} is equivalent to \ref{item:even_assoc_gr_fil}.
    We prove that \ref{item:even_transition_map_fil} is equivalent to \ref{item:X_is_even_filtered}.
    Write $\FilSp^\ev$ for the full subcategory of those filtered spectra $X$ satisfying \ref{item:even_transition_map_fil}, and write $\C$ for the full subcategory of $\FilSp$ generated by $\S^{n,\,2s}$ for all $n$ and $s$.
    Clearly $\FilSp^\ev$ is closed under colimits, so it follows that $\C \subseteq \FilSp^\ev$.
    To prove the other inclusion, by \cite[Corollary~2.5]{yanovski_monadic_tower}, it is enough to show that the functors $\pi_{n,2s}$ jointly detect isomorphisms on $\FilSp^\ev$.
    This follows from the fact that the functors $\pi_{n,s}$ jointly detect isomorphisms on $\FilSp$, and that when restricted to $\FilSp^\ev$, the functors $\pi_{n,\,2s}$ and $\pi_{n,\,2s-1}$ are naturally isomorphic (induced by $\tau$).
\end{proof}

\begin{definition}
\label{def:even_filtered_spectrum}
We say that a filtered spectrum is \defi{even} if it satisfies the equivalent conditions of \cref{lem:fil_characterisation_even}.
We write $\FilSp^\ev$ for the full subcategory of $\FilSp$ on the even filtered spectra.
\end{definition}

Evenness is a notion that is detectable on the spectral sequence: in first-page indexing on filtered spectra, it says that the first page vanishes in odd filtrations.

\begin{example}
    \leavevmode
    \begin{numberenum}
        \item A filtered sphere $\S^{n,s}$ is even if and only if $s$ is even.
        We may refer to a sphere of this form as an \emph{even filtered sphere}; note that this places no restrictions on the variable $n$.
        \item Let $X$ be a spectrum.
        Then its Whitehead filtration $\Wh X$ is even if and only if $\pi_*X$ is concentrated in even degrees, and likewise for $\Post X$.\qedhere
    \end{numberenum}
\end{example}

\begin{proposition}
    \label{prop:FilSp_ev_properties}
    \leavevmode
    \begin{numberenum}
        \item The $\infty$-category $\FilSp^\ev$ is a presentable stable $\infty$-category.
        \item The subcategory $\FilSp^\ev\subseteq \FilSp$ is closed under limits, colimits, and tensor products.
    \end{numberenum}
    Consequently, $\FilSp^\ev$ inherits a presentably symmetric monoidal structure from $\FilSp$, and the inclusion functor admits both a left and a right adjoint and is a symmetric monoidal functor.
\end{proposition}
\begin{proof}
    Presentability follows from characterisation~\ref{item:X_is_even_filtered} of \cref{lem:fil_characterisation_even}.
    Closure under limits and colimits follows from characterisation~\ref{item:even_transition_map_fil}, which in particular implies it is a stable subcategory.
    Closure under tensor products follows from characterisation~\ref{item:even_assoc_gr_fil}, using that $\Gr \colon \FilSp \to \grSp$ is symmetric monoidal (see \cref{rmk:Ctau_gives_Gr_symm_mon_str}), and that $\grSp^\ev \subseteq \grSp$ is obviously closed under tensor products.
    The inclusion admits a right adjoint because it preserves colimits.
    Because the even filtered spheres are compact, the inclusion also preserves compact objects, and hence admits a left adjoint.
\end{proof}

\begin{remark}
    \label{rmk:FilSp_ev_left_right_adjoint}
    One may check that the left adjoint to the inclusion $\FilSp^\ev \subseteq \FilSp$ sends a filtered spectrum $X$ to
    \[
        \begin{tikzcd}
            \dotsb \rar & X^1 \rar[equals] & X^1 \rar &  X^{-1} \rar[equals] & X^{-1} \rar & X^{-3} \rar[equals] & \dotsb
        \end{tikzcd}
    \]
    with $X^{-1}$ in filtrations $0$ and $-1$, and the right adjoint sends $X$ to
    \[
        \begin{tikzcd}
            \dotsb \rar & X^2 \rar[equals] & X^2 \rar &  X^0 \rar[equals] & X^{0} \rar & X^{-2} \rar[equals] & \dotsb
        \end{tikzcd}
    \]
    with $X^0$ in filtrations $0$ and $-1$.
\end{remark}

\begin{remark}
    \label{rmk:FilSp_ev_equiv_to_FilSp}
    There is a symmetric monoidal equivalence $p\colon \FilSp^\ev \simeqto \FilSp$ that `pinches' the transition maps together: it sends an even filtered spectrum $X$ to
    \[
        \begin{tikzcd}
            \dotsb \rar & X^{2} \rar & X^0 \rar & X^{-2} \rar & \dotsb.
        \end{tikzcd}
    \]
    For instance, this sends $\S^{n,\,2s}$ to $\S^{n,s}$.
    Formally, precomposition with $2\colon \Z^\op \to \Z^\op$ results in a symmetric monoidal functor $\FilSp \to \FilSp$.
    Restricting the domain to the even subcategory results in the claimed symmetic monoidal equivalence $p$.
    Indeed, this restricted functor has a two-sided inverse given by $d^*$, where $d\colon \Z^\op \to \Z^\op$ is the map of posets that sends $2s$ and $2s-1$ to~$s$.
\end{remark}

\begin{remark}
    \label{rmk:FilSp_ev_as_deformation}
    We can regard $\FilSp^\ev$ as a symmetric monoidal deformation in its own right, with deformation parameter $\tau^2$.
    Formally, we have a symmetric monoidal functor $2\colon \Z\to \Z$.
    Writing $i \colon \Z \to \FilSp$ for the functor from \cref{def:functor_Z_to_FilSp}, we see that $i\circ 2\colon \Z \to \FilSp$ is a symmetric monoidal functor that lands in the subcategory $\FilSp^\ev$.
    The underlying functor can be depicted by the diagram
    \[
        \begin{tikzcd}
            \dotsb \rar["\tau^2"] & \S^{0,-2} \rar["\tau^2"] & \S \rar["\tau^2"] & \S^{0,2} \rar["\tau^2"] & \dotsb.
        \end{tikzcd}
    \]
    Via \cref{not:rho_sigma_general_deformation}, this gives $\FilSp^\ev$ the structure of a symmetric monoidal deformation.
    The resulting right adjoint $\FilSp^\ev \to \FilSp$ is in fact the functor $p$ from \cref{rmk:FilSp_ev_equiv_to_FilSp}.
\end{remark}

\subsection{Even deformations}

We focus on the case of monoidal deformations.
We have a choice of using either characterisation \ref{item:even_assoc_gr_fil} or \ref{item:X_is_even_filtered} from \cref{lem:fil_characterisation_even} as the definition of evenness in a monoidal deformation.
In general these notions differ, but we will show that they agree under suitable hypotheses.
It is an arbitrary choice which of the two characterisations we regard as the definition.

\begin{definition}
\label{def:even_deformation}
    Let $\C$ be a monoidal deformation.
    The \defi{even subcategory} $\C^\ev$ of $\C$ is the smallest subcategory closed under colimits that contains $\rho(\S^{n,\,2s})$ for all~$n$ and~$s$.
    We say that an object $X \in \C$ is \defi{even} if it belongs to $\C^\ev$.
\end{definition}

Because our definition uses only the filtered spheres, it is only applicable in the cellular setting.
More precisely, the above definition implies that every even object is cellular.

\begin{proposition}
    Let $\C$ be a monoidal deformation.
    \begin{numberenum}
        \item The $\infty$-category $\C^\ev$ is a presentable stable $\infty$-category.
        \item The subcategory $\C^\ev \subseteq \C^\cell$ is closed under colimits and tensor products.
    \end{numberenum}
\end{proposition}

Only under the following hypotheses on $\C$ do we obtain a usable notion of evenness.

\begin{proposition}
\label{prop:deformation_even_characterisation}
    Let $\C$ be a cellular monoidal deformation with a compact unit.
    Suppose that the filtered spectrum $\sigma(\1_\C)$ is even.
    Then for every $X \in \C$, the following are equivalent.
    \begin{letterenum}
        \item The object $X$ belongs to $\C^\ev$.\label{item:X_is_even_deform}
        \item The homotopy groups $\oppi_{n,s}(X/\tau)$ vanish whenever $s$ is odd.\label{item:Xmodtau_even_deform}
        \item The filtered spectrum $\sigma X$ is even.\label{item:sigmaX_is_even_deform}
    \end{letterenum}
    Moreover, the adjunction $\rho \dashv \sigma$ restricts to an adjunction between $\FilSp^\ev$ and $\C^\ev$, which we denote by
    \[
        \begin{tikzcd}[column sep=3.5em]
            \FilSp^\ev \rar[shift left, "\rho^\ev"] & \C^\ev. \lar[shift left, "\sigma^\ev"]
        \end{tikzcd}
    \]
\end{proposition}

\begin{proof}
    Clearly \ref{item:Xmodtau_even_deform} is equivalent to \ref{item:sigmaX_is_even_deform}, because this is the case for filtered spectra by \cref{lem:fil_characterisation_even}.
    Using \cref{lem:projection_formula_deformation}\,\ref{item:proj_formula_dualisable}, for all $n$ and $s$ we have an isomorphism
    \[
        \sigma(\rho(\S^{n,s})) \cong \opSigma^{n,s}\sigma(\1_\C).
    \]
    It follows that if $\sigma(\1_\C)$ is even, then $\sigma(\rho(\S^{n,\,2s}))$ is even for all $n$ and $s$.
    Because the unit of $\C$ is compact, the functor $\sigma \colon \C \to \FilSp$ preserves colimits; see \cref{prop:deformation_sigma_tauinv}\,\ref{item:colims_sigma_condition}.
    Using this, we see that \ref{item:X_is_even_deform} implies \ref{item:sigmaX_is_even_deform}.
    To see that~\ref{item:sigmaX_is_even_deform} implies~\ref{item:X_is_even_deform}, one either reasons as in the proof of \cref{lem:fil_characterisation_even}, or one deduces it from the filtered case using \cref{thm:deformation_filtered_model}.
    Finally, these conditions show that when restricting $\sigma$ to $\C^\ev$, it lands in $\FilSp^\ev$, so that the adjunction restricts as indicated.
\end{proof}

In this case, evenness is preserved under limits too, which is not obvious from the definition of $\C^\ev$ from \cref{def:even_deformation}.

\begin{corollary}
    \label{cor:even_deformation_closed_limits}
    Let $\C$ be a monoidal deformation satisfying the conditions of \cref{prop:deformation_even_characterisation}.
    Then the subcategory $\C^\ev\subseteq \C$ is closed under limits.
\end{corollary}
\begin{proof}
    Because $\sigma$ preserves limits, this follows from \cref{prop:FilSp_ev_properties} and from characterisation \ref{item:sigmaX_is_even_deform} in \cref{prop:deformation_even_characterisation}.
\end{proof}

\begin{remark}
    \label{rmk:even_filtered_model}
    The filtered model of \cref{thm:deformation_filtered_model} can be adapted to the even case as well.
    If the conditions of \cref{prop:deformation_even_characterisation} hold, then $\C$ is equivalent to modules over $\sigma(\1_\C)$ by \cref{thm:deformation_filtered_model}, which by assumption is an $\E_\infty$-algebra in $\FilSp^\ev$.
    As $\FilSp^\ev$ is a monoidal subcategory of $\FilSp$, it follows that $\sigma$ induces an equivalence of $\infty$-categories
    \[
        \C^\ev \simeqto \Mod_{\sigma(\1_\C)}(\FilSp^\ev).
    \]
    Postcomposing this with the (symmetric monoidal) equivalence $p\colon \FilSp^\ev \simeqto \FilSp$ from \cref{rmk:FilSp_ev_equiv_to_FilSp}, we obtain an equivalence
    \[
        \C^\ev \simeqto \Mod_{p(\sigma(\1_\C))}(\FilSp).
    \]
    If $\C$ is symmetric monoidal, then these equivalences are of symmetric monoidal $\infty$-categories.
\end{remark}

\subsection{Even synthetic spectra}

We now specialise to the case of synthetic spectra.
We do not get a good notion of evenness in $\Syn_E$ for every $E$, because the signature of the synthetic sphere may not be even.
The following condition on $E$ will ensure this.

\begin{definition}[\cite{pstragowski_synthetic}, Definition~5.8]
    Let $E$ be a homotopy associative ring spectrum.
    \begin{numberenum}
        \item A finite spectrum $P$ is called \defi{finite even $E$-projective} if $E_*P$ is a projective $E_*$-module and is concentrated in even degrees.
        \item We say that $E$ is of \defi{even Adams type} if it is of Adams type (see \cref{def:Adams_type}) and can be written as a filtered colimit of finite even $E$-projective spectra.
    \end{numberenum}
\end{definition}

If $E$ is of even Adams type, then it follows that $E_*E$ is also concentrated in even degrees.
Beware that asking $E$ to be of even Adams type is stronger than asking for $E_*$ to be concentrated in even degrees and for $E$ to be of Adams type, as the following example illustrates.

\begin{example}
    \leavevmode
    \begin{numberenum}
        \item The sphere spectrum $\S$ is of even Adams type.
        \item The ring spectrum $\MU$ is of even Adams type; see \cite[Example~5.9]{pstragowski_synthetic}.
        \item The ring spectrum $\F_p$ is not of even Adams type, even though $\pi_*\F_p$ is even and $\F_p$ is of Adams type.
        Indeed, $\pi_*(\F_p \otimes \F_p)$ is the ($p$-primary) dual Steenrod algebra, which is not concentrated in even degrees.
        \item On the other hand, if $E$ is Landweber exact and $E_*$ is concentrated in even degrees, then $E$ is in fact of even Adams type; see \cite[Example~5.9]{pstragowski_synthetic}.\qedhere
    \end{numberenum}
\end{example}

\begin{lemma}
    \label{lem:split_comodule_even_odd}
    Let $(A,\Gamma)$ be a graded Hopf algebroid, with $A$ and $\Gamma$ concentrated in even degrees.
    Write
    \[
        \grComod_{(A,\Gamma)}^\ev \qquad \text{and} \qquad \grComod_{(A,\Gamma)}^\odd
    \]
    for the full subcategories of $\grComod_{(A,\Gamma)}$ on those graded comodules that are concentrated in even, respectively odd, degrees.
    Then we have a symmetric monoidal equivalence
    \[
        \grComod_{(A,\Gamma)} \simeq \grComod_{(A,\Gamma)}^\ev \times \grComod_{(A,\Gamma)}^\odd.
    \]
\end{lemma}
\begin{proof}
    This is immediate.
\end{proof}

The characterisation of evenness looks slightly different in the synthetic setting, due to the reindexing when passing from filtered to synthetic spheres from \cref{prop:synthetic_vs_filtered_htpy_groups}.
Namely, because of the relation
\[
    \rho(\S^{n,s}_\fil) = \S^{n,\,s-n}_\syn,
\]
we should think of a synthetic sphere $\S^{n,s}_\syn$ as being even when $n+s$ is even.

\begin{lemma}
    \label{lem:even_synthetic_spheres}
    Let $E$ be a homotopy-associative ring spectrum of even Adams type.
    Then $\sigma(\S_\syn)$ is even.
\end{lemma}
\begin{proof}
    Recall from \cref{prop:synthetic_vs_filtered_htpy_groups} that we have an isomorphism
    \[
        \oppi_{n,\,s}^\syn (\blank) \cong \oppi_{n,\,s+n}^\fil(\sigma (\blank)).
    \]
    It follows that $\sigma X$ is even if and only if $\oppi_{n,s}^\syn(X/\tau)$ vanishes whenever $n+s$ is odd.
    By \cref{ex:mod_tau_maps_between_nu}, we have an isomorphism
    \[
        \oppi_{n,s}(\S_\syn/\tau) \cong \Ext^{s,\,n+s}_{E_*E}(E_*,E_*) = \Ext^s_{E_*E}(E_*[n+s],\, E_*).
    \]
    Using the splitting of $\grComod_{E_*E}$ from \cref{lem:split_comodule_even_odd}, it follows that these groups vanish whenever $n+s$ is odd, proving that $\sigma(\S_\syn)$ is even.
\end{proof}

Translated to the synthetic setting, \cref{prop:deformation_even_characterisation} takes the following form.

\begin{proposition}
    \label{prop:even_synthetic_spectra}
    Let $E$ be a homotopy-associative ring spectrum of even Adams type.
    Let $X$ be a cellular $E$-synthetic spectrum.
    Then the following are equivalent.
    \begin{letterenum}
        \item The synthetic spectrum $X$ belongs to the smallest subcategory of $\Syn_E^\cell$ generated under colimits by the synthetic spheres $\S^{n,s}$ where $n+s$ is even.\label{item:syn_even_colimit_generation}
        \item The homotopy groups $\oppi_{n,s}(X/\tau)$ vanish whenever $n+s$ is odd; in other words, the second page of the signature spectral sequence of $X$ is concentrated in even total degree.\label{item:syn_even_modtau}
        \item The filtered spectrum $\sigma X$ is even in the sense of \cref{def:even_filtered_spectrum}.\label{item:syn_sigma_is_even}
    \end{letterenum}
    If this is the case, then we say that $X$ is \defi{even}.
    Moreover, the adjunction $\rho \dashv \sigma$ restricts to an adjunction
    \[
        \begin{tikzcd}[column sep=3.5em]
            \FilSp^\ev \rar[shift left, "\rho^\ev"] & (\Syn_E^\cell)^\ev. \lar[shift left, "\sigma^\ev"]
        \end{tikzcd}
    \]
\end{proposition}
\begin{proof}
    Using \cref{lem:even_synthetic_spheres}, it follows that \cref{prop:deformation_even_characterisation} applies.
\end{proof}

\begin{remark}
    The restriction to the cellular case is an artefact of our use of filtered spectra to define the notion of evenness.
    In the case of synthetic spectra, one can define a notion of evenness that does not require cellularity assumptions: see \cite[Section~5.2]{pstragowski_synthetic}.
    It follows from Theorem~5.13 of op.\ cit.\ that if $\Syn_E$ is cellular, then these two notions of evenness coincide.
    Moreover, \citeauthor{pstragowski_synthetic} proves that the homological t-structure on $\Syn_E$ restricts to a t-structure on $\Syn_E^\ev$ whose heart is equivalent to $\grComod_{E_*E}^\ev$; see Remark~5.11 of op.\ cit.
\end{remark}

\begin{example}
    To see why we have to restrict to $E$ of even Adams type, consider the case $E=\F_2$.
    Then for every $s\geq 0$, we have an isomorphism $\oppi_{0,s}(C\tau) \cong \F_2\cdot h_0^s$.
    It follows that $\opsigma(\S_\syn)$ is not an even filtered spectrum.
    As a result, the subcategory of $\Syn_{\F_2}$ generated under colimits by the spheres $\S^{n,s}$ for $n+s$ even is not related to $\FilSp^\ev$ in a natural way.
\end{example}

\begin{corollary}
    \label{cor:even_synthetic_analogue}
    Let $E$ be a homotopy-associative ring spectrum of even Adams type.
    Let $X$ be a spectrum.
    Then $(\nu X)^\cell$ is even if and only if $E_*X$ is concentrated in even degrees.
\end{corollary}
\begin{proof}
    This follows immediately from the isomorphism from \cref{ex:mod_tau_maps_between_nu}:
    \[
        \oppi_{n,s}(\nu X/\tau) \cong \Ext^{s,\,n+s}_{E_*E}(E_*,\, E_*X) = \Ext^{s}_{E_*E}(E_*[n+s],\, E_*X),
    \]
    combined with the splitting of \cref{lem:split_comodule_even_odd}.
\end{proof}

\begin{example}
    \label{ex:even_synthetic_spheres}
    By \cref{cor:even_synthetic_analogue}, the synthetic sphere $\S^{n,s}$ is even if and only if $n+s$ is even.
    In particular, the sphere $\S^{0,-1}$ is \emph{not} even, and as a result, the map $\tau$ does not live in $\Syn_E^\ev$.
    However, the map $\tau^2 \colon \S^{0,-2}\to\S$ does live in this subcategory.
    
    If we restrict our attention to only the even synthetic spectra, then it would be more useful to use the letter $\tau$ for what we would otherwise denote by $\tau^2$; in other words, to regard $\tau^2$ as the deformation parameter.
    (The resulting deformation adjunction would then factor as the deformation structure on $\FilSp^\ev$ from \cref{rmk:FilSp_ev_as_deformation} followed by the adjunction $\rho^\ev \dashv\sigma^\ev$.)
    From our perspective, this is what happens in the literature when \emph{motivic spectra} are used to study the stable stems; see \cref{ssec:motivic}, particularly \cref{thm:motivic_is_synthetic}, for more information.
\end{example}

Our approach naturally gives us limit-closure of even synthetic spectra, at least up to cellularisation.
Note that the limit-closure of even synthetic spectra is absent from the discussion in \cite[Section~5.2]{pstragowski_synthetic}.

\begin{corollary}
    \label{cor:even_synthetic_closed_limits}
    Let $E$ be a homotopy-associative ring spectrum of even Adams type.
    Then the subcategory $(\Syn_E^\cell)^\ev \subseteq \Syn_E^\cell$ is closed under limits.
    In particular, if $\Syn_E$ is cellular, then Pstr\k{a}gowski's subcategory $\Syn_E^\ev \subseteq \Syn_E$ is closed under limits.
\end{corollary}
\begin{proof}
    This follows from \cref{cor:even_deformation_closed_limits}.
\end{proof}

\begin{remark}
    \label{rmk:Syn_even_filtered_model}
    Via \cref{rmk:even_filtered_model}, the filtered models for $\Syn_E$ from \cref{sec:filtered_models_syn} also carry over to $\Syn_E^\ev$.
    For instance, we obtain a symmetric monoidal equivalence
    \[
        \Mod_{(\nu\S)_\tau^\wedge}(\Syn_E^\ev) \simeq \Mod_{\Tot(\tau_{\geq 2*}(E^{[\bullet]}))}(\FilSp).
    \]
\end{remark}

\begin{example}
    \label{ex:even_MU_synthetic}
    The most relevant case for us is where $E=\MU$.
    Since $\MU$-synthetic spectra are cellular, the above discussion applies.
    Recall that there is an equivalence
    \[
        \grComod_{\MU_*\MU}^\ev \simeq \QCoh(\M_\fg).
    \]
    As a result, the special fibres of $\Syn_\MU^\ev$ and $\Synhat{}_\MU^\ev$ are, respectively,
    \[
        \mathrm{IndCoh}(\M_\fg) \qquad \text{and} \qquad \D(\QCoh(\M_\fg)).
    \]
    From the geometric perspective, the category of all graded $\MU_*\MU$ is a little more awkward, where we have to work with two sheaves (the even and odd parts) separately, so that only even $\MU$-synthetic spectra have this more elegant geometric description.
\end{example}

\section{Motivic spectra}
\label{ssec:motivic}

We assume basic familiarity with motivic spectra.

\begin{notation}
    We write $\Sp_\bC$ for Morel--Voevodsky's symmetric monoidal $\infty$-category of \emph{motivic spectra} over $\Spec \bC$.
    We write $\Sp_\bC^\cell$ for the smallest stable subcategory of $\Sp_\bC$ that is closed under colimits and that contains both $\G_m$ and the (categorical) suspension of the unit.
\end{notation}

The objects $\G_m$ and the suspension of the unit are the two `circles' of motivic homotopy theory, making it have a notion of bigraded homotopy, and explaining the above notion of cellularity.


Motivic homotopy theory is by nature an algebro-geometric homotopy theory.
It is very surprising then that over $\bC$, it turns out to have a deep connection to the Adams--Novikov spectral sequence for ordinary spectra.
\begin{itemize}
    \item Levine \cite{levine_slice_ANSS} shows that the Betti realisation of the slice spectral sequence for the $\bC$-motivic sphere spectrum is isomorphic to the (d\'ecalage of the double speed) Adams--Novikov filtration for the sphere spectrum (more generally, Levine shows this over an algebraically closed field of characteristic zero).
    \item Hu--Kriz--Ormsby \cite{hu_kriz_ormsby_motivic_alg_closed} show that at the prime $2$, differentials in the $\bC$-motivic Adams--Novikov spectral sequence for the motivic sphere can be deduced formally from differentials in the ordinary Adams--Novikov spectral sequence for the sphere.
    (More generally, they work over an algebraically closed field of characteristic zero.)
    Stahn \cite{stahn_motivic_odd_primes} showed the analogous odd-primary version of this.
    \item There is a twisted endomorphism $\tau$ of the motivic sphere spectrum, and its cofibre admits an $\E_\infty$-structure.
    Gheorghe--Wang--Xu \cite{gheorghe_wang_xu_special_fibre} show that modules over $C\tau$ in $(\Sp_\bC^\cell)_p^\wedge$ is equivalent to the derived $\infty$-category of $\BP_*\BP$-comodules.
\end{itemize}

\begin{remark}
    Bachmann--Kong--Wang--Xu \cite{bachmann_kong_wang_xu_chow_tstructures} have generalised the results by Gheorghe--Wang--Xu on the structure of $\bC$-motivic spectra to an arbitrary base field.
    The answer is more complicated; very roughly speaking, the arithmetic of the base field starts to play a big role (which was not visible in the case of $\bC$ because it is algebraically closed).
\end{remark}

The above suggests that the $\bC$-motivic category should have a close relation to the $\MU$-synthetic (or $\BP$-synthetic) category.
This is true in a very strong sense: up to $p$-completion, $\MU$-synthetic spectra form a full subcategory of $\bC$-motivic spectra.

\begin{theorem}[Gheorghe--Isaksen--Krause--Ricka, \citeauthor{pstragowski_synthetic}]
    \label{thm:motivic_is_synthetic}
    Let $p$ be an arbitrary prime.
    There exists a symmetric monoidal equivalence
    \[
        (\Sp_\bC^\cell)_p^\wedge \simeqto (\Syn_\MU^\ev)_p^\wedge
    \]
    that sends $\tau$ to $\tau^2$.
    Under this equivalence, Betti realisation corresponds to $\tau$-inversion.
\end{theorem}


\begin{proof}
    See \cite[Theorem~6.12]{mmf} or \cite[Theorem~7.34]{pstragowski_synthetic}.
    These two equivalences are related via the equivalence of \cref{cor:filtered_model_Syn_MU} (or rather, the even version of it; see \cref{rmk:Syn_even_filtered_model}).
\end{proof}

The aforementioned connections between $\bC$-motivic spectra and Adams--Novikov spectral sequences now match up with synthetic notions we previously introduced.
For instance, the result of \cite{hu_kriz_ormsby_motivic_alg_closed} now corresponds to the Omnibus Theorem; see also \cite[Remark~A.2]{burklund_hahn_senger_manifolds}.

\begin{warning}
Unlike what the case of synthetic or ordinary spectra might suggest, the two modifications that we have to do on $\Sp_\bC$ are quite drastic.
\begin{itemize}
    \item Most varieties are not cellular, so that $\Sp_\bC^\cell$ does not see a lot of the algebraic geometry contained in that category.
    \item Unlike in the case of ordinary spectra, rational motivic spectra are highly nontrivial and contain a lot of interesting information.
    Rational objects are killed by $p$-completion, so we lose a lot of information by passing to $p$-complete objects.
\end{itemize}
\end{warning}

As explained in \cref{ch:intro_to_part1}, cellular motivic spectra were used in computational stable homotopy theory in various ways.
Anachronistically, the equivalence of \cref{thm:motivic_is_synthetic} explains why motivic spectra in particular were so useful in this regard.
More seriously, as explained in \cite[Section~1.1.2]{isaksen_wang_xu_dimension_90}, this provides a way to make these computational arguments using a much lighter technical setup, so that it does not logically depend on the setup of motivic homotopy theory.
This is useful as this usage of motivic homotopy theory relies on deep results, such as Voevodsky's computation of the motivic dual Steenrod algebra \cite{voevodsky_motivic_F2_cohomology,voevodsky_motivic_EM_spaces}.
The analogous computation in the synthetic setting is a more straightforward adaptation of the spectral one; see \cite[Section~6.2]{pstragowski_synthetic}.
(One cannot combine this with the equivalence of \cref{thm:motivic_is_synthetic} to obtain a new proof of these computations however, as Voevodsky's result is input to the proof of \cref{thm:motivic_is_synthetic}.)



\begin{example}
    \label{ex:GIKR_Gamma_is_nu}
    Gheorghe--Isaksen--Krause--Ricka \cite[Definition~3.2]{mmf} define a functor $\Gamma \colon \Sp \to \FilSp$ given by sending $X$ to
    \[
        \Tot(\tau_{\geq 2*}(\MU^{[\bullet]}\otimes X)).
    \]
    Observe that this is a cosimplicial d\'ecalage (\cref{def:cosimplicial_decalage}) using the double-speed Whitehead filtration instead of the ordinary one.
    We may identify $\Gamma$ with the functor~$\nu$, up to some slight caveats due to working with even objects and with cosimplicial objects (which can only capture completions).
    Specifically, under the equivalence
    \[
        \Mod_{\Tot(\optau_{\geq 2*}(\MU^{[\bullet]}))}(\FilSp) \simeq \Syn_\MU^\ev
    \]
    from \cref{rmk:Syn_even_filtered_model}, we see that on the full subcategory of $\Sp$ on those spectra with even $\MU$-homology, the functor $\Gamma$ can be identified with the $\tau$-completion of~$\nu$.
    On general spectra, it identifies $\Gamma$ with the right adjoint to the inclusion $\Syn_\MU^\ev \subseteq \Syn_\MU$ applied to $\tau$-completion of $\nu$.
    (One might refer to this right adjoint composed with $\nu$ as the \emph{even synthetic analogue}.)
\end{example}

Under the equivalence of \cref{thm:motivic_is_synthetic}, the functor $\nu$ (or equivalently $\Gamma$) is useful because, in general motivic homotopy theory, there is no `motivic analogue' functor.
At least over $\bC$, this allows for the construction of new (cellular) motivic spectra.


\begin{example}
    \label{ex:GIKR_mmf_is_smf}
    In \cite[Section~5]{mmf}, Gheorghe--Isaksen--Krause--Ricka use \cref{thm:motivic_is_synthetic} to define what they call \emph{motivic modular forms}.
    Namely, they consider $\Gamma(\tmf)$, which after $p$-completion defines a $\bC$-motivic spectrum.
    It plays an important part in the computation of Isaksen--Wang--Xu \cite{isaksen_wang_xu_dimension_90}.
    
    We prefer to think of this as \emph{synthetic modular forms}, as its construction is inherently synthetic.
    Phrased in synthetic terms via \cref{ex:GIKR_Gamma_is_nu}, this definition is given by
    \[
        \smf \defeq \opnu \tmf.
    \]
    This uses that the $\MU$-homology of $\tmf$ is even.
    (Note that, as $\tmf$ is connective, it is $\MU$-nilpotent complete, so its synthetic analogue is $\tau$-complete.)
    Note also that these synthetic modular forms are different from the nonconnective versions $\Smf$ and $\SMF$ introduced in joint work of Carrick, Davies and the author in \cite{CDvN_part1}.
\end{example}

Recall from \cref{rmk:motivic_grading_Syn} that there is an alternative grading convention for the synthetic spheres known as \emph{motivic grading}.
As the name suggests, this fits better with grading conventions of (stable) motivic homotopy theory.

\newcommand*{\mot}{\mathrm{mot}}

\begin{remark}[Indexing conventions]
    \label{rmk:synthetic_vs_motivic_sphere_indexing}
    Let us write $\S^{t,w}_\mot$ for the motivic bigraded $(t,w)$-sphere.
    As explained in \cite[Section~7.1]{pstragowski_synthetic}, this amounts to
    \[
        \S^{0,0}_\mot = \opSigma^\infty_+ \bA^0 \qquad \text{and} \qquad \S^{2,1}_\mot = \opSigma^\infty_+ \bP^1.
    \]
    Let us for the moment use the motivic grading on synthetic spectra, writing $\S^{t,w}_\syn$ for the synthetic $(t,w)$-sphere in the sense of \cref{rmk:motivic_grading_Syn}.
    Then $\S^{t,w}_\syn$ is even in the sense of \cref{prop:even_synthetic_spectra} if and only if $w$ is even; see \cref{ex:even_synthetic_spheres}.
    Under the equivalence of \cref{thm:motivic_is_synthetic}, we have the correspondence between
    \[
        \S^{t,w}_\mot \qquad \text{and} \qquad \S^{t,\,2w}_\syn.
    \]
    If we instead use Adams grading on synthetic spectra, then this correspondence is between
    \[
        \S^{t,w}_\mot \qquad \text{and} \qquad \S^{t,\,2w-t}_\syn.
    \]
\end{remark}

%% file: fil/sseq_informal.tex
\chapter{Informal introduction to spectral sequences}
\label{ch:informal_sseq}


This chapter is meant as an informal introduction to spectral sequences and the role of $\tau$.
Although it is certainly possible to take this chapter as a first introduction to spectral sequences, we particularly have in mind two kinds of readers: one who wants to look at spectral sequences for the second time, and one who is familiar with these, but wants to learn about the $\tau$-formalism specifically.

For the rest of this chapter, we fix a filtered spectrum $X \colon \Z^\op \to \Sp$.
Recall from \cref{ch:filtered_spectra_sseqs} that we regard $X$ as a tool to understand its colimit $X^{-\infty}$.
For simplicity, and since this covers most of our use cases, throughout this chapter we assume that $X$ is constant from degree $0$ onwards:
\[
    \dotsb \to X^2 \to X^1 \to X^0 \congto X^{-1} \congto \dotsb.
\]
As a result, we will simply ignore the spectra in negative filtration.
Our goal then is to understand $\pi_*X^0$.

We may do this one degree at a time, so henceforth we fix an integer $n$.
Hitting the above diagram with the functor $\pi_n$, we obtain a diagram of abelian groups:
\[
    \dotsb \to \pi_n X^2 \to \pi_n X^1 \to \pi_n X^0.
\]
This filtered abelian group induces a strict filtration on $\pi_n X^0$, and this is what we aim to understand.
Our first job then should be to understand when an element in $\pi_n X^0$ is in the image of $\pi_n X^1$; in other words, to determine which elements have filtration at least~$1$.

\section{The reconstruction problem}

We have a cofibre sequence
\[
    X^1 \to X^0 \to \Gr^0 X,
\]
leading to a long exact sequence
\[
    \dotsb \to \pi_n X^1 \to \pi_n X^0 \to \pi_n \Gr^0 X \to \dotsb.
\]
This allows us to test whether $\alpha \in \pi_n X^0$ has filtration at least $1$: this happens if and only if it goes to zero under the map $\pi_n X^0 \to \pi_n \Gr^0 X$.
This pattern continues: if $\alpha \in \pi_n X$ has filtration at least $1$, we then ask if it filtration is at least $2$.
Choosing a lift to $\pi_n X^1$, we look at the associated graded $\Gr^1 X$, whose homotopy sits in a long exact sequence
\[
    \dotsb \to \pi_n X^2 \to \pi_n X^1 \to \pi_n \Gr^1 X \to \dotsb,
\]
and we can iterate this procedure until $\alpha$ does not lift further, at which point we have determined the filtration of $\alpha$.

This way of thinking only goes so far: it presupposes that we understand the elements of $\pi_n X^s$, which we usually do not.
In practice, what is more understandable is the homotopy of the associated graded.
Instead of starting with the $\pi_n X^s$, we will start with the groups $\pi_n \Gr^s X$ for all $s$, and then try to piece the $\pi_n X^s$ back together from this data.
This presents two issues:
\begin{numberenum}
    \item \label{issue:fake_elements} not every element in $\pi_n \Gr^s X$ comes from $\pi_n X^s$ (in other words, there are ``fake elements''),
    \item \label{issue:kernels} even if an element in $\pi_n \Gr^s X$ lifts to $\pi_n X^s$ (in other words, it is not ``fake''), then it may map to zero in $\pi_n X^0$.
\end{numberenum}

We can solve both of these issues using the same mechanism.
We equip the homotopy of the associated graded with more information that will make it ``remember'' the homotopy of the filtered spectrum.
This additional information comes in the form of self-maps on the associated graded, known as the \emph{differentials}.
Concretely, a differential will connect a ``fake'' element to an element that maps to zero under (a composite of) the transition maps.
As a result, we see that the purpose of these ``fake'' elements is to introduce \emph{relations} in the homotopy of $\pi_*X^0$.

\section{Differentials: obstructions to lifting}

First, let us address issue~\ref{issue:fake_elements}.
For this, we use the long exact sequence
\[
    \dotsb \to \pi_n X^{s+1} \to \pi_n X^s \to \pi_n \Gr^s X \to \pi_{n-1} X^{s+1} \to \dotsb,
\]
which tells us that an element in $\pi_n \Gr^s X$ comes from $\pi_n X^s$ if and only if it maps to zero in $\pi_{n-1}X^{s+1}$.
The question, then, is how explicit we can make this condition, where `explicit' refers to describing it in terms of the associated graded as much as possible.
It would also be helpful to organise this information in a digestible way.

To make notation easier, we will focus on the case $s=0$.
Our situation is summarised by the diagram
\[
    \begin{tikzcd}
        \dotsb \ar[r] & X^2 \ar[r] & X^1 \ar[r] & X^0 \ar[d] \\
        & & & \Gr^0 X, \ar[ul,dashed,"\partial"]
    \end{tikzcd}
\]
where the dashed arrow indicates that the map is of degree $1$: it is the boundary map $\partial \colon \Gr^0 X \to \Sigma X^1$ of a cofibre sequence.
By exactness, an element $x \in \pi_n \Gr^0 X$ comes from $\pi_n X^0$ if and only if its image in $\pi_{n-1} X^1$ is zero.
However, as we said before, we usually do not know much about the homotopy groups of $X^1$, so this is not a helpful description.
To approximate the question of the image $\partial x \in \pi_{n-1} X^1$ being zero, we first ask if its image in the associated graded $X^1 \to \Gr^1 X$ is zero:
\[
    \begin{tikzcd}
        \dotsb \ar[r] & X^2 \ar[r] & X^1 \ar[d] \ar[r] & X^0 \ar[d] \\
        & & \Gr^1 X & \Gr^0 X. \ar[ul,dashed,"\partial"]
    \end{tikzcd}
\]
Write $d_1(x)$ for the image of $\partial x$ in $\pi_{n-1}\Gr^1X$.
If $d_1(x)\neq 0$, then $\partial x \neq 0$ as well, so in particular we learn that $x$ is not in the image of $\pi_nX^0$.
If on the other hand $d_1(x)=0$, then we are not yet done: all we learn is that $\partial x \in \pi_{n-1}X^1$ lifts to $\pi_{n-1}X^2$.
Choosing a lift, we can ask the same question, testing whether this element is zero by looking at its image in $\pi_{n-1} \Gr^2 X$:
\[
    \begin{tikzcd}
        \dotsb \ar[r] & X^2 \ar[d] \ar[r] & X^1 \ar[r] & X^0 \ar[d] \\
        & \Gr^2 X & & \Gr^0 X. \ar[ul,dashed,"\partial"]
    \end{tikzcd}
\]
This choice of lift will not be unique, and neither will the resulting class in $\pi_{n-1}\Gr^2 X$; the class in $\pi_{n-1}\Gr^2 X$ is only well defined up to the image of $d_1$.
We write $d_2(x)$ for this element in $(\pi_{n-1}\Gr^2 X)/ d_1$.
If $d_2(x)$ is nonzero, then $\partial x$ is also nonzero.
If $d_2(x)$ is zero, then we continue the story and define $d_3(x)$ in $\pi_{n-1}\Gr^3 X$ (only well defined up to $d_1$ and $d_2$), et cetera.

We obtain inductively defined elements $d_r(x)$ for $r\geq 1$.
If they all vanish, then our class $x$ lifts (possibly not uniquely) to an element of the limit $X^\infty = \lim_s X^s$.
This gets us into convergence issues.
In good situations, this limit vanishes; let us assume that this is the case.
This is good news: it means that we can detect whether $\partial x \in \pi_{n-1}X^1$ is zero by checking if the $d_r(x)$ are zero for all $r\geq 1$.
This, in turn, means that we can answer the question whether $x \in \pi_n \Gr^0X$ comes from $\pi_n X$.

In summary then: we have an inductively defined list of \emph{differentials} $d_r(x)$, which (in good cases) vanish if and only if $x$ comes from an element in $\pi_n X$.
While so far we only started with classes in $\pi_n \Gr^0 X$, the same applies when starting with an element of $\pi_n \Gr^s X$, which lifts to $\pi_n X^s$ if and only if all differentials on it vanish.

\section{Differentials: kernels of transition maps}

On to issue~\ref{issue:kernels}, which is asking what the kernel of $\pi_n X^s \to \pi_n X^0$ is.
Because the map $X^s \to X^0$ is a composite of $s$ maps, we can focus on the map $X^s \to X^{s-1}$ and iterate this procedure.
Here we will encounter some of the awkwardness of working solely in terms of the associated graded.
To illustrate this, we start with an element $\alpha \in \pi_n X^s$, and write $x$ for its image in $\pi_n \Gr^s X$.
Our aim is to understand whether $\alpha$ maps to zero in $\pi_n X^{s-1}$.
We have a long exact sequence
\[
    \dotsb \to \pi_{n+1} \Gr^{s-1} X \to \pi_n X^{s} \to \pi_n X^{s-1} \to \dotsb,
\]
so by exactness, $\alpha$ maps to zero in $\pi_n X^{s-1}$ if and only if it is in the image of the map $\pi_{n+1}\Gr^s X \to \pi_n X^s$.
Notice that in terms of $x$, this is equivalent to the existence of an element $y\in \pi_{n+1}\Gr^{s-1}$ such that $d_1(y) = x$.
Iterating this procedure, we find that the element $\alpha \in \pi_n X^s$ maps to a nonzero element in $\pi_n X$ if and only if $x$ is not in the image of $d_1,\dotsc,d_s$.

\begin{remark}
    It might appear there is an asymmetry in the above: to resolve issue~\ref{issue:fake_elements}, we had to check that infinitely many differentials on $x$ vanish, whereas for issue~\ref{issue:kernels} we only have to check a condition involving finitely many differentials.
    This is due to the simplifying assumption we made earlier that the filtered spectrum is constant after filtration~$0$.
    This is equivalent to the associated graded being zero in negative filtrations.
    In effect, this means that differentials originating in filtration below~$0$ automatically vanish, so that the condition of not being hit by them is vacuous.
\end{remark}

Phrasing the previous story solely in terms of the associated graded runs into some slightly delicate matters.
By this we mean that we do not start with a class $\alpha \in \pi_n X^s$, but only with a class $x \in \pi_n \Gr^s X$.
If all differentials on $x$ vanish, and moreover $x$ is not the target of a differential, then any lift of $x$ to $\pi_n X^s$ maps to a nonzero element in $\pi_n X^0$.
However, if $d_r(y)=x$ for some $r \leq s$ and some $y$, then we only learn that \emph{there exists} a lift of $x$ to $\pi_n X^s$ that will map to zero in $\pi_n X^{s-r}$.
It is not guaranteed that every lift will satisfy this: if $\alpha \in \pi_n X^s$ is a lift of $x$, then for any $\beta \in \pi_n X^s$ that comes from $\pi_n X^{s+1}$, the element $\alpha+\beta$ also lifts $x$.
But the associated graded has no control over $\beta$: it maps to zero in $\pi_n \Gr^s X$.
This is a matter we cannot ignore, since $\beta$ need not even map to zero in $\pi_n X^0$.
The summary then is that the associated graded $\pi_n \Gr^s X$ only sees phenomena \emph{up to higher filtration}.

This is a problem that we simply have to live with if all we understand is the associated graded.
It can be delicate matter to check that a lift of an element hit by a $d_r$-differential is the lift that dies $r$ filtrations down.
In practice, one might be able to bootstrap this together by comparing different spectral sequences: in one spectral sequence, there might be no elements of higher filtration, so that there are no problems choosing a lift.
This choice can then be transported to different spectral sequences where it is not clear how to choose this lift.


\section{Graphical depiction of spectral sequences}

At this point, we need a way to organise all of this information in a way to make it more approachable for humans.
Define
\[
    \uE_1^{n,s} \defeq \pi_n \Gr^s X,
\]
and define, for every $n,s$, the first differential
\[
    d_1 \colon \uE_1^{n,s} \to \uE_1^{n-1,\, s+1}
\]
as the boundary map $\uE^{n,s}_1 \to \pi_{n-1}X^{s+1}$ followed by the projection $\pi_{n-1}X^{s+1} \to \uE_1^{n-1,\, s+1}$.

We depict these by letting the horizontal axis correspond to the stem $n$, and the vertical axis correspond to the filtration $s$.
The differential $d_1$ goes one to the left, and one up.

The differential $d_r$ goes one to the left, and $r$ units up.
This map is however only well defined after taking homologies for the preceding differentials $d_1,\dotsc, d_{r-1}$.
We therefore inductively define, for $r \geq 2$,
\[
    \uE_r^{n,s} \defeq \uH^{n,s}(\uE_{r-1}^{*,*},\, d_{r-1}) = \frac{\ker (d_{r-1} \colon \uE_{r-1}^{n,\, s} \to \uE_{r-1}^{n-1,\, s+r-1} )}{\im (d_{r-1} \colon \uE_{r-1}^{n+1,\,s-r+1} \to \uE_{r-1}^{n,\, s}) }.
\]
Roughly speaking, doing this process infinitely many times results in page $\infty$, denoted $\uE_\infty^{n,s}$.
In good cases (the other part of convergence issues), this is isomorphic to the associated graded of the induced strict filtration on $\pi_n X^0$:
\[
    \uE_\infty^{n,s} \cong \frac{F^s\, \pi_nX^0}{F^{s+1}\, \pi_n X^0}.
\]
In general, without the simplifying assumption that $X^0 \cong X^{-\infty}$, this would instead be the associated graded of the filtration on $\pi_n X^{-\infty}$.

\begin{remark}
    The reason the differential goes $r$ units up is because we are using a \emph{cohomological} indexing on the filtration.
    If we instead indexed $X$ to be a functor $\Z \to \Sp$, which is a \emph{homological} indexing on filtration, then the differential would decrease the filtration.
\end{remark}

In summary: by passing from the filtered spectrum to the associated graded, we introduced ``fake'' elements.
These elements are responsible for recording which elements die under the transition maps $\pi_n X^s \to \pi_n X^{s-1}$.
Taking homology for a $d_r$-differential removes both the fake elements, and kills elements that die under $X^s \to X^{s-r}$.
Letting all differentials run brings us to the associated graded of the filtration we were trying to understand.

\section{Reformulation in terms of \texorpdfstring{$\tau$}{tau}}
\label{sec:informal_sseq_tau}

It would be useful to introduce some notation to make it easier to describe transitions maps.
This is what $\tau$ is designed to do.
If $X$ is a filtered spectrum, let us define
\[
    \oppi_{n,s}X \defeq \pi_n (X^s).
\]
Let $\Z[\tau]$ denote the bigraded ring where $\tau$ has bidegree $(0,-1)$.
We turn $\oppi_{*,*}X$ into a bigraded $\Z[\tau]$-module by letting $\tau$ act as the transition maps.

Practically, all this means is that if $\alpha \in \oppi_{n,s}X$, then we write $\tau\cdot \alpha$ for the image of $\alpha$ under the transition map $X^s \to X^{s-1}$.
This is helpful as we do not have to give each transition map its own name, which would become quite cumbersome when expressing more involved relations.

The previous discussion then takes the following form.
If $x \in \oppi_n \Gr^s X$ is a class that is hit by a $d_r$-differential, then there exists a lift $\alpha \in \oppi_{n,s}X$ of it such that
\[
    \tau^r \cdot \alpha = 0.
\]
Further, passing from $\oppi_{n,s}X$ to $\pi_n X^{-\infty}$ is given by the colimit along the transition maps, which in this language is given by inverting $\tau$.
In this way, we see the fundamental difference between $\oppi_{*,*}X$ and $\pi_*X^{-\infty}$: while differentials in the spectral sequence are responsible for killing elements in the latter, in the former they only kill $\tau$-power multiples of it, and thereby leave a trace of their existence.
